\documentclass[a4paper,12pt]{book}

\usepackage{theorem}
\usepackage[ngerman, english]{babel}
\usepackage{amssymb} 
\usepackage{enumerate} 
\usepackage{textcomp}
\usepackage{a4wide}
\usepackage[dvips]{graphicx}
\usepackage{pdfpages}
\usepackage[utf8]{inputenc}
\usepackage{amsmath}
\usepackage{caption}

\newcommand{\rem}[1]{}
\newcommand{\R}{\mathbb{R}}
\newcommand{\N}{\mathbb{N}}
\newcommand{\Z}{\mathbb{Z}}
\newcommand{\C}{\mathbb{C}}
\newcommand{\Hei}{\mathbb{H}}
\newcommand{\G}{\mathbb{G}}
\newcommand{\T}{\mathbb{T}}

\newcommand{\tr}{\mathrm{tr}}

\renewcommand{\ker}{\text{Ker}}

\newcommand{\id}{\text{id}}
\newcommand{\en}{\text{End}}
\newcommand{\spa}{\mathrm{span}}
\newcommand{\grad}{\mathrm{grad}}
\newcommand{\pr}{\mathrm{pr}}
\newcommand{\hor}{\mathrm{hor}}
\newcommand{\rank}{\text{Rank }}
\newcommand{\hook}{\ \lrcorner \ }
\newcommand{\dbar}{\ensuremath{\, \mathchar'26\mkern-12mu d}}

\theoremstyle{break}
\newtheorem{defin}{Definition}[section]
\newtheorem{prop}[defin]{Proposition}
\newtheorem{lemma}[defin]{Lemma}

\newtheorem{thm}[defin]{Theorem}
\newtheorem{cor}[defin]{Corollary}
\newtheorem{example}[defin]{Example}
\newtheorem{remark}[defin]{Remark}
\newtheorem{ass}[defin]{Assumption}
\newcommand{\B}{\normalfont \noindent \textbf{Proof: }}
\newcommand{\eB}{$\hfill \Box$}
\newcommand{\eBsp}{$\hfill \lhd$}
\newcommand{\Bem}{\normalfont \noindent \textbf{Remark: }}
\title{}
\author{}

\begin{document}

\frontmatter
\thispagestyle{empty}
\begin{centering}

{\Huge\textbf{Spectral Triples on Carnot Manifolds}

}

\vfill

{\large

Von der Fakultät für Mathematik und Physik \\
der Gottfried Wilhelm Leibniz Universität Hannover\\
zur Erlangung des Grades \\
Doktor der Naturwissenschaften\\
Dr. rer. nat.

\vfill

genehmigte Dissertation

\vfill

von

\vfill

{\Large Dipl.Math. Stefan Hasselmann}

\bigskip

geboren am 12.04.1981 in Ibbenbüren

\vfill

2014

}

\end{centering}

\newpage

\vspace*{\fill}

{\Large

Referent: Prof. Dr. Christian Bär \\[0.2 cm]
Koreferent: Prof. Dr. Elmar Schrohe \\[0.2 cm]
eingereicht am 14.05.2013 \\[0.2 cm] 
Tag der Promotion: 18.10.2013
}

\chapter*{Abstract}

We analyze whether one can construct a spectral triple for a Carnot manifold $M$, which detects its Carnot-Carath\'{e}odory metric and its graded dimension. Therefore we construct self-adjoint horizontal Dirac operators $D^H$ and show that each horizontal Dirac operator detects the metric via Connes' formula, but we also find that in no case these operators are hypoelliptic, which means they fail to have a compact resolvent.\smallskip

First we consider an example on compact Carnot nilmanifolds in detail, where we present a construction for a horizontal Dirac operator arising via pullback from the Dirac operator on the torus. Following an approach by Christian Bär to decompose the horizontal Clifford bundle, we detect that this operator has an infinite dimensional kernel. But in spite of this, in the case of Heisenberg nilmanifolds we will be able to discover the graded dimension from the asymptotic behavior of the eigenvalues of this horizontal Dirac operator. Afterwards we turn to the general case, showing that any horizontal Dirac operator fails to be hypoelliptic. Doing this, we develop a criterion from which hypoellipticity of certain graded differential operators can be excluded by considering the situation on a Heisenberg manifold, for which a complete characterization of hypoellipticity in known by the Rockland condition.\smallskip

Finally, we show how spectral triples can be constructed from horizontal Laplacians via the Heisenberg pseudodifferential calculus developed by Richard Beals and Peter Greiner. We suggest a few of these constructions, and discuss under which assumptions it may be possible to get an equivalent metric to the Carnot-Carath\'{e}odory metric from these operators. In addition, we mention a formula by which the Carnot-Carath\'{e}odory metric can be detected from arbitrary horizontal Laplacians.\bigskip

\textbf{Keywords:} Spectral triple, Carnot-Carath\'{e}odory metric, Hypoellipticity.

\chapter*{Zusammenfassung}

Wir untersuchen, inwiefern man auf einer Carnot-Mannigfaltigkeit $M$  ein spektrales Tripel konstruieren kann, welches die Carnot-Carath\'{e}odory Metrik und die gradierte Dimension von $M$ erkennen soll. Zu diesem Zweck konstruieren wir selbst-adjungierte horizontale Dirac Operatoren $D^H$ und zeigen, dass zwar jeder horizontale Dirac Operator über Connes' Formel die Metrik erkennt, allerdings in keinem Fall hypoelliptisch ist und somit keine kompakte Resolvente besitzen kann.\smallskip

Zunächst betrachten wir ein Beispiel auf kompakten Carnot Nilmannigfaltigkeit detailliert, wobei wir eine Konstruktion für einen horizontalen Dirac Operator über den Pull-back des Dirac Operators auf dem Torus durchführen. Einer Methode von Christian Bär folgend können wir das horizontale Clifford Bündel dieses Operators zerlegen und erkennen, dass der Operator einen unendlich dimensionalen Kern besitzt. Dennoch können wir im Fall von Heisenberg Nilmannigfaltigkeiten die gradierte Dimension aus dem asymptotischen Verhalten der Eigenwerte dieses horizontalen Dirac-Operators erkennen. Anschließend wenden wir uns dem allgemeinen Fall zu, indem wir zeigen dass ein beliebiger horizontaler Dirac Operator nicht hypoelliptisch ist. Dazu entwickeln wir ein Kriterium mit dem man die Hypoelliptizität von bestimmten gradierten Differentialoperatoren ausschließen kann indem man die Situation auf einer Heisenberg Mannigfaltigkeit betrachtet, für welche eine vollständige Charakterisierung der Hypoelliptizität durch die Rockland-Bedingung gegeben ist.\smallskip

Schließlich zeigen wir wie man spektrale Tripel aus horizontalen Laplace Operatoren mit Hilfe des Heisenberg Pseudodifferentialkalküls, das von Richard Beals und Peter Greiner entwickelt wurde, konstruieren kann. Wir stellen ein paar explizite Konstruktionen vor und diskutieren, unter welchen Voraussetzungen es möglich sein kann aus diesen Operatoren eine zu der Carnot-Carath\'{e}odory Metrik äquivalente Metrik zu erhalten. Zusätzlich erwähnen wir eine Formel, mit der die Carnot-Carath\'{e}odory Metrik aus beliebigen horizontalen Laplace Operatoren erkannt werden kann.\bigskip

\textbf{Schlüsselwörter:} Spektrales Tripel, Carnot-Carath\'{e}odory Metrik, Hypoelliptizität.

\cleardoublepage

\begin{figure}
  \centering
  \includegraphics[scale=.6]{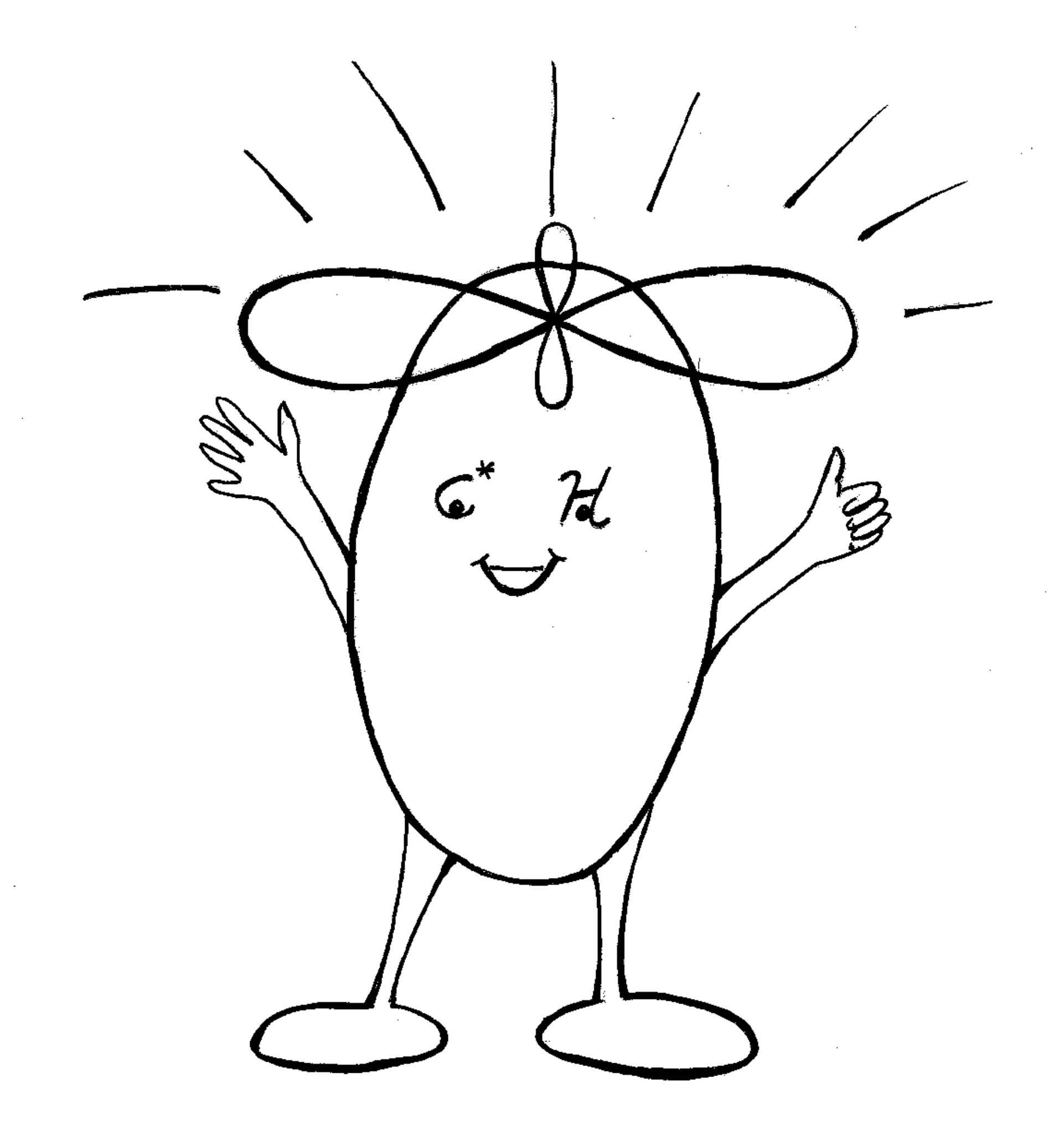}
  \caption*{\emph{A spectral triple.} By S. Wiencierz.}
  \label{img:grafik-dummy}
\end{figure}

\clearpage

\vspace*{\fill}

\begin{flushright}
 \emph{Du sollst dich zu einer Stadt begeben, die den Namen Tsnips-Eg'N-Rih trägt!} \smallskip

  (Walter Moers)
\end{flushright}

\chapter*{Introduction}

In the 1980s, Alain Connes presented the concept of non-commutative geometry as an extension of the usual notion of a topological space (see e.g. \cite{ConnesNCG}, \cite{Con1}, \cite{Con3}). The idea goes back to the 1940s, when Israel Gelfand and Mark Naimark showed that every commutative $C^\ast$-Algebra is isomorphic to the $C^\ast$-algebra $C_0(X)$ of continuous functions on a locally compact Hausdorff space $X$ vanishing at infinity (see \cite{GN}). From this starting point, one sees that a lot of properties of the space $X$ can be translated into properties of its $C^\ast$-algebra $C_0(X)$. Motivated by this, in non-commutative geometry one considers a general $C^\ast$-algebra as a \textgravedbl non-commutative space\textacutedbl.\smallskip

To describe geometry on a non-commutative space, Connes introduced so-called \emph{spectral triples}. The definition of a spectral triple is suggested by the fact that many geometric properties of a compact connected Riemannian spin manifold $M$ without boundary can be obtained from the Dirac operator $D$ acting on a Clifford bundle $\Sigma M$ over $M$. For example, one can reproduce the dimension of $M$ by the asymptotic growth of the eigenvalues of $D$ via
$$\dim M = \inf \left\{ p \in \R: \left(D^2 + I\right)^{-\frac{p}{2}} \ \text{is trace class} \right\},$$
and one can detect the geodesic distance on $M$ by the formula
$$d_{geo}(x,y) = \sup \left\{ |f(x)-f(y)|: f \in C^\infty(M), \left\|[D,f]\right\|_{L^2(\Sigma M)} \leq 1 \right\}.$$
These properties can be transported to the picture of $C^\ast$-algebras: A spectral triple is a triple $(\mathcal{A}, \mathcal{H}, \mathcal{D})$ where $\mathcal{A}$ is a $C^\ast$-algebra, $\mathcal{H}$ is a Hilbert space carrying a faithful action of $\mathcal{A}$ on $\mathcal{B}(\mathcal{H})$ and $\mathcal{D}$ is a self-adjoint operator on $\mathcal{H}$ such that $[\mathcal{D},a]$ is bounded for $a$ belonging to a dense sub-algebra of $\mathcal{A}$ and such that the resolvent $(\mathcal{D}^2+I)^{-1/2}$ of $\mathcal{D}$ is compact. For a general spectral triple, one can define notions of dimension and metric in analogy to the formulas above. Hence on a closed Riemannian spin manifold with Dirac operator $D$, a spectral triple which detects the dimension and the geodesic metric on $M$ is given by the triple $(C(M), L^2(\Sigma M), D)$.\smallskip

During the last decades there have been several approaches to construct spectral triples for more examples than the one of a closed Riemannian spin manifold. In particular there are some constructions for certain fractals which have a non-integer Hausdorff dimension, done for example by Erik Christensen, Christina Ivan and their collaborators (see e.g. \cite{CIL}, \cite{CIS}): For some fractals, it is possible to detect the Hausdorff dimension as well as the geodesic distance of the space from a spectral triple. A more general approach has been suggested by Ian Palmer (\cite{Pal}) and by John Pearson and Jean Bellissard \cite{PB}). Palmer shows that under mild conditions the Hausdorff dimension of every compact metric space can be discovered by a spectral triple. For these triples one can also find an estimate from above for the metric. But the constructions Palmer is considering lead away from the geometry of the space, since these constructions make use of an approximation of the space by a discrete subset.\smallskip

From the aspect of dimension, it is also interesting to note that there has been a construction by Erik Christensen and Christina Ivan showing that to any given positive real number $s$, one can construct a spectral triple of dimension $s$, where $\mathcal{A}$ is a limit of finite-dimensional $C^\ast$-algebras and $\mathcal{D}$ is a limit of finite-dimensional operators (\cite{CI}). Hence at least theoretically it is possible to define spectral triples of arbitrary dimension on certain spaces.\smallskip

The question we are dealing with in this thesis is whether it is possible to define two different spectral triples on one space, which both give reasonable geometries (in terms of dimension and metric). Therefore we consider so-called \emph{sub-Riemannian manifolds} (or in a more specialized setting \emph{Carnot manifolds}), which are Riemannian manifolds $M$ equipped with a bracket-generating horizontal sub-bundle $HM$ of their tangent bundle $TM$. It has been detected by Wei-Liang Chow (\cite{Cho}) in 1939 that in this case and if $M$ is connected, any two points on $M$ can be connected by a curve which is tangent to the horizontal distribution $HM$. This means that we obtain a metric on $M$ via the formula
$$d_{CC}(x,y) := \inf \left\{ \int_0^1 \left\| \dot{\gamma}(t) \right\| dt: \gamma \text{ horizontal path with } \gamma(0) = x \text{ and } \gamma(1) = y \right\},$$
which differs from the metric induced by the geodesic distance on $M$. In addition, the Hausdorff dimension of the metric space $(M,d_{CC})$ turns out to be strictly greater than the Hausdorff dimension of the metric space $(M,d_{geo})$. This result is due to John Mitchell (\cite{Mit}) and is also known as \emph{Mitchell's Measure Theorem}.\smallskip

The most important example for such a sub-Riemannian manifold is the $(2m+1)$-dimensional Heisenberg group $\Hei^{2m+1}$. It can be represented as the matrix group consisting of matrices of the form
$$\Hei^{2m+1} = \left\{ \begin{pmatrix} 1&x^t&z \\ 0&1&y \\ 0&0&1 \end{pmatrix}: x,y \in \R^m, z \in \R \right\},$$
where the group composition is given by matrix multiplication. Note that as a point set $\Hei^{2m+1}$ is isomorphic to $\R^{2m+1}$, and that $\Hei^{2m+1}$ has the structure of a graded nilpotent Lie group. The grading is induced by the Lie algebra $\mathfrak{h}_{2m+1}$ of $\Hei^{2m+1}$, which is of the form $\mathfrak{h}_{2m+1} = V_1 \oplus V_2$ such that $\dim V_1 = 2m$, $\dim V_2 = 1$ and $[V_1,V_1] = V_2$. More generally, we will consider Carnot groups. A Carnot group is a nilpotent Lie group $\G$ whose Lie algebra is carrying a grading such that $\mathfrak{g} = V_1 \oplus \ldots \oplus V_R$ with $[V_S,V_1] = V_{S+1}$ for $S \leq R-1$ and $[V_S,V_R] = 0$ for all $1 \leq S \leq R$. For our work, we will consider Carnot manifolds which will be defined to be Riemannian manifolds which carry such a grading structure on their tangent bundle. \smallskip

In this thesis, we will use the geometric approach by Dirac operators to discuss whether one can construct spectral triples on compact Carnot manifolds, which furnish the Hausdorff dimension of $(M,d_{CC})$ and the Carnot-Carath\'{e}odory metric of $d_{CC}$. We will indeed be able to construct so-called \emph{horizontal Dirac operators} $D^H$ acting on a horizontal Clifford bundle $\Sigma^H M$ over a Carnot manifold $M$ in analogy to classical Dirac operators in a quite general setting. In a frame $\{X_1, \ldots, X_d\}$ of the horizontal distribution of a Carnot manifold such an operator will have the form
$$D^H = \sum_{j=1}^d c^H(X_j) \nabla^{\Sigma^H}_{X_j} + \gamma,$$
where $c^H$ denotes the Clifford action of $HM$ and $\gamma$ is an endomorphism on the bundle. Then indeed the Carnot-Carath\'{e}odory metric can be detected from $D^H$, as we will show in detail, and it does not matter if one uses the Lipschitz functions with respect to the Carnot-Carath\'{e}odory metric or the smooth functions for the dense sub-algebra of $C(M)$ appearing in Connes' metric formula. But we will see that these operators are not hypoelliptic, which means that they do not have a compact resolvent (arising from pseudodifferential calculus). Hence they do not furnish a spectral triple, providing us the unexpected result that the theory of spectral triples does not apply to the Carnot manifolds $(M,d_{CC})$ in the way one would expect.\smallskip

After introducing horizontal Dirac operators and proving the metric formula, we will construct a class of examples. To this end we consider local homogeneous spaces of Carnot groups, which arise from the left-action of a lattice sub-group $\Gamma$ of a Carnot group $\G$. We call these closed Carnot manifold $M = \Gamma \backslash \G$ compact Carnot nilmanifolds. Then the idea is to consider the submersion
$$\pi: \Gamma \backslash M \rightarrow \T^d$$
and to define a horizontal Clifford module and a horizontal Dirac operator $D^H$ by pulling back the spinor bundle and the Dirac operator from a spin structure on $\T^d$. For the case where $\G \cong \Hei^{2m+1} \times \R^{n}$, we will be able to calculate the spectrum of this horizontal Dirac operator completely. To do this, we adapt an argument by Christian Bär to our situation(\cite{Bae}): Bär decomposes the spinor bundle belonging to a classical Dirac operator on a Heisenberg nilmanifold $\Gamma \backslash \Hei^{2m+1}$ into its irreducible components under which the operator is invariant. From these irreducible components, he is able to calculate the spectrum of the Dirac operators by means of the representation theory of the Heisenberg group.\smallskip

In our case, the horizontal Clifford bundle decomposes in the same way. We will present these calculations in detail, and afterwards we will deduce that our horizontal pull-back Dirac operator $D^H$ on a compact Carnot nilmanifold possesses at least one infinite dimensional eigenspace, hence it cannot have a compact resolvent. In spite of that we will be able to detect the Hausdorff dimension of $(M,d_{CC})$ in the Heisenberg group case from the asymptotic behavior of the non-degenerate eigenvalues of $D^H$. We will also extend this approach to the setting of the compact nilmanifold of a general Carnot group $\G$, where we will show that in the spectral decomposition of such a group there is at least one infinite dimensional eigenspace. The reason for this is that in the spectral decomposition of $D^H$ there is at least one subspace isomorphic to the space of $L^2$-sections on a horizontal Clifford bundle belonging to a group of nilpotency step $2$. On this subspace, we will be able to detect an infinite dimensional eigenspace of the horizontal pull-back Dirac operator from the Heisenberg case.\smallskip

Further we will see that the problem we detected in this concrete example is not due to a bad choice of the horizontal Dirac operator. In fact it is a general phenomenon: We will argue for that using techniques from pseudodifferential calculus. There is a calculus invented by Richard Beals and Peter Greiner (see \cite{BG}) on Heisenberg manifolds, from which it can be derived that hypoellipticity of a (self-adjoint) operator of positive order implies that it has a compact resolvent and that furthermore the Hausdorff dimension of the Heisenberg manifold can be detected by the eigenvalue asymptotics of this operator. We will prove in this thesis that any horizontal Dirac operator on an arbitrary Carnot manifold cannot be hypoelliptic. This is a big difference to the classical case, where any Dirac operator is elliptic.\smallskip

On the way, we will develop a criterion for the non-hypoellipticity of an arbitrary graded differential operator of the from
$$D = D(X_1, \ldots, X_n) \in \mathcal{U}(\mathfrak{g}),$$
where $\{X_1, \ldots, X_n\}$ is a frame of $\R^n$ which is forming a graded nilpotent Lie algebra $\mathfrak{g}$. As soon as there is a graded differential operator
$$\tilde{D} = \tilde{D}(\tilde{X}_1, \ldots \tilde{X}_m) \in \mathcal{U}(\tilde{\mathfrak{g}}) \ \text{for} \ m < n, $$
induced by a projection $\pr: X_j \mapsto \tilde{X}_j$ of $\mathfrak{g}$ onto a lower dimensional Lie algebra $\tilde{\mathfrak{g}}$, which is not hypoelliptic, then $D$ cannot be hypoelliptic. This criterion will serve us well, since we will be able to reduce the problem of showing general non-hypoellipticity of a horizontal Dirac operator to the Heisenberg case. In the Heisenberg case we have a complete characterization of hypoellipticity of horizontal Laplacians arising from the representation theory of the Heisenberg group, from which we will be able to exclude that a horizontal Dirac operator is hypoelliptic.\smallskip

Finally, once having introduced the Heisenberg pseudodifferential calculus, we make use of this calculus and show how hypoelliptic Heisenberg pseudodifferential operators furnishing a spectral triple and detecting in addition the Hausdorff dimension of the Heisenberg manifold can be constructed. We will suggest a few concrete operators, but it remains unclear whether one can detect or at least estimate the Carnot-Carath\'{e}odory metric from them. But we will show that the Carnot-Carath\'{e}odory metric can be detected by horizontal Laplacians instead of horizontal Dirac operators via the formula
$$  d_{CC}(x,y) = \sup \left\{ |f(x) - f(y)|: f \in C^\infty(M), \left\| \frac{1}{2} \left[\left[\Delta^{\mathrm{hor}},f\right],f\right] \right\| \leq 1 \right\}.$$ 
Hence, maybe in this case the second order operators are the right operators to look at when we want to do non-commutative geometry on Carnot manifolds.\smallskip

We conclude the thesis with some approaches from which one may be able to estimate or to approximate the Connes metric by first order Heisenberg pseudodifferential operators, and we present some criteria from which such an estimate would follow.\smallskip

The structure of this thesis will be the following.

\begin{itemize}
 \item In Chapter 1, we introduce the notion of spectral triples, their metric dimension and the Connes metric, and state a few well-known examples. Then we will turn to the more general approach by Mark Rieffel of compact quantum order-unit spaces and prove a few criteria to detect convergence of a family of spectral triples to such a compact quantum metric space (which does not necessarily have to be a spectral triple).

 \item In Chapter 2, we give an overview of sub-Riemannian geometry and Carnot manifolds. We introduce the Carnot-Carath\'{e}odory metric and state some important theorems in this context. The Sections 2.3 and 2.4 serve to introduce the concept of a Levi form and of certain submersions between Carnot groups which will be of importance later.

 \item In Chapter 3, we construct the horizontal Dirac operator. We start by analyzing horizontal connections, which will be the connections we want a horizontal Clifford bundle to be compatible with. Then we construct a self-adjoint horizontal Dirac operator on these bundles. In preparation of what we need later we will calculate its square locally and we will state a proposition about the eigenvalues of certain sums of Clifford matrices. Finally, in Section 3.3, we will show that any horizontal Dirac operator on a Carnot manifold $M$ detects the Carnot-Carath\'{e}odory metric via Connes' metric formula, since the norm of the commutator $[D^H,f]$ coincides with the Lip-norm of $f$ with respect to the Carnot-Carath\'{e}odory metric. We will also show that the metric is already detected by the smooth functions on $M$.

 \item In Chapter 4, we treat in detail the example of nilmanifolds $M = \Gamma \backslash \G$ from Carnot groups, arising from the left action of the standard lattice. First we construct a horizontal Dirac operator on $M$ by pulling back the spinor bundle of the horizontal torus. Then it is our aim to show that this horizontal Dirac operator possesses infinite dimensional eigenspaces. Therefore we use an approach which was presented by Christian Bär and Bernd Amman (\cite{Bae}, \cite{AB}) for the case of the classical Dirac operator on Heisenberg nilmanifolds to find a spectral decomposition of the horizontal Clifford bundle $\Sigma^H M$ which is invariant under $D^H$. For the case of a general Carnot group, we find that one part of this decomposition is isomorphic to the horizontal Clifford bundle belonging to a horizontal pull-back operator of a Carnot nilmanifold of lower commutator step. For the case of Heisenberg nilmanifold (where the horizontal distribution is of co-dimension 1) we will be able to calculate all eigenvalues of $D^H$ from the approach by Bär and Ammann, and we will detect that $D^H$ has an infinite-dimensional kernel. But in spite of that we will show that the Hausdorff dimension of $M$ can still be detected by these eigenvalues for the case $\G \cong \Hei^{2m+1}$. In section 4.4 we will put all these results together to show that on any Carnot nilmanifold arising from the left-action of the standard lattice sub-group of $\G$ the horizontal Dirac operator we constructed has an infinite dimensional kernel.

 \item In Chapter 5, we introduce the Heisenberg pseudodifferential calculus while referring to Richard Beals and Peter Greiner (\cite{BG}) and to Rapha\"{e}l Ponge. We will discuss symbol classes, composition of operators, parametrices and the role of hypoellipticity. Finally we mention the well-known results that hypoellipticity in the Heisenberg calculus implies the existence of complex powers and certain eigenvalue asymptotics.

 \item In Chapter 6, we show that the absence of a compact resolvent detected in Chapter 4 for a specific class of examples is of general nature. In particular we show that any horizontal Dirac operator is not hypoelliptic. Therefore we first review some classical hypoellipticity theorems with special attention to the case of Heisenberg manifolds. Then in Section 6.2 we prove a theorem from which hypoellipticity can be excluded for certain graded differential operators by going back to a graded differential operator acting on a lower dimensional Carnot group. We will present some consequences which arise from going back to the Heisenberg case via this reduction: This reduction criterion provides us with the possibility to exclude hypoellipticity of an arbitrary graded differential operator by looking at the Levi form of its underlying graded Lie algebra. Finally, we will apply the reduction criterion to show that an arbitrary horizontal Dirac operator cannot be hypoelliptic.

 \item In Chapter 7, we discuss the possibility of constructing spectral triples by taking square roots of horizontal Laplacians and the question if one can get any metric information from these spectral triples. First we show how the application of Heisenberg calculus furnishes spectral triples, and we suggest a few operators being not too far away from the horizontal Dirac operator which produce spectral triples. In Section 7.2, we show that the Carnot-Carath\'{e}odory metric can be detected by any horizontal Laplacian, while the business is much more difficult if we want to detect the metric from an arbitrary first-order hypoelliptic and self-adjoint Heisenberg operator. Unfortunately, all we can present concerning the last question are some criteria and some ideas which might approximate the Carnot-Carath\'{e}odory metric by metrics arising from spectral triples, but so far we have not been able to prove such an approximation or estimate completely.
\end{itemize}

\chapter*{Acknowledgements}

Before starting with the text of the thesis, there must be some time to say thank you to all the people who supported me during the last three years. First of all, I want to express my gratitude to my advisor Elmar Schrohe who was always interested in my work and found time for discussion, no matter how busy he was. He supported me all the time and motivated me, especially during the times when I got stuck with my work. I am deeply grateful to the members of the Institute of Analysis and the Graduiertenkolleg 1643 for a warm and friendly atmosphere at work and many interesting discussions, not only on mathematics. I found some wonderful friends there! In particular I have to mention Magnus Goffeng in this context, a great mathematician, who supported my work with many discussions and hints. He found the time to proofread this thesis and made some really helpful remarks. Further more I say thank you to Erik Christensen for hosting me in Copenhagen for two months, and to Christian Bär for accepting the task of being the second referee for this thesis. And it is an important matter to me to send some special thanks to Natascha Krienen and Susanne Rudolph, not only for taking care of all the administrative stuff but also for being the ``good souls'' of the Graduiertenkolleg and the institute.\smallskip

Thanks for financial supports are adressed to the DFG, since for most of the time of my PhD studies I was a member of the Graduiertenkolleg 1643 - Analysis, Geometry and String Theory. But I think I would never have finished this thesis without all the support from my friends and my parents. Finally my gratitude towards my dear wife Susanne cannot be expressend in words: She really had a hard time with me during the last months listening to all of my lamentation, but nonetheless she always kept on encouraging me and giving me amounts of love.

\tableofcontents

\mainmatter

\chapter{Non-commutative Metric Spaces}

The intention of the first chapter of this thesis is to introduce the concepts of non-commutative geometry we want to consider. In Section 1.1 we will present the concept of a spectral triple introduced by Alain Connes, which describes the geometry of a space in non-commutative geometry. This will only be a rough introduction covering only the notion of dimension and metric arising from a spectral triple, since these are the objects we are interested in. After defining these objects, we will give a few classical examples for commutative spectral triples to motivate the work of this thesis. The Carnot-Carath\'{e}odory spaces which we want to examine from the point of view of non-commutative geometry will be treated in Chapter 2.\smallskip

Since it will turn out that we do not get spectral triples in the sense of Connes from our constructions, we will state a modified definition of a degenerate spectral triple, to denote an object being a spectral triple except for a degeneration of certain eigenspaces of its Dirac operator. Later in Chapter 4 we will see that this definition fits in our situation. Besides that, in Section 1.2 we will present the more general concept of compact quantum metric spaces, which was introduced by Mark Rieffel and which only provides metric information in the abstract setting of Lip-norms. A Lip-norm can arise from an operator which does not provide a spectral triple, as it will be the case for our situation.\smallskip

We will close Section 1.2 by some simple observations we made leading to a sufficient criterion for the convergence of a sequence of compact quantum metric spaces to a certain given space. In detail, we show that the metrics $\rho_\theta$ arising from a family of Lip-norms $L_\theta$ converge to a given metric $\rho$, if the corresponding Lip-norms converge to the Lip-norm $L$ corresponding to $\rho$ in a uniform way. We will refer to these observations in the final chapter of this thesis, where we will be in the situation that we have a family of spectral triples providing metrics close to the metric we want to discover. But the desired metric itself will be detected by an operator which does not furnish a spectral triple. \bigskip

\section{Spectral Triples}

We will give a rough overview about the theory of spectral triples now, mentioning only the concepts we will use in this thesis. Most of the definitions and examples we give are well known and can be found in any textbook on non-commutative geometry, see e.g. \cite{ConnesNCG}, \cite{GVF} or \cite{Landi}. First of all we define what a spectral triple is.\medskip

\begin{defin} \label{def spectral triple} 
 Let $\mathcal{A}$ be a unital $C^\ast$-algebra and let $\mathcal{H}$ be a Hilbert space which carries an injective unitary representation $\pi: \mathcal{A} \rightarrow \mathcal{B}(\mathcal{H})$. Furthermore let $\mathcal{D}$ be an unbounded self-adjoint operator on $\mathcal{H}$ such that
\begin{enumerate}[(i)]
 \item The algebra
\begin{equation} \label{eq dense subalgebra of spectral triple}
 \mathcal{A}' := \left\{ a \in \mathcal{A}: [\mathcal{D},\pi(a)] \text{ is densely defined and bounded } \right\}
\end{equation}
is a dense sub-algebra of $\mathcal{A}$.
 \item For any number $\lambda \notin \mathrm{spec}(\mathcal{D})$ the operator $(\mathcal{D}-\lambda I)^{-1}$ is compact.
\end{enumerate}
Then the triple $(\mathcal{A}, \mathcal{H}, D)$ is called a \emph{spectral triple}. \eBsp \smallskip

\Bem Considering condition (ii), it follows from the Hilbert identity that the property of $(\mathcal{D} -\lambda I)^{-1}$ being compact for one $\lambda \notin \mathrm{spec}(\mathcal{D})$ implies that $(\mathcal{D}-\lambda I)^{-1}$ is compact for every $\lambda \notin \mathrm{spec}(\mathcal{D})$. Hence one often reformulates condition (ii) in the sense that
$$\left( \mathcal{D}^2 + I \right)^{-\frac{1}{2}} \in \mathcal{K}(\mathcal{H}),$$
or by simply demanding that \emph{$\mathcal{D}$ has a compact resolvent}. \eBsp
\end{defin}\medskip

For a spectral triple we can define a notion of dimension.\medskip

\begin{defin} \label{def metric dimension} 
 A spectral triple $(\mathcal{A}, \mathcal{H}, \mathcal{D})$ is called $s$-summable for some $s \in \R$, if we have
\begin{equation} \label{eq condition s-summable}
 \left(\mathcal{D}^2 + I\right)^{-\frac{s}{2}} \in \mathcal{L}_1(\mathcal{H}),
\end{equation}
 where $\mathcal{L}_1(\mathcal{H}) \subset \mathcal{K}(\mathcal{H})$ denotes the ideal of trace-class operators on $\mathcal{H}$. The number
\begin{equation} \label{eq condition metric dimension}
 s_0 = \inf \left\{s \in \R: \left(\mathcal{D}^2 + I\right)^{-\frac{s}{2}} \in \mathcal{L}_1(\mathcal{H}) \right\}
\end{equation}
is called the \emph{metric dimension} of $(\mathcal{A}, \mathcal{H}, \mathcal{D})$. \eBsp \smallskip

\Bem Like in Definition \ref{def spectral triple}, one can replace the operator $(\mathcal{D}^2 + I)^{-s/2}$ in \eqref{eq condition s-summable} and \eqref{eq condition metric dimension} by any operator of the form $(\mathcal{D}-\lambda I)^{-s}$ for $\lambda \notin \mathrm{spec}(\mathcal{D})$. \eBsp
\end{defin}\medskip

One can also define a metric from a spectral triple $(\mathcal{A}, \mathcal{H}, \mathcal{D})$ on the state space $\mathcal{S}(\mathcal{A})$ of its $C^\ast$-algebra $\mathcal{A}$, equipped with the weak $^\ast$-topology. It is not hard to check that the expression
\begin{equation} \label{eq Connes metric}
 d_\mathcal{D}(\phi,\psi) := \sup \left\{ |\phi(a) - \psi(a)|: a \in \mathcal{A}', \left\| [\mathcal{D},\pi(a)] \right\| \leq 1 \right\},
\end{equation}
where $\mathcal{A}'$ is the dense sub-algebra \eqref{eq dense subalgebra of spectral triple} from condition (i) of the definition of a spectral triple, gives a metric on $\mathcal{S}(\mathcal{A})$.

\begin{defin} \label{def metric on state space} 
 Let $(\mathcal{A}, \mathcal{H}, \mathcal{D})$ be a spectral triple, $\mathcal{S}(\mathcal{A})$ be the state space of its $C^\ast$-algebra $\mathcal{A}$. Then the metric $d_\mathcal{D}$ on $\mathcal{S}(\mathcal{A})$ from \eqref{eq Connes metric} is called the \emph{Connes metric} arising from $(\mathcal{A}, \mathcal{H}, \mathcal{D})$. \eBsp
\end{defin}\medskip

It turns out that it is sufficient to consider the expression \eqref{eq Connes metric} only for the positive elements of $\mathcal{A}$. We cite this proposition here, for a proof we refer to \cite{IKM}, Section 2, Lemma 1.\medskip

\begin{prop} \label{prop metric from positive elements}
 Let $(\mathcal{A}, \mathcal{H}, \mathcal{D})$ be a spectral triple, and let $\mathcal{A}_+$ denote the subset of positive elements of the $C^\ast$-algebra $\mathcal{A}$. Then the Connes metric \eqref{eq Connes metric} on the state space $\mathcal{S}(\mathcal{A})$ is given by
$$d_\mathcal{D}(\phi,\psi) = \sup \left\{ |\phi(a) - \psi(a)|: a \in \mathcal{A}_+, \left\| [\mathcal{D},\pi(a)] \right\| \leq 1 \right\}$$
for all $\phi, \psi \in \mathcal{S}(\mathcal{A})$. \eB 
\end{prop}\medskip

Spectral triples can be defined on arbitrary $C^\ast$-algebras. In non-commutative geometry, the notion of a $C^\ast$-algebra replaces in a way the notion of a space. This is motivated by the Gelfand-Naimark theory: The famous theorem of Israel Gelfand and Mark Naimark (see \cite{GN}) states that every unital commutative $C^\ast$-algebra $\mathcal{A}$ is isometrically isomorphic to the $C^\ast$-algebra $C(X)$ of the continuous functions on some compact Hausdorff space $X$. In this case the state space of $\mathcal{A}$ is exactly the space $X$, which turns \eqref{eq Connes metric} into the expression
\begin{equation} \label{eq Connes metric commutative version}
 d_\mathcal{D}(x,y) := \sup \left\{ |f(x) - f(y)|: f \in \mathcal{A}' \subset C(X), \left\| [\mathcal{D},\pi(a)] \right\| \leq 1 \right\},
\end{equation}
defining a metric on $X$. In this thesis, we will only consider this commutative situation.\smallskip

We will now state the classical (and most important) example of a commutative spectral triple, where $X$ is a compact Riemannian manifold. The central point of this thesis is to construct another commutative spectral triple on this manifold $X$ which describes a different geometry, meaning that one gets a different dimension and a different metric from it.\medskip

\begin{example} \label{ex canonical spectral triple on manifold} \normalfont
 Let $M$ be an $n$-dimensional closed Riemannian manifold equipped with a Clifford bundle $\Sigma M$ and a Dirac operator $D$. We further consider the algebra $C(M)$ of continuous functions on $M$ and the Hilbert space $L^2(\Sigma M)$ of $L^2$-sections in the spinor bundle, on which $C(M)$ has a representation by left multiplication. Then we have the following (see e.g. \cite{GVF}, \cite{Landi} or \cite{ConnesNCG}):
\begin{itemize}
 \item The triple $(C(M), L^2(\Sigma M), D)$ is a spectral triple. A dense sub-algebra of $C(M)$, which furnishes a bounded commutator with $D$, is given by $C^\infty(M)$ or also by the Lipschitz functions $\mathrm{Lip}(M)$ on $M$. And since $D$ is an elliptic and self-adjoint differential operator of order $1$ with discrete spectrum, it has a compact resolvent.

 \item The metric dimension of this spectral triple is exactly the dimension $n$ of the manifold: By Weyl asymptotics we have $\lambda_k \sim k^{1/n}$ for the eigenvalues of $D$ (since it is an elliptic and self-adjoint differential operator of order $1$), which implies that for $\lambda \notin \mathrm{spec}(D)$ the operator $(D - \lambda I)^{-s}$ is trace class if and only if $s>n$.

 \item The geodesic distance $d_{geo}$ on $M$ can be detected by the formula
 $$d_{geo}(x,y) = d_D(x,y) = \sup \left\{ |f(x) - f(y)|: f \in \mathrm{Lip}(M), \left\| [D,\pi(f)] \right\| \leq 1 \right\}.$$
This is the case, since $\left\| [D,\pi(f)] \right\|$ is exactly the essential supremum of the gradient of $f$: One can show 
$$|f(x)-f(y)| = \int_\gamma df \leq \mathrm{ess} \sup \left\| df \right\| \int_0^1 \left\| \dot{\gamma}(t) \right\| dt$$
for every geodesic curve connecting $x$ and $y$, which gives  $d_{geo}(x,y) \geq d_D(x,y)$. On the other hand we have $d_D(x,y) \geq d_{geo}(x,y)$ since for a fixed $y \in M$ the essential supremum norm of the gradient of the function $g(x) := d_{geo}(x,y)$ is bounded by $1$.
\end{itemize} \eBsp
\end{example}\medskip

In the commutative world, one can detect further examples for spectral triples which give a meaningful geometry by considering the Hausdorff dimension of a metric space. We briefly sketch its construction. \medskip

\begin{defin} \label{def Hausdorff dimension}
 Let $(X,d)$ be a metric space, $\Omega \subset X$.
\begin{enumerate}[(i)]
 \item Let $s>0$. For $\varepsilon > 0$ we set
$$\mu^s_\varepsilon(\Omega) := \inf \left\{ \sum_{\alpha} \left(\mathrm{diam}(U_\alpha)\right)^s: \ \mathcal{U} = \{U_\alpha\} \text{ open cover of } \Omega, \ \sup_{\alpha} \left(\mathrm{diam}(U_\alpha)\right) < \varepsilon \right\}.$$
Then the \emph{$s$-dimensional Hausdorff measure} of $\Omega$ is given by the number
$$\mu^s(\Omega) = \lim_{\varepsilon \rightarrow 0} \mu^s_\varepsilon(\Omega).$$
 \item The unique number $0 \leq \dim_H(\Omega) \leq \infty$ such that $\mu^s(\Omega) = \infty$ for all $s<\dim_H(\Omega)$ and $\mu^s(\Omega) = 0$ for all $s>\dim_H(\Omega)$ is called the \emph{Hausdorff dimension} of $\Omega$.
\end{enumerate} \eBsp
\end{defin}\medskip

The Hausdorff dimension is in particular used in fractal geometry, where one gets metric spaces of non-integer Hausdorff dimension. One classical example for a spectral triple by Alain Connes is a spectral triple for the Cantor set (see \cite{ConnesNCG}, Chapter 3.4.$\varepsilon$), whose metric dimension coincides with the Hausdorff dimension of this set. During the last years, spectral triples were constructed for fractals, detecting the Hausdorff dimension and the geodesic distance on this fractal. We refer for example to the work of Erik Christensen, Christina Ivan and their collaborators for that (see e.g. \cite{CIL} or \cite{CIS}).\smallskip

Under mild conditions, it is indeed possible to construct a spectral triple for an arbitrary compact metric space $X$. The idea is to approximate the space by finite sets of points. Since $X$ is assumed to be compact, one can find a sequence of open covers of $X$ whose diameter decreases to zero. Then one chooses two points in any open set belonging to one of the covers. Using this sequence, it is possible to define a Dirac operator of a spectral triple on the Hilbert space of $l^2$-sequences indexed by the sets belonging to the covers with values in $\C^2$. This construction can be found in detail in the PhD Thesis of Ian Christian Palmer (see \cite{Pal}) and has also been used by John Pearson and Jean Bellissard (see e.g. \cite{PB}). It is shown that there exists such a sequence such that the corresponding spectral triple gives back the Hausdorff dimension of the space (\cite{Pal}, Theorem 4.2.2), and that the Connes metric provided by this spectral triple dominates the original metric on $X$ (\cite{Pal}, Proposition 5.2.1).\smallskip

In this thesis, we will consider compact Riemannian manifolds which carry a second metric besides the geodesic one: These so-called \emph{Carnot manifolds} will be introduced in Chapter 2. Indeed we will be able to construct a differential operator detecting the second metric, but it will turn out that this operator does not provide a spectral triple. But we can modify the definition of a spectral triple a little bit, such that it will fit to our situation.\medskip

\begin{defin} \label{def degenerate spectral tripel}
 Let $\mathcal{A}$ be a unital $C^\ast$-algebra and let $\mathcal{H}$ be a Hilbert space which carries an injective unitary representation $\pi_ \mathcal{A} \rightarrow \mathcal{B}(\mathcal{H})$. Furthermore let $\mathcal{D}$ be an unbounded self-adjoint operator on $\mathcal{H}$ such that
\begin{enumerate}[(i)]
 \item The algebra
\begin{equation*}
 \mathcal{A}' := \left\{ a \in \mathcal{A}: [\mathcal{D},\pi(a)] \text{ is densely defined and bounded } \right\}
\end{equation*}
is a dense sub-algebra of $\mathcal{A}$.
 \item The spectrum of the operator $\mathcal{D}$ is discrete. If in addition $\Lambda$ denotes the set of all eigenvalues $\lambda$ of $\mathcal{D}$ which have an infinite dimensional eigenspace $(\lambda) \subset \mathcal{H}$, then for any number $\mu \notin \mathrm{spec}(D)$ the operator 
$$\left( (D-\mu I) \bigg{|}_{\left(\bigoplus_{\lambda \in \Lambda} E(\lambda)\right)^\perp} \right)^{-1}$$
 is compact.
\end{enumerate}
Then we call the triple $(\mathcal{A}, \mathcal{H}, \mathcal{D})$ a \emph{degenerate spectral triple}. \eBsp
\end{defin}\bigskip

\section{Compact Quantum Order Unit Spaces}

In this section we consider a more general version of a non-commutative metric space. The idea is that for any compact metric space $(X,d)$, on can associate a Lipschitz semi-norm to the space $C(X)$ of continuous functions on $X$ via
$$L(f) := \sup_{x \neq y} \frac{|f(x)-f(y)|}{d(x,y)}.$$
Now since a non-commutative space can be viewed as a generalization of the (commutative) $C^\ast$-algebra $C(X)$, the idea is to introduce the notion of a Lip-norm on a general unital $C^\ast$-algebra, from which one gets a metric on the state space of $\mathcal{A}$. This approach has been introduced by Mark Rieffel in \cite{Rie2} and \cite{Rie1}. In his approach, instead of $C^\ast$-algebras, Rieffel works in the more abstract setting of \emph{order-unit spaces} (see \cite{Rie}).\medskip

\begin{defin} \label{def order-unit space}
 An \emph{order-unit space} is a real partially ordered vector space $A$ with a distinguished element $e$ (the \emph{order unit}), which satisfies
\begin{enumerate} [(i)]
 \item For each $a \in A$ there is an $r \in \R$ such that $a \leq re$.
 \item If $a \in A$ and if $a \leq re$ for all $r \in \R$ with $r>0$, then $a \leq 0$.
\end{enumerate}
Furthermore, the norm of an order-unit space is given by
$$\left\| a \right\| := \inf \left\{ r\in \R: -re \leq a \leq re \right\}.$$
\eBsp
\end{defin}\medskip

Note that the self-adjoint elements of a unital $C^\ast$-algebra $\mathcal{A}$ form an order unit space, such that the non-commutative spaces in the sense of Alain Connes are included in this setting. Now we introduce a semi-norm on these spaces. \medskip

\begin{defin} \label{def compact quantum metric space}
 Let $(A,e)$ be an order-unit space. Then a \emph{Lip-norm} on $A$ is a semi-norm $L$ on $A$ with the following properties. 
\begin{enumerate}[(i)]
 \item We have $L(a) = 0$ $\Leftrightarrow$ $a \in \R e$.
 \item The topology on the state space $\mathcal{S}(A)$ of $A$ from the metric
\begin{equation} \label{eq metric from Lip-norm}
 \rho_L(\phi, \psi) := \sup \left\{ |\phi(a) - \psi(a)|: L(a) \leq 1 \right\}
\end{equation}
is the $w^\ast$-topology.
\end{enumerate}
We call a pair $(A,L)$, consisting of an order-unit space $A$ and a Lip-norm $L$ on $A$ a \emph{compact quantum metric space}. \eBsp
\end{defin}\medskip

We recognize that the metric \eqref{eq metric from Lip-norm} looks similar to the metric \eqref{eq Connes metric} arising from a spectral triple. Indeed, if $(\mathcal{A},\mathcal{H},\mathcal{D})$ is a spectral triple, where the dense sub-algebra of $\mathcal{A}$ for which $[\mathcal{D},a]$ is bounded is denoted by $\mathcal{A}'$, then the term
$$L(a) := \left\| [\mathcal{D},a] \right\|$$
gives a Lip-norm on the order unit space $A$ consisting of self-adjoint elements of $\mathcal{A}'$. By Proposition \ref{prop metric from positive elements}, the Connes metric \eqref{eq Connes metric} does only depend on the positive elements of $\mathcal{A}'$, and hence it is in particular determined by the self-adjoint elements of $\mathcal{A}'$ forming the order-unit space $A$.\medskip

\begin{remark} \normalfont
 It is shown in \cite{Rie1}, Section 11, that for any essentially self-adjoint operator $D$ on a Hilbert space $\mathcal{H}$, which carries a faithful representation $\pi$ of an order-unit space $A$ such that $[D,\pi(a)]$ is bounded for a dense subspace $A'$ of $A$, furnishes a Lip-norm via 
$$L(a) := \|[D,\pi(a)]\|.$$
 On the other hand, in \cite{Rie}, Appendix 2, it is shown that for any lower semi-continuous Lip-norm $L$ on an order-unit space $A$ one can define such an operator $D$ describing $L$ via $L(f) = \|[D,f]\|$ on the space $C(\mathcal{S}(A))$ of continuous functions on the state space of $A$.\smallskip

 Note that in both cases the operator $D$ does not need to have a compact resolvent. In particular, the operator $D$ constructed to a given semi-norm $L$ is in general far away from having a compact resolvent. \eBsp
\end{remark}\medskip

We will see that for the compact spaces we study in this thesis, the metric will be detected by a Dirac operator which does not have a compact resolvent. But it will fit into the setting of a compact quantum metric space. And in addition, we will detect that we can repair the lack of not being a spectral triple by disturbing the operator a little bit. Then the question arises whether the metrics from the disturbed spectral triples are equivalent to the original metric or whether they even converge towards this metric as the perturbation goes to zero. Sadly we have not been able to give a satisfying answer to this in the later chapters of this thesis. But nonetheless we have made some simple observations, which give criteria for the convergence of a series of compact quantum metric spaces, which we will present now.\smallskip 

We will see that the metric convergence of a family of compact quantum metric spaces will follow from the condition that the corresponding family of Lip-norms is uniformly continuous. Our version of uniform continuity of Lip-norms is characterized by the following definition:\medskip

\begin{defin} \label{def uniformly continuous families of cqms}
 Let $A$ be a commutative order unit space (which means that it can be realized as the set of continuous real-valued functions on a compact metric space), and for $\theta \in [0,1]$ let $L_\theta$ be a family of Lip-norms on $A$. We set
\begin{equation} \label{eq unit sphere cqms}
 \Sigma_0 := \left\{ f \in A: L_0(f) = 1\right\}.
\end{equation}
 Then we call the family $\theta \mapsto (A,L_\theta)$ of compact quantum metric spaces \emph{uniformly continuous} towards $(A,L_0)$, if
\begin{equation} \label{eq condition unifomly continuous family of cqms}
 \forall \varepsilon > 0 \ \exists \delta >0: \ 0 < \theta < \delta \ \Rightarrow \ \left| L_\theta(f) - L_0(f) \right| < \varepsilon \ \forall f \in \Sigma_0.
\end{equation} \eBsp
\end{defin} \medskip

We will see now that the property of uniform continuity of a family of Lip-norms implies the fact that the metrics $\rho_\theta$ corresponding to $(A,L_\theta)$ converge to $\rho_0$, and that they are equivalent in case $\theta$ is small enough. Remember that $\rho_\theta$ is defined on the state space $\mathcal{S}(A)$ of $A$, which is isometrically isomorphic to a compact metric space $X$ with $A \cong C(X)$. In this case the metric is given via the formula
$$\rho_\theta (x,y) = \sup \left\{|f(x)-f(y)|: L_\theta(f) \leq 1 \right\}.$$
Then we have the following proposition. \medskip

\begin{prop} \label{prop uniformly continuity implies convergence of metrices}
 Let $A$ be a commutative order unit space and for $\theta \in [0,1]$ let $L_\theta$ be a family of Lip-norms on $A$, such that the family $\theta \mapsto (A,L_\theta)$ of compact quantum metric spaces is uniformly continuous in the sense of \eqref{eq condition unifomly continuous family of cqms}.\smallskip

Then for every $\varepsilon > 0$ with $\varepsilon < 1$ there is a $\delta > 0$ such that
\begin{equation} \label{eq equivalence and convergence of metricies}
 \left(1-\varepsilon\right) \rho_\theta(x,y) \leq \rho_0(x,y) \leq \left(1+\varepsilon\right) \rho_\theta(x,y) \ \ \ \forall x,y \in \mathcal{S}(A)
\end{equation}
 for every $0 < \theta < \delta$. \smallskip

\B Since $A$ is commutative, we consider 
$$A = \left\{f \in C(X): f = f^\ast\right\} \subset C(X)$$
for a compact metric space $X$ such that $\mathcal{S}(A) \cong X$. Note that we can forget about the constant functions when calculating $\rho_\theta$ since in this case we have $|f(x)-f(y)|=0$ for any two points $x$ and $y$.\smallskip

Now assume $\theta \mapsto (A,L_\theta)$ is uniformly continuous, and let $\varepsilon > 0$. Then we find a $\delta >0$ such that
$$\left(1-\varepsilon\right) \leq L_\theta(\tilde{f}) \leq \left(1+\varepsilon\right) \ \forall \theta < \delta$$
for every $\tilde{f} \in A$ with $L_0(\tilde{f}) = 1$. For an arbitrary function $f$ which is not constant on $X$ we have $L_0(f) > 0$. Therefore we can set 
$$\tilde{f} := \frac{f}{L_0(f)} \in \Sigma_0, $$
such that $L_\theta(\tilde{f}) = L_\theta(f) / L_0(f)$ and the above estimate becomes
\begin{equation} \label{eq equivalence lip-norms}
 L_0(f) \left(1-\varepsilon\right) \leq L_\theta(f) \leq \left(1+\varepsilon\right) L_0(f)
\end{equation}
uniformly for all $f \in A$ which are not constant.\smallskip

Now for given points $x,y \in M$, \eqref{eq equivalence lip-norms} leads to the estimates
\begin{eqnarray*}
 \rho_\theta(x,y) &=& \sup \left\{|f(x)-f(y)|: L_\theta(f) \leq 1 \right\} \\
 && \begin{cases} \leq \sup \left\{|f(x)-f(y)|: (1-\varepsilon) L_0(f) \leq 1 \right\}\\
\geq \sup \left\{|f(x)-f(y)|: (1+\varepsilon) L_0(f) \leq 1 \right\}. \end{cases}
\end{eqnarray*}
Since for every constant $C := 1 \pm \varepsilon >0$ one has
\begin{eqnarray*}
 \sup \left\{|f(x)-f(y)|: C L_0(f) \leq 1 \right\} &=& \sup \left\{\frac{|f(x)-f(y)|}{C}: L_0(f) \leq 1 \right\}\\
&=& \frac{\rho_0(x,y)}{C},
\end{eqnarray*}
this shows \eqref{eq equivalence and convergence of metricies}, and therefore the proposition is proved. \eB
\end{prop} \medskip

We note a similar observation, which gives equivalence of the metrics in case one can estimate the difference of two Lip-norms against one of them.\medskip

\begin{prop} \label{prop equivalence condition for metricies}
 Let $A$ be a commutative order unit space and let $L_1$, $L_0$ be two Lip-norms on $A$. Assume that there is a $C<1$ such that we have an estimate
$$ \left| L_0(f) - L_1(f) \right| \leq C L_0(f)$$ 
for all $f \in A$. Then the metrics $\rho_0$ and $\rho_1$, arising from $L_0$ and $L_1$, are equivalent.\smallskip

\B On the one hand, we have the estimate
$$  L_1(f) \geq L_0(f) - \left| L_0(f) - L_1(f) \right| \geq (1 - C) L_0(f), $$
and on the other hand we have
$$ L_1(f) \leq L_0(f) + \left| L_1(f) - L_0(f) \right| \leq (1+C) L_0(f). $$
Now the equivalence of the metrics follows by an argument analogous to the argument given in the proof of Proposition \ref{prop uniformly continuity implies convergence of metrices}, using $C$ instead of $\varepsilon$ in \eqref{eq equivalence lip-norms}.\eB 
\end{prop}\medskip

Now one can keep on playing the game and show that the statement of Proposition \ref{prop uniformly continuity implies convergence of metrices} (and therefore the uniform continuity of the family of Lip-norms) implies that the sequence of compact quantum metric spaces $(C(X),L_\theta)$ converges to the compact quantum metric space $(C(X),L_0)$ in the so-called \emph{quantum Gromov-Hausdorff convergence}. The concept of quantum Gromov-Hausdorff convergence was introduced by Mark Rieffel (see \cite{Rie}) as an quantum analogy to the Gromov-Hausdorff convergence of classical metric spaces. To prove this, one can adopt the arguments given by Rieffel in \cite{Rie}, Section 11, where this convergence is proved for a field of Lip-norms on finite-dimensional vector spaces.\smallskip

In our situation, we have to assume the uniform continuity condition from Definition \ref{def uniformly continuous families of cqms}, which can be shown to be fulfilled in the finite dimensional case. But since this does not affect the questions we are considering in this thesis, we will not write down the argumentation here.

\chapter{A Review of Sub-Riemannian Geometry}

This chapter is devoted to the introduction of the spaces for which we want to construct spectral triples. On some Riemannian manifolds $M$ it is possible to establish another metric besides the geodesic one. One can define the so-called Carnot-Carath\'{e}odory distance, determined by shortest paths tangent to a certain sub-bundle of the tangent bundle. Whenever this distance gives a metric on $M$ (meaning that any two points of $M$ can be connected by such a path), the Hausdorff dimension of this metric space differs from the usual dimension of $M$.\smallskip

In the first section we will introduce these sub-Riemannian spaces and review the properties mentioned above. Then, we will turn to the example of Carnot groups, which are nilpotent Lie groups carrying a certain grading inside their Lie algebras, and define Carnot manifolds. These will be the central objects to construct spectral triples on in this thesis; and they are fundamental within sub-Riemannian geometry since the tangent space of Riemannian geometry is generalized to a Carnot group in the sub-Riemannian case. The last two sections of this chapter are meant to provide some techniques we will take advantage of later: Section 2.3 deals with the case of Heisenberg manifolds, which have a horizontal distribution of co-dimension $1$ and occur at different points in mathematics and physics, and introduces the tool of the Levi form to describe the Lie group structure. Finally, in Section 2.4, we will construct submersions between Carnot groups which will allow us to reduce certain problems for general Carnot groups to lower dimensional cases.\smallskip

Throughout this section, $(M,g)$ will denote a Riemannian manifold with tangent bundle $TM$. Since most parts of this chapter are intended to sum up things which are already known, we will refer to the basic literature, which includes in our case the books and monographs by Richard Montgomery \cite{Mon}, by Luca Capogna, Donatella Danielli, Scott D. Pauls and Jeremy T. Tyson \cite{CDPT} and by Mikhael Gromov \cite{Gro}. In addition, we refer to books by Ovidiu Calin and Der-Chen Chang \cite{CC} and by A. Bonfiglioli, E. Lanconelli, F. Uguzzoni \cite{BLU}.\bigskip

\section{Sub-Riemannian Manifolds}

We start this chapter by giving the definition of a horizontal distribution, on which a sub-Riemannian geometry is based on.\medskip

\begin{defin} \label{def horizontal distribution}
Let $M$ be a Riemannian manifold with a Riemannian metric $g \in \Gamma^\infty(T^\ast M \otimes T^\ast M)$.
 \begin{enumerate}[(i)]

  \item A \emph{horizontal distribution} of the tangent bundle $TM$ is given by a sub-bundle $HM \subset TM$ together with a fiber wise inner product $\langle \cdot,\cdot \rangle_H$, such that at each point $x_0 \in M$ we have
$$\left. g(X,Y) \right|_{x_0} = \left\langle X(x_0), Y(x_0) \right\rangle_H$$
  for all $X,Y \in HM$.\smallskip 

  A vector field $X \in TM$ is called \emph{horizontal} if $X \in HM$.

  \item A (smooth) path $\gamma: [0,1] \rightarrow M$ is called \emph{horizontal} if at any point $t \in [0,1]$ we have $\dot{\gamma}(t) \in H_{\gamma(t)}M$.

  \item If, for $d \in \N$, $HM = \spa \{X_1, \ldots, X_d\}$ and $X_1, \ldots, X_d$ are linearly independent at each point $x \in M$, the integer $d$ is called the \emph{Rank} of $HM$. In this case we call $\{X_1, \ldots, X_d\}$ a \emph{frame} for the horizontal distribution $HM$.

  \item Let $\{X_1, \ldots, X_d\}$ be a frame for $HM$. Then for $f \in C^1(M)$, the vector field
$$\grad^H f := \sum_{j=1}^d X_j(f) X_j$$
  is called the \emph{horizontal gradient} of $f$.

  \item A horizontal distribution $HM$ is called \emph{involutive} if $[X,Y] \in HM$ for any $X,Y \in HM$.

  \item A horizontal distribution $HM$ is called \emph{bracket-generating} if the Lie hull of $HM$, which is the collection of all vector fields of $HM$ and their (multi-step) commutators, generates the tangent bundle $TM$. The smallest number $R \in \N$, such that $HM$ together with all its $R$-step commutators generates $TM$, is called the \emph{step} of a bracket generating distribution.
 \end{enumerate}
\end{defin}\smallskip

\Bem A horizontal distribution can also be defined by a set of $1$-forms: If $\{X_1, \ldots, X_n\}$ is a frame of $TM$ such that for $d<n$ $\{X_1, \ldots, X_d\}$ spans $HM$, we denote by $\{d\omega^1, \ldots, d\omega^n\}$ the corresponding dual frame of $T^\ast M$. Then we have 
$$HM = \ker \left( \spa \{d\omega^{d+1}, \ldots, d\omega^n\}\right).$$
\eBsp \medskip 

The Frobenius theorem (see for example \cite{CC}, Theorem 1.3.1) asserts that $HM$ is involutive if and only if it is integrable, which means that the set of all horizontal paths through a fixed point $x \in M$ sweeps out a smooth immersed sub-manifold $N$ of $M$ with $\dim N = \rank HM$. In this thesis we are interested in non-involutive, but bracket generating horizontal distributions: One can define a distance between two points $x,y \in M$ by considering the length of the shortest horizontal part connecting these points.\medskip

\begin{defin} \label{def carnot-caratheodory distance}
Let $(M,g)$ be a Riemannian manifold which is equipped with a horizontal distribution $HM = \spa \{X_1, \ldots, X_d\}$.
\begin{enumerate}[(i)]
 \item The \emph{horizontal length} of a smooth horizontal path $\gamma: [0,1] \rightarrow M$ is given by the number
$$L_{CC}(\gamma) := \int_0^1 \left( \sum_{j=1}^d g\left\langle \gamma'(t), X_{j}(\gamma(t)) \right\rangle_H \right)^\frac{1}{2} dt.$$
Here, $\langle \cdot,\cdot \rangle_H$ denotes the fiber wise inner product of $HM$ induced by the Riemannian metric of $M$, see Definition \ref{def horizontal distribution}.
 \item The \emph{Carnot-Carath\'{e}odory distance} between to points $x,y \in M$ is given by the (not necessarily finite) number
$$d_{CC}(x,y) := \inf \left\{ L_{CC}(\gamma): \gamma \text{ horizontal path with } \gamma(0) = x \text{ and } \gamma(1) = y \right\}.$$
 \item If $d_{CC}(x,y)$ is finite for any $x,y \in M$, we call $(M, d_{CC})$ a \emph{sub-Riemannian manifold}.
 \end{enumerate}
\end{defin}\medskip

In general, this distance needs of course not to be finite, but a famous theorem by Chow says that this distance is indeed finite for arbitrary $x,y \in M$ if $HM$ is bracket-generating (see also \cite{Mon}, Theorem 1.17):\medskip

\begin{thm} \label{thm Chows theorem}
 If $HM$ is a bracket-generating horizontal distribution on a connected manifold $M$, then any two points $x,y \in M$ can be joined by a horizontal path. \eB
\end{thm}\medskip

In particular this means that if the horizontal distribution $HM$ is bracket-generating, then $(M,d_{CC})$ is a metric space. This leads us to the question whether one can compare the metric spaces $(M,d_{CC})$ and $(M,d_{geo})$, where $d_{geo}$ denotes the usual geodesic distance. First of all, it is obvious that for any $x,y \in M$ we have $d_{geo}(x,y) \leq d_{CC}(x,y)$, but in general the metrics are not equivalent. But they do indeed induce the same topologies on $M$, as the following theorem shows (see \cite{Mon}, Theorem 2.3).\medskip

\begin{thm} \label{thm theorem on topologies}
 If $HM$ is a bracket-generating horizontal distribution on $M$, then the topology on $M$ induced by the Carnot-Carath\'{e}odory distance $d_{CC}$ is the same as the usual manifold topology induced by $d_{geo}$. \eB
\end{thm}

An unexpected observation is that the metric spaces $(M,d_{CC})$ and $(M,d_{geo})$ have a different Hausdorff dimension (see Definition \ref{def Hausdorff dimension} for the construction of this measure theoretic dimension). While the Hausdorff dimension of $(M,d_{geo})$ coincides with the (usual) topological dimension $n$ of $M$, the Hausdorff dimension of $(M,d_{CC})$ is in general strictly greater than $n$: We will see that the Hausdorff dimension of $(M,d_{CC})$ is exactly the so-called \emph{graded dimension} of the horizontal distribution.\medskip

To define the graded dimension of a sub-Riemannian manifold, we need a little bit of preparation. Since the horizontal distribution $HM$ is bracket generating of step $R$, we have a sequence of vector bundles
$$HM \subset H^2M \subset \ldots \subset H^RM = TM,$$
where
$$H^{S+1}M := H^SM + [HM,H^SM]$$
with $[HM,H^SM] := \spa \left\{[X,Y]: X \in HM, Y \in H^SM \right\}$ for $1 \leq S \leq R-1$ (using the convention that $H^1M = HM$). Now we assume that for a given $1 \leq S \leq R$ the dimension of the space $H^S_xM$ is the same for every point $x \in M$ (this will be an assumption throughout this thesis). In this situation, we are able to make the following definition.\medskip

\begin{defin} \label{def graded dimension}
 Let $M$ be a sub-Riemannian manifold with horizontal distribution $HM$ of rank $d$ and step $R$ as above. For $1 \leq S \leq R-1$ we denote by $d_1 := d := \rank HM$ and by
$$d_{S+1} := \rank H^{S+1}M - \rank H^{S}M$$
 the ranks of the spaces $H^{S+1}M \left/ H^SM \right.$. Then the \emph{graded dimension} (or \emph{homogeneous dimension}) of $M$ is the number
$$\dim_G(M) := \sum_{S=1}^R S \cdot d_S.$$\smallskip

\Bem Note that the topological dimension of $M$ is given by the number $\sum_{S=1}^R d_S$ in this context. \eBsp
\end{defin}\medskip

Now we are ready to formulate the so-called \emph{Mitchell`s Measure Theorem}, from which it follows that the Hausdorff dimension of $(M,d_{CC})$ differs from its topological dimension (see also \cite{Mon}, Theorem 2.17).\medskip

\begin{thm} \label{thm Mitchells measure theorem}
 The Hausdorff dimension of a sub-Riemannian manifold $M$ is equal to its graded dimension, i.e. under the notations of Definition \ref{def graded dimension} we have
$$\dim_H(M) = \sum_{S=1}^R S \cdot d_S.$$
This means in particular, that in general the Hausdorff dimension of the metric space $(M,d_{CC})$ is strictly greater than the Hausdorff dimension of the metric space $(M,d_{geo})$. \eB 
\end{thm}\medskip

We now finish this basic section by giving the most important example for a sub-Riemannian manifold: The $(2m+1)$-dimensional Heisenberg group. In dimension $3$ this is the easiest example where we have one space on which two different (with respect to metric and dimension) geometries can be established, and it will also be an example for a Carnot group, which we will introduce in the next section.\medskip

\begin{example} \label{ex Heisenberg group} \normalfont
For $m \in \N$, we consider the space $\R^{2m+1}$ equipped with a Riemannian metric $g$ such that the vector fields
$$X_j = \partial_{x_j} - \frac{1}{2}x_{m+j} \partial_{x_{2m+1}}, \ \ \ X_{m+1} = \partial_{x_{m+j}} + \frac{1}{2}x_j\partial_{x_{2m+1}}, \ \ \ X_{2m+1} = \partial_{x_{2m+1}}$$
for $1 \leq j \leq m$ form an orthonormal frame of $T\R^{2m+1}$. Note that $[X_j, X_{m+j}] = X_{2m+1}$, such that $H \R^{2m+1} = \spa \{X_1, \ldots, X_{2m}\}$ forms a non-involutive, $2$-step bracket generating horizontal sub-bundle of $T\R^{2m+1}$. Therefore $(\R^{2m+1},d_{CC})$ is a metric space of Hausdorff dimension $2m+2 = 2m \cdot 1 + 1 \cdot 2$.\smallskip

Since at each point, the tangent space generated by $\{X_1, \ldots, X_{2m+1} \}$ has the structure of a $2$-step nilpotent Lie algebra $\mathfrak{h}_{2m+1}$ the exponential mapping furnishes a simply-connected nilpotent Lie group $\Hei^{2m+1} = \exp \mathfrak{h}_{2m+1}$. $\Hei^{2m+1}$ is called the \emph{$2m+1$-dimensional Heisenberg group}. It is a well known fact that $\Hei^{2m+1}$ can be realized as the space of upper triangle $((m+2) \times (m+2))$-matrices which have the form
 $$\Hei^{2m+1} = \left\{ \begin{pmatrix} 1&x^t&z \\ 0&1_{(m \times m)}&y \\ 0&0&1 \end{pmatrix}: x,y \in \R^m; z \in \R \right\},$$
where the group rule is given by matrix multiplication. On $\R^{2m+1}$, this multiplication can be written as
\begin{equation} \label{eq composition rule Heisenberg group polarized coordinates}
(x,y,z) \cdot (\tilde{x},\tilde{y},\tilde{z}) = (x + \tilde{x}, y + \tilde{y}, z + \tilde{z} + \sum_{j=1}^m x_j \tilde{y}_j).
\end{equation}
The Lie algebra $\mathfrak{h}_{2m+1}$ of $\Hei^{2m+1}$ can be realized by strictly upper triangle matrices
$$\mathfrak{h}_{2m+1} = \left\{ \begin{pmatrix} 0&x^t&z \\ 0&0_{(m \times m)}&y \\ 0&0&0 \end{pmatrix}: x,y \in \R^m; z \in \R \right\},$$
where the algebra multiplication is once again the matrix multiplication. \eBsp \smallskip

\Bem The composition rule \eqref{eq composition rule Heisenberg group polarized coordinates} of the Heisenberg group is written down in the so-called \emph{polarized coordinates}. These are exactly the coordinates derived from the matrix model of this Lie group. On the other hand, one can describe points of $\Hei^{2m+1}$ by their \emph{exponential coordinates}, which are the coordinates of the corresponding Lie algebra having a one-to-one correspondence with the coordinates of $\Hei^{2m+1}$ by the exponential mapping. We will introduce both kinds of coordinates in the next section in the context of Carnot groups. \eBsp
\end{example}\bigskip

\section{Carnot Groups and Carnot Manifolds}

The Heisenberg group $\Hei^{2m+1}$ introduced in Example \ref{ex Heisenberg group} is an example for a greater class of nilpotent Lie groups, which are called Carnot groups.\medskip 

\begin{defin} \label{def Carnot group}
 A \emph{Carnot group} of step $R \in \N$ is a simply connected Lie group $\G$ whose Lie algebra $\mathfrak{g}$ has a stratification (or grading)
$$\mathfrak{g} = \bigoplus_{S=1}^R V_S,$$
such that $V_1, \ldots, V_R$ are vector spaces satisfying the conditions
\begin{enumerate}[(i)]
 \item $[V_1,V_S] = V_{S+1}$ for $S=1, \ldots, R-1$
 \item $[V_S,V_R] = 0$ for $S = 1, \ldots, R$.
\end{enumerate}
We call the number $d_1 := \dim V_1$ the \emph{bracket-generating dimension} of $\G$ and the number $\dim \G - d_1$ the \emph{bracket-generating co-dimension} of $\G$. \eBsp
\end{defin}\medskip

Carnot groups are the canonical generalization of the Euclidean space $\R^n$ in sub-Riemannian geometry, since the tangent space (or better the tangent cone) of a Carnot manifold has the structure of a Carnot group. This can be seen as a generalization of the Riemannian case, where the tangent space at any point is isomorphic to the ($1$-step nilpotent) Carnot group $\R^n$. Without going into detail, the situation is as follows: Let $X$ be any metric space, $x_0 \in X$. Then we define the \emph{tangent cone} as the pointed Gromov-Hausdorff limit of the family $(\lambda X,x_0)$ of pointed metric spaces for $\lambda \rightarrow \infty$, if it exists. It has been proved by John Mitchell \cite{Mit} that in the case where $(M,d_{CC})$ is a Carnot manifold, this limit exists and has the structure of a Carnot group (see also \cite{Mon}, Theorem 8.8):\medskip

\begin{thm} \label{thm Mitchells theorem on tangent cone}
 Let $x$ be a regular point of a Carnot manifold $(M,d_{CC})$. Then, at every $x \in M$, the tangent cone exists and is a Carnot group, which is arising from the nilpotentization of the horizontal distribution at $x$. \eB \smallskip 

\Bem We will not explain the details of the process of nilpotentization here and refer to \cite{Mon} or \cite{Bel} instead, since the technical details will not affect this thesis. \eBsp
\end{thm}\medskip

Now getting back to the objects we are considering, we finally define a Carnot manifold as a sub-Riemannian manifold, which has globally a Carnot group structure on its tangent bundle.\medskip

\begin{defin} \label{def Carnot manifold}
 Let $(M,g)$ be an $n$-dimensional Riemannian manifold, whose tangent bundle $TM$ carries a grading 
\begin{equation} \label{eq grading Carnot manifold}
 TM = \bigoplus_{S=1}^R V_SM
\end{equation}
such that $V_1M, \ldots, V_RM$ are vector bundles satisfying the conditions
\begin{enumerate}[(i)]
 \item $[V_1M,V_SM] = V_{S+1}M$ for $S=1, \ldots, R-1$
 \item $[V_RM,V_RM] = 0$ for $S = 1, \ldots, R$.
\end{enumerate}
Then we call $M$ a \emph{Carnot manifold} of step $R \in \N$. We further call the number $d_1 = \rank V_1M$ the \emph{bracket-generating dimension} of $M$ and the number $n - d_1$ the \emph{bracket-generating co-dimension} of $M$. \eBsp \smallskip

\Bem It is clear that every Carnot group is a Carnot manifold. It is also clear by the conditions (i) and (ii) on the brackets of the vector fields that every Carnot manifold contains the structure of a sub-Riemannian manifold, where the horizontal distribution is given by $HM = V_1M$. It should be mentioned that every Carnot manifold has locally the structure of a Carnot group, while for a general sub-Riemannian manifold this is a bit more involved (see Theorem \ref{thm Mitchells theorem on tangent cone}).\smallskip

Note that we have not fixed the geometry we are considering in the above definition. So by Section 2.1 we can establish at least two different geometries on a Carnot manifold $M$: The Riemannian one (equipped with the geodesic distance $d_{geo}$) and the sub-Riemannian one (equipped with the Carnot-Carath\'{e}odory distance $d_{CC}$). \eBsp
\end{defin}\medskip

When working on a Carnot manifold, we will always assume that we have a Riemannian metric $g$ on $M$ such that the spaces $V_SM$ appearing in the grading \eqref{eq grading Carnot manifold} are mutually orthogonal. For $S=1, \ldots, R$ we set $d_S := \rank V_SM$ and denote by $\left\{ X_{S,j}: j = 1, \ldots, d_S \right\}$ an orthonormal frame of $V_SM$. Hence, a frame of $TM$ is given by
$$ \left\{X_{1,1}, \ldots, X_{1,d_1}, X_{2,1}, \ldots, X_{2,d_2}, \ldots, X_{R,1}, \ldots, X_{R,d_R} \right\}.$$
We will sometimes use the abbreviation $X^{(S)} = (X_{S,1}, \ldots, X_{S,d_S})$ to denote the frame of $V_SM$.\smallskip

From the grading of its Lie algebra $\mathfrak{g}$, one can introduce coordinates on a Carnot group $\G$. Thereby we make use of the fact that because of the nilpotency the exponential mapping $\exp: \mathfrak{g} \rightarrow \G$ is a diffeomorphism. We present two different types of coordinates here, where it will be depending on the situation which type is more comfortable to work with.\medskip

\begin{defin} \label{def coordinates of Carnot groups}
Let $\G$ be a Carnot group with grading $\mathfrak{g} = V_1 \oplus \ldots \oplus V_R$ of its Lie algebra, such that for $1 \leq S \leq R$, $X^{(S)}$ denotes a basis of $V_S$.
\begin{enumerate}[(a)]
 \item The coordinates 
$$x = \left(x^{(1)}, \ldots, x^{(R)} \right) = \left(x_{1,1}, \ldots, x_{1,d_1}, \ldots, x_{R,1}, \ldots, x_{R,d_R}\right),$$
given by
\begin{equation} \label{eq exponential coordinates of a Carnot group}
 x \leftrightarrow \exp \left( \sum_{S=1}^R \sum_{j=1}^{d_S} x_{S,j}X_{S,j} \right),
\end{equation}
are called \emph{exponential coordinates} or \emph{canonical coordinates of the first kind} of $\G$.
\item The coordinates 
$$y = \left(y^{(1)}, \ldots, y^{(R)} \right) = \left(y_{1,1}, \ldots, y_{1,d_1}, \ldots, y_{R,1}, \ldots, y_{R,d_R}\right),$$
given by
\begin{equation} \label{eq polarized coordinates of a Carnot group}
 y \leftrightarrow \prod_{j=1}^{d_1} \exp \left( y_{1,j}X_{1,j} \right) . \prod_{j=1}^{d_2} \exp \left( y_{2,j}X_{2,j} \right) . \  \ldots \ . \prod_{j=1}^{d_R} \exp \left(y_{R,j}X_{R,j} \right)
\end{equation}
are called \emph{polarized coordinates} or \emph{canonical coordinates of the second kind} of $\G$.
\end{enumerate}
Thereby, $x.y$ denotes the group composition on $\G$, and the products are taken with respect to this group composition. \eBsp \smallskip

\Bem It is known that there is an isomorphism between the exponential and the polarized coordinates of a Carnot group. For example, in the case of the $(2m+1)$-dimensional Heisenberg group $\Hei^{2m+1}$ this isomorphism is given via
$$\phi(x_1, \ldots, x_{2m}, x_{2m+1}) = \left(x_1, \ldots, x_{2m}, x_{2m+1} + \frac{1}{2} \sum_{j=1}^m x_j x_{m+j} \right),$$
where $(x_1, \ldots, x_{2m}, x_{2m+1})$ denote the exponential coordinates (see \cite{Fol}, Section 1.2). \eBsp
\end{defin} \medskip

When calculating a group composition on $\G$, one uses the Baker-Campbell-Hausdorff formula. This formula is given on $\mathfrak{g}$ by
\begin{equation} \label{eq Baker-Campbell-Hausdorff formula}
 \exp X . \exp X = \exp \left(X+Y + B(X,Y) \right),
\end{equation}
where $B(X,Y)$ is a sum of multi-step commutators of order $2, 3, \ldots, R$, see \cite{Knapp}. Therefore, $B(X,Y)$ is a polynomial of degree smaller or equal to the step of $\G$, which does not depend on vectors belonging to $V_R$ since $V_R$ commutes with every $X \in \mathfrak{g}$. In the case of a $2$-step Carnot group, we have
$$\exp X . \exp Y = \exp \left(X+Y + \frac{1}{2}[X,Y]\right).$$
Using this expression on the Lie algebra, one can derive the composition rule in exponential or polarized coordinates on the Carnot group $\G$. \smallskip

Maybe the most important property of a Carnot group is its homogeneity, which is expressed by a weighted dilation of it. Once again there is a one-to-one relation between the dilation on a Carnot group and the corresponding weighted dilation on its graded Lie algebra, given by the exponential map. \medskip

\begin{defin} \label{def delation on Carnot group} \normalfont
\begin{enumerate}[(i)]
 \item  Let $\G$ be a Carnot group with coordinates $x = (x^{(1)}, \ldots, x^{(R)})$ (exponential or polarized), $\lambda > 0$. Then we define by
\begin{equation} \label{eq dilations Carnot group}
 \delta_\lambda: \G \rightarrow \G, \ \ \ x = \left(x^{(1)}, \ldots, x^{(R)}\right) \mapsto \lambda.x := \left(\lambda x^{(1)}, \lambda^2 x^{(2)}, \ldots, \lambda^R x^{(R)}\right)
\end{equation}
the \emph{(weighted) dilation} on $\G$ by $\lambda$.\smallskip

A function $f: \G \rightarrow \R$ is called \emph{homogeneous} of degree $\mu \in \R$ with respect to $\delta_\lambda$ if we have
$$f(\delta_\lambda(x)) = \lambda^\mu \cdot f(x)$$
for all $x \in \G$, $\lambda > 0$.

 \item Let $\mathfrak{g}$ be a graded nilpotent Lie algebra with grading $\mathfrak{g} \cong V_1 \oplus \ldots \oplus V_R$, $\lambda > 0$. Then we define by
\begin{equation} \label{eq dilations graded Lie algebra}
 \hat{\delta}_\lambda: \mathfrak{g} \rightarrow \mathfrak{g}, \ \ \ X = \sum_{S=1}^R \sum_{j=1}^{d_S} x_{S,j} X_{S,j} \mapsto \hat{\delta}_\lambda.X := \sum_{S=1}^R \sum_{j=1}^{d_S} \lambda^S x_{S,j} X_{S,j}
\end{equation}
the \emph{(weighted) dilation} on $\mathfrak{g}$ by $\lambda$. \smallskip

A function $f: \mathfrak{g} \rightarrow \R$ is called \emph{homogeneous} of degree $\mu \in \R$ with respect to $\hat{\delta}_\lambda$ of degree $\mu$ if we have
$$f(\hat{\delta}_\lambda(X)) = \lambda^\mu \cdot f(X)$$
for all $X \in \mathfrak{g}$, $\lambda > 0$.

\item Let $M$ be a Carnot manifold, $\delta_\lambda$ and $\hat{\delta}_\lambda$ the weighted dilations by $\lambda > 0$ defined point-wise on $TM$. A vector field $X \in \Gamma(TM)$ is called \emph{homogeneous} of degree $\mu \in \R$ if we have
\begin{equation} \label{eq homogenicity of vector fields}
 \hat{\delta}_\lambda X = \lambda^\mu \cdot X
\end{equation}
for every $\lambda > 0$.
\end{enumerate} \smallskip

\Bem Note that if the grading of $TM$ is given by $TM \cong V_1M \oplus \ldots \oplus V_RM$, we have that $X \in TM$ is homogeneous of degree $S$ if and only if $X \in V_SM$. \eBsp
\end{defin}\medskip

We remark that the weight $S$ of $\delta_\lambda$ (or $\hat{\delta}_\lambda$) belonging to certain coordinates reflects exactly the commutator step of the space $V_S$ of $\mathfrak{g}$ which belongs to $x^{(S)}$ (or $X^{(S)}$). These weighted dilations play an important role when one wants to describe a symbol calculus with respect to a Carnot group structure: Asymptotic expansions are given in terms of homogeneous functions with respect to weighted dilations. Furthermore one can use weighted dilations to define a Carnot group, see for example \cite{BLU}, Definition 2.2.1, where a (homogeneous) Carnot group is defined to be a Lie group structure on $\R^n$ which respects these dilations and has a bracket generating structure. In \cite{BLU} it is also shown that these two definitions are equivalent (up to isomorphisms).\smallskip

Let us mention now how every Carnot group can be realized as a certain Lie group structure on $\R^n$, using certain vector fields to represent the basis of its Lie algebra. We have already seen this in Example \ref{ex Heisenberg group} for the case of the Heisenberg group. The idea is to find a frame of vector fields $X_1, \ldots X_n$ which satisfies the grading conditions of Definition \ref{def Carnot group} and which has the property that $X_j(0) = \left. \partial_{x_j} \right|_{x=0}$ for each $j$. This is indeed possible, and the requested vector fields $X_j$ have polynomial coefficients, as the following proposition shows (see \cite{BLU}, Remark 1.4.6).\medskip

\begin{prop} \label{prop Carnot Lie algebra as vector fields}
 Let $\G$ be a Carnot group with Lie algebra $\mathfrak{g}$, where the grading of $\mathfrak{g}$ is given by $\mathfrak{g} = \bigoplus_{S=1}^R V_S$ with $d_S = \dim V_S$. Then for each $V_S$ there exists a frame of vector fields $\left\{X_{S,1}, \ldots, X_{S,d_S}\right\}$ such that 
\begin{equation} \label{eq graded vector fields}
 X_{S,j} = \partial_{x_{S,j}} + \sum_{L=S+1}^R \sum_{k=1}^{d_L} p_{j,k}^{(S,L)} \left(x^{(1)}, \ldots, x^{(L-S)} \right) \partial_{x_{L,k}}. 
\end{equation}
Here, $p_{j,k}^{(S,L)}$ is a polynomial which is homogeneous with respect to the dilations $\delta_\lambda$ from \eqref{eq dilations Carnot group} of degree $L-S$, which means
$$ p_{j,k}^{(S,L)} \left(\delta_\lambda \left(x^{(1)}, \ldots, x^{(L-S)} \right) \right) = \lambda^{L-S} p_{j,k}^{(S,L)} \left(x^{(1)}, \ldots, x^{(L-S)} \right).$$
\eB
\end{prop}\medskip

Now one can consider the homogeneous vector fields from \eqref{eq graded vector fields} as homogeneous differential operators. Thus any polynomial of these vector fields is a differential operator. These so-called \emph{graded differential operators} will play a big role in this thesis: The horizontal Dirac operators we will construct in Chapter $3$ fall into this category, and the horizontal Laplacians of (homogeneous) degree $2$ will also turn out to be very important.\medskip

\begin{defin} \label{def graded differential operator}
Let $\mathfrak{g} \cong V_1 \oplus \ldots \oplus V_R$ be a graded Lie algebra with $d_S = \rank V_S$, which is represented by vector fields on $\R^{d_1 + \ldots + d_R}$ as in Proposition \ref{prop Carnot Lie algebra as vector fields}. A frame for $V_S$ shall consist of the vector fields $\{X_{S,1}, \ldots, X_{S,d_S}\}$, where we write $\{X_1, \ldots, X_d\}$ for the frame of $V_1$ (with $d=d_1$).\smallskip

Then a \emph{graded differential operator} is a differential operator of the from
$$L = p\left(X_1, \ldots, X_d, X_{2,1}, \ldots, X_{R,d_R} \right),$$
where $p$ is a polynomial with matrix-valued $C^\infty$ coefficients. If $p$ is homogeneous of degree $\mu \in \R$ we call $L$ \emph{homogeneous} of degree $\mu$. \eBsp \smallskip

In particular, a \emph{horizontal Laplacian} is a graded differential operator of order $2$, which means it is an operator of the form
$$\Delta^{\mathrm{hor}} = - \sum_{j=1}^d X_j^2 + \sum_{j=1}^{d_2} b_{2,j} X_{2,j} + \sum_{j=1}^d b_{1,j} X_j + b_0,$$
where all the $b_{S,j}$ and $b_0$ are smooth (matrix valued) functions on $\C^n$. \eBsp
\end{defin}\medskip 

The above definition suggests to consider a graded differential operator $L$ as an element of the universal enveloping algebra $\mathcal{U}(\mathfrak{g})$ of the graded Lie algebra $\mathfrak{g}$ generated by $\{X_1, \ldots, X_d\}$. This interpretation will be used in Chapter $6$ when we analyze graded differential operators using the representation theory of their underlying Lie algebras $\mathfrak{g}$.\smallskip

We return once again to the dilations from Definition \ref{def delation on Carnot group}: For our purposes we should mention an additional another property of $\delta_\lambda$: For arbitrary points $x,y \in \G$ we have
$$d_{CC}(\delta_\lambda(x), \delta_\lambda(y)) = \lambda \cdot d_{CC}(x,y),$$
so the Carnot-Carath\'{e}odory metric on $\G$ respects these dilations. Now we will introduce a quasi-norm (or gauge norm) $\| \cdot \|_\G$ on a Carnot group $\G$, which is homogeneous with respect to these dilations and provides us with a quasi-metric on $\G$ which is equivalent to $d_{CC}$. It will rather be a quasi-norm than a norm since the triangle inequality on $\G$ must be replaced by the condition
$$\left\|x.y\right\|_\G \leq C \left\|x \right\|_\G \left\| y \right\|_\G.$$
This quasi-norm will play an important role when we define a symbol calculus according to Carnot groups. \medskip

\begin{defin} \label{def Koranyi gauge}
 Let $\G$ be a Carnot group of step $R$. The quasi-norm $\left\| \cdot \right\|_\G$ on $\G$, defined by
\begin{equation} \label{eq Koranyi gauge}
 \left\| x \right\|_\G^{2R!} := \sum_{S=1}^R \sum_{j=1}^{d_S} \left| x_{S,j} \right|^\frac{2R!}{S},
\end{equation}
is called the \emph{Koranyi gauge} of $\G$. The quasi-metric
$$d_\G(x,y) := \left\| y^{-1}.x \right\|_\G$$
arising from the Koranyi gauge will be called the \emph{Koranyi (quasi-)metric} on $\G$. \eBsp\smallskip

\Bem Note that in the case $\G = \Hei^{2m+1}$ the Koranyi gauge is given by the formula
\begin{equation} \label{eq Koranyi gauge on Heisenberg group}
 \left\| x \right\|_{\Hei^{2m+1}} = \left( \sum_{j=1}^{2m} \left| x_j \right|^4 + \left| x_{2m+1} \right|^2 \right)^\frac{1}{4}.
\end{equation}
We also note that there are several definitions of the Koranyi gauge on a Carnot group which are equivalent, for example in \cite{Stein} the Koranyi gauge on $\Hei^{2m+1}$ is defined via
$$\rho(x) := \max \left\{ \left\| (x_1, \ldots, x_{2m}) \right\|, \left| x_{2m+1} \right|^{1/2} \right\},$$
where $\| \cdot \|$ denotes the Euclidean norm of the corresponding vector. \eBsp
\end{defin}\medskip
 
As noted in \cite{CDPT}, the Koranyi quasi-metric on a Carnot group $\G$ is equivalent to the Carnot-Carath\'{e}odory metric on $\G$. \medskip

\begin{prop} \label{prop equivalence Koranyi metric Carnot-Caratheodory metric}
 Let $\G$ be a Carnot group, and let $d_\G$ denote the Koranyi quasi-metric and $d_{CC}$ the Carnot-Carath\'{e}odory metric on $\G$. Then there are constants $c > 0$ and $C > 0$ such that for all $x,y \in \G$ we have
$$c \cdot d_{CC}(x,y) \leq d_\G(x,y) \leq C \cdot d_{CC}(x,y).$$
\eB
\end{prop}\medskip

Now having finished a rough review about the theory on Carnot groups and Carnot manifolds, we want to finish this section by giving some simple examples for compact Carnot manifolds. We consider the local homogeneous space of a discrete lattice subgroup of a Carnot group $\G$. In a way, this is the analogous object to the torus $\T^n$, which arises as the quotient by $\Z^n$ on $\R^n$. Since it is quite comfortable to do calculations on these objects, they will serve as the main example for the considerations of this thesis.\medskip

\begin{example} \label{ex compact Carnot nilmanifolds} \normalfont
 Let $\G$ be a Carnot group, equipped with a (left-invariant) Riemannian metric, with grading $\mathfrak{g} = V_1 \oplus \ldots \oplus V_R$, where $\{X_1, \ldots, X_d\}$ is a frame for $V_1$. We consider the discrete subgroup $\Gamma$ of $\G$ generated by the basis vectors of $V_1$, that is
\begin{equation} \label{eq standard lattice}
 \Gamma := \left\langle \ \{\gamma_j = \exp(X_j): 1 \leq j \leq d \} \ \right\rangle_{\G}.
\end{equation}
$\Gamma$ acts on $\G$ via the left multiplication $(\gamma,x) \mapsto \gamma.x$ for all $x \in \G$. It follows easily that $\Gamma$ is a lattice in $\G$ (see \cite{Mon}, Section 9.3), which we will call the \emph{standard lattice} of $\G$.\smallskip

 Since $\Gamma$ is a lattice, the local homogeneous space $\Gamma \backslash \G$, consisting of the orbits of this group action, is a compact Riemannian manifold which is locally isometric to $\G$. Hence $M := \Gamma \backslash \G$ is a compact Carnot manifold, where the grading of the tangent bundle $TM$ comes from the grading of $\mathfrak{g}$. We will call $M := \Gamma \backslash \G$ the (standard) \emph{compact Carnot nilmanifold} of $\G$.\eBsp
\end{example}\medskip

\begin{remark}\normalfont
 The group $\Gamma$ defined via \eqref{eq standard lattice} can also be viewed as a discrete nilpotent group. From the grading of $\G$ one can detect a so-called \emph{central descending series} 
$$0 = \Gamma^{R+1} \subset \Gamma^R \subset \ldots \subset \Gamma^2 \subset \Gamma^1 = \Gamma,$$
which means we have $\Gamma^{S+1} = [\Gamma, \Gamma^S]$, where $\Gamma^S$ is the standard lattice of the Carnot group $\exp(V_S \oplus \ldots \oplus V_R)$. One can check that for the graded dimension (and therefore by Theorem \ref{thm Mitchells measure theorem} also for the Hausdorff dimension) of $\G$ we have the identity 
\begin{equation} \label{eq graded dimension from discrete groups}
 \dim_G(\G) = \sum_{S=1}^R S \cdot \rank \left(\Gamma^S / \Gamma^{S+1} \right).
\end{equation}
Now there is a famous theorem by Bass, Milnor and Wolf (see \cite{Mon}, Theorem 9.3) which states that the right hand side of \eqref{eq graded dimension from discrete groups} is equal to the polynomial growth of $\Gamma$. So, altogether, the Hausdorff dimension of $\G$ coincides with the polynomial growth of its standard lattice subgroup. See \cite{Mon}, Sections 9.2 and 9.3, for details on this. \eBsp
\end{remark}\bigskip

\section{Heisenberg Manifolds and Levi Forms}

We will now pay attention to the most important class of Carnot manifolds: It is the case where the graded co-dimension is equal to $1$. Those manifolds are also known as Heisenberg manifolds. In the context of non-commutative geometry, they were treated in details in the work of Rapha\"{e}l Ponge (see for example \cite{Pon1} and the references there), and recently they have also been a big research area in index theory, for which we refer to Erik van Erp and his collaborators (see e.g. \cite{BE}). For most of the following definitions and proposition, we refer to Ponge (\cite{Pon1}).\medskip

\begin{defin} \label{def Heisenberg manifold}
 A Heisenberg manifold is a smooth Riemannian manifold $(M,g)$ of dimension $n=d+1$ equipped with a hyperplane bundle $HM \subset TM$ of rank $d$, which is bracket-generating. \eBsp \smallskip

\Bem One can formulate the definition of a Heisenberg manifold a little bit more general by dropping the assumption that $HM$ has to be bracket generating, see e.g. \cite{Pon1}. In this case, objects like contact manifolds or CR manifolds are included in the class of Heisenberg manifolds. But since we are following more or less theoretical aspects and want to consider manifolds on which we have the Carnot-Carath\'{e}odory geometry, we make this further assumption. \eBsp
\end{defin}\medskip

An important tool to handle Heisenberg manifolds is the so called \emph{Levi form}. In case $M$ is a Carnot manifold of graded co-dimension $1$, this is the $1$-form describing the Lie algebra structure of the (graded) tangent bundle. \medskip

\begin{defin} \label{defin Levi form}
 Let $M$ be a Carnot manifold of co-dimension $1$, where the grading of its tangent bundle is given by $TM = V_1M \oplus V_2M$, with $\rank V_2M = 1$ and $V_2M = [V_1M,V_1M]$. Let $\{X_1, \ldots, X_d\}$ be a frame of $V_1M$ and let $\{X_{d+1}\}$ be a frame of $V_2M$. Then the \emph{Levi form} $\mathcal{L}$ of $M$ is the (antisymmetric) bilinear form 
\begin{equation} \label{eq Levi form}
\mathcal{L}: V_1M \times V_1M \rightarrow V_2M, \ \ \ (Y_1,Y_2) \mapsto [Y_1,Y_2] \mod V_1M.
\end{equation}
For $\mathcal{L}(X_j,X_k) = L_{ik} X_{d+1}$ with $L_{ik} \in \R$ for $j,k \in \{1, \ldots d\}$, we denote by $L = \left(L_{jk}\right)$ the antisymmetric $(d \times d)$-matrix describing $\mathcal{L}$ according to the basis $\{X_1, \ldots, X_d\}$.
\end{defin}\medskip

Now one can define a Levi form on any Heisenberg manifold $M$ the same way: On can see rather easily that from the Lie bracket of vector fields on $HM$ one gets a $2$-form
$$ \mathcal{L}: HM \times HM \rightarrow TM / HM ,$$
such that for any section $X,Y \in HM$ we have 
$$\mathcal{L}_{x_0}\left(X(x_0), Y(x_0)\right) = [X,Y](x_0) \mod H_{x_0}M$$
near a point $x_0 \in M$ (see \cite{Pon1}). But this allows us to define a bundle of graded $2$-step nilpotent Lie algebras $\mathfrak{g}M \cong HM \oplus (TM/HM)$, which gives rise to a bundle of Carnot groups $\G M$ of step $2$ and graded co-dimension $1$ over $M$ via the exponential mapping. This bundle $\G M$ is also called the \emph{tangent Lie group bundle} of $M$.\medskip

We can even say how this tangent Lie group bundle looks like:\medskip

\begin{prop} \label{prop tangent Lie group bundle and rank of Levi form}
 Let $(M,g)$ be a Heisenberg manifold. Then for any point $x_0 \in M$, the point-wise Levi form $\mathcal{L}_{x_0}$ has rank $2m$ if and only if $\G_{x_0}M \cong \Hei^{2m+1} \times \R^{d-2m}$. In especially, this means $\mathcal{L}$ has constant rank $2m$ if and only if $\G M$ is a fiber bundle with typical fiber $\Hei^{2m+1} \times \R^{d-2m}$.\smallskip

In addition, if $\rank \mathcal{L} = 2m$, it is always possible to find an orthonormal basis $\{X_1, \ldots, X_d\}$ of $V_1M$ such that the matrix representation $L = (L_{jk})$ of $\mathcal{L}$ becomes
\begin{equation} \label{eq Levi form matrix normal form}
 L = \begin{pmatrix} 0&D&0 \\ -D&0&0 \\ 0&0&0 \end{pmatrix},
\end{equation}
where $D \in \mathrm{Mat}_{m\times m}(\R)$ is a diagonal matrix carrying the absolute values $\lambda_1, \ldots, \lambda_m$ of the non-zero eigenvalues on its diagonal. \smallskip

\B The first statement is exactly the statement of \cite{Pon1}, Proposition 2.1.6. The second statement follows by linear algebra, since $L$ is a skew-symmetric matrix, so it is known to have the non-zero eigenvalues $\pm i \lambda_1, \ldots, \pm i \lambda_m$. Then for any orthonormal frame of $HM$ the form \eqref{eq Levi form matrix normal form} can be reached by an orthonormal basis transformation at every point of $HM$. \eB \smallskip

\Bem By our assumption of $HM$ being bracket-generating, we always have $\rank \mathcal{L} \geq 2$ and therefore $m \geq 1$. Of course the case $m=0$ also fits into the theory: In this case the hyperplane bundle $HM$ induces a foliation on $M$. \eBsp
\end{prop}\medskip

We can further introduce Levi forms on arbitrary Carnot manifolds: Since $[V_1M, V_1M] = V_2M$ for the first two summands appearing in the grading of the tangent space of any Carnot manifold, the following definition is well-defined.\medskip

\begin{defin}\label{def generalized Levi form}
 Let $M$ be a Carnot manifold with grading $TM = \bigoplus_{S=1}^R V_SM$ of its tangent bundle. We denote by $\{X_{1,1}, \ldots, X_{1,d}\}$ an orthonormal frame of $V_1M$ and by $\{X_{2,1}, \ldots, X_{2,d_2}\}$ an orthonormal frame of $V_2M$. Then for $\nu \in \{1, \ldots, d_2\}$, the \emph{$\nu$-Levi form} of $M$ is given by the bilinear form
$$\mathcal{L}_\nu: V_1 \times V_1 \rightarrow \spa\{X_{2,\nu}\}, \ \ \ (Y_1,Y_2) \mapsto [Y_1,Y_2] \mod \left( \spa \{X_{2,\nu}\} \right)^\perp.$$
For $\mathcal{L}_\nu(X_j,X_k) = L_{ik}^{(\nu)} X_\nu$ with $L_{ik}^{(\nu)} \in \R$, we denote by $L^{(\nu)} = \left(L_{jk}\right)$ the antisymmetric matrix describing $\mathcal{L}_\nu$. \eBsp
\end{defin}\medskip

In other words, the collection of $\nu$-Levi forms $\mathcal{L}_\nu$ describes the structure of the first step commutators of a Carnot manifold. We will use this notation later when we look for structures of Heisenberg manifolds inside a Carnot manifold.\bigskip

\section{Submersions of Carnot Groups}

In the final section of this chapter we will introduce a technique which will play an important role in our later considerations: We will show how one gets a submersion from a given Carnot group $\G_1$ onto a lower dimensional Carnot group $\G_2$. This will provide us the possibility to pull back objects defined on $\G_2$ to objects on $\G_1$, such that statements about $\G_2$ can be transported to $\G_1$.\smallskip

Once again we consider a Carnot group $\G$ whose Lie algebra has the grading $\mathfrak{g} = \bigoplus_{S=1}^R V_S$. Let $\tilde{V} \subset \mathfrak{g}$ be a linear subspace of (the vector space) $\mathfrak{g}$ which has the structure
\begin{equation} \label{eq projection Lie algebra subspace}
 \tilde{V} := \bigoplus_{S=1}^{M-1} V_S \oplus \tilde{V}_{M}, \ \ \ \text{where } \tilde{V}_M \subset V_M \ \text{is a linear subspace,}
\end{equation}
for some $1 \leq M \leq R$. We take a look at the canonical orthogonal projection
\begin{equation} \label{eq projection Lie algebra}
 \mathrm{pr}: \mathfrak{g} \rightarrow \tilde{V}, \ \ \ v \mapsto v \mod \tilde{V}^\perp. 
\end{equation}
We show now that $\mathrm{pr}$ induces a homomorphism of Lie algebras, which gives rise to a submersion (and also to a homomorphism) of Lie groups.\medskip

\begin{prop} \label{prop properties projection lie algebra}
 Let $\G$, $\mathfrak{g}$, $\tilde{V}$ and $\pr$ be as above. Then the vector space $\tilde{\mathfrak{g}} := \pr(\mathfrak{g}) \cong \tilde{V}$ can be equipped with a Lie algebra structure via
\begin{equation} \label{eq Lie brackets projected}
 [X,Y]_\mathrm{pr} := \mathrm{pr} \left([X,Y]\right) \ \ \ \forall X,Y \in \tilde{V}
\end{equation}
such that $\pr: (\mathfrak{g}, [\cdot,\cdot]) \rightarrow (\mathfrak{\tilde{g}}, [\cdot,\cdot]_\mathrm{pr})$ is an homomorphism of graded Lie algebras.\smallskip

Further, if we set $\mathfrak{n} := \ker(\mathrm{pr})$, the spaces $N := \exp(\mathfrak{n})$ and $\tilde{\G} := \exp(\tilde{\mathfrak{g}})$ are Carnot groups such that $\tilde{\G} \cong \G \left/ N \right.$ and hence $\G \cong \tilde{\G} \times N$. The resulting map
\begin{equation} \label{eq projection Lie groups}
 \psi := \exp_{\tilde{\G}} \circ \ \mathrm{pr} \ \circ \exp_{\G^{-1}}: \G \rightarrow \tilde{\G}
\end{equation}
is a homomorphism of Carnot groups and a submersion of Riemannian manifolds. \smallskip

\B First of all, one checks that the Lie brackets defined via \eqref{eq Lie brackets projected} are indeed Lie brackets: The bilinearity of $[\cdot,\cdot]_\mathrm{pr}$ is obvious because of the bilinearity of $[\cdot,\cdot]$ and the linearity of $\mathrm{pr}$; and the Jacobi identity also follows from these properties in connection with the Jacobi identity of $[\cdot,\cdot]$.\smallskip

To check that $\mathrm{pr}$ is a homomorphism of Lie algebras, one has to check that
\begin{equation} \label{eq condition Lie algebra homomorphism}
 [\mathrm{pr}(X),\mathrm{pr}(Y)]_\mathrm{pr} = \mathrm{pr} \left([X,Y]\right) \ \ \ \forall \ X,Y \in \mathfrak{g}.
\end{equation}
This is by definition true for all $X,Y \in \tilde{V}$. Now let (without loss of generality) $X \in \tilde{V}^\perp$. In this case, $\mathrm{pr}(X) = 0$ and therefore the left hand side of \eqref{eq condition Lie algebra homomorphism} is zero. But the right hand side is also zero, since $X \in \tilde{V}^\perp \subset \bigoplus_{S=M}^R V_S$ and therefore $[X,Y] \in \bigoplus_{S=M+1}^R V_S$ because of the graded structure of $\mathfrak{g}$. Therefore \eqref{eq condition Lie algebra homomorphism} is true, which shows that the linear map $\mathrm{pr}$ is a Lie algebra homomorphism. But this means that $\mathfrak{n} := \ker (\pr)$ is an ideal in $\mathfrak{g}$ and therefore also a graded Lie algebra. The grading structures of $\tilde{g}$ follows immediately from the grading structure of $\mathfrak{g}$ since 
$$\tilde{\mathfrak{g}} \cong \bigoplus_{S=1}^{M-1} V_S \oplus \tilde{V}_M.$$
Thus $\tilde{\mathfrak{g}}$ is obviously a graded Lie algebra of step $M$. \smallskip

Now we consider the map $\psi$ from \eqref{eq projection Lie groups}. Because of the nilpotency the exponential maps from the Lie algebras $\tilde{\mathfrak{g}}$ and $\mathfrak{g}$ to their Lie groups $\tilde{\G}$ and $\G$ are diffeomorphisms, and since $\mathrm{pr}$ is a surjective linear map (which is smooth because of its linearity) we have the result that $\psi$ is a submersion. To check that it is a group homomorphism, we calculate, using the Baker-Campbell-Hausdorff formula \eqref{eq Baker-Campbell-Hausdorff formula}:
\begin{eqnarray*}
 \psi(x._\G y) &=& \exp_{\tilde{\G}} \circ \pr \left(\exp_\G^{-1} (\exp_\G\left(\sum_{S=1}^R \sum_{j=1}^{d_S} x_{S,j}X_{S,j} \right) ._\G \exp_\G\left( \sum_{S=1}^R \sum_{j=1}^{d_S} y_{S,j}X_{S,j} \right) ) \right)\\
&=& \exp_{\tilde{\G}} \circ \pr \left( \sum_{S,j} (x_{S,j} + y_{S,j}) X_{S,j} + B(\sum_{S,j} x_{S,j}X_{S,j},\sum_{S,j} y_{S,j}X_{S,j}) \right)\\
&=& \exp_{\tilde{\G}} \left(\sum_{S,j} x_{S,j}\pr(X_{S,j}) + \sum_{S,j} y_{S,j}\pr(X_{S,j}) + \tilde{B}(\sum_{S,j} x_{S,j} \pr(X_{S,j}),\sum_{S,j} y_{S,j} \pr(X_{S,j})) \right)\\
&=& \exp_{\tilde{\G}}\left(\sum_{S,j} x_{S,j}\pr(X_{S,j})\right) ._{\tilde{\G}} \exp_{\tilde{\G}}\left(\sum_{S,j} y_{S,j}\pr(X_{S,j})\right)\\
&=& \psi(x) ._{\tilde{\G}} \psi(y).
\end{eqnarray*}
In this calculation we have used the fact that $\pr$ is a Lie algebra homomorphism. In particular, the third equation is true since $B$ is a sum of (multi-step) commutators of the vector fields $X_{S,j}$, such that applying $\mathrm{pr}$ to $B$ furnishes the polynomial $\tilde{B}$ in the Baker-Campbell-Hausdorff formula on $\tilde{\G}$.\smallskip

Finally we have 
$$N = \exp \mathfrak{n} = \ker \psi,$$
such that $N$ is a normal Lie subgroup of $\G$. But this shows in addition that $\tilde{\G} \cong \G / N$ and that $\G \cong \tilde{\G} \times N$. Now every statement of the proposition is proved. \eB
\end{prop} \medskip

We now show briefly that such a submersion can be lifted to the compact nilmanifolds arising from $\G$ from Example \ref{ex compact Carnot nilmanifolds}. Let $\Gamma$ be generated by the images of the basis vector fields of $V_1$ under the exponential mapping, and let $M = \Gamma \backslash \G$ be the local homogeneous space of the left action of $\Gamma$ on $\G$. Since the submersion $\psi$ from \eqref{eq projection Lie groups} is a Lie group homomorphism, the image $\psi(\Gamma)$ under $\psi$ is a discrete subgroup of $\psi(\G)$, and by the definition of $\psi$ it is clear that $\psi(\Gamma)$ is generated by the image of the basis vectors of $V_1$ under the corresponding projection $\pr: \exp_\G^{-1} (\G) \rightarrow \exp_{\psi(\G)}^{-1} (\psi(\G))$ of the Lie algebras. Therefore, we get a compact Carnot nilmanifold $\tilde{M} := \psi(\Gamma) \backslash \psi(\G)$.\smallskip

For the action of $\psi(\Gamma)$ on $\psi(\G)$ we have
$$\psi(\gamma.x) = \psi(\gamma).\psi(x)$$
for any $\gamma \in \Gamma$ and $x \in \G$, therefore any orbit of the action of $\Gamma$ on $\G$ is mapped to a orbit of the action of $\psi(\Gamma)$ on $\psi(\G)$. Let $[x]_\Gamma$ denote the orbit belonging to an element $x \in \G$ under the action of $\Gamma$. Because $\psi$ is a submersion by Proposition \ref{prop properties projection lie algebra}, this means that the mapping
\begin{equation} \label{eq projection Carnot nilmanifolds}
 \pi: M = \Gamma \backslash \G \rightarrow \tilde{M} = \psi(\Gamma) \backslash \psi(\G), \ \ \ [x]_\Gamma \mapsto [\psi(x)]_{\psi(\Gamma)}
\end{equation}
is also a submersion. We summarize the above argumentations in the following corollary:\medskip

\begin{cor} \label{cor projection Carnot nilmanifolds}
 In the above situation, the map $\pi$ from \eqref{eq projection Carnot nilmanifolds} is a submersion of Riemannian manifolds. Locally, $\pi$ coincides with the submersion $\psi: \G \rightarrow \psi(\G)$ from \eqref{eq projection Lie groups} of the corresponding Carnot groups.\smallskip

\B The fact that $\pi$ is a submersion has already been deduced in the discussion ahead of this corollary. Since the nilmanifolds $\Gamma \backslash \G$ and $\psi(\Gamma) \backslash \psi(\G)$ are locally isometric to the Carnot groups $\G$ and $\tilde{\G}$ (see Example \ref{ex compact Carnot nilmanifolds}), the second statement is obvious by the construction of $\pi$ from $\psi$. \eB
\end{cor}\bigskip

\chapter{Horizontal Dirac-Operators}

The intention of this chapter is to find a first order differential operator which detects the Carnot-Carath\'{e}odory metric via Connes' formula. Therefore we follow an approach analogous to the standard example for a spectral triple on a Riemannian manifold: We construct a so-called \emph{horizontal Dirac operator}, arising from the Clifford action of a horizontal distribution of a Carnot manifold.\smallskip

To really cover the horizontal geometry, we want our horizontal Dirac operator to be compatible with the horizontal part of the Levi-Civita connection on a Carnot manifold $M$. We thus define and discuss a (partial) connection $\nabla^{H}$ according to the sub-bundle $HM$ of $TM$ in Section 7.1 from the Levi-Civita connection on $M$. Then in section 7.2, we introduce horizontal Clifford bundles arising from the Clifford action of the bundle $HM$, which carry a connection compatible with $\nabla^H$. On these bundles, we will be able to define horizontal Dirac operators $D^H$ in a general sense (in analogy to classical Dirac operators on Clifford bundles), and we will be able to modify these operators such that they are self-adjoint. In the end of Section 7.2, we calculate the square of $D^H$ and proof a technical proposition concerning the eigenvalues of a sum of certain Clifford matrices in preparation for future arguments.\smallskip

Finally in Section 3.3, we will show that the horizontal Dirac operators we constructed are indeed the right operators to detect the Carnot-Carath\'{e}odory metric on a Carnot manifold $M$: We show that the norm of $[D^H,f]$ coincides with the Lipschitz norm of $f$ with respect to the Carnot-Carath\'{e}odory metric, such that we can apply Connes' metric formula to the sub-algebra of $C(M)$ consisting of these Lipschitz functions. In addition we show that this metric is already detected by the sub-algebra of $C^\infty$-functions. In the following chapters, we will seen that in spite of all these good properties this operator does not define a spectral triple.\smallskip

Throughout this section, $M$ will be a compact Carnot of step $R$ with Carnot-Carath\'{e}odory metric $d_{CC}$ of dimension $n$. We will use the notation $\{ X_1, \ldots, X_d \}$ for a frame of the horizontal distribution $HM$ of $M$, and if not stated otherwise the Riemannian metric $g$ on $M$ will be chosen such that this frame is orthonormal. The horizontal Dirac operator will be defined on closed (meaning compact without boundary) Carnot manifolds. \bigskip

\section{Horizontal Connections}

To define a horizontal Dirac operator using the Clifford module arising from the horizontal distribution $HM$ we need to introduce a connection on $HM$. This happens straight forward and can be found at various parts of the literature (see e.g. \cite{DGN}, but a horizontal connection is also mentioned in the textbooks \cite{CDPT} and \cite{CC}). The idea is simply to start with the Levi-Civita connection on $M$ and to project it onto the horizontal bundle. This approach is justified by the following proposition (see \cite{CDPT}, Proposition 4.2; not that the case of Carnot manifolds follows directly from the case of Carnot groups formulated there). \medskip

\begin{prop} \label{prop extension of horizontal metric}
 Let $M$ be a Carnot manifold with horizontal distribution $HM \subset TM$ which carries a smoothly varying inner product $\langle \cdot,\cdot \rangle_H$. Let $VM \subset TM$ be a sub-bundle which is complementary to $HM$. If $g_1$ and $g_2$ are Riemannian metrics which make $VM$ and $HM$ orthogonal such that for $j = 1,2$
$$\left. g_j(X,Y) \right|_{x_0} = \left\langle X(x_0),Y(x_0) \right\rangle_H \ \ \ \forall X,Y \in HM \ \forall x_0 \in M,$$
then the associated Levi-Civita connections $\nabla_1$ and $\nabla_2$ coincide when projected to $HM$: We have
\begin{equation} \label{eq equality of extended Levi-Civita connections}
 g_1 \left( \nabla_{1X}Y,Z \right) = g_2 \left(\nabla_{2X}Y,Z \right)
\end{equation}
for all sections $X,Y,Z$ of $HM$. \eB
\end{prop}\medskip

This proposition allows us to define a connection only depending on the horizontal bundle, and it shows that it is in fact well defined when constructing it from any Levi-Civita connection of the extension of the horizontal inner product. We follow \cite{CDPT}, where this is done for the Heisenberg group.\medskip

\begin{defin} \label{def horizontal connection} \normalfont
 Let $M$ be a Carnot manifold with horizontal distribution $HM \subset TM$, which carries a smoothly varying inner product $\langle \cdot,\cdot \rangle_H$. Let $\pi_H: \Gamma^\infty(TM) \rightarrow \Gamma^\infty(HM)$ denote the projection of a tangent vector field onto its horizontal component. We then define
\begin{equation} \label{eq horizontal connection}
 \nabla^H: \Gamma^\infty(HM) \times \Gamma^\infty(HM) \rightarrow \Gamma^\infty(HM), \ \ \ (X,Y) \mapsto \nabla^H_X(Y) := \pi^H \nabla_X(Y),
\end{equation}
where $\nabla$ denotes the Levi-Civita connection of any extension of the horizontal inner product. $\nabla^H$ is called the \emph{horizontal (Levi-Civita-)connection} of $M$. \eBsp
\end{defin}\medskip

\begin{prop} \label{prop properties horizontal connection}
 The map $\nabla^H$ defined in Definition \ref{def horizontal connection} is indeed a (partial) connection defined on the horizontal bundle, which means we have
\begin{enumerate}[(i)]
 \item For all $X,X',Y \in \Gamma^\infty(HM)$ and for all $f,g \in C^\infty(M)$:
$$\nabla^H_{fX+gX'} Y = f \cdot \nabla^H_X Y + g \cdot \nabla^H_{X'} Y.$$
 \item For all $X,Y,Y' \in \Gamma^\infty(HM)$ and for all $f,g \in C^\infty(M)$:
$$\nabla^H_X (fY + gY') = f \cdot \nabla^H_X Y + g \cdot \nabla^H_X Y' + (Xf) \cdot Y + (Xg) \cdot Y'.$$ 
\end{enumerate}
$\nabla^H$ is \emph{metric} with respect to the point-wise inner product $\langle \cdot,\cdot \rangle_H$ of $HM$, which means that we have
\begin{enumerate}[(i)]
 \item[(iii)] For all $X,Y,Y' \in \Gamma^\infty(HM)$:
$$X \left\langle Y,Y' \right\rangle_H = \left\langle \nabla^H_X Y,Y' \right\rangle_H + \left\langle Y,\nabla^H_X Y' \right\rangle_H.$$
\end{enumerate}
$\nabla^H$ is \emph{torsion free} in the horizontal direction, which means that we have
\begin{enumerate}[(i)]
 \item [(iv)] for all $X,Y \in \Gamma(HM)$:
$$\pi^H \left( \nabla^H_X Y - \nabla^H_Y X \right) = \pi^H ([X,Y]). $$
\end{enumerate}
\smallskip

\B Let $\nabla$ be the Levi-Civita connection corresponding to any Riemannian metric $g$ on $TM$ which is an extension of $\langle \cdot,\cdot \rangle_H$. The conditions (i) and (ii) which show that $\nabla^H$ is indeed a connection follow immediately from the corresponding properties of the Levi-Civita connection after projection onto the horizontal distribution. The metric property (iii) follows immediately from the metric property of $\nabla$ since $\langle \cdot,\cdot \rangle_H$ is just a restriction of the Riemannian metric $g$.\smallskip

Finally, the torsion freeness into the horizontal direction follows because
\begin{eqnarray*}
 \pi^H \left( \nabla^H_X Y - \nabla^H_Y X \right) &=& \pi^H \left( \pi^H \left(\nabla_X Y - \nabla_Y X\right) \right)\\
&=& \pi^H \left(\nabla_X Y - \nabla_Y X\right)\\
&=& \pi^H \left( [X,Y]\right)
\end{eqnarray*}
for all $X,Y \in \Gamma^\infty(HM)$, since $\pi_H^2 = \pi_H$ because $\pi_H$ is a projection and since the Levi-Civita connection $\nabla$ is torsion free. Hence the proposition is proved. \eB \smallskip

\Bem We note that we cannot expect to have torsion freeness in the sense that $\nabla^H_X Y - \nabla^H_Y X = [X,Y]$ for all $X,Y \in \Gamma^\infty(HM)$, as it is the case for the Levi-Civita connection: Since $HM$ is not involutive, we have $[X,Y] \notin HM$ for some $X,Y \in HM$, but by definition of the horizontal connection the vector field $\nabla^H_X Y - \nabla^H_Y X$ must be horizontal again. \eBsp 
\end{prop}\medskip

\begin{remark} \normalfont
 We have only defined our horizontal connection to be a partial connection, which means we only allow differentiation into horizontal directions. There are several possibilities to extend $\nabla^H$ to a connection
$$\tilde{\nabla}^H : \Gamma^\infty(TM) \times \Gamma^\infty(HM) \rightarrow \Gamma^\infty(HM):$$
For instance, the expression \eqref{eq horizontal connection} makes sense for any $X \in \Gamma^\infty(TM)$, and Proposition \ref{prop properties horizontal connection} would also hold for this case with the same proof. Another possibility of defining a connection is given by setting $\nabla^H_X := 0$ for all $X \in VM$. \eBsp
\end{remark}\medskip

\begin{remark} \normalfont
We cannot get rid of the vertical directions completely as soon as we want to apply horizontal covariant derivatives more than one time: By the tensorial property (ii) of Proposition \ref{prop properties horizontal connection} we have for any $X_1, X_2, Y \in HM$ and for any $f \in C^\infty(M)$:
\begin{eqnarray*}
 \left( \nabla^H_{X_1} \nabla^H_{X_2} - \nabla^H_{X_2} \nabla^H_{X_1} \right) f \cdot Y &=& \nabla^H_{X_1} \left( f \cdot \nabla^H_{X_2} Y + \left(X_2 f\right) \cdot Y \right)\\
& & - \nabla^H_{X_2} \left( f \cdot \nabla^H_{X_1} Y + \left(X_1 f\right) \cdot Y \right) \\
&=& f \cdot \nabla^H_{X_1} \nabla^H_{X_2} Y + \left(X_1 f\right) \cdot \nabla^H_{X_2} Y + \left(X_2 f\right) \cdot \nabla^H_{X_1}Y\\
& & + \left(X_1X_2 f\right) \cdot Y - f \cdot \nabla^H_{X_2} \nabla^H_{X_1} Y - \left(X_2 f\right) \cdot \nabla^H_{X_1} Y \\
& & - \left(X_1 f\right) \cdot \nabla^H_{X_2}Y - \left(X_2X_1 f\right) \cdot Y \\
&=& f \cdot \left(\nabla^H_{X_1} \nabla^H_{X_2} - \nabla^H_{X_2} \nabla^H_{X_1} \right) Y + \left( [X_1,X_2]f \right) \cdot Y,
\end{eqnarray*}
 and since $HM$ is not involutive we have $[X_1,X_2] \notin \Gamma^\infty(HM)$.\smallskip

This has consequences when one wants to define some kind of curvature belonging to the horizontal connection: It turns out that this curvature shows into the transversal direction of $TM$. But we will not discuss this aspect any further and refer to $\cite{CC}$ and $\cite{Mon}$ instead. \eBsp
\end{remark}\medskip

We close this short section with a proposition which shows us how the calculate a horizontal covariant derivative in local coordinates respecting the grading structure of our Carnot manifold $M$. This is more or less trivial, since after calculating the covariant derivative with respect to the Levi-Civita connection we simply project onto the horizontal distribution.\medskip

\begin{prop} \label{prop horizontal connection in local coordinates}
 Let $M$ be a Carnot manifold with a horizontal distribution $HM$ of rank $d$. Let $\{X_1, \ldots, X_n\}$ denote an orthonormal frame for $TM$ such that for $d \leq n$  $\{X_1, \ldots, X_d\}$ is an orthonormal frame for $HM$. Then we have locally
$$\nabla^H_{X_j} X_k = \sum_{l=1}^d \Gamma_{jk}^l X_l,$$
where the $\Gamma_{jk}^l$ are the Christoffel symbols of the Levi-Civita connection $\nabla$ of $TM$ with respect to the frame $\{X_1, \ldots, X_n\}$. \eB
\end{prop}\bigskip

\section{Construction of Horizontal Dirac Operators}

We want to construct a horizontal Dirac operator in analogy to the (classical) Dirac operator. In the classical case this construction is outlined for example in \cite{Roe} or \cite{LawMich}.\smallskip

We start with a review about the definition of a Clifford algebra and Clifford action. Remember that for each vector space $V$ equipped with a symmetric bilinear form $\langle \cdot,\cdot \rangle$ there exists an (up to isomorphisms) unique \emph{Clifford algebra} $A = \mathrm{Cl}(V)$ which is a unital algebra, equipped with a map $\varphi: V \rightarrow A$ such that
\begin{equation} \label{eq clifford-algebra square}
 \varphi(v)^2 = - \langle v,v \rangle \cdot 1.
\end{equation}
 In addition, $\varphi$ is supposed to fulfill the universal property in the sense that for any other unital algebra $A'$ equipped with map $\varphi': V \rightarrow A'$ satisfying \eqref{eq clifford-algebra square}, there is a unique algebra homomorphism $\alpha: A \rightarrow A'$ such that $\varphi' = \alpha \circ \varphi$ (see \cite{Roe}, Definition 3.1 and Proposition 3.2). Multiplication in this unique algebra is determined by the rule
\begin{equation} \label{eq cliffor-algebra multiplication}
 \varphi(v_1) \varphi(v_2) + \varphi(v_2) \varphi(v_1) = -2\langle v_1,v_2 \rangle,
\end{equation}
 and we know that for $\dim V = n$ we have $\dim \mathrm{Cl}(V) = 2^n$. Using this map $\varphi$ to define multiplication, we can construct a left module over the complex algebra $\mathrm{Cl}(V) \otimes_\R \C$ which will be called a \emph{Clifford module}. In other words, a \emph{Clifford module} $S$ for a real inner product space $V$ is a complex vector space $S$ equipped with an $\R$-linear map 
\begin{equation} \label{eq Clifford algebra representation}
 c: V \rightarrow \en_\C(S) 
\end{equation}
such that $c(v)^2 = - \langle v,v \rangle \cdot 1$ for all $v \in V$. The map $c$ is called the \emph{Clifford action}. \smallskip

Now the question arises what the minimal possible dimension of such a Clifford module $S$ is, depending on the dimension of $V$; in other words, we are looking for irreducible representations of $A$. Such an irreducible representation is given by the so-called \emph{spin representation} $\Sigma$, and it is known that for $n = \dim V$ we have
$$\dim \Sigma = 2^{\left[\frac{n}{2}\right]},$$
where $[ n/2 ]$ denotes the Gaussian bracket which gives the greatest integer smaller or equal then $m/2$ (see \cite{Roe}, Chapter 4, or \cite{LawMich}).\smallskip

The idea for our situation is now to take a look at the Clifford algebras generated by the fibers of the horizontal distribution $HM = \spa \{X_1, \ldots, X_d\}$. Since the Riemannian metric on $M$ is chosen in a way that $\{X_1, \ldots, X_d\}$ forms an orthonormal frame at each point, we immediately get the fundamental properties of this horizontal Clifford action, which we will denote by $c^H$, from the above discussion. \medskip

\begin{prop} \label{prop clifford module over HM}
 For $x \in M$, let $S_x$ be a Clifford module for $H_xM$ with (horizontal) Clifford action $c^H: H_xM \rightarrow \en_\C(S_x)$. Then we have:
\begin{enumerate}[(i)]
 \item $c^H(X_j)^2 = -\mathrm{Id}$ for all $j \in \{1, \ldots, d\}$.
 \item $c^H(X_j)c^H(X_k) + c^H(X_k)c^H(X_j) = 0$ for all $j,k \in \{1, \ldots, d\}$, $j \neq k$. \eB
\end{enumerate}
\end{prop}\medskip

Now let $S^HM$ be a bundle of Clifford modules for the horizontal distribution $HM$. We need to equip $S^HM$ with an point-wise (hermitian) inner product and with a connection, for which we claim certain compatibility conditions. For a classical Dirac operator, we ask the connection on $SM$ to be compatible with the Levi-Civita connection $\nabla$ on $TM$, see e.g. \cite{Roe}, Definition 3.4. Hence to define a suitable bundle and connection where a horizontal Dirac operator can act on, we would like to have compatibility with the horizontal connection $\nabla^H$ on $HM$ defined in Section 3.1 via
\begin{equation} \label{eq horizontal Levi-Civita connection Section 3.2}
  \nabla^H: \Gamma^\infty(HM) \times \Gamma^\infty(HM) \rightarrow \Gamma^\infty(HM), \ \ \ (X,Y) \mapsto \nabla^H_X Y := \pi_H \nabla_X Y,
\end{equation}
where $\pi_H$ is the orthogonal projection onto the horizontal distribution. For this horizontal connection we formulate the compatibility conditions in the following way:\medskip

\begin{defin} \label{def horizontal clifford bundle}
 Let $S^H M$ be a bundle of horizontal Clifford moduls over $M$ which is equipped with a fiber-wise Hermitian metric $\left( \cdot,\cdot \right)_H$ and a metric connection $\nabla^{S^H}$. $S^H M$ is called a \emph{horizontal Clifford bundle} if
  \begin{enumerate}[(i)]
   \item For each $x \in M$ we have 
$$\left( c^H(X_x) \sigma_1(x) , \sigma_2(x) \right)_H + \left( \sigma_1(x), c^H(X_x) \sigma_2(x) \right)_H = 0$$
 for all $X_x \in H_xM$, $\sigma_1, \sigma_2 \in \Gamma^\infty(S^HM)$.
   \item If $\nabla^H$ is the horizontal Levi-Civita connection \eqref{eq horizontal Levi-Civita connection Section 3.2} on $HM$, we have 
$$\nabla^{S^H}_X (c^H(Y) \sigma) = c^H(\nabla^H_X Y) \sigma + c^H(Y) \nabla^{S^H}_X \sigma$$
 for all $X,Y \in \Gamma^\infty(HM)$ and for all sections $\sigma \in \Gamma^\infty(S^HM)$.
  \end{enumerate}
\eBsp
\end{defin} \medskip

Now, still in analogy to the classical Dirac operator, we can define a horizontal Dirac operator. Later we will see that this operator has to be modified by adding an endomorphism in order to be symmetric. \medskip

\begin{defin} \label{def formal horizontal dirac operator}
 The \emph{formal horizontal Dirac operator} $\tilde{D}^H$ of $S^HM$ is the first order differential operator on $\Gamma^\infty(S^HM)$ defined by the composition
 $$ \Gamma^\infty(S^HM) \rightarrow \Gamma^\infty(H^\ast M \otimes S^HM) \rightarrow \Gamma^\infty(HM \otimes S^HM) \rightarrow \Gamma^\infty(S^HM). $$
 Here the first arrow is given by the connection $\nabla^{S^H}$, the second arrow is given by the identification of $H^\ast M$ and $HM$ via the horizontal metric, and the third arrow is given by the Clifford action.\smallskip

 If we choose a local orthonormal frame $\{X_1, \ldots, X_d\}$ of sections of $HM$, we can write
 \begin{equation}  \label{eq formal D^H local}
  \tilde{D}^H \sigma = \sum_{j=1}^d c^H(X_j) \nabla_{X_j}^{S^H} \sigma. 
 \end{equation}
\eBsp
\end{defin}\medskip

Before we go on with our construction to get a self-adjoint operator from $\tilde{D}^H$, we have a look at the most natural example.\medskip

\begin{example} \label{ex exterior bundle} \normalfont
Consider the exterior bundle $\Omega M := \bigwedge T^\ast M \otimes \C$ on a compact Riemannian manifold $M$ without boundary equipped with its natural metric and connection. It is well known (see for example \cite{Roe}, 49ff.) that (in the classical sense) this is a Clifford bundle with multiplication given by the wedge product and Clifford action of $e \in T^\ast M \cong TM$ given by
$$c(e) \sigma = e \wedge \sigma + e \hook \sigma.$$
The Dirac operator of this Clifford bundle is given by $d+d^\ast$. Remember that the interior product is defined by
$$e \hook \sigma = (-1)^{nk+n+1} \ast (e \wedge \ast \sigma)$$
via the Hodge star operator $\ast \sigma$, which (for a $k$-form $\sigma$) is defined to be the unique $(n-k)$-form $\alpha$ such that
$$(\sigma, \alpha) \mathrm{vol} = \alpha \wedge \ast \sigma.$$
Using this Hodge star operator, the operator $d^\ast$ is defined by
$$d^\ast \sigma = \left(-1\right)^{nk+n+1} \ast d \ast \sigma.$$

Now a horizontal Clifford bundle can be constructed analogously. We assume $M$ to be a Carnot manifold with horizontal distribution $HM = \spa \{X_1, \ldots, X_d\}$ such that the frame $\{X_1, \ldots, X_d, X_{d+1}, \ldots, X_n\}$ is orthonormal, Hence we can identify this frame with its dual frame $\{d\omega^1, \ldots, d\omega^n\}$. We further identify $H^\ast M$ with the sub-bundle of $T^\ast M$ annihilating $(HM)^\perp$, which means $H^\ast M = \spa \{d\omega^1, \ldots, d\omega^d\}$. Then $\Omega^HM := \bigwedge H^\ast M \otimes \C$ is a horizontal Clifford bundle: The Clifford action of $e \in H^\ast M \cong HM$ is again given by
$$c^H(e) \sigma := e \wedge \sigma + e \hook \sigma;$$
note that in case $\sigma \in S^HM$ and $e \in H^\ast M$ the right side of this equation lies still in $S^HM$. If we use the horizontal exterior derivative
$$d^H(\sigma) := \pi_H(d \sigma),$$
where $\pi_H$ is the orthonormal projection onto $H^\ast M$, we get a connection which is compatible with the horizontal Levi-Civita connection is defined on $H^\ast M$. Using this, we can define the formal horizontal Dirac operator on $S^H M$ and see that it is given by
$$D^H \sigma = \sum_{j=1}^d d^H\omega^j \wedge \sigma + d^H\omega^j \hook \omega = d^H \sigma + \left(d^H\right)^\ast \sigma.$$
\eBsp
\end{example}\medskip

Our approach of constructing a horizontal Clifford bundle and a horizontal Dirac operator can be found in the literature in a greater generality. Igor Prokhorenkov and Ken Richardson recently introduced a class of so-called \emph{transversally Dirac operators} in \cite{PR} by considering a distribution $QM \subset TM$ of the tangent bundle (which may be integrable or not), which furnishes a $\mathrm{Cl}(QM)$-module structure on a complex Hermitian vector bundle $EM$ over $M$ with Clifford action $c: QM \rightarrow \en_\C(EM)$ and a connection $\nabla^E$ on $EM$ fulfilling the requirements of Definition \ref{def horizontal clifford bundle}. Hence our construction of a horizontal Clifford bundle fits into this setting. \smallskip

It turns out that the gap of this construction so far is that the operator $\tilde{D}^H$ is not symmetric: When one wants to calculate the $L^2$-adjoint of $\tilde{D}^H$, Stokes' theorem causes an additional term of mean curvature into the vertical direction, projected to $HM$ via the horizontal connection. Hence to get a symmetric operator, we have to add the (horizontal) Clifford action by the mean curvature of the orthogonal complement of the horizontal distribution. This was shown by Prokhorenkov and Richardson for the general case of a distribution of the tangent bundle (\cite{PR}), such that we refer to their work formulating the following theorem. It shows in addition that the resulting operator is essential self-adjoint, which follows by a theorem of Paul Chernoff about the self-adjointness of certain differential operators (\cite{Che}).\smallskip

Altogether, we have the following theorem (see \cite{PR}, Theorem 3.1), which will finally provide us with the possibility to define (essential self-adjoint) horizontal Dirac operators. \medskip

\begin{thm} \label{thm D^H selfadjoint}
 Let $\tilde{D}^H$ be a formal horizontal Dirac operator, acting on the smooth sections of a horizontal Clifford bundle $S^HM$ over a compact Carnot manifold $M$ without boundary. Then we have
\begin{enumerate}[(i)]
 \item The formal $L^2$-adjoint of $\tilde{D}^H$ is given by
\begin{equation} \label{eq adjoint of formal D^H}
 \left(\tilde{D}^H\right)^\ast = \tilde{D}^H - c^H\left( \sum_{j=d+1}^n \pi^H \nabla_{X_j} X_j \right),
\end{equation}
where $\pi^H: TM \rightarrow HM$ is the orthogonal projection onto the horizontal distribution and $\nabla$ is the Levi-Civita connection of $TM$.
\item The operator
\begin{equation} \label{eq D^H}
 D^H := \sum_{j=1}^k c^H(X_j) \nabla_{X_j}^{S^H} - \frac{1}{2} c^H\left( \sum_{j=d+1}^n \pi^H \nabla_{X_j} X_j \right)
\end{equation}
is essentially self adjoint on $L^2(S^HM)$.
\end{enumerate}
\eB 
\end{thm}\medskip

\begin{defin} \label{def horizontal dirac operator}
 Let $S^H M$ be a horizontal Clifford bundle over a closed Carnot manifold $M$. Then the operator defined by \eqref{eq D^H}, acting on $\Gamma^\infty(S^H M)$, is called the \emph{horizontal Dirac operator} of $S^HM$. \eBsp
\end{defin} \medskip

At this point, we want to state another example which is due to \cite{PR}. It states that on any (classical) Clifford bundle $EM$ over a closed Carnot manifold $M$ one can implement the structure of a horizontal Clifford bundle by adjusting the bundle connection. We refer to \cite{PR}, Section 2, for the following proposition. \medskip

\begin{prop} \label{prop horizontal Clifford bundle from spinor bundle}
 Let $M$ be a Carnot manifold with horizontal distribution $HM$, such that $\{X_1, \ldots, X_d\}$ is an orthonormal frame for $HM$ and $\{X_{d+1}, \ldots, X_n\}$ is an orthonormal frame for $HM^\perp$. Let in addition $EM$ be a Clifford bundle over $M$ with Clifford action $c$ and bundle connection $\nabla^E$, and let $c^H$ denote the restriction of $c$ to the horizontal distribution $HM$.\smallskip

 Then there is a connection $\tilde{\nabla}^E$ on $E$ such that $E$ equipped with $c^H$ and $\tilde{\nabla}^E$ is a horizontal Clifford module. This connection $\tilde{\nabla}^E$ is given by
\begin{equation}
 \tilde{\nabla}^E_X = \nabla^E_X + \frac{1}{2} \sum_{j=d+1}^n c^H\left(\pi^H \nabla_X X_j \right) c(X_j),
\end{equation}
where $\pi^H$ denotes the orthogonal projection of $TM$ onto $HM$ and $\nabla$ denotes the Levi-Civita connection on $TM$. \eB
\end{prop}\medskip

In Chapter 4, we will consider another example in detail: On the local homogeneous space of a Carnot group $\G$ of rank $d$ under the action of a lattice subgroup, there is a submersion $\pi: \Gamma \backslash \G \rightarrow \T^d$ onto the $d$-dimensional torus by Section 2.4. Then one can define a horizontal Clifford bundle and a horizontal Dirac operator by pulling back the objects from $\T^d$, which does not depend on the non-horizontal directions of $TM$. We will discuss this example extensively in Chapter 4, since it serves well as a toy model to the question whether one can define spectral triples from horizontal Dirac operators.\smallskip

For further considerations it will be important to work with the square of the horizontal Dirac operator. Especially it will follow in Section 6.3 that this operator, considered as an operator on $L^2(S^H M)$, is not hypoelliptic, and we will conclude from that that a horizontal Dirac operator in not hypoelliptic. Hence we calculate $\left(D^H\right)^2$, whose principal term will be a horizontal Laplacian, in a local expression.\medskip

\begin{prop} \label{prop D^H squared}
Let $M$ be a Carnot manifold, and let $\{X_1, \ldots, X_d\}$ be an orthonormal frame for its horizontal distribution $HM$. Let $D^H$ be a horizontal Dirac operator acting on a horizontal Clifford bundle $S^H M$ over $M$. Then we have locally
\begin{equation} \label{eq D^H squared locally}
 \left(D^H\right)^2 = - \sum_{j=1}^d X_j^2 + \sum_{j<k} c^H(X_j)c^H(X_k) \left[X_j,X_k\right] + O_H(1),
\end{equation}
where $X_j$ is to be understood as a component-wise directional derivative in a local chart and $O_H(1)$ denotes a graded differential operator of order smaller or equal to $1$ (which only depends on the differential operators $X_1, \ldots, X_d$ and endomorphisms on the bundle).\smallskip

\B We work with the local expression of $D^H$, given by \eqref{eq D^H}. From this we get for $\sigma \in S^HM$
\begin{equation} \label{eq D^H squared first equation}
\begin{split}
  \left(D^H\right)^2 \sigma = & \sum_{j=1}^d c^H(X_j) \nabla^{S^H}_{X_j} \left( \sum_{k=1}^d c^H(X_k) \nabla^{S^H}_{X_k} \sigma - \frac{1}{2} c^H\left(Z \right) \sigma \right)\\
& - \ \frac{1}{2} c^H(Z) \left( \sum_{k=1}^d c^H(X_k) \nabla^{S^H}_{X_k} \sigma - \frac{1}{2} c^H(Z) \sigma \right)
\end{split}
\end{equation}
where
$$Z = \sum_{j=d+1}^n \pi^H \nabla_{X_j} X_j \in HM$$
is the mean curvature of the bundle $HM^\perp$ like in Theorem \ref{thm D^H selfadjoint}. Note that we have $c^H(Z) \in \en_\C (S^HM)$. We observe that the second summand on the right hand side of \eqref{eq D^H squared first equation} consists only of differential operators of graded order $1$ or $0$ applied to $\sigma$, so it is contained in the expression $O_H(1)$. Further we have
$$\nabla^{S^H}_{X_j} \left(c^H\left(Z\right) \sigma \right) = c^H\left(\nabla^H_{X_j} Z\right) \sigma + c^H\left( Z \right) \nabla^{S^H}_{X_j} \sigma$$
by Definition \ref{def horizontal clifford bundle}, since $S^HM$ is a horizontal Clifford bundle, so this expression is also contained in $O_H(1)$. \smallskip

Now we cannot assume $\nabla^H_{X_j} X_k = 0$ as one does when calculating the square of an ordinary Dirac operator in this case, since our frame is fixed and therefore cannot be chosen to be synchronous at a point $x \in M$. We calculate for $j,k \in \{1, \ldots, d\}$, using once again the properties of a horizontal Clifford bundle:
\begin{eqnarray*}
 c^H(X_j) \nabla^{S^H}_{X_j} \left(c^H(X_k) \nabla^{S^H}_{X_k} \sigma \right) &=& c^H(X_j) \left( c^H\left(\nabla^H_{X_j} X_k\right) \nabla^{S^H}_{X_k} \sigma + c^H(X_k) \nabla^{S^H}_{X_j} \nabla^{S^H}_{X_k} \sigma \right)\\
&=& c^H(X_j)c^H(X_k) \nabla^{S^H}_{X_j} \nabla^{S^H}_{X_k} \sigma + O_H(1) \sigma.
\end{eqnarray*}
Since we are working in local charts, $\nabla^{S^H}_{X_j}$ is locally given by the expression $X_j + \Gamma_j$ with $\Gamma_j \in \en (S^HM)$. Applying this to the last expression, we find
\begin{eqnarray*}
 \nabla^{S^H}_{X_j} \nabla^{S^H}_{X_k} &=& \left(X_j + \Gamma_j\right) \left(X_k + \Gamma_k\right)\sigma\\
&=& X_jX_k \left(\sigma\right) + X_j \left(\Gamma_k \sigma\right) + \Gamma_j X_k(\sigma) + \Gamma_j \Gamma_k \sigma\\
 &=& X_jX_k \left(\sigma\right) + O_H(1) \sigma
\end{eqnarray*}
and therefore
\begin{eqnarray*}
& & \sum_{j=1}^d \sum_{k=1}^d c^H(X_j) \nabla^{S^H}_{X_j} \left(c^H(X_k) \nabla^{S^H}_{X_k} \sigma \right)\\
&=& \sum_{j=1}^d \sum_{k=1}^d c^H(X_j) c^H(X_k) X_jX_k \left(\sigma\right) + O_H(1) \sigma\\
&=& \sum_{j=k} -X_j^2 \left(\sigma\right) + \sum_{j \neq k} c^H(X_j)c^H(X_k) X_jX_k \left(\sigma\right) + O_H(1) \sigma\\
&=& -\sum_{j=1}^d X_j^2 \left(\sigma\right) + \sum_{j<k} c^H(X_j)c^H(X_k) \left(X_jX_k - X_kX_j\right) \left(\sigma\right) + O_H(1) \sigma,
\end{eqnarray*}
since $c^H(X_j)^2 = -1$ and $c^H(X_j)c^H(X_k) = - c^H(X_k)c^H(X_j)$ by the characterization of the Clifford algebra. Finally we plug this into \eqref{eq D^H squared first equation} and find together with the discussion above
\begin{eqnarray*}
  \left(D^H\right)^2 \sigma &=& \sum_{j=1}^d c^H(X_j) \nabla^{S^H}_{X_j} \left( \sum_{k=1}^d c^H(X_k) \nabla^{S^H}_{X_k} \sigma \right) + O_H(1) \sigma\\
&=& -\sum_{j=1}^d X_j^2 \left(\sigma\right) + \sum_{j<k} c^H(X_j)c^H(X_k) \left[X_j,X_k\right] \left(\sigma\right) + O_H(1)\sigma.
\end{eqnarray*}
This proves the proposition. \eB \smallskip

\Bem The proof of this proposition shows that the leading term 
$$- \sum_{j=1}^d X_j^2 + \sum_{j<k} c^H(X_j)c^H(X_k) \left[X_j,X_k\right]$$
 in the local expression \eqref{eq D^H squared locally} of $(D^H)^2$ will not change if we modify $D^H$ by an term of graded order zero. For an operator which is locally of the type
$$\tilde{D}^H = D^H + \gamma(x)$$
with $\gamma \in C^\infty(M, \en_\C (S^HM) )$ we still have 
$$ \left(\tilde{D}^H\right)^2 = - \sum_{j=1}^d X_j^2 + \sum_{j<k} c^H(X_j)c^H(X_k) \left[X_j,X_k\right] + O_H(1).$$
In particular this shows that for a modified connection
$$\tilde{\nabla}^{S^H}_{X_j} = \nabla^{S^H}_{X_j} + \gamma_0(x)$$
with $\gamma_0 \in C^\infty(M, \en_\C (S^HM))$ the leading term of the local expression of $(D^H)^2$ according to this connection will not change. We will refer to this remark when we check $D^H$ for hypoellipticity in Section 6.3, since the property of being hypoelliptic only depends of the leading term of a horizontal Laplacian. \eBsp
\end{prop}\medskip

Before we conclude this section, we want to prove a rather technical lemma about the eigenvalues a certain combination of Clifford matrices can have. This will be needed in the following chapters, and since it will be needed more than one time we decided to put it here. Its statement is quite general, since it can be applied to the Clifford action arising from any vector bundle, horizontal or not. \medskip

\begin{prop} \label{prop eigenvalues sum c(X_j)c(X_k)}
 Let $V$ be a vector space of dimension $d$, and let $A = \mathrm{Cl}(V) \otimes \C$ denote its (complexified) Clifford algebra, which is represented by matrices on a complex vector space $S$ via \eqref{eq Clifford algebra representation}. We assume that we have an $m \in \N$ such that $2m \leq d$ and such that $\{e_1, \ldots, e_{2m}, \ldots, e_d\}$ is a basis of $V$.\smallskip

Then we have the following:
\begin{enumerate}[(i)]
 \item For any $k \neq l$ with $1 \leq k,l \leq d$, all the eigenvalues of the matrix $c(e_k)c(e_l)$ are given by the numbers $\pm i$, where both eigenvalues have the same multiplicity.

 \item The eigenvalues of the matrix 
\begin{equation} \label{eq sum Clifford matrices}
 \sum_{j=1}^m c(e_j)c(e_{m+j})
\end{equation}
are exactly the numbers
$$\mu_l = i \left(-m + 2l\right) \ \ \ \text{with} \ l = 0, \ldots, m .$$
In the case where $d=2m$ and the Clifford action of $V$ on $S$ is irreducible, which means we have $\dim S = 2^m$, each eigenvalue $\mu_l$ has the multiplicity $\tbinom{m}{l}$.

 \item For $\lambda_j \in \R$ with $\lambda_j > 0$ (with $j \in \{1, \ldots, m\}$), all the eigenvalues of the matrix
\begin{equation} \label{eq sum Clifford matrices weighted}
 \sum_{j=1}^m \lambda_j c(e_j)c(e_{m+j}).
\end{equation}
are included in the interval 
$$ \left[-i \sum_{j=1}^m \lambda_j, i \sum_{j=1}^m \lambda_j \right] \subset i \R $$
on the imaginary line. Thereby, the numbers $i \sum_{j=1}^m \lambda_j$ and $-i \sum_{j=1}^m \lambda_j$ are eigenvalues of the matrix \eqref{eq sum Clifford matrices weighted}.
\end{enumerate} \smallskip

\Bem Although this proposition has been formulated for vector spaces, it can be transferred to Clifford bundles immediately, since point-wise we have the eigenvalues noted above. In particular, the situation of the proposition occurs for the horizontal Clifford action on a Clifford bundle on a Carnot manifold $M$, where $M$ is locally diffeomorph to a Carnot group of the type $\Hei^{2m+1} \times \R^{d-2m}$. In this case we will identify the horizontal frame $\{X_1, \ldots, X_d\}$ at each point $x \in M$ with the vector space spanned by $\{e_1, \ldots, e_d\}$, for which the proposition is formulated. \eBsp

\B First of all, since
\begin{eqnarray*}
 \left(c(e_k)c(e_l)\right)^2 &=& c(e_k)c(e_l)c(e_k)c(e_l)\\
& =&  -\left(c(e_k)\right)^2 \left(c(e_l)\right)^2 \\
& =&  -\id \ \ \ \forall 1 \leq k,l \leq d,
\end{eqnarray*}
each of the products $c(e_k)c(e_l)$ has eigenvalues $i$ and $-i$. Furthermore, the matrices $c(e_k)c(e_l)$ and $c(e_l)c(e_k)$ commute since
$$c(e_k)c(e_l)c(e_l)c(e_k) = c(e_l)c(e_k)c(e_k)c(e_l) = \id,$$
and are thus simultaneously diagonalizable. Because of the relation $c(e_k)c(e_l) = -c(e_l)c(e_k)$ this means that we have the same number of $i$- and $-i$-eigenvalues, and therefore statement (i) is proved. \smallskip

The next observation is that for $j \neq k$ and $1 \leq j,k \leq m$, the matrices $c(e_j)c(e_{m+j})$ and $c(e_k)c(e_{m+k})$ commute: We have
\begin{eqnarray*}
 & & [c(e_j)c(e_{m+j}),c(e_k)c(e_{m+k})]\\
&=& c(e_j)c(e_{m+j})c(e_k)c(e_{m+k}) - c(e_k)c(e_{m+k})c(e_j)c(e_{m+j})\\
&=& c(e_k)c(e_{m+k})c(e_j)c(e_{m+j}) - c(e_k)c(e_{m+k})c(e_j)c(e_{m+j})\\
&=& 0
\end{eqnarray*}
by the commutation rules of the horizontal Clifford action, see Proposition \ref{prop clifford module over HM}. But this means that all the summands of the matrix $\sum_{j=1}^m \lambda_j c(e_j)c(e_{m+j})$ are simultaneously diagonalizable. Let $S$ be a matrix diagonalizing these summands simultaneously, this means its diagonal matrix is given by
\begin{equation} \label{eq diagonal matrix non-hypoellipticity D^H squared}
  S^{-1} \left(\sum_{j=1}^m \lambda_j c(e_j)c(e_{m+j}) \right) S = \sum_{j=1}^m \lambda_j S^{-1} c(e_j)c(e_{m+j}) S.
\end{equation}
Using this, we can prove the rest of the statements by induction over $m$.\smallskip

From this point on, we will assume for a moment that the representation of our Clifford algebra is irreducible. This means that we have a map
$$c: V \rightarrow \en_\C(S)$$
for a vector space $S$ with $\dim S = 2^{[d/2]}$, $[ \cdot ]$ denoting the Gaussian bracket, such that for every $j \in \{1, \ldots, d\}$ we have 
$$c(e_j) \in \mathrm{Mat}_{2^{[d/2]} \times 2^{[d/2]}}(\C)$$
(see the discussion in the beginning of this section for this). We further denote the (unique) complexified Clifford algebra arising from a (complex) vector space of dimension $n$ by $\mathrm{Cl}_\C(n)$. For the proof by induction, we will use the isomorphism
\begin{equation} \label{eq rekursive formula for Clifford algebra}
 \mathrm{Cl}_\C(n+2) \cong \mathrm{Cl}_\C(n) \otimes \mathrm{Cl}_\C(2),
\end{equation}
see \cite{LawMich}, Theorem I.4.3.\smallskip

Now we start with the induction argument: For the case $m=1$ the statements (ii) and (iii) already follow from statement (i). Now we assume the statements to be true for $m-1$. In detail, this means the following: For a vector space $\tilde{V} \cong \R^{d-2}$ with basis $\left\{\tilde{e}_1, \ldots, \tilde{e}_{d-2}\right\}$, we have a Clifford algebra $\mathrm{Cl}_{m-1}(\tilde{V}) \otimes \C \cong \mathrm{Cl}_\C(d-2)$. Since we assumed irreducibility of the Clifford action, the elements $c_{m-1}(\tilde{e}_k)$ are given by $(2^{[d/2] - 1} \times 2^{[d/2] - 1})$-matrices, such that the statements (ii) and (iii) are true for sums of the form
$$\sum_{j=1}^{m-1} \lambda_j c(\tilde{e}_j) c(\tilde{e}_{m-1+j}),$$
considering $m-1$ instead of $m$.\smallskip

Now for a vector space $V \cong \R^d$ we have $\mathrm{Cl}_d(V) \otimes \C \cong \mathrm{Cl}_\C(d)$ and consider a basis $\left\{e_1, \ldots, e_{d} \right\}$. Hence after using the isomorphism from \eqref{eq rekursive formula for Clifford algebra}, we can work with the following representations of the elements $e_1, \ldots, e_d$:
\begin{equation} \label{eq representations Cl(2m)}
 c_m(e_j) = \begin{cases}
             c_{m-1}(\tilde{e}_j) \otimes \begin{pmatrix} i&0 \\ 0&-i \end{pmatrix} & \text{for} \ 1 \leq j \leq m-1 \\
             \id_{2^{[d/2] - 1}} \otimes \begin{pmatrix} 0&i \\ i&0 \end{pmatrix} & \text{for} \ j=m \\
             c_{m-1}(\tilde{e}_{j-1}) \otimes \begin{pmatrix} i&0 \\ 0&-i \end{pmatrix} & \text{for} \ m+1 \leq j \leq 2m-1 \\
             \id_{2^{[d/2] - 1}} \otimes \begin{pmatrix} 0&-1 \\ 1&0 \end{pmatrix} & \text{for} \ j=2m \\
             c_{m-1}(\tilde{e}_{j-2}) \otimes \begin{pmatrix} i&0 \\ 0&-i \end{pmatrix} & \text{for} \ 2m+1 \leq j \leq d,\\
            \end{cases}
\end{equation}
where $c_{m-1}$ denotes the Clifford action by $\mathrm{Cl}_{m-1}(\tilde{V}) \otimes \C$ and $c_m$ denotes the Clifford action by $\mathrm{Cl}_{m}(V) \otimes \C$ on the corresponding vector spaces. Note that the matrices
$$E_1 = \begin{pmatrix} i&0 \\ 0&-i \end{pmatrix}, \ E_2 = \begin{pmatrix} 0&i \\ i&0 \end{pmatrix}, \ E_3 = \begin{pmatrix} 0&-1 \\ 1&0 \end{pmatrix},$$
together with the $(2 \times 2)$-unit matrix, form a basis for the Clifford algebra $\mathrm{Cl}_\C(2)$. For the representations from \eqref{eq representations Cl(2m)} we get by the rules of the tensor product for matrices:
\begin{eqnarray*}
 & & \sum_{j=1}^m \lambda_j c_m(e_j)c_m(e_{m+j}) \\
&=& \sum_{j=1}^{m-1} \lambda_j c_m(e_j)c_m(e_{m+j}) + \lambda_m c_m(e_m)c_m(e_{2m})\\
&=& \sum_{j=1}^{m-1} \lambda_j \left(c_{m-1}(\tilde{e}_j) \otimes \begin{pmatrix} i&0 \\ 0&-i \end{pmatrix} \right) \cdot \left( c_{m-1}(\tilde{e}_{m-1+j}) \otimes \begin{pmatrix} i&0 \\ 0&-i \end{pmatrix} \right) \\
& & + \lambda_m \left( \id_{2^{[d/2] - 1}} \otimes \begin{pmatrix} 0&i \\ i&0 \end{pmatrix} \right) \cdot \left( \id_{2^{[d/2] - 1}} \otimes \begin{pmatrix} 0&-1 \\ 1&0 \end{pmatrix} \right) \\
&=& \sum_{j=1}^{m-1} \lambda_j \left(c_{m-1}(\tilde{e}_j) c_{m-1}(\tilde{e}_{m-1+j}) \otimes \begin{pmatrix} -1&0 \\ 0&-1 \end{pmatrix} \right) + \lambda_m \cdot \id_{2^{[d/2] - 1}} \otimes \begin{pmatrix} i&0 \\ 0&-i \end{pmatrix}\\
&=& \left( \sum_{j=1}^{m-1} \lambda_j c_{m-1}(\tilde{e}_j) c_{m-1}(\tilde{e}_{m-1+j}) \right) \otimes \begin{pmatrix} -1&0 \\ 0&-1 \end{pmatrix} + \lambda_m \cdot \id_{2^{[d/2] - 1}} \otimes \begin{pmatrix} i&0 \\ 0&-i \end{pmatrix}.
\end{eqnarray*}
Now let $\tilde{S}$ be the matrix diagonalizing the matrices $c_{m-1}(\tilde{e}_j)c_{m-1}(\tilde{e}_{m+j})$ simultaneously, such that we have \eqref{eq diagonal matrix non-hypoellipticity D^H squared} in this situation. But then we see from the above calculation that the matrix 
$$S := \tilde{S} \otimes \begin{pmatrix} 1&0 \\ 0&1 \end{pmatrix}$$
diagonalizes $\sum_{j=1}^m \lambda_j c_m(e_j)c_m(e_{m+j})$ simultaneously, which means we have because of the above calculation
\begin{equation} \label{eq diagonalized sum of Cl(2m)}
\begin{split}
  & S^{-1} \sum_{j=1}^m \lambda_j c_m(e_j)c_m(e_{m+j}) S  \\
= & \left( \tilde{S}^{-1} \otimes \begin{pmatrix} 1&0 \\ 0&1 \end{pmatrix} \right) \cdot \left( \sum_{j=1}^m \lambda_j c_m(e_j)c_m(e_{m+j}) \right) \cdot \left( \tilde{S} \otimes \begin{pmatrix} 1&0 \\ 0&1 \end{pmatrix} \right)\\
= & \tilde{S}^{-1} \left( \sum_{j=1}^{m-1} \lambda_j c_{m-1}(\tilde{e}_j) c_{m-1}(\tilde{e}_{m-1+j}) \right) \tilde{S} \otimes \begin{pmatrix} -1&0 \\ 0&-1 \end{pmatrix} + \lambda_m \cdot \id_{2^{[d/2] - 1}} \otimes \begin{pmatrix} i&0 \\ 0&-i \end{pmatrix}.
\end{split}
\end{equation}
Since all these matrices are diagonal matrices, we immediately see that for each eigenvalue $\tilde{\mu}_l$, $l \in \{0, \ldots, m-1\}$, of $\sum_{j=1}^{m-1} \lambda_j c_{m-1}(\tilde{e}_j) c_{m-1}(\tilde{e}_{m-1+j})$ with multiplicity $\tilde{\nu}_l$ the numbers 
\begin{equation} \label{eq eigenvalues mu_l}
 \mu_l^+ := -\tilde{\mu}_l + i\lambda_m \ \ \ \text{and} \ \ \ \mu_l^- := -\tilde{\mu}_l-i\lambda_m 
\end{equation}
are eigenvalues of $ \sum_{j=1}^m \lambda_j c_m(e_j)c_m(e_{m+j})$. But from \eqref{eq eigenvalues mu_l} we can prove the statements (ii) and (iii) by induction:
\begin{itemize}
 \item Statement (iii) can be seen immediately from \eqref{eq eigenvalues mu_l}. If $\tilde{\mu}_0 := i\sum_{j=1}^{m-1} \lambda_j $ is the greatest and $\tilde{\mu}_{m-1} := -i\sum_{j=1}^{m-1} \lambda_j$ is the lowest eigenvalue of 
$\sum_{j=1}^{m-1} \lambda_j c(\tilde{e}_j) c(\tilde{e}_{m-1+j})$, then the greatest and lowest eigenvalue of \eqref{eq sum Clifford matrices weighted} are given by the numbers 
$$\mu_{m+1}^+ = \tilde{\mu}_{m-1} + i\lambda_m \ \ \ \text{and} \ \ \ \mu_0^- = \tilde{\mu}_0 - i \lambda_m,$$
which proves both statements of (iii).

 \item Statement (ii) follows after choosing $\lambda_1 = \ldots = \lambda_m = 1$ in \eqref{eq eigenvalues mu_l}. By the assumption that (ii) is true for $m-1$, we have
$$\mu_l^+ = i((m-1) + 1 -2l) = i(m-2l)$$
and
$$\mu_l^- = i((m-1) -1 -2l) = i(m-2(l+1))$$
for $0 \leq l \leq m-1$, where each $\mu^\pm_l$ has multiplicity $\binom{m-1}{l}$. We observe for $l \leq m-2$ that 
$$\mu_{m-l+1} := \mu^+_{l+1} = \mu^-_l,$$
hence each of these eigenvalues has multiplicity $\binom{m-1}{l} + \binom{m-1}{l+1} = \binom{m}{l+1}$. In addition we have the eigenvalues $\mu_0 := \mu^-_{m-1} = -im$ and $\mu_m := \mu^+_0 = im$, which are both of multiplicity $1$.
\end{itemize}

We finally drop the restriction that the representation of the Clifford algebra is irreducible: We do not get any new eigenvalues, because every (reducible) representation is a direct sum of irreducible ones, and hence (ii) and (iii) are also true for this case. But we cannot make a general statement about the multiplicity of the eigenvalues in the reducible case. \smallskip

Altogether, every statement of this proposition is proved.\eB
\end{prop}\bigskip

\section{Detection of the Carnot-Carath\'{e}odory Metric}

Throughout this section we will assume that $M$ is a compact Carnot manifold without boundary, such that the algebra $C(M)$ is a unital $C^\ast$-algebra which can be represented on $L^2(M)$ via left multiplication. Our intention is to show that the operator $D^H$ from the previous section detects the Carnot-Carath\'{e}odory metric via the Connes metric formula
\begin{equation} \label{eq Connes metric formula for D^H}
 d_{CC}(x,y) = \sup \left\{|f(x)-f(y): f \in \mathcal{A}', \left\|[D^H,f]\right\| \leq 1\right\},
\end{equation}
where $\mathcal{A}'$ is a dense sub-algebra of $C(M)$. This means that although $D^H$ does not furnish a spectral triple (which we will see in general in Chapter 6), we can consider the triple $(C(M), L^2(M), D^H)$ as a compact quantum metric space in the sense of Mark Rieffel, see Definition \ref{def compact quantum metric space}, with the corresponding Lip-norm $L(f) := \left\|[D^H,f]\right\|$, whose metric is exactly the Carnot-Carath\'{e}odory metric on $M$. \smallskip

The key observation is that, in analogy to the classical case from the standard example for a spectral triple, the commutator $[D^H,f]$ acts as Clifford action by the horizontal gradient of a function $f$. Remember that on a Carnot manifold $M$ the horizontal gradient of a function $f \in C^1(M)$ is given by the vector field
$$\grad^H(f) = \sum_{j=1}^d X_j(f) \cdot X_j,$$ 
where, like before, $HM = \mathrm{span} \{X_1, \ldots, X_d\} \subset TM$ is a horizontal distribution of $M$.\medskip

\begin{prop} \label{prop D^H,f = c grad^Hf}
 Let $D^H$ be a horizontal Dirac operator acting on a horizontal Clifford bundle $S^HM$ with horizontal Clifford action $c^H$ over a closed Carnot manifold $M$. Then, for any function $f \in C^1(M)$, we have
$$\left[D^H, f\right] = c^H(\grad^H f).$$

\B Let us first show that the horizontal Dirac operator fulfills the Leibniz rule. This is just an easy calculation using the properties of the connection and the Clifford action: For any $\sigma \in \Gamma^\infty(S^HM)$ we have
\begin{eqnarray*}
 D^H (f \cdot \sigma) &=& \sum_{j=1}^d c^H\left(X_j\right) \nabla_{X_j}^{S^H} (f \cdot \sigma)\\
&=& \sum_{j=1}^d c^H\left(X_j\right) \left(f \cdot \nabla_{X_j}^{S^H} \sigma + X_j(f) \cdot \sigma \right)\\
&=& f \cdot D^H \sigma + \left(\sum_{j=1}^d c^H\left(X_j\right) X_j(f) \right) \sigma\\
&=& f \cdot D^H \sigma + c^H\left(\sum_{j=1}^d X_j(f) \cdot X_j\right) \sigma\\
&=& f \cdot D^H \sigma + c^H\left(\grad^H f\right) \sigma.
\end{eqnarray*}
Now the statement follows immediately:
$$\left[D^H, f\right] \sigma = D^H(f \cdot \sigma) - f \cdot D^H \sigma = c^H\left(\grad^H f\right) \sigma,$$
and the proposition is proved. \eB
\end{prop}\medskip

Using this proposition, we have to show that the Lip-norm defined by the Connes metric coming from $D^H$ coincides with the supremum of the horizontal gradient of a function $f$. But yet it is not clear how the sub-algebra $\mathcal{A}'$ of $C(M)$ has to look like such that the Connes metric formula \eqref{eq Connes metric formula for D^H} is true. In the classical case, this is exactly the algebra of Lipschitz functions, and it will turn out that in the Carnot case this will be the algebra of functions which are Lipschitz with respect to the Carnot-Carath\'{e}odory metric.\medskip

\begin{defin} \label{def Lip_CC}
  Let $(M,d_{CC})$ be a Carnot manifold, $f \in C(M)$. Then we call the number
$$\mathrm{Lip}_{CC}(f) := \sup \left\{\frac{|f(x)-f(y)|}{d_{CC}(x,y)}: x,y \in M, x \neq y\right\}$$
the \emph{Carnot-Carath\'{e}odory-Lipschitz constant} of $f$. If $\mathrm{Lip}_{CC}(f)$ is finite, we call $f$ a \emph{Carnot-Carath\'{e}odory-Lipschitz function}.\smallskip

We denote the algebra of all Carnot-Carath\'{e}odory-Lipschitz functions on $M$, equipped with the semi-norm $\mathrm{Lip}_{CC}(f)$, by $\mathrm{Lip}_{CC}(M)$. \eBsp
\end{defin}\medskip

The proof that the horizontal Dirac-operator $D^H$ detects the Carnot-Carath\'{e}odory metric works similar to the classical one that a Dirac operator detects the geodesic metric on a Riemannian spin manifold (see e.g. \cite{ConnesNCG}, \cite{GVF} or \cite{Landi}): If we assume that for $f \in \mathrm{Lip}_{CC}(M)$ the horizontal gradient $\mathrm{grad}^H f(x)$ exists almost everywhere, we show that the number $\mathrm{Lip}_{CC}(f)$ coincides with the essential supremum of $\mathrm{grad}^H f$. After that we will show that the Carnot-Carath\'{e}odory metric can be described via the Carnot-Carath\'{e}odory-Lipschitz constant, such that we can apply Proposition \ref{prop D^H,f = c grad^Hf} to get the result.\smallskip

We need to show that the assumption that $\mathrm{grad}^H f(x)$ exists almost everywhere for $f \in \mathrm{Lip}_{CC}(M)$. Note that the analogous statement for classical Lipschitz functions is well known (see e.g. \cite{Fed}). The horizontal case is shown in \cite{CDPT} for the case of the Heisenberg group (see \cite{CDPT}, Proposition 6.12), and this result can be generalized easily to arbitrary Carnot manifolds as we show in the following proposition. In a greater generality, this is also a consequence of the Pansu-Rademacher theorem (\cite{Pan}; see also for example \cite{CDPT}, Theorem 6.4), which states that a Lipschitz map between two Carnot groups, the so-called \emph{Pansu differential} (see \cite{Pan}), exists. \medskip

\begin{prop} \label{prop Lip_CC functions are differentialble almost everywhere}
 Let $M$ be a closed Carnot manifold with horizontal distribution $HM = \spa \{X_1, \ldots, X_d\}$ and let $f: M \rightarrow \R$ be a Carnot-Carath\'{e}odory-Lipschitz function. Then the horizontal gradient
$$\grad^H f = \sum_{j=1}^d X_j(f) \cdot X_j$$
exists almost everywhere on $M$.\smallskip

\B We consider the case where $M = \Omega \subset \G$ is an open subset of a Carnot group $\G$; then the general case follows after restricting ourselves to local coordinates (on which we have a Carnot group structure) since the coordinate changes are smooth and therefore do not affect the regularity of $f$. Thus we choose an arbitrary point $x_0 \in \Omega$. We fix a $j \in \{1, \ldots, d\}$. After an affine change of coordinates on $\Omega$ (such that $0 \in \Omega$) we can assume 
\begin{equation} \label{eq coordinate description Omega_0}
 x_0 \in \left\{ \left((x_{1,1}, \ldots, x_{1,d}, x^{(2)}, \ldots, x^{(R)}\right) \in \Omega: x_{1,j} = 0 \right\} =: \Omega_{0,j},
\end{equation}
where the coordinates on $\Omega \subset \G$ are meant to be either exponential or polarized coordinates, see Section 2.2.\smallskip

The idea is to consider the integral curves arising from $x_0$ into the direction of the vector field $X_j$ via the exponential map, which means we have a curve
\begin{equation} \label{eq horizontal exponential curve from hyperplane}
 \gamma_{x_0}: [0,a] \rightarrow \G, \ \ \ t \mapsto x_0 . \exp{t X_j} 
\end{equation}
for any $a > 0$. We denote by $L_{x_0} := \gamma_{x_0}([0,a]) \cap \Omega$ the path of $\gamma_{x_0}$. Since this path is horizontal, the map
$$f_{x_0}: \R \rightarrow \R, \ \ \ t \mapsto f\left(x_0 \cdot \exp{t X_j} \right)$$
is a (classical) Lipschitz function on $\R$ because of the Carnot-Carath\'{e}odory Lipschitz property of $f$. Since $f_{x_0}$ is Lipschitz, it is differentiable almost everywhere (see e.g. \cite{Fed}), and its derivative is given by
$$f_{x_0}'(t) = X_j f \left(x_0 \cdot \exp{t X_j} \right).$$
But this shows that $X_jf$ exists almost everywhere on $L_{x_0}$. This is true for any starting point $x_0 \in \Omega_{0,j}$, such that we can conclude that $X_j f$ exists almost everywhere on $\Omega$ because any point of a Carnot group $\G$ can be reached from the hyperplane $\Omega_{0,j}$ by a horizontal curve of the type \eqref{eq horizontal exponential curve from hyperplane}. \smallskip

Now the above argument is true for any $j \in \{1, \ldots, d\}$, and hence any horizontal partial derivative $X_j f$ and therefore also the horizontal gradient exists almost everywhere on $\Omega$. \eB
\end{prop}\medskip

To show that the Carnot-Carath\'{e}odory-Lipschitz constant coincides with the supremum norm of the horizontal gradient, we work with the corresponding object from the cotangent bundle of $M$. We assume we have a Riemannian metric $g$ on $M$ such that $\{X_1, \ldots, X_n\}$ is an orthonormal frame with respect to $g$, which respects the grading structure of $TM$. Remember from Section 2.1 that we can choose a basis $\{d\omega^1, \ldots, d\omega^n\}$ of $T^\ast M$ such that 
\begin{equation} \label{eq HM as kernel}
 HM = \spa \{X_1, \ldots, X_d\} = \ker \left( \mathrm{span} \left\{d\omega^{d+1}, \ldots, d\omega^n\right\} \right).
\end{equation}
According to this basis, the horizontal differential of $f$ is given by
$$d^H f = \sum_{j=1}^d X_j(f) d\omega^j.$$
Obviously we have $\sup_{x \in M} \|\mathrm{grad}^H f(x)\| = \sup_{x \in M} \|d^Hf(x)\|$ for any function $f \in C^1(M)$, where $\| \cdot \|$ denotes the (Euclidean) norm of a horizontal tangent (or cotangent) vector on $M$, coming from the Riemannian metric $g$.\medskip

\begin{lemma} \label{lemma Lip_CC = grad^H}
Let $M$ be a closed Carnot manifold with Riemannian metric $g$ as above, and let $f \in \mathrm{Lip}_{CC}(M)$ such that, by Proposition \ref{prop Lip_CC functions are differentialble almost everywhere}, the horizontal gradient $\mathrm{grad}^H f$ of $f$ exists almost everywhere. Then we have
$$\mathrm{Lip}_{CC}(f) = \mathrm{ess} \sup_{x \in M} \left\| \grad^H f(x) \right\|.$$
\smallskip

\B For $f \in \mathrm{Lip}_{CC}(M)$ we show that  $\mathrm{Lip}_{CC}(f) \leq \mathrm{ess} \sup_{x \in M} \left\| \grad^H f(x) \right\|$ and that $\mathrm{ess} \sup_{x \in M} \left\| \grad^H f(x) \right\| \leq \mathrm{Lip}_{CC}(f)$.\smallskip

Let $x,y \in M$, and let $\gamma: [0,1] \rightarrow M$ be a smooth horizontal curve connecting $x$ and $y$, which means $\gamma(0) = x$, $\gamma(1) = y$ and $\dot{\gamma}(t) \in H_{\gamma(t)}M$ for all $t \in [0,1]$. Note that because of \eqref{eq HM as kernel} and the characterization of the horizontal differential this means $df(\dot{\gamma}(t)) = d^Hf(\dot{\gamma}(t))$ for all $t$, and we have for $f \in C^1(M)$
\begin{eqnarray*}
 f(x)-f(y) & = & f(\gamma(1)) - f(\gamma(0))\\
&=& \int_0^1 \frac{d}{dt} f(\gamma(t)) dt\\
&=& \int_0^1 df \left(\dot{\gamma}(t)\right) dt \\
&=& \int_0^1 d^Hf \left(\dot{\gamma}(t)\right) dt\\
&=& \int_0^1 g\left(\mathrm{grad}^H f (\gamma(t)), \dot{\gamma}(t) \right) dt.\\
& \leq & \int_0^1 \left\| \mathrm{\grad}^H f (\gamma(t)) \right\| \cdot \left\| \dot{\gamma}(t) \right\| dt \\
&\leq& \sup_{x \in M} \left\| \grad^H f(x) \right\| \cdot \int_0^1 \left\|\dot{\gamma}(t) \right\| dt,
\end{eqnarray*}
where we have used the Cauchy-Schwarz inequality. Taking the infimum over all horizontal curves connecting $x$ and $y$, we find
\begin{equation} \label{eq estimate for |f(x)-f(y)| using d_CC}
 \left|f(x)-f(y)\right| \leq \sup_{x \in M} \left\| \mathrm{grad}^H f (x) \right\| \cdot d_{CC}(x,y).
\end{equation}
Now the above calculation holds not only for $C^1$-functions $f$, but for all functions which have a horizontal gradient almost everywhere, such that $\grad^H f$ is defined as an essentially bounded vector field on $M$. Hence for every $f$ which fulfills the assumptions of the Lemma, \eqref{eq estimate for |f(x)-f(y)| using d_CC} becomes
$$ \left|f(x)-f(y)\right| \leq \mathrm{ess} \sup_{x \in M} \left\| \mathrm{grad}^H f (x) \right\| \cdot d_{CC}(x,y),$$
and since $x$ and $y$ can be chosen arbitrarily we see
\begin{equation} \label{eq Lip_CC <= sup d^Hf}
  \mathrm{Lip}_{CC}(f) = \sup_{\{x,y \in M: \ x \neq y\}} \frac{\left|f(x)-f(y)\right|}{d_{CC}(x,y)} \leq \mathrm{ess} \sup_{x\in M} \left\|\mathrm{grad}^Hf(x) \right\|.
\end{equation}\smallskip

On the other hand, we choose an $x_0 \in M$ such that $\mathrm{grad}^H f(x_0)$ exists and consider the integral curve along the vector field $\mathrm{grad}^H f (x_0)$ arising from the exponential map $\exp_{x_0}: T_{x_0}M \rightarrow M$ from Riemannian geometry, that is
$$\gamma_{x_0}(t) := x_0 . \exp \left(t \cdot \grad^H f(x_0) \right).$$
For a small $\varepsilon > 0$ we set $x := \gamma_{x_0}(\varepsilon)$ and denote the length of a horizontal curve $\gamma_{x_0}([0,\varepsilon])$ which is connecting $x_0$ and $x$ by $L_{\gamma_{x_0}}(x,x_0)$. We observe that we have
\begin{equation} \label{eq diffquotient horizontal curve}
\begin{split}
  \lim_{x \rightarrow x_0} \left| \frac{f(x)-f(x_0)}{L_{\gamma_{x_0}}(x,x_0)} \right| & = \lim_{\varepsilon \rightarrow 0} \left| \frac{f\left(x_0 . \exp \left(\varepsilon \cdot \grad^H f(x_0) \right)\right) - f(x_0)}{\varepsilon} \right|\\
& = \left\| \grad^H f(x_0) \right\|.
\end{split}
\end{equation}
Now, since $\gamma_{x_0}$ is a horizontal curve connecting $x$ and $x_0$ we have $d_{CC}(x_0,x) \leq L_{\gamma_{x_0}}(x,x_0)$, and hence we get from \eqref{eq diffquotient horizontal curve}
\begin{eqnarray*}
 \left\| \grad^H f(x_0) \right\| &\leq& \lim_{x \rightarrow x_0} \left| \frac{f(x)-f(x_0)}{d_{CC}(x,x_0)} \right|\\
& \leq & \lim_{x \rightarrow x_0} \frac{\mathrm{Lip}_{CC}(f) \cdot d_{CC}(x,x_0)}{d_{CC}(x,x_0)} = \mathrm{Lip}_{CC}(f)
\end{eqnarray*}
by definition of the Carnot-Carath\'{e}odory-Lipschitz constant. But since this works for every $x_0 \in M$ where $\grad^H f (x_0)$ exists (which is almost every $x_0 \in M$), we have 
\begin{equation} \label{eq Lip_CC >= sup d^Hf}
 \mathrm{ess} \sup_{x \in M} \left\| \grad^Hf(x) \right\| \leq \mathrm{Lip}_{CC}(f).
\end{equation}
Altogether \eqref{eq Lip_CC <= sup d^Hf} and \eqref{eq Lip_CC >= sup d^Hf} prove the statement of the lemma. \eB 
\end{lemma}\medskip

To make use of the above Lemma, one has to show that the Carnot-Carath\'{e}odory distance can be expressed using the Carnot-Carath\'{e}odory-Lipschitz constant. This should be obvious, since it is just the Lipschitz semi-norm belonging to the compact metric space $(M, d_{CC})$, but for completeness we write down the proof. From our point of view it is important that the function describing the Carnot-Carath\'{e}odory distance from a fixed point $x_0 \in M$ is a Carnot-Carath\'{e}odory-Lipschitz function, which has Lip-norm bounded by $1$ and is differentiable into horizontal directions almost everywhere. But note that this is not a $C^1$-function.\medskip

\begin{lemma} \label{lemma d_CC from Lip_CC}
  On a closed Carnot manifold $M$ the Carnot-Carath\'{e}odory distance between two points $x,y \in M$ is given by
$$d_{CC}(x,y) = \sup \left\{ |f(x)-f(y)|: f \in \mathrm{Lip}_{CC}(M), \mathrm{Lip}_{CC}(f) \leq 1 \right\}.$$ \smallskip

\B For every $f \in C(M)$ such that $\mathrm{Lip}_{CC}(f) \leq 1$ we have for all $x,y \in M$:
\begin{equation} \label{eq d_CC from Lip_CC <=}
 \left| f(x) - f(y) \right| \leq \mathrm{Lip}_{CC}(f) \cdot d_{CC}(x,y) \leq d_{CC}(x,y).
\end{equation}
On the other hand, we define a function $h: M \rightarrow \R$ via $h(y) := d_{CC}(x,y)$ for a given $x \in M$. Obviously $h$ is continuous, and since $d_{CC}$ is a metric on $M$ we have for any $z \in M$
$$\left|h(y)-h(z)\right| = \left|d_{CC}(x,y) - d_{CC}(x,z)\right| \leq d_{CC}(y,z).$$
This shows $\mathrm{Lip}_{CC}(h) \leq 1$, and since $|h(x)-h(y)| = d_{CC}(x,y)$ we get
\begin{equation} \label{eq d_CC from Lip_CC >=}
 d_{CC}(x,y) \leq \sup \left\{|f(x)-f(y)|: \mathrm{Lip}_{CC}(f) \leq 1 \right\}.
\end{equation}
From \eqref{eq d_CC from Lip_CC <=} and \eqref{eq d_CC from Lip_CC >=} the statement of the lemma follows. \eB \smallskip

\Bem The fact that the function $h$ appearing in the proof is differentiable almost everywhere can also be deduced from the fact that any two points $x,y \in M$ can be joint by a so-called minimizing geodesic, which is a part from $x$ to $y$ realizing the Carnot-Carath\'{e}odory distance; fulfilling the additional property that it has a derivative for almost all $t$ whose components are measurable functions (see \cite{Mon}, Theorem 1.19). We do not need this argument, since the proof above shows that $h$ is a Carnot-Carath\'{e}odory-Lipschitz function which has the property mentioned above by Proposition \ref{prop Lip_CC functions are differentialble almost everywhere}. \eBsp
\end{lemma}\medskip

Now we simply have to put everything together to get the identity of the metrics.\medskip

\begin{thm} \label{thm d_CC from D^H}
  Let $M$ be a closed Carnot manifold and $D^H$ the horizontal Dirac operator acting on a horizontal Clifford bundle $S^HM$ over $M$. Then the Carnot-Carath\'{e}odory metric of $M$ can be detected via the formula
\begin{equation} \label{eq detect ccmetric}
 d_{CC}(x,y) = \sup \left\{ |f(x) - f(y)|: f \in \mathrm{Lip}_{CC}(M), \ \left\| [D^H, f] \right\| \leq 1\right\}.
\end{equation}\smallskip

\B By Lemma \ref{lemma Lip_CC = grad^H} and Lemma \ref{lemma d_CC from Lip_CC} we have
$$d_{CC}(x,y) = \sup \left\{ |f(x)-f(y)|: f \in \mathrm{Lip}_{CC}(M), \mathrm{ess} \sup_{x \in M} \left\| \grad^H f(x) \right\| \leq 1 \right\}.$$
Now Proposition \ref{prop D^H,f = c grad^Hf} tells us that
\begin{equation} \label{eq D^H,f = c grad^Hf in proof}
 [D^H,f] = c^H(\grad^H f)
\end{equation}
for every $f \in C^1(M)$, where $c^H: HM \rightarrow \en_\C(S^HM)$ denotes the horizontal Clifford action, and the norm of the operator $c^H(\grad^H f)$ is given by
$$\left\| c^H\left(\grad^H f\right) \right\| = \sup_{x \in M} \left\| c^H\left(\grad^H f(x)\right) \right\|.$$
But since because of $(c^H)^2 = -\id$ the map $c^H: H_xM \rightarrow S^H_xM$ is an isometry for any $x \in M$, this shows
$$\left\| c^H\left(\grad^H f\right) \right\| = \sup_{x \in M} \left\| \grad^H f(x) \right\|.$$
Hence together with \eqref{eq D^H,f = c grad^Hf in proof} we see that the identity
\begin{equation} \label{eq norm commutator = norm horizontal gradient}
 \sup_{x \in M} \left\| \grad^H f(x) \right\| = \left\| \left[D^H,f\right] \right\|
\end{equation}
is true for any $f \in C^1(M)$. Since for $f \in \mathrm{Lip}_{CC}(M)$ the horizontal gradient $\grad^H f$ exists almost everywhere by Proposition \ref{prop Lip_CC functions are differentialble almost everywhere}, \eqref{eq norm commutator = norm horizontal gradient} implies
$$\mathrm{ess} \sup_{x \in M} \left\| \grad^H f(x) \right\| = \left\| \left[D^H,f\right] \right\| \ \ \ \forall f \in \mathrm{Lip}_{CC}(M),$$
and therefore the theorem is proved. \eB
\end{thm}\medskip

In Theorem \ref{thm d_CC from D^H} we have seen that any horizontal Dirac operator detects the Carnot-Carath\'{e}odory metric via Connes metric formula, where the supremum is taken over the Carnot-Carath\'{e}odory-Lipschitz functions on $M$. We will see now that it suffices to take the supremum over all $C^\infty$-functions, since each $f \in \mathrm{Lip}_{CC}(M)$ can be approximated by functions $f_\varepsilon \in C^\infty(M)$ with smaller $\mathrm{Lip}_{CC}$-norm by a standard approximation argument.\medskip

\begin{cor} \label{cor d_CC from D^H via smooth functions}
   Let $M$ be a closed Carnot manifold and $D^H$ the horizontal Dirac operator acting on a horizontal Clifford bundle $S^HM$ over $M$. Then the Carnot-Carath\'{e}odory metric of $M$ can be detected via the formula
\begin{equation} \label{eq detect ccmetric smooth}
 d_{CC}(x,y) = \sup \left\{ |f(x) - f(y)|: f \in C^\infty(M), \ \left\| [D^H, f] \right\| \leq 1\right\}.
\end{equation}\smallskip

\B Let $\G$ be the tangent Carnot group $\G$ of the Carnot manifold $M$, where like before $R$ is the nilpotency step of $\G$ and for $1 \leq S \leq R$ the number $d_S$ denotes the dimension of the vector space $V_S$ belonging to the grading $\mathfrak{g} = V_1 \oplus \ldots \oplus V_R$ of $\G$. We show that each function $f \in \mathrm{Lip}_{CC}(\G)$ can be approximated uniformly by a sequence of $C^\infty$-functions $f_\varepsilon$ such that $\mathrm{Lip}_{CC}(f_\varepsilon) \leq \mathrm{Lip}_{CC}(f)$ for all $\varepsilon$. Then the statement for $\mathrm{Lip}_{CC}$-functions on the manifold $M$ follows, since $M$ is compact, by considering local charts, and hence the corollary follows immediately from Theorem \ref{thm d_CC from D^H}.\smallskip

Let $f \in \mathrm{Lip}_{CC}(\G)$. We use the Koranyi gauge
$$\left\| x \right\|_{\G} = \left( \sum_{S=1}^R \sum_{j=1}^{d_S} \left|x_{S,j}\right|^\frac{2R!}{S} \right)^\frac{1}{2R!}, $$
see Definition \ref{def Koranyi gauge}, and consider the unit ball
$$B_\G(0,1) = \left\{x \in \G: \left\|x\right\|_\G \leq 1\right\}$$
with respect to this semi-norm. Then we can consider the smooth function
$$u: \R \rightarrow \R, \ \ \ u(t) = \begin{cases} e^{-\frac{1}{t}}, & t>0 \\ 0 & \text{otherwise} \end{cases}.$$
and set 
$$\varphi(x) := c \cdot u \left( 1 - \left\| x \right\|_\G^{2R!} \right),$$
where we choose the constant $c>0$ such that $\int_\G \varphi(x) dx = \int_{B_\G(0,1)} \varphi(x) dx = 1$. But this means that we have for all $\varepsilon > 0$
\begin{equation} \label{eq integral varphi(epsilon x)}
 \int_\G \varphi\left(\delta_{\varepsilon^{-1}}(x) \right) dx = \left(\varepsilon^{\dim_H(\G)}\right)^{-1},
\end{equation}
where $\delta_{\varepsilon^{-1}}$ denotes the weighted dilation on $\G$ by $\varepsilon^{-1}$ (see Definition \ref{def delation on Carnot group}), and 
$$\dim_H(\G) = \sum_{S=1}^R S \cdot \dim V_S$$
is the Hausdorff dimension of $\G$ (see Theorem \ref{thm Mitchells measure theorem}). Note that because of the smoothness of $u$ the function $\varphi$ is a smooth function on $\G$, which is compactly supported in the unit ball with respect to the Koranyi gauge on $\G$.\smallskip

We finally define
\begin{equation} \label{eq smooth family varphi_epsilon}
 \varphi_\varepsilon (x) := \varepsilon^{\dim_H(\G)} \varphi\left(\delta_{\varepsilon^{-1}}(x) \right),
\end{equation}
which provides us a family of functions $\varphi_\varepsilon \in C^\infty_c(\G)$ with $\varphi_\varepsilon \geq 0$ such that  $\int_{\G} \varphi_\varepsilon(x) dx = 1$ for all $\varepsilon > 0$ (because of \eqref{eq integral varphi(epsilon x)}) and $\lim_{\varepsilon \rightarrow 0} \varphi_\varepsilon = \delta_0$ in the sense of distributions. Because of the equivalence of the Koranyi gauge and the Carnot-Carath\'{e}odory metric (see Proposition \ref{prop equivalence Koranyi metric Carnot-Caratheodory metric}) there is a $C > 0$ such that each of these functions $\varphi_\varepsilon$ is supported in $C \cdot B_{CC}(0,\varepsilon)$, where
$$B_{CC}(0,\varepsilon) := \left\{ x \in \G: d_{CC}(0,x) < \varepsilon\right\}$$
denotes the Carnot-Carath\'{e}odory ball with radius $\varepsilon$.\smallskip

Using the compactly supported smooth functions $\varphi_\varepsilon$ from \ref{eq smooth family varphi_epsilon}, we consider
$$f_\varepsilon(x) := f \ast_\G \varphi_\varepsilon (x) := \int_{\G} f(y^{-1}.x) \varphi_\varepsilon(y) dy,$$
where $.$ denotes the composition on $\G$. It is clear from the rules of convolution that we have $f_\varepsilon \in C^\infty(\G)$. Further, since $f$ in continuous (with respect to the Carnot-Carath\'{e}odory metric on $\G$), we know that for any compact subset $K \subset \G$ and any $\delta > 0$ there is an $\varepsilon' > 0$ such that $|f(x)-f(z)| < \delta$ for all $x \in K$ and for all $z \in C \cdot B_{CC}(x,\varepsilon')$, where
$$B_{CC}(x,\varepsilon') := \left\{ y \in \G: d_{CC}(x,y) < \varepsilon'\right\}$$
and the constant $C>0$ is chosen such that $\varphi_{\varepsilon'}$ is supported in $C \cdot B_{CC}(0,\varepsilon')$. But this shows that for all $\varepsilon < \varepsilon'$ we have because of the translation invariance of the Carnot-Carath\'{e}odory metric with respect to the composition on $\G$ (see \cite{CDPT}) and the properties of $\varphi_\varepsilon$
\begin{eqnarray*}
 \left|f(x) -f_\varepsilon(x)\right| &=& \left| C \cdot \int_{B_{CC}(0,\varepsilon)} \left( f(x) - f(y^{-1}.x) \right) \varphi_\varepsilon(y) dy \right| \\
&\leq& \int_{C \cdot B_{CC}(0,\varepsilon)} \left| f(x) - f(y^{-1}.x) \right| \varphi_\varepsilon(y) dy\\
&<& \delta \int_{C \cdot B_{CC}(0,\varepsilon)} \varphi_\varepsilon(y) dy\\
&=& \delta,
\end{eqnarray*}
for all $x \in K$. Since $M$ is a compact manifold (hence we can restrict ourselves to compact subsets in any chart), this shows that any $\mathrm{Lip}_{CC}$-function $f$ on $M$ can be approximated by these $C^\infty$-functions $f_\varepsilon$ in the supremum norm.\smallskip

Finally, we show that for any $\varepsilon > 0$ we have $\mathrm{Lip}_{CC}(f_\varepsilon) \leq \mathrm{Lip}_{CC}(f)$. This follows once again because the Carnot-Carath\'{e}odory metric is translation invariant with respect to the composition on the Carnot group $\G$ (see e.g. \cite{CDPT}). From the definition of the Carnot-Carath\'{e}odory-Lipschitz constant this invariance in connection with the properties of $\varphi_\varepsilon$ leads to the estimate
\begin{eqnarray*}
 \mathrm{Lip}_{CC} (f_\varepsilon) &=& \sup_{x_1 \neq x_2} \left| \frac{f_\varepsilon(x_1) - f_\varepsilon(x_2)}{d_{CC}(x_1,x_2)} \right|\\
&=& \sup_{x_1 \neq x_2} \left| \int_\G \frac{f(y^{-1} . x_1)\varphi_\varepsilon(y) - f(y^{-1}.x_2) \varphi_\varepsilon(y)}{d_{CC}(x_1,x_2)} dy \right|\\
&\leq& \sup_{x_1 \neq x_2} \int_\G \left| \frac{f(y^{-1} . x_1) - f(x_2.y^{-1})}{d_{CC}(y^{-1}.x_1,y^{-1}.x_2)} \right| \varphi_\varepsilon(y) dy\\
&\leq& \sup_{x_1 \neq x_2} \int_\G \mathrm{Lip}_{CC}(f) \varphi_\varepsilon(y) dy \\
&=& \mathrm{Lip}_{CC}(f).
\end{eqnarray*}
Altogether we have proved that any function $f \in \mathrm{Lip}_{CC}(\R^n)$ can be approximated uniformly by a sequence of $C^\infty$-functions $f_\varepsilon$ such that $\mathrm{Lip}_{CC}(f_\varepsilon) \leq \mathrm{Lip}_{CC}(f)$ for all $\varepsilon$, and the statement of the corollary follows. \eB
\end{cor}\medskip

Since $\mathrm{Lip}_{CC}(M)$ (or $C^\infty(M)$) is a dense sub-algebra of $C(M)$, Proposition \ref{prop D^H,f = c grad^Hf} shows that $[D^H,f]$ is bounded for a dense sub-algebra of $C(M)$ and therefore $D^H$ fulfills condition (i) for a spectral triple (see Definition \ref{def spectral triple}. But as we have already mentioned (and will see in the next chapter for a concrete class of examples and in Chapter 6 in general), $D^H$ fails to have a compact resolvent, and therefore $(C(M), L^2(S^HM), D^H)$ is not a spectral triple. But on the other hand, Theorem \ref{thm d_CC from D^H} and the preceding lemmas suggest that $D^H$ seems to be the logical candidate for a first order differential operator to detect the Carnot-Carath\'{e}odory metric.\smallskip

We will now show that at least $(C(M), L^2(\Sigma M), D^H)$ gives rise to a compact quantum metric space in the sense of Mark Rieffel (see Definition \ref{def compact quantum metric space}) if we consider the real-valued functions of $C(M)$ as an order-unit space.\medskip

\begin{cor} \label{cor D^H provides compact quantum metric space}
 The pair $(A,L)$, where
$$A := \left\{f \in C(M): f^\ast = f \right\} \ \ \ \text{and} \ \ \ L(f) := \left\| [D^H,f] \right\|,$$
is a compact quantum metric space which detects the Carnot-Carath\'{e}odory metric on a closed Carnot manifold $M$.\smallskip

\B It is obvious that $A$ is an order-unit space with norm
$$\|f\| := \sup_{x \in M} |f(x)|.$$
Now for calculating the Connes metric from $L$ it suffices to consider only the self-adjoint elements of $C(M)$ (see Proposition \ref{prop metric from positive elements}, which is exactly the space $A$. We have to show that $L$ is a Lip-norm on $A$, i.e.
\begin{enumerate}[(i)]
 \item For every $f \in A$ we have $L(f) = 0$ $\Leftrightarrow$ $f \in \R \cdot 1$.
 \item The topology on the state space $\mathcal{S}(A)$ of $A$ from the Connes metric defined by $L$ is the $w^\ast$-topology.
\end{enumerate}
The non-trivial part of condition (i) follows from Proposition \ref{prop D^H,f = c grad^Hf}: Since 
$$[D^H,f] = c^H(\left( \grad^H f \right),$$
$[D^H,f] = 0$ implies $\grad^H f = 0$ almost everywhere, which implies $X_j(f) = 0$ for all $j \in \{1, \ldots, d\}$ almost everywhere. But this also means $X_k(f) = 0$ almost everywhere for every $k \in \{d+1, \ldots, n\}$ because every vector field of the $X_k$'s can be written as a commutator consisting of the $X_j$'s. Since $\{X_1, \ldots, X_n\}$ spans the tangent space of $M$ and $f$ is continuous by assumption, this implies that $f$ must have been a constant. Therefore (i) is proved.\smallskip

For condition (ii), note that $\mathcal{S}(A) = \mathcal{S}(C(M)) \cong M$ by Gelfand-Naimark theory, where the $w^\ast$-topology on $M$ is exactly the usual manifold topology. Now, by the sub-Riemannian theorem on topologies (see Theorem \ref{thm theorem on topologies}), this topology coincides with the topology induced by the Carnot-Carath\'{e}odory metric $d_{CC}$. But by Theorem \ref{thm d_CC from D^H}, 
$$d_{CC}(x,y) = \rho_L(x,y) := \sup \left\{ |f(x)-f(y)|: L(f) \leq 1 \right\},$$
which shows condition (ii). \smallskip

Hence $(A,L)$ is a compact quantum metric space, and the statement about the metric is just the statement of Theorem \ref{thm d_CC from D^H}. \eB \smallskip

\Bem Note that the Carnot-Carath\'{e}odory-Lipschitz constant $L$ also provides a Lip-norm on $A$. As we have seen in Lemma \ref{lemma d_CC from Lip_CC}, this is exactly the Lip-norm which belongs to the Carnot-Carath\'{e}odory metric $d_{CC}$. In this context, Corollary \ref{cor D^H provides compact quantum metric space} shows that the compact quantum metric spaces $\left( A,\mathrm{Lip}_{CC}(\cdot )\right)$ and $\left(A, \left\|[D^H, \cdot ]\right\| \right)$ are identical. \eBsp
\end{cor}\bigskip

\chapter{Degenerate Spectral Triples on Nilmanifolds}

In the previous chapter we constructed horizontal Dirac operators more or less analogous to classical Dirac operators and we saw that they detect the Carnot-Carath\'{e}odory metric. Therefore they seem to be a natural candidate to construct a spectral triple which covers the horizontal geometry of a Carnot manifold.\smallskip

In this chapter we will do a concrete construction of horizontal Dirac operators $D^H$ on local homogeneous spaces of Carnot groups which arise from the action of a lattice subgroup, namely compact Carnot nilmanifolds $M = \Gamma \backslash \G$. This can be seen as a generalization of the torus in the non-abelian case. We will make use of the spin structures of the horizontal torus, arising as the image of a submersion from $M$, and observe that we obtain a horizontal Clifford structure via pullback where the representation of the horizontal Clifford algebra is irreducible. Afterwards we will use an approach developed by Christian Bär and Bernd Ammann (see \cite{Bae} and \cite{AB}) to get information about the spectrum of $D^H$. It will follow that our horizontal Dirac operator has an infinite dimensional eigenspace. In particular its resolvent in not compact and thus it does not furnish a spectral triple. The strategy is to decompose the $L^2$-space of horizontal Clifford sections for the case where the horizontal distribution has co-dimension $1$: In this case we can use well-known results from the representation theory of Heisenberg groups to calculate the eigenvalues of $D^H$. Using this, we will be able to deduce the statement for the general case. These calculations on Carnot nilmanifolds can be seen as an example of a more general statement: In later chapters, we will use this idea of transferring the problem to the co-dimension $1$ case to show that on any Carnot manifold a horizontal Dirac operator is not hypoelliptic.\smallskip

Despite this lack, as an additional result we are still able to extract the Hausdorff dimension of $(M,d_{CC})$ from the asymptotic behavior of the non-degenerate eigenvalues of $D^H$ in the case where $\G = \Hei^{2m+1}$ is a Heisenberg group.\smallskip

We will use the notion of a compact Carnot nilmanifold introduced in Example \ref{ex compact Carnot nilmanifolds}. Throughout this chapter, we assume that $\Gamma$ is the standard lattice of a Carnot group $\G$. We denote the resulting local homogeneous space by $M = \Gamma \backslash \G$.\bigskip

\section{The Pullback Construction}

Let $\G$ be a Carnot group with horizontal distribution of rank $d$ and nilpotency step $R$, and let $\Gamma \vartriangleleft \G$ be its standard lattice given by
$$  \Gamma := \left\langle \ \{\gamma_j = \exp(X_j): 1 \leq j \leq d \} \ \right\rangle_{\G},$$
see Example \ref{ex compact Carnot nilmanifolds}. We consider the local homogeneous space $M = \Gamma \backslash \G$, where the action of $\Gamma$ on $\G$ is given by the group composition from the left. We equip $M$ with a Riemannian metric $g$ such that the vector bundles $V_1M \oplus \ldots \oplus V_RM$, forming the grading of $TM$, are pairwise orthogonal at each point.\smallskip

To construct an example for a horizontal Clifford connection on $M$, we start by calculating the Christoffel symbols for the horizontal connection on $M$ arising from the Levi-Civita connection (see Section 3.1). \medskip

\begin{prop} \label{prop Christoffel symbols homogeneous space}
Let $\nabla^H$ be the horizontal connection arising from the Levi-Civita connection $\nabla$ on the compact Carnot nilmanifold $M = \Gamma \backslash \G$. If $\{X_1, \ldots, X_d\}$ is an orthonormal frame for the (bracket generating) horizontal distribution $HM = V_1M$, then we have
$$\Gamma_{jk}^l = 0$$
for all Christoffel symbols $\Gamma_{jk}^l$ of $\nabla^H$, $j,k,l \in \{1, \ldots, d\}$, corresponding to this basis.\smallskip

In addition, if we extend the horizontal frame $\{X_1, \ldots, X_d\}$ to an orthonormal tangent frame $\{X_1, \ldots, X_n\}$ of $M$ which respects the grading of $TM$, then all Christoffel symbols of $\nabla$ belonging to this frame satisfy 
$$\Gamma_{jj}^l = 0$$
for $j = d+1, \ldots, n$ and $l = 1, \ldots, d$.\smallskip

\B By the construction of the horizontal connection the horizontal Christoffel symbols are exactly the horizontal Christoffel symbols of the Levi-Civita connection on $M$, see Proposition \ref{prop horizontal connection in local coordinates}. If $g$ is a Riemannian metric on $M$ such that $\{X_1, \ldots, X_d\}$ forms an orthonormal frame at every point, they can be calculated using the properties of the Levi-Civita connection. We have locally $\nabla_{X_j}X_K = \sum_{l=1}^n \Gamma_{jk}^l X_l$ (with $n = \dim M)$, and hence we get by the Koszul formula for the Levi-Civita connection because of the orthonormality of the frame $\{X_1, \ldots, X_d\}$
\begin{equation}\label{eq Koszul formula}
\begin{split}
 \Gamma_{jk}^l = & g(\nabla_{X_j}X_k,X_l)\\
  = & \frac{1}{2} \left( \partial_{X_j} g(X_k,X_l) + \partial_{X_k} g(X_l,X_j) - \partial_{X_l} g(X_j,X_k) \right) \\
   &+\frac{1}{2} \left(-g(X_k,[X_j,X_l]) - g(X_l,[X_k,X_j]) + g(X_j,[X_k,X_l])\right)\\
= & \frac{1}{2} \left(-g(X_k,[X_j,X_l]) - g(X_l,[X_k,X_j]) + g(X_j,[X_k,X_l])\right)
\end{split}
\end{equation}
for all $1 \leq j,k,l \leq d$. Since for any choice of vector fields $X,Y,Z \in V_1M$ the vector fields $X \in V_1M$ and $[Y,Z] \in V_2M$ are orthogonal with respect to $g$, the right hand side of \eqref{eq Koszul formula} vanishes. Therefore all the horizontal Christoffel symbols are $0$.\smallskip

The fact that $\Gamma_{jj}^l = 0$ for all $d+1 \leq j \leq n$ and for all $1 \leq l \leq d$ can also be seen immediately from \eqref{eq Koszul formula}: If $X_j \in V_SM$ for $S \geq 2$, we have $[X_j,X_j] = 0$ and $[X_j,X_l] \in V_{S+1}M$ for $l = 1, \ldots, d$ (which means $X_l \in V_1M$) by the grading structure of $TM$. But this means that we also have $g(X_j,[X_j,X_l]) = 0$ since $X_j \perp V_{S+1}$. Hence every term on the right hand side of \eqref{eq Koszul formula} is $0$, and the additional statement of the proposition follows. \eB \smallskip

\Bem Using formula \eqref{eq Koszul formula}, we can also calculate all the other Christoffel symbols belonging to the Levi-Civita connection of $M$  as soon as we know the commutator relations of the Lie algebra generated by the vector fields $\{X_1, \ldots, X_d\}$. \eBsp
\end{prop}\medskip

We now construct an irreducible horizontal Clifford bundle over $M$ only involving the horizontal distribution of $M$. The idea is to consider the submersion $\psi: \G \rightarrow \R^d$ from Section 2.4, which is given in exponential coordinates via
$$\psi: \G \rightarrow \R^d, \ \ \ \psi\left(x^{(1)}, \ldots, x^{(R)}\right) = x^{(1)}.$$
We have seen in Corollary \ref{cor projection Carnot nilmanifolds} that $\psi$ can be lifted to the nilmanifold given by the action of $\Gamma$ on $\G$. This means we have a submersion
\begin{equation}\label{eq submersion homogeneous Carnot space onto torus}
 \pi: M \rightarrow \T^d \cong \Z^d \backslash \R^d
\end{equation}
of $M$ onto the $d$-dimensional torus, which (by Corollary \ref{cor projection Carnot nilmanifolds}) coincides locally with the submersion $\psi$ of the Carnot groups. If $\{\tilde{X}_1, \ldots, \tilde{X_d}\}$ is a local frame for $T\T^d$ such that we have $\psi(\exp_\G(X_j)) = \exp_{\tilde{\R}^d}(\tilde{X}_j)$ locally on the Carnot groups, then the differential of the submersion $\pi$ applied to our orthonormal frame $\{X_1, \ldots, X_d, \ldots, X_n\}$ is given by
\begin{equation}\label{eq differential submersion homogeneous Carnot space onto torus}
 D\pi: TM \rightarrow T\T^d, \ \ \ D\pi(X_j) = \begin{cases}
              \tilde{X_j}, & 1 \leq j \leq d\\
              0 & \text{otherwise}
             \end{cases}.
\end{equation}
We choose the Riemannian metric $g^{\T^d}$ on $\T^d$ such that $\{\tilde{X}_1, \ldots, \tilde{X}_d\}$ forms an orthonormal frame. \smallskip

The idea for our construction of a horizontal Dirac operator is now to exploit the fact that $\T^d$ is a spin manifold. The starting point for our constructions is the following theorem, which summarizes the well known facts about the spin structures on $\T^d$ and the realization of their corresponding spinor bundles.\medskip

\begin{thm} \label{thm spin structures on T^d}
There are $2^d$ different spin structures $\Sigma^{\T^d}_\delta$ on $\T^d$, which are in one-to-one correspondence to the group homomorphisms
\begin{equation} \label{eq spin structures torus}
 \varepsilon: \Z^d \rightarrow \Z / 2\Z.
\end{equation}
They are indexed by 
$$\delta = (\delta_1, \ldots, \delta_d) \in \left( \Z / 2\Z \right)^d,$$
where $\delta_j$ is the image of the generator $e_j$ of $\Z^d$ under $\varepsilon$.\smallskip

Furthermore, for each spin structure the spinor bundle $\Sigma_\delta \T^d$ is a (complex) vector bundle of rank $2^{[d/2]}$, where $[ \cdot ]$ denotes the Gaussian bracket, and the space of sections of $\Gamma^\infty(\Sigma_\delta \T^d)$ can be identified with functions $\tilde{\sigma} \in C^\infty(\R^d, \C^{2^{[d/2]}})$ such that
\begin{equation} \label{eq periodicity spinors on torus}
 \tilde{\sigma}(a + x) = \varepsilon(a) \tilde{\sigma}(x) \ \ \ \text{for all} \ x \in \R^d, a \in \Z^d.
\end{equation} \smallskip

\B The spin structures of a connected Riemannian manifold $M$ (if they exist) are characterized by group homomorphisms from the fundamental group of $M$ to $\Z / 2\Z$ (see e.g. \cite{LawMich}, Theorem II.2.1). In the case of the $d$-dimensional torus $\T^d$, the fundamental group is $\Z^d$, and since a group homomorphism $\varepsilon: \Z^d \rightarrow \Z / 2\Z$ is uniquely determined by the images $\delta_j = \varepsilon(e_j)$ of the $d$ generators $e_j$ of $\Z^d$, there are $2^d$ possibilities for such a homomorphism.\smallskip

The second statement of the theorem follows from the construction of the spinor bundle corresponding to a spin structure (see e.g. \cite{LawMich} or \cite{Roe}).\eB\smallskip

\Bem Note that \eqref{eq periodicity spinors on torus} is also true for $\sigma \in L^2(\Sigma_\delta \T^d$ to be an $L^2$-spinor, since $\Gamma^\infty(\Sigma_\delta \T^d)$ is dense in $L^2(\Sigma_\delta \T^d)$.
\end{thm}\medskip

For a given spin structure $\Sigma^{\T^d}_\delta$ on the torus $\T^d$, we have a corresponding spinor bundle $\Sigma_\delta \T^d$, equipped with a spinor connection $\nabla^{\Sigma_\delta \T^d}$. Since $\pi$ is a submersion, there exist unique pullbacks of these objects on $M$ (see e.g. \cite{AH}). In detail, we have the following:

\begin{itemize}
 \item The sections of the pullback $\pi^\ast \Sigma_\delta \T^d$ of the spinor bundle have the form
\begin{equation} \label{eq charakterization pullback spinors}
  \pi^\ast \Sigma_\delta \T^d = \left\{\sum_{j=1}^{2^{[d/2]}} f_j \pi^\ast\varphi_j: \ f_j \in C^\infty(M), \ \varphi_1, \ldots, \varphi_{2^{[d/2]}} \ \text{basis sections of} \ \Sigma_\delta \T^d\right\}.
\end{equation}

\item If $\langle \cdot,\cdot \rangle_{\Sigma_\delta \T^d}$ is the bundle metric on $\Sigma_\delta \T^d$, a metric on $\pi^\ast \Sigma_\delta \T^d$ is given via the pullback 
\begin{equation} \label{eq horizontal spinor metric from pullback}
 \left\langle \pi^\ast \varphi_1, \pi^\ast \varphi_2 \right\rangle_{\pi^\ast \Sigma_\delta \T^d} := \left\langle \varphi_1, \varphi_2 \right\rangle_{\Sigma_\delta \T^d}.
\end{equation}

 \item The pullback $\pi^\ast \nabla^{\Sigma_\delta \T^d}$ of the spinor connection to $\pi^\ast \Sigma_\delta \T^d$ has the form
\begin{equation} \label{eq horizontal connection from pullback}
 \pi^\ast \nabla^{\Sigma_\delta \T^d}_{X_j} \left(\pi^\ast \varphi \right) = \pi^\ast \left(\nabla^{\Sigma_\delta \T^d}_{DX_j} \varphi\right) = 
\begin{cases}
 \pi^\ast \left(\nabla^{\Sigma_\delta \T^d}_{\tilde{X}_j} \varphi\right) & \text{for } 1 \leq j \leq d\\
 0 & \text{otherwise}
\end{cases}
\end{equation}
 on the pull-backs of sections of $\Sigma_\delta \T^d$. For an arbitrary element of $\Gamma^\infty(\pi^\ast \Sigma_\delta \T^d)$ described by \eqref{eq charakterization pullback spinors} $\pi^\ast \nabla^{\Sigma_\delta \T^d}$ is defined by using the linear and tensorial behavior of a connection.
\end{itemize}

Considering the Clifford action $c^\T: T\T^d \rightarrow \en_\C(\Sigma_\delta \T^d)$ on $\Sigma_\delta \T^d$, we can use the identification $HM \cong \pi^\ast T\T^d$ given by the differential $D\pi$ from \eqref{eq differential submersion homogeneous Carnot space onto torus} to define a horizontal Clifford action on $\pi^\ast \Sigma_\delta \T^d$. This is done by the pull-back of the endomorphism bundle: For $X \in HM$ we define
\begin{equation*} 
 c^H (X) := \pi^\ast \left(c^{\T^d}(D\pi(X))\right),
\end{equation*}
which is an endomorphism on $\pi^\ast \Sigma_\delta \T^d$ since $c^{\T^d}$ is an endomorphism on $\Sigma_\delta \T^d$. In detail, point-wise we have by the definition of the pull-back of an endomorphism
\begin{equation} \label{eq horizontal Clifford action from pullback}
c^H (X) (\pi^\ast \varphi) = \pi^\ast \left( c^{\T^d}(D\pi(X)) \right) (\pi^\ast \varphi) = \pi^\ast \left( c^{\T^d}(D\pi(X)) \varphi \right)
\end{equation}
for any basis section $\varphi$ of $\Sigma_\delta \T^d$, which extends to the whole bundle via linearity. From the identification $HM \cong \pi^\ast T\T^d$ via $D\pi$ we can conclude that the restriction of the Riemannian metric $g$ on $M$ to $HM$ is exactly the pull-back of the Riemannian metric $g^{\T^d}$ on $\T^d$, which was chosen such that $\{D\pi(X_1), \ldots, D\pi(X_d)\}$ forms an orthonormal frame for $T\T^d$. \smallskip

By writing
$$\Sigma_\delta^H M := \pi^\ast \Sigma_\delta \T^d \ \ \ \text{and} \ \ \ \nabla^{\Sigma_\delta^H} := \pi^\ast \nabla^{\Sigma_\delta \T^d},$$
we will now show that this structure indeed gives a horizontal Clifford bundle. Furthermore, we can write down the horizontal Dirac operator for this bundle. \medskip

\begin{thm} \label{thm horizontal Cilfford bundle from submersion}
 Let $M = \Gamma \backslash \G$ be the nilmanifold of a Carnot group $\G$, and let $\Sigma_\delta^H M$, $\nabla^{\Sigma_\delta^H}$ and $c^H$ be as above. Then $\Sigma_\delta^H M$ equipped with the connection $\nabla^{\Sigma_\delta^H}$ and the horizontal Clifford multiplication $c^H$ is a horizontal Clifford bundle over $M$.\smallskip

$L^2$-sections $\sigma \in L^2(\Sigma_\delta^H M)$ can be identified with $\C^{2^{[d/2]}}$-valued functions $\sigma$ on $\G$ such that 
\begin{equation} \label{eq periodicity pullback spinors}
 \sigma\left( \left(a^{(1)}, \ldots, a^{(R)}\right).x\right) = \varepsilon\left(a^{(1)}\right) \sigma(x)
\end{equation}
for all $a \in \Gamma$ and $x \in \G$, where $\varepsilon$ is the group homomorphism \eqref{eq spin structures torus} describing the spin structure of $\T^d$ and $a.x$ denotes the group operation on $\G$.\smallskip

The horizontal Dirac operator $D^H$ acting on $\Gamma^\infty(\Sigma_\delta^H M)$ is given in local coordinates by
\begin{equation} \label{eq horizontal Dirac on homogeneous space}
 D^H \sigma = \sum_{j=1}^d c^H(X_j) \partial_{X_j},
\end{equation}
 where $\partial_{X_j}$ denotes the partial derivative belonging to a local coordinate chart of $M$.\smallskip

\B First of all, it is clear that the action of $HM$ on $\Sigma^H_\delta M$ via $c^H$ furnishes a Clifford module over $M$: Since $c^{\T^d}$ furnishes a Clifford module structure on $\Sigma^{\T^d}_\delta M$, we have for any $X \in HM$
\begin{eqnarray*}
 \left(c^H (X) \right)^2 (\pi^\ast \varphi) &=& \pi^\ast \left( c^{\T^d}(D\pi(X)) \right)^2 (\varphi)\\
&=& - g^{\T^d}(D\pi(X),D\pi(X)) \pi^\ast \varphi\\
& =& -g(X,X) \pi^\ast \varphi
\end{eqnarray*}
for the basis sections $\pi^\ast \varphi$ of $\Sigma^H_\delta$ by \eqref{eq horizontal Clifford action from pullback}. Note thereby that the restriction of the Riemannian metric $g$ on $M$ to the horizontal distribution $HM$ is exactly the pull-back of the Riemannian metric $g^{\T^d}$ we chose above.\smallskip

Analogously, the condition (i) and (ii) of Definition \ref{def horizontal clifford bundle} follow because the (classical) Clifford bundle $\Sigma_\delta \T^d$ on the torus fulfills these conditions: For the metric \eqref{eq horizontal spinor metric from pullback} on $\Sigma^H_\delta$ we have point-wise for any basis sections $\pi^\ast \varphi_1$, $\pi^\ast \varphi_2$ and for any horizontal vector field $X \in HM$
\begin{eqnarray*}
 \left\langle c^H(X) \pi^\ast \varphi_1, \pi^\ast \varphi_2 \right\rangle_{\Sigma_\delta^H M} &=& \left\langle c^{\T^d}(D\pi(X)) \varphi_1, \varphi_2 \right\rangle_{\Sigma_\delta \T^d}\\
&=& \left\langle \varphi_1,  c^{\T^d}(D\pi(X)) \varphi_2 \right\rangle_{\Sigma_\delta \T^d}\\
&=& \left\langle \pi^\ast \varphi_1, c^H(X) \pi^\ast \varphi_2 \right\rangle_{\Sigma_\delta^H M}
\end{eqnarray*}
by \eqref{eq horizontal spinor metric from pullback} and \eqref{eq horizontal Clifford action from pullback}, and hence the metric compatibility (i) follows.\smallskip

To show the compatibility of $\nabla^{\Sigma_\delta^H}$ with the horizontal connection $\nabla^H$, we first calculate this for any basis section $\pi^\ast \varphi$ of $\Sigma^H_\delta M$. For any $X,Y \in HM$ we have 
\begin{eqnarray*}
 \nabla^{\Sigma^H_\delta}_X (c^H(Y) \pi^\ast \varphi) &=& \nabla^{\Sigma^H_\delta}_X (\pi^\ast c^{\T^d}(D\pi(Y)) \varphi) \\
&=& \pi^\ast \nabla^{\Sigma_\delta \T^d}_{D\pi(X)} (c^{\T^d}(D\pi(Y)) \varphi)\\
&=& \pi^\ast \left(c^{\T^d} \left(\nabla^{\T^d}_{D\pi(X)} D\pi(Y) \right) \varphi + c^{\T^d}(D\pi(Y)) \nabla^{\Sigma_\delta \T^d}_{D\pi(X)} \varphi \right)\\
&=& c^H \left(\nabla^H_X Y\right) \pi^\ast \varphi + \pi^\ast \left( c^{\T^d}(D\pi(Y)) \nabla^{\Sigma_\delta \T^d}_{D\pi(X)} \varphi \right)\\
&=& c^H \left(\nabla^H_X Y\right) \pi^\ast \varphi + c^H(Y) \pi^\ast \left( \nabla^{\Sigma_\delta \T^d}_{D\pi(X)} \varphi \right)\\
&=& c^H \left(\nabla^H_X Y\right) \pi^\ast \varphi + c^H(Y) \nabla^{\Sigma^H_\delta}_X \pi^\ast \varphi
\end{eqnarray*}
by \eqref{eq horizontal connection from pullback}, \eqref{eq horizontal Clifford action from pullback} and since the compatibility condition is true on the spinor bundle $\Sigma_\delta \T^d$, using the Levi-Civita connection $\nabla^{\T^d}$ over $\T^d$. Thereby the fourth equation needs a bit more explanation. The crucial point is that for all $X,Y \in HM$ we have
\begin{equation} \label{eq identification nabla^H with pullback nabla^Td}
 \nabla^H_X Y = \pi^\ast \nabla_{D\pi(X)} D\pi(Y)
\end{equation}
under the identification $HM \cong \pi^\ast T\T^d$ via $D\pi$ from \eqref{eq differential submersion homogeneous Carnot space onto torus}: If we write both sides of \eqref{eq identification nabla^H with pullback nabla^Td} in local coordinates, we see that they coincide since the Christoffel symbols of $\T^d$ with respect to frame $\{D\pi(X_1), \ldots, D\pi(X_d)\}$ vanish as well as the Christoffel symbols of the horizontal Levi-Civita connection on $M$ with respect to the frame $\{X_1, \ldots, X_d\}$ do. But this means by \eqref{eq horizontal Clifford action from pullback} 
$$\pi^\ast \left( c^{\T^d} \left(\nabla^{\T^d}_{D\pi(X)} D\pi(Y)\right) \varphi \right) = c^H \left(\nabla^H_X Y\right) \pi^\ast \varphi $$
for all $\varphi \in \Gamma^\infty(\Sigma_\delta \T^d)$
by linearity of $c^{\T^d}$ and $c^H$ and by \eqref{eq horizontal Clifford action from pullback}. \smallskip

We still need to show the compatibility with the connection for an arbitrary element $\sigma \in \Gamma^\infty(\Sigma^H_\delta M)$. By \eqref{eq charakterization pullback spinors}, such a $\sigma$ can be written in the form
\begin{equation} \label{eq horizontal spinors aslinear combination}
 \sigma(x) = \sum_{j=1}^d f_j(x) \pi^\ast \varphi_j(x),
\end{equation}
where $\varphi_1, \ldots, \varphi_d$ is a local spinor basis for $\Sigma_\delta \T^d$ and $f_j \in C^\infty(M)$ for $1 \leq j \leq d$. But then the general compatibility follows from the above calculation and by the linearity and the tensorial behavior of the connection $\nabla^{\Sigma^H}$, since for each $j$ we have
\begin{eqnarray*}
\nabla^{\Sigma^H_\delta}_X (c^H(Y) f_j \pi^\ast \varphi_j) &=& X(f_j) c^H(Y) \pi^\ast \varphi_j + f_j \nabla^{\Sigma^H_\delta}_X c^H(Y) \pi^\ast \varphi_j \\
&=&  X(f_j) c^H(Y) \pi^\ast \varphi_j + f_j c^H(\nabla^H_X Y) \pi^\ast \varphi_j + f_j c^H(Y) \nabla_X^{\Sigma^H_\delta} \pi^\ast \varphi_j\\
&=& f_j c^H(\nabla^H_X Y) \pi^\ast \varphi_j + c^H(Y) \nabla_X^{\Sigma^H_\delta} f_j \pi^\ast \varphi_j\\
\end{eqnarray*}
Altogether we have shown that $\Sigma^H_\delta M$ is a horizontal Clifford bundle over $M$.\smallskip

Turning to the representation \eqref{eq periodicity pullback spinors} of our horizontal spinors, we use \eqref{eq horizontal spinors aslinear combination} to describe a general section $\sigma \in \Gamma^\infty(\Sigma_\delta^H M)$. Since $M = \Gamma \backslash \G$ is a local homogeneous space, each $f_j \in C^\infty(M)$ can be considered as a function $f \in C^\infty(\G)$ such that $f(a.x) = f(x)$ for all $a \in \Gamma$ and $x \in \G$. And since the submersion $\pi$ is a Lie group homomorphism according to the group operations on $\G$ and $\R^d$, we have for each $\pi^\ast \varphi_j$
\begin{eqnarray*}
 \pi^\ast \varphi_j (a.x) &=& \varphi_j (\pi(a.x)) = \varphi_j (a^{(1)} + x^{(1)}) = \varepsilon\left(a^{(1)}\right) \varphi_j (x^{(1)})\\
& = & \varepsilon\left(a^{(1)}\right) \pi^\ast \varphi_j(x)
\end{eqnarray*}
by Equation \eqref{eq periodicity spinors on torus} from Theorem \ref{thm spin structures on T^d}. From this, \eqref{eq periodicity pullback spinors} follows for sections $\sigma \in \Gamma^\infty(\Sigma^H_\delta M)$, and since $\Gamma^\infty(\Sigma^H_\delta M)$ is dense in the Hilbert space $L^2(\Sigma^H_\delta M)$ of $L^2$-sections of this vector bundle, the second statement of the theorem is proved.\smallskip

Finally we prove the local description \eqref{eq horizontal Dirac on homogeneous space} of the horizontal Dirac operator from this Clifford module. Let $\tilde{\Gamma}_{jk}^l$ denote the Christoffel symbols of $T\T^d$ according to the frame $\{ \tilde{X}_1, \ldots, \tilde{X}_d\} = \{ D\pi(X_1), \ldots, D\pi(X_d)\}$, which are all zero. Then the formula \eqref{eq horizontal Dirac on homogeneous space} follows, because it is known that we have locally
\begin{equation} \label{eq spinor derivative on torus}
 \nabla^{\Sigma_\delta \T^d}_{\tilde{X}_j} \sigma = \partial_{\tilde{X}_j} \sigma - \frac{1}{4} \sum_{k,l=1}^d \tilde{\Gamma}_{jk}^l c^{\T^d}(\tilde{X}_k)c^{\T^d}(\tilde{X}_l) \sigma = \partial_{\tilde{X}_j} \sigma
\end{equation}
for the spinor connection in the torus (see e.g. \cite{Bae}, Lemma 4.1): Using once again the representation \eqref{eq horizontal spinors aslinear combination} for elements of $\Gamma(\Sigma_\delta^H M)$, we get from the definition \eqref{eq horizontal connection from pullback} of $\nabla^{\Sigma^H_\delta}$ for every $j \in \{1, \ldots, d\}$:
\begin{eqnarray*}
 \nabla_{X_j}^{\Sigma_\delta^H} \left( \sum_{j=1}^d f_j \pi^\ast \varphi\right) &=& \sum_{j=1}^d \left( X_j(f) \cdot \pi^\ast \varphi + f_j \cdot \pi^\ast \left(\nabla_{\tilde{X}_j}^{\Sigma_\delta \T^d} \varphi_j \right) \right)\\
&=& \sum_{j=1}^d \left( X_j(f) \cdot \pi^\ast \varphi + f_j \cdot \pi^\ast \left(\partial_{\tilde{X}_j} \varphi_j \right) \right)\\
&=& \sum_{j=1}^d \left( X_j(f) \cdot \pi^\ast \varphi + f_j \cdot \partial_{X_j} \pi^\ast \varphi_j \right)\\
&=& \partial_{X_j} \left( \sum_{j=1}^d f_j \pi^\ast \varphi\right).
\end{eqnarray*}

Now by Definition \ref{def horizontal dirac operator}, the horizontal Dirac operator is given locally via
\begin{equation} \label{eq horizontal Dirac operator local in proof} 
 D^H = \sum_{j=1}^d c^H(X_j) \nabla^{\Sigma^H_\delta}_{X_j} - \frac{1}{2} c^H\left( \sum_{j=d+1}^n \pi^H \nabla_{X_j} X_j \right),
\end{equation}
where $\nabla$ is the Levi-Civita connection on $TM$ and $\pi^H$ denotes the orthonormal projection of $TM$ onto the horizontal distribution $HM$. For the first term on the right hand side of \eqref{eq horizontal Dirac operator local in proof} we can plug in the result from the above calculation, and for the second term we can use the additional statement of Proposition \ref{prop Christoffel symbols homogeneous space} which states that all the Christoffel symbols $\Gamma_{jj}^l$ for $j = d+1, \ldots, n$ and $l = 1, \ldots, d$ are zero: By definition of the Levi-Civita connection this implies that we have locally
$$ \pi^H \nabla_{X_j} X_j = \pi^H \left( \sum_{l=1}^n \Gamma_{jj}^l X_l \right) = \sum_{l=1}^d \Gamma_{jj}^l X_l = 0 $$
for any $j \in \{d+1, \ldots, n\}$. Altogether, the local expression \eqref{eq horizontal Dirac on homogeneous space} for the horizontal Dirac operator we constructed follows. \eB
\end{thm}\medskip

\begin{defin}
 Let $M = \Gamma \backslash \G$ be a compact Carnot nilmanifold, equipped with the objects $\Sigma_\delta^H M$, $\nabla^{\Sigma_\delta^H}$, $c^H$ and $D^H$ from Theorem \ref{thm horizontal Cilfford bundle from submersion}. Then we call 
$$D^H: \Gamma^\infty( \Sigma_\delta^H M) \rightarrow \Gamma^\infty(\Sigma_\delta^H M)$$
the \emph{horizontal pull-back Dirac operator} on $M$. Further, we call $\Sigma_\delta^H M$ the \emph{horizontal spinor bundle} and $\nabla^{\Sigma_\delta^H}$ the \emph{horizontal spinor connection} on $M$. \eBsp
\end{defin}\medskip

Let us summarize what we have done in this section: We have constructed a horizontal Clifford bundle together with a horizontal Dirac operator on an arbitrary compact Carnot nilmanifold $M = \Gamma \backslash \G$ only depend on the horizontal distribution of this manifold. For dimensional reasons this representation of the (bundle of) Clifford algebras over $HM$ is irreducible, since our Clifford bundle $\Sigma^HM$ has rank $2^{[d/2]}$ whenever $d$ is the rank of $HM$. Therefore we can claim we have constructed a natural candidate for a horizontal Dirac operator on $M$.\smallskip

We remark that one can also construct horizontal Dirac operators from already existing (classical) Clifford or spinor bundles of $M$. For this we take a look at Proposition \ref{prop horizontal Clifford bundle from spinor bundle} and note that the horizontal connection in this case depends on the Clifford action of the whole tangent bundle $TM$. But nonetheless, whenever we choose an orthonormal frame $\{X_1, \ldots, X_n\}$ of $TM$ we will get a local expression similar to \eqref{eq horizontal Dirac on homogeneous space} for the resulting horizontal Dirac operator, which only differs by a matrix term arising from the Clifford action. Therefore, methods similar to the ones described in the following sections can also be applied to this situation, and we expect similar results.\bigskip

\section{Spectral Decomposition from the Center}

Our aim is to show that the horizontal pull-back Dirac operator constructed in the preceding section, which detects the Carnot-Carath\'{e}odory metric by Section 3.3, does not have a compact resolvent. This will be the case if we find an eigenvalue of $D^H$ which possesses an infinite dimensional eigenspace. Therefore we are interested in getting information about the spectrum of this operator.\smallskip

We intend to use the local expression \eqref{eq horizontal Dirac on homogeneous space} of the horizontal pull-back operator $D^H$ we constructed in the last section. From this expression, we will be able to use techniques involving the representation theory of its underlying Carnot group: The following proposition shows that $D^H$ can be expressed using the right regular representation of $\G$ (see also \cite{Bae} for the case of the classical Dirac operator).\medskip

\begin{prop} \label{prop expression of horizontal Dirac via representation}
 Let $M = \Gamma \backslash \G$ be a compact Carnot nilmanifold, and let $D^H$ be the horizontal pull-back Dirac operator defined on the horizontal spinor bundle $\Sigma^H_\delta M$, arising from a spin structure $\Sigma^{\T^d}_\delta$ on $\T^d$. \smallskip

We denote by $R: \G \rightarrow L^2(\Sigma_\delta^H M)$ the right regular representation of the Carnot group $\G$ on the Hilbert space $L^2(\Sigma_\delta^H M)$, which is defined by
\begin{equation} \label{eq right regular representation}
 \left(R(x_0) \sigma\right) (x) := \sigma \left( x.x_0 \right)
\end{equation}
for all $x_0 \in \G$. Then $D^H$ can be expressed locally using $R$ via
\begin{equation} \label{eq D^H via right regular representation}
 D^H \sigma (x) = \sum_{j=1}^d c^H(X_j) \left. \frac{d}{dt} \left(R(\exp tX_j) \sigma\right) (x) \right|_{t=0}
\end{equation}
for any $\sigma \in \Gamma^\infty(\Sigma^H_\delta M)$.\smallskip

As a consequence, a closed subspace $\mathcal{H}' \subset L^2(\Sigma_\delta^H M)$ is invariant under $D^H$ if it is invariant under $R$ and under the Clifford action of $HM$.\smallskip

\Bem The expression \eqref{eq right regular representation} is well-defined since elements of $L^2(\Sigma_\delta^H M)$ can be viewed as periodic functions on $\G$ by Theorem \ref{thm horizontal Cilfford bundle from submersion}. \eBsp \smallskip

\Bem For an element $X \in \mathfrak{g}$, where $\mathfrak{g}$ is the Lie algebra of the Carnot group $\G$, we can define
$$R_\ast(X) \sigma (x) := \left. \frac{d}{dt} \left(R(\exp tX) \sigma\right) (x) \right|_{t=0}.$$
This is the so-called \emph{right-regular representation} of the Lie algebra $\mathfrak{g}$ on $L^2(\Sigma^H_\delta M)$ adopted from $R$. Using this representation $R_\ast$, we can rewrite \eqref{eq D^H via right regular representation} in the form
$$D^H = \sum_{j=1}^d R_\ast(X_j) \otimes c^H(X_j).$$
\eBsp \smallskip

\B The local expression of the horizontal pull-back Dirac operator follows immediately from the expression of the directional differentiation along a vector field $X$ on a Lie group, which is
$$\partial_X \sigma (x) = \left. \frac{d}{dt} \sigma \left(x.\exp(tX)\right) \right|_{t=0}.$$
The statement about the invariance is straight forward. \eB
\end{prop}\medskip

The idea, which has been used by Christian Bär and Bernd Amman for the classical Dirac operator on compact nilmanifolds of Heisenberg groups and which we will adopt to our situation, is now to find a direct sum decomposition of the horizontal Clifford bundle $L^2(\Sigma_\delta^H M)$, which is invariant under the horizontal Dirac operator. Thus the determination of the spectrum of $D^H$ splits into parts which are easier to handle. Remember that by \eqref{eq periodicity pullback spinors} from Theorem \ref{thm horizontal Cilfford bundle from submersion} we can consider elements of $L^2(\Sigma_\delta^H M)$ as functions on $\G$ which have certain periodicity properties. We intend to use these properties to decompose $L^2(\Sigma_\delta^H M)$. In the language of representation theory and in view of Proposition \ref{prop expression of horizontal Dirac via representation} this is just the decomposition of the unitary right-regular representation $R$, acting on the Hilbert space $L^2(\Sigma_\delta^H M)$, into its irreducible components.\smallskip

Let $\mathfrak{g} = V_1 \oplus \ldots \oplus V_R$ be the grading of the Lie algebra belonging to $\G$. We start with the decomposition of $L^2(\Sigma_\delta^H M)$ using the periodicities arising from the subgroup $\exp_\G(V_R)$ of $\G$.\medskip

\begin{defin} \label{def Carnot center}
 For the grading $\mathfrak{g} = V_1 \oplus \ldots \oplus V_R$ of the Carnot group $\G = \exp_\G \mathfrak{g}$, we call 
$$Z(\G) := \exp_\G (V_R) \subset \G$$
the \emph{Carnot center} of $\G$. \eBsp \smallskip

\Bem Note that the center of a general Lie group $G$ is the maximal (normal) subgroup of $G$ commuting with every element of $G$. In the case $\G$ is a Carnot group, the center of $\G$ needs not to coincide with the Carnot center of $\G$ defined above. Consider for example the Carnot group $\G = \Hei^{2m+1} \times \R^{d-2m}$ for $d>2m$, where the components belonging to $\R^{d-2m}$ belong to the center but not to the Carnot center of $\G$. \eBsp
\end{defin}\medskip

Now the decomposition of $L^2(\Sigma^H_\delta M)$ from the Carnot center works as follows. In addition, we will detect that the horizontal spinor space of a compact Carnot nilmanifold of lower step can be found in this decomposition.\medskip

\begin{thm} \label{thm decomposition from center}
 Let $M = \Gamma \backslash \G$ be the compact nilmanifold of a Carnot group $\G$ with center $Z(\G) \cong \R^{d_R}$, and let $\Sigma^H_\delta M$ be the horizontal Clifford bundle arising via pull-back from a spin structure $\Sigma^H_\delta$ on $\T^d$.\smallskip

Then there is a decomposition
\begin{equation} \label{eq first decomposition general case}
 L^2 \left(\Sigma_\delta^H M\right) = \bigoplus_{\tau \in \Z^{d_R}} \mathcal{H}_\tau
\end{equation}
into Hilbert spaces $\mathcal{H}_\tau$, which are invariant under the pull-back horizontal Dirac operator $D^H$ acting on $\Sigma^H_\delta M$. The elements $\sigma \in \mathcal{H}_\tau$ are exactly those elements $\sigma \in L^2(\Sigma_\delta^H M)$ fulfilling
\begin{equation} \label{eq characterization elements of H_tau}
 \sigma \left(x^{(1)}, \ldots, x^{(R)} \right) = e^{2\pi i \langle \tau, x^{(R)} \rangle} \cdot \sigma\left(x^{(1)}, \ldots, x^{(R-1)}, 0\right).
\end{equation}
In addition, the space $\mathcal{H}_0$ is isomorphic to the space $L^2(\Sigma_\delta^H \tilde{M})$, where 
$$\tilde{M} = \left( \Gamma / Z(\Gamma)\right) \left\backslash \left( \G / Z(\G) \right) \right.$$
is the compact nilmanifold of the Carnot group $\tilde{\G} := \G / Z(\G)$ of step $R-1$. Under this isomorphism, the restriction of the pull-back horizontal Dirac operator $D^H$ to $\mathcal{H}_0$ can be identified with the pull-back horizontal Dirac operator $\tilde{D}^H$ acting on $\Sigma^H_\delta \tilde{M}$.\smallskip

\B Throughout this proof we will use the characterization of a horizontal spinor $\sigma \in L^2(\Sigma^H_\delta M)$ as a map $\sigma: \G \rightarrow \C^{2^{[d/2]}}$ such that
\begin{equation} \label{eq periodicity pullback spinors in proof}
 \sigma\left( \left(a^{(1)}, \ldots, a^{(R)}\right).x\right) = \varepsilon\left(a^{(1)}\right) \sigma(x)
\end{equation}
for every $a = (a^{(1)}, \ldots, a^{(R)}) \in \Gamma$ given by Theorem \ref{thm horizontal Cilfford bundle from submersion}. Note that this $\sigma$, considered as a periodic function on $\G$, is an $L^2$-function on any fundamental domain for the action of $\Gamma$ on $\G$ by left translation.\smallskip

We will work with the projection $\pi_{\tilde{\G}}: \G \rightarrow \G / Z(\G)$. The Carnot center of $\G$ can be described via (exponential or polarized) coordinates by
$$Z(\G) = \exp(V_R) = \left\{ \left(0, \ldots, 0, x^{(R)} \right) \in \G: x^{(R)} \in \R^{d_R} \right\} \cong \R^{d_R},$$
where $\mathfrak{g} \cong V_1 \oplus \ldots \oplus V_R$ is the grading of the Lie algebra $\mathfrak{g}$ of $\G$ with $d_R = \dim V_R$. The image $\tilde{\G} \cong \G / Z(\G)$ of $\pi_{\tilde{\G}}$ has the structure of a Carnot group of step $R-1$: This follows from Proposition \ref{prop properties projection lie algebra}, using the inherited group composition from $\G$. We can realize $\tilde{\G}$ the subset
\begin{equation} \label{eq tilde G as subset of G}
 \tilde{\G} = \left\{ \left(x^{(1)}, \ldots, x^{(R)} \right) \in \G: x^{(R)} = 0 \right\}
\end{equation}
of $\G$, which leads to the coordinate expression
\begin{equation} \label{eq projection onto complement of center}
 \pi_{\tilde{\G}}: \G \rightarrow \tilde{\G}, \ \ \ \left(x^{(1)}, \ldots, x^{(R-1)}, x^{(R)}\right) \mapsto \left(x^{(1)}, \ldots, x^{(R-1)}, 0 \right)
\end{equation}
of $\pi_{\tilde{\G}}$. In this description, the group composition on $\tilde{\G}$ is given by executing the group composition on $\G$ first and projecting the result to $\tilde{\G}$ afterwards, i.e.
$$x ._{\tilde{\G}} y := \pi_{\tilde{\G}} \left(x ._\G y \right).$$
For the rest of the proof we will work with this coordinate description.\smallskip

In addition, the Carnot center of the standard lattice $\Gamma \vartriangleleft \G$ of $\G$ is given by 
$$Z(\Gamma) := Z(\G) \cap \Gamma \cong \Z^{d_R}.$$
Hence, and since the generators of the Lie algebra of $\G$ can be identified with the generators of the Lie algebra of $\tilde{\G}$ via $\pi_{\tilde{\G}}$, the discrete Carnot group $\tilde{\Gamma} := \Gamma / Z(\Gamma)$ is the standard lattice of $\tilde{\G} := \G / Z(\G)$. This means that we can form the local homogeneous space $\tilde{M} = \tilde{\Gamma} \backslash \tilde{\G}$, which is a compact Carnot nilmanifold of $\tilde{\G}$. From the projection $\pi_{\tilde{\G}}: \G \rightarrow \G / Z(\G)$ we see that $M$ has the structure of a principle $\T^{d_R}$-bundle over $\tilde{M}$, since $\T^{d_R}$ is the local homogeneous space of $\R^{d_R} \cong Z(\G)$ under the action of $\Z^{d_R} \cong Z(\Gamma)$.\smallskip

After this preparation, we can finally prove the statements of the theorem. We fix a point $y = (y^{(1)}, \ldots, y^{(R)}) \in \G$. Then for any $\sigma \in L^2(\Sigma_\delta^HM)$, we can define a map $\varphi_y$ via
\begin{equation} \label{eq definition phi_y}
 \varphi_y: Z(\G) \cong \R^{d_R} \rightarrow \C^{2^{[d/2]}}, \ \ \ z \mapsto \sigma\left(y ._\G (0, \ldots, 0 ,z) \right).
\end{equation}
As one can calculate using the Baker-Campbell-Hausdorff formula, we have for each $z \in \R^{d_R}$
$$\left(y^{(1)}, \ldots, y^{(R)}\right) ._\G \left(0, \ldots, 0, z\right)  = \left(y^{(1)}, \ldots, y^{(R-1)}, y^{(R)} + z\right)$$
since $(0, \ldots, 0, z)$ lies in the Carnot center of $\G$, which means it commutes with any other element of $\G$. Thus, since $\varepsilon(0) = 1$ for any homomorphism $\varepsilon$ characterizing a spin structure $\Sigma^H_\delta$ on $\T^d$ (see Theorem \ref{thm spin structures on T^d}), \eqref{eq periodicity pullback spinors in proof} provides us the periodicity
\begin{equation} \label{eq periodicity varphi_y}
 \begin{split}
   \varphi_y(z+a) &= \sigma \left(y.(0, \ldots, 0, z).(0, \ldots, 0,a)\right) \\
   &= \varepsilon(0) \sigma \left(y.(0, \ldots, 0, z)\right)\\
   &= \varphi_y(z)
 \end{split}
\end{equation}
for each $a \in \Z^{d_R}$. Furthermore, $\varphi_y$ is an $L^2$-section on $\T^{d_R} \cong \Z^{d_R} \backslash \R^{d_R}$ (or on a corresponding fundamental domain of $\Z^{d_R}$ on $\R^{d_R}$, which is any cube of edge length $1$): This follows from the fact that $M$ is a $\T^{d_R}$-principle bundle over $\tilde{M}$ and hence integration over $M$ is integration over the fiber $\T^{d_R}$ followed by integration over the base space $\tilde{M}$. But since $\sigma$ is an $L^2$ section on $M$, it must also be an $L^2$-section on $\T^{d_{R}}$. Together with the periodicity \eqref{eq periodicity varphi_y} this means that we can develop $\varphi_y$ into a Fourier series
\begin{equation} \label{eq Fourier series phi_y}
 \varphi_y(z) = \sum_{\tau \in \Z^{d_R}} \varphi_\tau(y) \cdot e^{2\pi i \langle \tau,z \rangle}
\end{equation}
with Fourier coefficients $\varphi_\tau$. If we choose $z=0$ in \eqref{eq definition phi_y} and \eqref{eq Fourier series phi_y}, we find
\begin{equation} \label{eq desomposition from center using varphi_tau}
 \sigma(y) = \varphi_y(0) = \sum_{\tau \in \Z^{d_R}} \varphi_\tau(y)
\end{equation}
for any $y \in \G$, which proves the decomposition \eqref{eq first decomposition general case}.\smallskip

To find an explicit description of the elements of the spaces $\mathcal{H}_\tau$, we use the expression of the Fourier coefficients from \eqref{eq Fourier series phi_y} via the integral
\begin{equation} \label{eq Fourier coefficiants phi_y}
 \varphi_\tau(y) = \int_{[0,1]^{d_R}} \sigma\left(y.(0, \ldots, 0, t)\right) e^{-2\pi i \langle \tau,t \rangle} \dbar^{d_R}t.
\end{equation}
Assume $\sigma \in \mathcal{H}_\tau$, which means $\sigma(y) = \varphi_\tau(y)$ by \eqref{eq desomposition from center using varphi_tau}. Let the coordinates of $\G$ belonging to the Carnot center be denoted by $z \in \R^{d_R}$, which means we set $z := y^{(R)}$. Since we have 
$$(y^{(1)}, \ldots, y^{(R-1)}, z) . (0, \ldots, 0, t) = (y^{(1)}, \ldots, y^{(R-1)}, 0) . (0, \ldots, 0, z+t)$$
from the group rule in $\G$, we have the following calculation using the substitution $u := t + z$ in the integral and the periodicity \eqref{eq periodicity pullback spinors in proof} of $\sigma$:
\begin{eqnarray*}
 \varphi_\tau(y) &=& \int_{[0,1]^{d_R}} \sigma\left((y^{(1)}, \ldots, y^{(R-1)}, 0) . (0, \ldots, 0, z+t)\right) e^{-2\pi i \langle \tau,t \rangle} \dbar^{d_R}t\\
&=& e^{2\pi i \langle \tau, z \rangle} \int_{[z_1,z_1+1] \times \ldots \times [z_{d_R}, z_{d_R}+1]} \sigma \left((y^{(1)}, \ldots, y^{(R-1)}, 0).(0, \ldots, 0, u) \right)e^{-2 \pi i \langle \tau,u \rangle} \dbar^{d_R}u \\
&=& e^{2\pi i \langle \tau, z \rangle} \int_{[0,1]^{d_R}} \sigma \left((y^{(1)}, \ldots, y^{(R-1)}, 0).(0, \ldots, 0, u) \right)e^{-2 \pi i \langle \tau,u \rangle} \dbar^{d_R}u\\
&=& e^{2\pi i \langle \tau, z \rangle} \cdot \sigma \left(y^{(1)}, \ldots, y^{(R-1)}, 0\right) ,
\end{eqnarray*}
where we used the identity \eqref{eq Fourier coefficiants phi_y} once again in the last line. But from this we get the description \eqref{eq characterization elements of H_tau} for any element $\sigma \in \mathcal{H}_\tau$.\smallskip

On the other hand, for any $\sigma \in L^2(\Sigma^H_\delta M)$ fulfilling
$$\sigma\left(x^{(1)}, \ldots, x^{(R)} \right) = e^{2\pi i \left\langle \tau, x^{(R)} \right\rangle} \cdot \sigma \left(x^{(1)}, \ldots, x^{(R-1)}, 0\right)$$
we fix $y \in \G$ and write down the Fourier series of the corresponding function $\varphi_y(z)$ from \eqref{eq definition phi_y}. For any $\tau' \in \Z^{d_R}$ we calculate the Fourier coefficients $\varphi_{\tau'}(y)$ from \eqref{eq Fourier coefficiants phi_y}:
\begin{eqnarray*}
 \varphi_{\tau'}(y) &=& \int_{[0,1]^{d_R}} \sigma\left(y.(0, \ldots, 0, t)\right) e^{-2\pi i \langle \tau',t \rangle} \dbar^{d_R}t \\
 &=& \int_{[0,1]^{d_R}} e^{2\pi i \left\langle \tau, y^{(R)} + t \right\rangle} \cdot \sigma \left(y^{(1)}, \ldots, y^{(R-1)}, 0\right) e^{-2\pi i \langle \tau',t \rangle} \dbar^{d_R}t \\
&=& e^{2\pi i \left\langle \tau, y^{(R)} \right\rangle} \cdot \sigma \left(y^{(1)}, \ldots, y^{(R-1)}, 0\right) \int_{[0,1]^{d_R}} e^{2\pi i \left\langle \tau - \tau', t \right\rangle} \dbar^{d_R}t \\
&=& \begin{cases}
      \sigma(y) & \text{for} \ \tau = \tau' \\
      0 & \text{otherwise}
    \end{cases}.
\end{eqnarray*}
This shows immediately that we must have $\sigma \in \mathcal{H}_\tau$.
\smallskip

Using this description, we also see the invariance under $D^H$: By Proposition \ref{prop expression of horizontal Dirac via representation} we have to check the invariance of $\sigma \in \mathcal{H}_\tau$ under the Clifford action of $HM$ and the right-regular representation $R$. While the first invariance is trivial, the invariance under $R$ follows from a small calculation: Let $\sigma \in \mathcal{H}_\tau$. For an $x_0 \in \G$ we set 
$$\tilde{\sigma}(x) := (R(x_0) \sigma) (x) = \sigma(x.x_0).$$
Because of the Baker-Campbell-Hausdorff formula we have in (exponential or polarized) coordinates
\begin{equation} \label{eq BCH formula last component}
(x.x_0)^{(R)} = x^{(R)} + x_0^{(R)} + B(x,x_0),
\end{equation}
where $B(x,x_0)$ is a polynomial in the coordinates of $x^{(1)}, \ldots, x^{(R-1)}$ and $x_0^{(1)}, \ldots, x_0^{(R-1)}$ (see Section 2.2). Then we have because of \eqref{eq characterization elements of H_tau} and \eqref{eq BCH formula last component}
\begin{eqnarray*}
 \tilde{\sigma}(x) &=& e^{2\pi i \left\langle \tau, (x.x_0)^{(R)} \right\rangle} \sigma \left( (x.x_0)^{(1)}, \ldots, (x.x_0)^{(R-1)}, 0 \right) \\
&=&   e^{2\pi i \left\langle \tau, x^{(R)} \right\rangle} \sigma \left( (x.x_0)^{(1)}, \ldots, (x.x_0)^{(R-1)}, x_0^{(R)} + B(x,x_0) \right) \\
&=& e^{2\pi i \left\langle \tau, x^{(R)} \right\rangle} \sigma\left( (x^{(1)}, \ldots, x^{(R-1)}, 0) . (x_0^{(1)}, \ldots, x_0^{(R-1)}, x_0^{(R)}) \right) \\
&=& e^{2\pi i \left\langle \tau, x^{(R)} \right\rangle} \tilde{\sigma}(x^{(1)}, \ldots, x^{(R-1)}, 0),
\end{eqnarray*}
since $B(x,x_0)$ does not depend on the components $x^{(R)}$ and $x_0^{(R)}$. This shows $\tilde{\sigma} \in \mathcal{H}_\tau$ and therefore the invariance of $\mathcal{H}_\tau$ under $R$. \smallskip

The next step is to show that $\mathcal{H}_0 \cong L^2(\Sigma^H_\delta \tilde{M})$: First of all, the horizontal distributions of the compact Carnot nilmanifolds $M = \Gamma \backslash \G $ and $\tilde{M} = (\Gamma / Z(\Gamma)) \backslash (\G / Z(\G))$ can be identified using the Carnot group homomorphism $\pi_{\tilde{\G}}: \G \rightarrow \G / Z(\G)$ from \eqref{eq projection onto complement of center}, which means that for a given spin structure $\Sigma^{\T^d}_\delta$ on the horizontal torus $\T^d$ we find horizontal spinor bundles $\Sigma^H_\delta M$ and $\Sigma^H_\delta \tilde{M}$ on both manifolds. Now for any $\tilde{\sigma} \in L^2(\Sigma^H_\delta \tilde{M})$ we find can define
$$\sigma(x) := \tilde{\sigma}(\pi_{\tilde{\G}}(x)),$$
which is an element of $\mathcal{H}_0$: If we consider $\pi_{\tilde{\G}}(\G)$ as the subset \eqref{eq tilde G as subset of G} of $\G$, we have $\sigma(\pi_{\tilde{\G}}(x)) = \sigma(x)$, which is exactly the characterization of the space $\mathcal{H}_0$ from \eqref{eq characterization elements of H_tau}. And since elements of $\mathcal{H}_0$ are uniquely determined by its values on $\pi_{\tilde{\G}}(\G)$, the map
\begin{equation} \label{eq isomorphism H_0 smaller horizontal spinor bundle}
\varphi: L^2(\Sigma^H_\delta \tilde{M}) \rightarrow \mathcal{H}_0, \ \ \ \tilde{\sigma} \mapsto \sigma = \tilde{\sigma} \circ \pi_{\tilde{\G}}
\end{equation}
is an isomorphism. \smallskip

Finally, we show that the restriction of the pull-back horizontal Dirac operator $D^H$ on $L^2(\Sigma^H_\delta M)$ to $\mathcal{H}_0$ can be identified with the pull-back horizontal Dirac operator $\tilde{D}^H$ on $L^2(\Sigma^H_\delta \tilde{M})$, coming from the same spin structure $\Sigma^{\T^d}_\delta$ on the horizontal torus $\T^d$. Note that by the characterization \eqref{eq tilde G as subset of G} of $\tilde{\G}$ as a subset of $\G$, a horizontal frame $\{X_1, \ldots, X_d\}$ of $HM$ is also a horizontal frame of $H\tilde{M}$, and since the horizontal spinor bundles on $M$ and $\tilde{M}$ are constructed from the same horizontal torus $\T^d$ the horizontal Clifford action from elements of this frame coincides in both cases.\ Then the identification of the operators follows immediately from the isomorphism \eqref{eq isomorphism H_0 smaller horizontal spinor bundle} when we use the local expression \eqref{eq D^H via right regular representation} of the pull-back horizontal Dirac operator from Proposition \ref{prop expression of horizontal Dirac via representation}: We have for every $\sigma \in \mathcal{H}_0$
\begin{eqnarray*}
 D^H \sigma (x) &=& \sum_{j=1}^d c^H(X_d) \left. \frac{d}{dt} \sigma \left(x . \exp_{\G}(tX_j) \right) \right|_{t=0} \\
&=& \sum_{j=1}^d c^H(X_d) \left. \frac{d}{dt} \tilde{\sigma} \left( \pi_{\tilde{\G}} \left(x.\exp_{\G}(tX_j)\right) \right) \right|_{t=0}\\
&=& \sum_{j=1}^d c^H(X_d) \left. \frac{d}{dt} \tilde{\sigma} \left( \pi_{\tilde{\G}} (x) . \exp_{\tilde{\G}}(tX_j) \right) \right|_{t=0}\\
&=& \tilde{D}^H \tilde{\sigma} \left( \pi_{\tilde{\G}} (x) \right),
\end{eqnarray*}
and since $\mathcal{H}_0$ is invariant under $D^H$ this shows $\varphi (\tilde{D}^H \tilde{\sigma}) = D^H \sigma$ for the isomorphism $\varphi$ from \eqref{eq isomorphism H_0 smaller horizontal spinor bundle}.\smallskip

Altogether every statement of the theorem is proved. \eB \smallskip

\Bem Note that the coordinate expressions we used in the proof do not depend on whether we choose exponential or polarized coordinates on the Carnot group $\G$. This is the case since only actions by the Carnot center are involved. In the next section, it will turn our to be more comfortable to work with polarized coordinates. \eBsp
\end{thm} \medskip

The above theorem reduces the problem of showing that the horizontal Dirac operator has an infinite dimensional eigenspace to the corresponding horizontal Dirac operator on a Carnot group of a lower step. We are interested in using this theorem inductively until $\tilde{M}$ is the compact nilmanifold of a Carnot group of step $2$. As we will see in the next section, for a compact Heisenberg nilmanifold $\tilde{M}$ we will be able to do a complete spectral decomposition of $L^2(\tilde{M})$ and detect infinite dimensional eigenspaces from this. Then in Section 4.4, we will use Theorem \ref{thm decomposition from center} to lift these infinite dimensional eigenspace to compact nilmanifolds of general Carnot groups.\bigskip

\section{The Case of Compact Heisenberg Nilmanifolds}

We will now do the complete spectral decomposition of the horizontal pull-back Dirac operator on compact nilmanifolds $M = \Gamma \backslash \G$, where $\G \cong \Hei^{2m+1} \times \R^{d-2m}$ is a Carnot group of step $2$ and rank $d$ for some integer $1 \leq m \leq d/2$. Hence the horizontal distribution is of co-dimension $1$.\smallskip

For our spectral decomposition we will follow an argument of Christian Bär (see Section II.2 of his PhD thesis \cite{Bae}), which was used to determine the spectrum of Dirac operators on nilmanifolds from Heisenberg groups. This method was also used by Christian Bär and Bernd Ammann in \cite{AB} for the case of the Dirac operator on the $3$-dimensional Heisenberg group. Although the situation in case of our horizontal Dirac operator differs a bit, the space $L^2(\Sigma^H_\delta M)$ can be decomposed in exactly the same way since our local expression of $D^H$ from \eqref{eq horizontal Dirac on homogeneous space} is quite similar to the local expression of the Dirac operator Bär is considering. In this thesis, we have to consider a slightly more general case where $\G \cong \Hei^{2m+1} \times \R^{d-2m}$, but the commutative part of the group will not effect the general strategy too much. After doing the decomposition, we will be able to detect infinite-dimensional eigenspaces of $D^H$; and as an additional result we will show that the asymptotic behaviors of the non-degenerate eigenvalues gives back the homogeneous dimension of $M$.\smallskip

We start by introducing some notation. Let the grading of $TM$ be given by $TM = V_1M \oplus V_2M$, where $\{X_1, \ldots, X_d\}$ is an orthonormal frame of $V_1M$ such that the Levi form according to this frame is given by
$$L = \begin{pmatrix} 0 & D & \\ -D & 0 & \\ & & 0_{d-2m} \end{pmatrix}, \ \ \ \text{where} \ D = \begin{pmatrix} \lambda_1 & & \\ & \ddots & \\ & & \lambda_m \end{pmatrix} \ \text{with} \lambda_j > 0. $$ 
Note that this can always be achieved by Proposition \ref{prop tangent Lie group bundle and rank of Levi form} and that the numbers $\lambda_1, \ldots, \lambda_m$ are exactly the absolute values of the non-zero eigenvalues $\pm i\lambda_j$ of the Levi form of $M$. If $V_2M$ is spanned by $X_{d+1}$, this means that we have the commutator relations
\begin{equation} \label{eq commutator relations codimension 1} 
 [X_j,X_k] = \begin{cases}
                \lambda_j X_{d+1} & \text{for} \ 1 \leq j \leq m, \ k = m+j\\
                -\lambda_j X_{d+1} & \text{for} \ 1 \leq k \leq m, \ j = k+m \\
                0 &\text{otherwise}
               \end{cases}.
\end{equation}
In what follows, it will be more comfortable to use the polarized coordinates instead of the exponential coordinates of $\G$ (see Definition \ref{def coordinates of Carnot groups}). We will use polarized coordinates with respect to the frame $\{\tilde{X}_1, \ldots, \tilde{X}_{d+1}\}$ of $TM$, where 
\begin{equation} \label{eq frame tilde X}
 \tilde{X}_j = \frac{1}{\sqrt{\lambda_j}} X_j \ \ \text{and} \ \ \tilde{X}_{m+j} = \frac{1}{\sqrt{\lambda_j}} X_{m+j} \ \ \text{for} \ 1 \leq j \leq m
\end{equation}
and $\tilde{X}_k = X_k$ otherwise. We denote these polarized coordinates on $\G$ by $(x,y,z,t)$ with $x \in \R^m$, $y \in \R^m$, $z \in \R^{d-2m}$ and $t \in \R$, where
$$(x,y,z,t) = \prod_{j=1}^m \exp (x_jX_j) \ . \ \prod_{j=1}^m \exp(y_jX_{m+j}) \ . \ \prod_{k=1}^{d-2m} \exp (z_kX_{2m+k}) \ . \ \exp (tX_{d+1}).$$
Obviously, $\prod$ denotes the product according to the group composition $.$ of $\G$. In other words, the coordinates $(x,y,t)$ describe the $\Hei^{2m+1}$-part of $\G$ and $z$ describes the $\R^{d-2m}$-part of $\G$.\smallskip

Now let us calculate the composition rule on $\G$ in these coordinates. In exponential coordinates corresponding to the frame $\{\tilde{X}_1, \ldots, \tilde{X}_{d+1}\}$ we set for abbreviation
$$V := \sum_{j=1}^m a_j \tilde{X}_j + \sum_{j=1}^{m} b_j \tilde{X}_{m+j} + \sum_{j=1}^{d-2m} c_j \tilde{X}_{2m+j} + s\tilde{X}_{d+1}$$
and
$$W := \sum_{j=1}^m x_j \tilde{X}_j + \sum_{j=1}^{m} y_j \tilde{X}_{m+j} + \sum_{j=1}^{d-2m} z_j \tilde{X}_{2m+j} + t\tilde{X}_{d+1},$$
such that we have because of \eqref{eq commutator relations codimension 1}
\begin{eqnarray*}
[V,W] &=& \sum_{j=1}^m \left( a_jy_j \cdot \frac{1}{\lambda_j} [X_j,X_{m+j}] + b_jx_j \cdot \frac{1}{\lambda_j} [X_{m+j},X_j] \right) \\
&=& \left(\sum_{j=1}^m a_jy_j - b_jx_j \right) \tilde{X}_{d+1}.
\end{eqnarray*}
By the Baker-Campbell-Hausdorff formula (see Section 2.2), we calculate
\begin{eqnarray*}
 (a,b,c,s).(x,y,z,t) &=& \exp V . \exp W \\
&=& \exp \left(V+W + \frac{1}{2}[V,W] \right) \\
&=& \exp \left( \sum_{j=1}^m (a_j+x_j) \tilde{X}_j + \sum_{j=1}^{m} (b_j+y_j) \tilde{X}_{m+j} + \sum_{j=1}^{d-2m} (c_j+z_j) \tilde{X}_{2m+j} \right. \times \\
& & \ \ \ \ \ \ \ \times \left. + \frac{1}{2}\left(s+t + \sum_{k=1}^m (a_ky_k - b_kx_k) \right) {X}_{d+1} \right)\\
&=& \left( a+x, \ b+y, \ c+z, \ s+t + \frac{1}{2}\sum_{k=1}^m (a_ky_k - b_kx_k) \right)
\end{eqnarray*}
for the composition rule on $\G$ in exponential coordinates. We can now use the isomorphism between the exponential and the polarized coordinates on a Heisenberg group (see the remark after Definition \ref{def coordinates of Carnot groups}) and conclude that for the composition rule in coordinates according to the frame $\{\tilde{X}_1, \ldots, \tilde{X}_{d+1}\}$ we have
\begin{equation} \label{eq group composition codimension 1 polarized coordinates}
 (a,b,c,s).(x,y,z,t) = \left(a+x,b+y,c+z,s+t+\sum_{j=1}^m a_jy_{j+m} \right).
\end{equation} \smallskip 

Using this composition rule, we can start with the spectral decomposition of $L^2(\Sigma^H_\delta M)$. The first decomposition of $L^2(\Sigma_\delta^H M)$ is already given by Theorem \ref{thm decomposition from center}. In the situation $\G \cong \Hei^{2m+1} \times \R^{d-2m}$ we have
\begin{equation} \label{eq decomposition from center codimension 1}
 L^2(\Sigma_\delta^H M) \cong \bigoplus_{\tau \in \Z} \mathcal{H}_\tau,
\end{equation}
which is invariant under $D^H$ according to the theorem, where elements of $\mathcal{H}_\tau$ are identified by the relation
\begin{equation} \label{eq characterization elements of H_tau codim 1}
 \sigma(x,y,z,t) = e^{2\pi i \tau t} \cdot \sigma(x,y,z,0).
\end{equation}
We set
\begin{equation} \label{eq definition f_sigma}
 f_\sigma (x,y,z) := \sigma(x,y,z,0),
\end{equation}
considered as a function on $\R^d$. Note that $f_\sigma$ must be an $L^2$-function on the torus $\T^d$ since $\sigma$ is an $L^2$-function on $M$, and that any $\sigma \in \mathcal{H}_\tau$ is fully determined by $f_\sigma$.\smallskip

The strategy is to look for periodicities of $f_\sigma$ in order to find a further decomposition of the spaces $\mathcal{H}_\tau$. These periodicities can be detected from the lemma below.\smallskip

\begin{lemma} \label{lemma periodicity f_sigma}
 Let $\varepsilon: \Z^d \rightarrow \Z / 2\Z$ be the group homomorphism characterizing the spin structure $\Sigma^{\T^d}_\delta$ on $\T^d$ from Theorem \ref{thm spin structures on T^d}, which induces the horizontal spinor bundle $\Sigma^H_\delta M$ on the compact Heisenberg nilmanifold $M$. For $\tau \in \Z$, we assume $\sigma \in \mathcal{H}_\tau$ according to the decomposition \eqref{eq decomposition from center codimension 1} of $L^2(\Sigma^H_\delta M)$. Then for any $c \in \Z^d$ we have
\begin{equation} \label{eq periodicity f_sigma}
 f_\sigma\left((x,y,z)\right) = \varepsilon(c) e^{2\pi i \left\langle \sum_{j=1}^m c_jy_j, \tau \right\rangle} \cdot f_\sigma \left((x,y,z) + c\right).
\end{equation} \smallskip

\B We use the periodicities of $\sigma$ from \eqref{eq periodicity pullback spinors} of Theorem \ref{thm horizontal Cilfford bundle from submersion} and the composition rule \eqref{eq group composition codimension 1 polarized coordinates} to calculate for every $c \in \Z^d$:
\begin{eqnarray*}
 f_\sigma(x,y,z) &=& \varepsilon(c) \cdot \sigma( (c,0). (x,y,z,0) )\\
 &=& \varepsilon(c) \sigma( \left( (x,y,z) + c, \sum_{j=1}^m c_jy_j \right) )\\
&=& \varepsilon(c) e^{2\pi i \tau \left(\sum_{j=1}^m c_jy_j\right)} \cdot \sigma((x,y,z) + c, 0 )\\
&=& \varepsilon(c) e^{2\pi i \tau \left(\sum_{j=1}^m c_jy_j\right)} \cdot f_\sigma( (x,y,z) + c).
\end{eqnarray*}
Thereby, in the third equation we have used the characterization \eqref{eq characterization elements of H_tau codim 1} of $\sigma \in H_\tau$. \eB
\end{lemma}\medskip

From this lemma, we can detect further decompositions of the spaces $H_\tau$. We start with the case $\tau = 0$:\medskip

\begin{lemma} \label{lemma spectral decomposition H_0}
 For the space $\mathcal{H}_0$ from the decomposition \eqref{eq decomposition from center codimension 1} of $L^2(\Sigma^H_\delta M)$ we have the further decomposition
$$\mathcal{H}_0 \cong \bigoplus_{ \left\{ \alpha \in \frac{1}{2}\Z^d: \ e^{2\pi i \alpha_j} = \delta_j \right\}} \mathcal{H}_0^\alpha$$
 with $\mathcal{H}_0^\alpha \cong \C^{2^{[d/2]}}$ for every $\alpha$, where each of the spaces $\mathcal{H}_0^\alpha$ is invariant under the horizontal pull-back Dirac operator $D^H$ on $M$.\smallskip

\B Let $\sigma \in \mathcal{H}_0$, which means we have by \eqref{eq characterization elements of H_tau codim 1} and \eqref{eq definition f_sigma}
$$\sigma(x,y,z,t) = f_\sigma(x,y,z)$$
with $f_\sigma$ as in Lemma \ref{lemma periodicity f_sigma}. If we set $\tau = 0$ in \eqref{eq periodicity f_sigma}, we immediately see that $f_\sigma$ is $2\Z^d$-periodic, and therefore can be developed into a Fourier series
$$f_\sigma(x,y,z) = \sum_{\alpha \in \frac{1}{2} \Z^d} a_\alpha e^{2\pi i \langle \alpha, (x,y,z) \rangle}$$
with $a_\alpha \in \C^{2^{[d/2]}}$. Now the homomorphism $\varepsilon: \Z^d \rightarrow \Z / 2\Z$ from \eqref{eq periodicity f_sigma} is given by 
$$\varepsilon(c) = \delta^{c_1}_1 \cdot \ldots \cdot \delta^{c_d}_d,$$
see Theorem \ref{thm spin structures on T^d}, such that Lemma \ref{lemma periodicity f_sigma} provides us further
\begin{eqnarray*}
 \sum_{\alpha \in \frac{1}{2} \Z^d} a_\alpha e^{2\pi i \langle \alpha, (x,y,z) \rangle} &=& \varepsilon(c) \sum_{\alpha \in \frac{1}{2} \Z^d} a_\alpha e^{2\pi i \langle \alpha, (x,y,z) + c \rangle}\\
&=& \sum_{\alpha \in \frac{1}{2} \Z^d} a_\alpha e^{2\pi i \langle \alpha, (x,y,z) \rangle} \cdot \delta_1^{c_1} \ldots \delta_d^{c_d} e^{2\pi i \langle \alpha, c \rangle}
\end{eqnarray*}
for every $c \in \Z^d$. But this leads to the constraint that for every $j \in \{1, \ldots, d\}$ we must have $\delta_j e^{2\pi i \alpha_j} = 1$ and therefore $e^{2\pi i \alpha_j} = \delta_j$. This means that we have 
$$\alpha_j \in \Z \ \Leftrightarrow \ \delta_j = 1 \ \ \ \text{and} \ \ \ \alpha_j \in \Z + \frac{1}{2} \ \Leftrightarrow \ \delta_j = -1,$$
and so the Fourier series of $f_\sigma$ becomes
\begin{equation} \label{eq fourier series f_sigma on H_0}
 f_{\sigma}(x,y,z) = \sum_{\left\{\alpha \in \frac{1}{2} \Z^d: \ e^{2\pi i\alpha_j} = \delta_j \right\}} a_\alpha e^{2\pi i \langle \alpha, (x,y,z) \rangle},
\end{equation}
with $a_\alpha \in \C^{[d/2]}$. Therefore we have the decomposition
\begin{equation} \label{eq decomposition H_0 codim 1}
 \mathcal{H}_0 \cong \bigoplus_{\left\{\alpha \in \frac{1}{2} \Z^d: \ e^{2\pi i\alpha_j} = \delta_j \right\}} \mathcal{H}_0^\alpha,
\end{equation}
where each $\mathcal{H}_0^\alpha$ is spanned by $ e^{2\pi i \langle \alpha, (x,y,z) \rangle}$ and hence isomorphic to $\C^{2^{[d/2]}}$.\smallskip

Now we immediately see that each $\mathcal{H}_0^\alpha$ is invariant under the Clifford action of any $X_j \in HM$ (which only acts on the coefficient $a_\alpha$) and also under the right regular representation of $\G$: Since $\sigma \in \mathcal{H}_0^\alpha$, we have for every given $(x_0,y_0,z_0,t_0) \in \G$ 
\begin{eqnarray*}
 \sigma\left((x,y,z,0).(x_0,y_0.z_0,t_0)\right) &=& \sigma \left(x+x_0,y+y_0,z+z_0,0\right) \\
&=& a_\alpha e^{2\pi i \langle \alpha, (x_0,y_0,z_0,0) \rangle} \cdot e^{2\pi i \langle \alpha, (x,y,z,0) \rangle},
\end{eqnarray*}
with $a_\alpha \in \C^{[d/2]}$. But this shows that $R(x_0,y_0,z_0,0) \sigma \in \mathcal{H}_0^\alpha$, and altogether every statement of the lemma is proved. \eB \smallskip

\Bem Since $\G / Z(\G) \cong \R^d$, we already know by Theorem \ref{thm decomposition from center} that the space $\mathcal{H}_0$ is isomorphic to the space $L^2(\Sigma_\delta^H \T^d)$, which is by construction exactly the space $L^2(\Sigma_\delta \T^d)$, and the horizontal pull-back Dirac operator acting on $L^2(\Sigma_\delta \T^d)$ is exactly the classical Dirac operator $D^{\T^d}$ on the torus. Since the spectral decomposition of the Dirac operator on the torus is well-known, we could have deduced the statement of this lemma directly from Theorem \ref{thm decomposition from center}. \eBsp
\end{lemma}\medskip

For the cases $\tau \neq 0$ the spectral decomposition is a bit more involved.\medskip

\begin{lemma} \label{lemma spectral decomposition H_tau, tau neq 0}
 Let $\tau \neq 0$. Then for the space $\mathcal{H}_\tau$ from the decomposition \eqref{eq decomposition from center codimension 1} of $L^2(\Sigma^H_\delta M)$ we have the decomposition
\begin{equation} \label{eq decomposition H_tau codim 1}
\mathcal{H}_\tau \cong \bigoplus_{\left\{\gamma \in \frac{1}{2}\Z^{d-2m}: \ e^{2\pi i \gamma_j} = \delta_{2m+j}\right\}} \bigoplus_{J=1}^{|\tau|^m} \mathcal{H}_{\tau, \gamma}^J,
\end{equation}
with $\mathcal{H}_{\tau, \gamma}^J \cong L^2(\R^m, \C^{2^{[d/2]}})$, where each of the spaces $\mathcal{H}_{\tau, \gamma}^J$ is invariant under the horizontal pull-back Dirac operator $D^H$ on $M$.\smallskip

\B Once again we use the identity \eqref{eq periodicity f_sigma} for $f_\sigma$ from Lemma \ref{lemma periodicity f_sigma}, looking for periodicities in the case $\tau \neq 0$. We write $b = (\xi, \beta, \gamma) \in \Z^d$ with $\xi \in \Z^m$, $\beta \in \Z^m$ and $\gamma \in \Z^{d-2m}$ to distinguish the periodicities which belong to $x$, $y$ and $z$. Remember that we have $\G \cong \Hei^{2m+1} \times \R^{d-2m}$ for our Carnot group $\G$, such that the coordinates $x$ and $y$ (together with the coordinate from the center) belong to the Heisenberg part and $z$ forms the (commutative) $\R^{d-2m}$-part of $\G$.\smallskip

From \eqref{eq periodicity f_sigma} we see immediately that for any $\gamma \in \Z^{d-2m}$ we have
$$f_\sigma(x,y,z) = \varepsilon(0,0,\gamma) \cdot f_\sigma(x,y,z+\gamma),$$
with $\varepsilon: \Z^{d} \rightarrow \Z / 2\Z$ describing the spin structure of $\T^d$ from which $D^H$ is constructed. Therefore $f_\sigma$ is $2\Z^{d-2m}$-periodic in $z$ and after fixing $x$ and $y$ we can develop $f_\sigma$ in a Fourier series
\begin{equation} \label{eq periodicity f_sigma in z}
f_\sigma(x,y,z) = \sum_{\gamma \in \frac{1}{2} \Z^{d-2m}} a_\gamma^\sigma(x,y) e^{2\pi i \langle \gamma,z \rangle},
\end{equation}
where $a_\gamma^\sigma(x,y)$ is a $\C^{2^{[d/2]}}$-valued function. We note that $a_\gamma^\sigma$ must be an $L^2$-function on $\T^{2m}$, since $f_\sigma$ is an $L^2$-function on $\T^d \cong \T^{2m} \times \T^{d-2m}$.\smallskip

As we did before we can put further conditions on the indices $\gamma$ over which we sum up in \eqref{eq periodicity f_sigma in z}: For any $c \in \Z^{d-2m}$ we have 
$$\varepsilon(0,0,c) = \delta^{c_1}_{2m+1} \cdot \ldots \cdot \delta^{c_{d-2m}}_d,$$
and thus from Lemma \ref{lemma periodicity f_sigma} we can deduce the restriction
\begin{eqnarray*}
  \sum_{\gamma \in \frac{1}{2} \Z^{d-2m}} a^\sigma_\gamma(x,y) e^{2\pi i \langle \gamma,z \rangle} &=& \varepsilon (0,0,\gamma)  \sum_{\gamma \in \frac{1}{2} \Z^{d-2m}} a^\sigma_\gamma(x,y) e^{2\pi i \langle \gamma,z+c \rangle}\\
&=&  \sum_{\gamma \in \frac{1}{2} \Z^{d-2m}} a_\gamma(x,y)^\sigma e^{2\pi i \langle \gamma,z \rangle} \cdot \delta_{2m+1}^{c_1} \ldots \delta_{d}^{c_{d-2m}} e^{2\pi i \langle \gamma,c \rangle}
\end{eqnarray*}
for every $c \in \Z^{d-2m}$ by \eqref{eq periodicity f_sigma}, which leads to the constraint that for all $j \in \{1, \ldots, d-2m\}$ we must have $\delta_{2m+j} e^{2\pi i \gamma_j} = 1$ and therefore $e^{2\pi i \gamma_j} = \delta_{2m+j}$. But this means that we have 
$$\gamma_j \in \Z \ \Leftrightarrow \ \delta_{2m+j} = 1 \ \ \ \text{and} \ \ \ \gamma_j \in \Z + \frac{1}{2} \ \Leftrightarrow \ \delta_{2m+j} = -1,$$
and so the Fourier series of $f_\sigma$ becomes
\begin{equation} \label{eq fourier series f_sigma in z}
 f_{\sigma}(x,y,z) = \sum_{\left\{\gamma \in \frac{1}{2} \Z^{d-2m}: \ e^{2\pi i\gamma_j} = \delta_{2m+j} \right\}} a_\gamma^\sigma(x,y) e^{2\pi i \langle \gamma,z \rangle}.
\end{equation}
In Equation \eqref{eq periodicity f_sigma} we also see that we have the same periodicity properties for $f_\sigma$ in the variable $y$ (since the $\beta$-components of $b = (\xi,\beta,\gamma) \in \Z^d$ do not appear in the exponent of $e$), so we can repeat the above argument to show that
$$f_\sigma(x,y,z) = \varepsilon(0,\beta,0) f_\sigma(x,y+\beta,z),$$
which gives us
$$a_\gamma^\sigma(x,y) = \varepsilon(0,\beta,0) a_\gamma^\sigma(x,y+\beta)$$
together with \eqref{eq fourier series f_sigma in z}. But therefore we can develop $a_\gamma^\sigma$ into a Fourier series with respect to $y$, and after taking care of the restrictions arising from the components $\delta_{m+1}, \ldots, \delta_{2m}$ of the spin structure similarly to the way we did above we find that
$$a_\gamma^\sigma(x,y) = \sum_{\left\{\beta \in \frac{1}{2} \Z^{m}: \ e^{2\pi i\beta_j} = \delta_{m+j} \right\}} b_{\gamma,\beta}^\sigma(x) e^{2\pi i \langle \beta,y \rangle}$$
for a $\C^{2^{[d/2]}}$-valued function $b_{\gamma,\beta}^\sigma$ on $\R^m$. Plugging this into \eqref{eq fourier series f_sigma in z}, we can write
\begin{equation} \label{eq fourier series f_sigma in z and y}
 f_{\sigma}(x,y,z) = \sum_{\left\{\gamma \in \frac{1}{2} \Z^{d-2m}: \ e^{2\pi i\gamma_j} = \delta_{2m+j} \right\}} \sum_{\left\{\beta \in \frac{1}{2} \Z^{m}: \ e^{2\pi i\beta_j} = \delta_{m+j} \right\}} b_{\gamma,\beta}^\sigma(x) e^{2\pi i \langle \beta,y \rangle} e^{2\pi i \langle \gamma,z \rangle}.
\end{equation}\smallskip

We continue to find periodicities for $f_\sigma$ by considering dilations by $\xi \in \Z^m$ of the $x$ variable. This will lead to a further characterization of the functions $b_{\gamma,\beta}^\sigma$ appearing in \eqref{eq fourier series f_sigma in z and y}. Let $\sigma \in H_\tau$. We use the the composition rule \eqref{eq group composition codimension 1 polarized coordinates} on $\G$, the characterization \eqref{eq characterization elements of H_tau codim 1} of elements belonging to $H_\tau$ and the Fourier series development \eqref{eq fourier series f_sigma in z and y} to calculate for any $\xi \in \Z^m$
\begin{eqnarray*}
 \sigma( (\xi,0,0,0).(x,y,z,t)) &=& \sigma (\left(\xi+x,y,z,t + \sum_{j=1}^m \xi_jy_j \right) )\\
&=& e^{2\pi i \tau \left(t + \sum_{j=1}^m \xi_jy_j \right)} \cdot f_\sigma(\xi+x,y,z)\\
&=& e^{2\pi i \tau \left(t + \langle \xi,y \rangle \right)} \cdot \sum_{\gamma} \sum_{\beta} b_{\gamma,\beta}^\sigma(x+\xi) e^{2\pi i \langle \beta,y \rangle} e^{2\pi i \langle \gamma,z \rangle}\\
&=& e^{2\pi i \tau t} \cdot \sum_{\gamma} \sum_{\beta} b_{\gamma,\beta}^\sigma(x+\xi) e^{2\pi i \langle \beta + \tau \xi,y \rangle} e^{2\pi i \langle \gamma,z \rangle},
\end{eqnarray*}
where the summation over $\gamma$ and $\beta$ is given as in \eqref{eq fourier series f_sigma in z and y}. Now, since $\tau \xi \in \Z^m$ and therefore $e^{2 \pi i \beta} = e^{2 \pi i (\beta + \tau \xi)}$, we can use the substitution $\eta := \beta + \tau \xi$ in the last equation of the above calculation and find
\begin{equation*}
 \sigma( (\xi,0,0,0).(x,y,z,t)) = e^{2\pi i \tau t} \cdot \sum_{\gamma} \sum_{\left\{\eta \in \frac{1}{2}\Z^m: \ e^{2\pi i\eta_j} = \delta_{m+j} \right\}} b_{\gamma, \eta-\tau\xi}^\sigma(x + \xi) e^{2\pi i \langle \eta,y \rangle} e^{2\pi i \langle \gamma,z \rangle}.
\end{equation*}
Because of $\sigma( (\xi,0,0,0).(x,y,z,t)) = \varepsilon(\xi,0,0) \sigma(x,y,z,t)$ for every element $\sigma \in L^2(\Sigma_\delta^H M)$ this means
\begin{equation} \label{eq derivation periodicity b_gamma,beta}
 \varepsilon(\xi,0,0) \sigma(x,y,z,t) = e^{2\pi i \tau t} \cdot \sum_{\gamma} \sum_{\left\{\eta \in \frac{1}{2}\Z^m: \ e^{2\pi i\eta_j} = \delta_{m+j} \right\}} b_{\gamma, \eta-\tau\xi}^\sigma(x + \xi) e^{2\pi i \langle \eta,y \rangle} e^{2\pi i \langle \gamma,z \rangle}.
\end{equation}
On the other hand we have because of \eqref{eq characterization elements of H_tau codim 1} and \eqref{eq fourier series f_sigma in z and y}
\begin{equation} \label{eq characterization elements of H_tau via b_gamma,beta}
\begin{split}
 \sigma(x,y,z,t) &= e^{2\pi i \tau t} f_\sigma(x,y,z)\\
&= e^{2\pi i \tau t} \sum_{\gamma} \sum_{\left\{\beta \in \frac{1}{2} \Z^{m}: \ e^{2\pi i\beta_j} = \delta_{m+j} \right\}} b_{\gamma,\beta}^\sigma(x) e^{2\pi i \langle \beta,y \rangle} e^{2\pi i \langle \gamma,z \rangle},
\end{split}
\end{equation}
and together with \eqref{eq derivation periodicity b_gamma,beta} this leads to the identity $b_{\gamma,\beta}^\sigma(x) = \varepsilon(\xi,0,0) \cdot b_{\gamma, \beta -\tau\xi}^\sigma(x + \xi)$ or, equivalently,
\begin{equation} \label{eq periodicity b_gamma,beta}
 b_{\gamma,\beta}^\sigma(x+\xi) = \varepsilon(\xi,0,0) b_{\gamma,\beta + \tau\xi}^\sigma(x)
\end{equation}
for all $\xi \in \Z^m$. But this means that for a given $\gamma$ we have altogether $|\tau|^m$ independent functions $b_{\gamma,\beta}^\sigma$ which characterize $H_\tau$.\smallskip

Now we can choose $|\tau|^m$ independent functions which determine $\sigma \in \mathcal{H}_\tau$ completely. After fixing a $\beta_1 \in \{\beta \in \frac{1}{2} \Z^{m}: \ e^{2\pi i\beta_j} = \delta_{m+j}\}$, we consider the set
$$B := \left\{ \beta_1 + \sum_{j=1}^m b_j e_j: b_j \in \{0, \ldots, |\tau|-1\} \right\},$$ 
where $e_1, \ldots, e_m$ are the generators of $\Z^m$. Obviously we have $\# B = |\tau|^m$, and after enumerating the elements of $B$ by $\beta_1, \ldots, \beta_{|\tau|^m}$ we have the identity 
$$\left\{\beta \in \frac{1}{2} \Z^{m}: \ e^{2\pi i\beta_j} = \delta_{m+j} \right\} = \dot{\bigcup}_{J=1}^{|\tau|^m} \left\{ \beta_J + \tau \xi: \xi \in \Z^d \right\} $$
for the index set over which the $\beta$'s are summated. Note that the union on the right hand side of this identity is disjoint, 
and after using the property \eqref{eq periodicity b_gamma,beta} we see that for every $\xi \in \Z^m$ we have
$$ b_{\gamma,\beta_J + \tau\xi}^\sigma(x) = \varepsilon(\xi,0,0) b_{\gamma,\beta_J}^\sigma(x+\xi).$$
Plugging this into the expression \eqref{eq characterization elements of H_tau via b_gamma,beta} for elements $\sigma \in \mathcal{H}_\tau$, we see that every $\sigma \in \mathcal{H}_\tau$ is uniquely determined by the functions $b_{\gamma,\beta_1}^\sigma, \ldots, b_{\gamma,\beta_{|\tau|^m}}^\sigma$ via
\begin{equation} \label{eq characterization elements of H_tau via b_gamma,J}
\begin{split}
 \sigma(x,y,z,t) & = e^{2\pi i \tau t} \cdot \sum_\gamma a_\gamma^\sigma (x,y) e^{2\pi i \langle \gamma,z \rangle}\\
& = e^{2\pi i \tau t} \cdot \sum_{\gamma} \sum_{J=1}^{|\tau|^m} \sum_{\xi \in \Z^m} \varepsilon(\xi,0,0) b_{\gamma,\beta_J}^\sigma(x + \xi) e^{2\pi i \langle \beta_J + \tau \xi,y \rangle} e^{2\pi i \langle \gamma,z \rangle}\\
\end{split}
\end{equation}
We further make use of the fact that for any $\gamma$ $a_\gamma^\sigma(x,y)$ is an $L^2$-function on $\T^{2m}$ to derive that $b_{\gamma,\beta_J}^\sigma \in L^2(\R^m, \C^{2^{[d/2]}})$. We can argue using the expression \eqref{eq characterization elements of H_tau via b_gamma,J} as follows: For any $\gamma$ we have
\begin{eqnarray*}
 \int_{[0,1]^{2m}} \left| a_\gamma^\sigma (x,y) \right|^2 dy dx &=& \int_{[0,1]^m} \int_{[0,1]^m} \left| \sum_{J=1}^{|\tau|^m} \sum_{\xi \in \Z^m} \varepsilon(\xi,0,0)  e^{2\pi i \langle \beta_J + \tau \xi,y \rangle} \right|^2 dy dx\\
&=& \int_{[0,1]^m} \int_{[0,1]^m} \sum_{J=1}^{|\tau|^m} \sum_{\xi \in \Z^m} \left| b_{\gamma,\beta_J}^\sigma(x + \xi) \right|^2 dy dx\\
&=& \sum_{J=1}^{|\tau|^m} \int_{[0,1]^m} \sum_{\xi \in \Z^m} \left| b_{\gamma,\beta_J}^\sigma(x + \xi) \right|^2 dx\\
&=& \sum_{J=1}^{|\tau|^m} \int_{\R^m} \left| b_{\gamma,\beta_J}^\sigma(x) \right|^2 dx,
\end{eqnarray*}
and from the last line we can conclude $b_{\gamma,\beta_J}^\sigma \in L^2(\R^m, \C^{2^{[d/2]}})$.\smallskip

Since elements of $\mathcal{H}_\tau$ are fully characterized by the values of
$$\sigma(x,y,z,0) = \sum_{\gamma} \sum_{J=1}^{|\tau|^m} \sum_{\xi \in \Z^m} \varepsilon(\xi,0,0) b_{\gamma,\beta_J}^\sigma(x + \xi) e^{2\pi i \langle \beta_J + \tau \xi,y \rangle} e^{2\pi i \langle \gamma,z \rangle},$$
from the above calculation it furthermore follows that there is an isometry
\begin{equation*}
 H_\tau \cong \bigoplus_{\left\{\gamma \in \frac{1}{2}\Z^{d-2m}, \ e^{2\pi i \gamma_j} = \delta_{2m+j}\right\}} \bigoplus_{J=1}^{|\tau|^m} L^2(\R^m, \C^{2^{[d/2]}}).
\end{equation*}
But this shows the decomposition \eqref{eq decomposition H_tau codim 1} stated in the theorem.\smallskip

To finish the proof, we have to show that this decomposition is invariant under $D^H$. While once again the invariance under the Clifford action is trivial, the invariance under the right regular representation $R$ follows after a simple calculation. Remember that in the polarized coordinates we are using on $\G$, for a given point $\bar{x}_0 = (x_0,y_0,z_0,t_0) \in \G$ the right regular action on $L^2(\Sigma^H_\delta M)$ is given via
$$\left( R(\bar{x}_0) \sigma \right)(x,y,z,t) = \sigma ( \left(x+x_0, y+y_0, z+z_0, t+t_0 + \langle x,y_0 \rangle \right) ),$$
see \eqref{eq composition rule Heisenberg group polarized coordinates}. Now since $\sigma \in \mathcal{H}_\tau$, we can use the characterization \eqref{eq characterization elements of H_tau via b_gamma,J} of $\sigma$ and find that
\begin{eqnarray*}
 \left( R(\bar{x}_0) \sigma \right)(x,y,z,t) &=& e^{2\pi i \tau \left(t+t_0 + \langle x,y_0 \rangle\right)} \cdot \ \times \\
& & \ \times \sum_{\gamma} \sum_{J=1}^{|\tau|} \sum_{\xi \in \Z^m} \varepsilon(\xi,0,0) b_{\gamma,\beta_J}^\sigma(x + x_0 + \xi) e^{2\pi i \langle \beta,y+y_0 \rangle} e^{2\pi i \langle \gamma,z+z_0 \rangle}\\ 
&=& C(x_0,y_0,z_0,t_0) e^{2\pi i \tau t} \cdot \ \times \\
& & \ \times \sum_{\gamma} \sum_{J=1}^{|\tau|} \sum_{\xi \in \Z^m} \varepsilon(\xi,0,0) e^{2\pi i \tau \langle x,y_0 \rangle} b_{\gamma,\beta_J}^\sigma(x + x_0 + \xi) e^{2\pi i \langle \beta,y \rangle} e^{2\pi i \langle \gamma,z \rangle}.
\end{eqnarray*}
But since for $b_{\beta,\gamma_J}^\sigma \in L^2(\R^m, \C^{2^{[d/2]}})$ we also have $e^{2\pi i \tau \langle \cdot,y_0 \rangle} b_{\gamma,\beta_J}^\sigma(\cdot + x_0) \in L^2(\R^m, \C^{2^{[d/2]}})$, this proves the invariance of this decomposition \eqref{eq decomposition H_tau codim 1} under $R$, which means that it is also invariant under $D^H$. Hence the lemma is proved.\eB
\end{lemma}\medskip

The arguments for the decomposition of $L^2(\Sigma^H_\delta M)$ used so far (except for the small modification of including the commutative part $\R^{d-2m}$) are exactly the arguments Bär uses to decompose the (not horizontal) spinor bundle on $M$ in \cite{Bae}. But now we are ready to make some interesting observation for our situation by looking at the eigenvalues resulting from this decomposition. We can calculate all eigenvalues of $D^H$, and first of all we will see that $D^H$ does not have a compact resolvent, and therefore does not provide a spectral triple.\medskip

\begin{thm} \label{thm spectral decomposition D^H on Heisenberg group}
Let $\G \cong \Hei^{2m+1} \times \R^{d-2m}$ be a Carnot group of step $2$, rank $d$ and co-rank $1$ of its horizontal distribution with compact nilmanifold $M = \Gamma \backslash \G $. The non-zero eigenvalues of the Levi form of $\G$ shall be given by the numbers $\pm \lambda_j$ with $\lambda_j > 0$ for $1 \leq j \leq m$. Let $D^H$ be a horizontal pull-back Dirac operator on $M$ arising via pullback from a spin structure $\Sigma^{\T^d}_\delta$ with $\delta = (\delta_1, \ldots, \delta_d) \in (\Z / 2\Z)^d$ on the torus $\T^d$.\smallskip

Then we have the following statements about the spectrum of $D^H$: 
\begin{enumerate}[(a)]
 \item In any case, the spectrum of $D^H$ is discrete. The absolute values of the eigenvalues of $D^H$ are of the types
\begin{enumerate}[(i)]
\item $\left| \mu_{0,\alpha} \right| = 2 \pi \sqrt{\sum_{j=1}^m \lambda_j (\alpha_j^2 + \alpha_{m+j}^2) + \sum_{j=2m+1}^d \alpha_j^2}$ \ \ \ of multiplicity $2^{[d/2]}$ for $\alpha = (\alpha_1, \ldots, \alpha_d) \in  \frac{1}{2} \Z^d$ such that $e^{2\pi i \alpha_j} = \delta_{j}$.
\item $\left| \mu_{\tau,\gamma,\kappa} \right| \sim \sqrt{2 \pi |\tau| \sum_{j=1}^m \lambda_j (2\kappa_j + 1) + 4 \pi^2 \sum_{j=2m+1}^d \gamma_j^2}$ \ \ \ of multiplicity $2^{[d/2]} |\tau|^m$ for $\kappa \in \N^m$, $\tau \in \Z \setminus \{0\}$ and $\gamma \in \frac{1}{2} \Z^{d-2m}$ such that $e^{2\pi i \gamma_j} = \delta_{2m+j}$.
\end{enumerate}
 \item In the case $2m=d$, $D^H$ has an infinite dimensional kernel, and any eigenspace belonging to an eigenvalue $\mu \neq 0$ of $D^H$ is finite dimensional.
 \item In the case $2m<d$, there are infinitely many eigenvalues of $D^H$ which have an infinite dimensional eigenspace.
\end{enumerate}
In particular, since $D^H$ has in any case an infinite-dimensional eigenspace, it does not have a compact resolvent.\smallskip

\B The strategy is to use the decomposition of $L^2(\Sigma_\delta^H M)$  we get from Theorem \ref{thm decomposition from center}, Lemma \ref{lemma spectral decomposition H_0} and Lemma \ref{lemma spectral decomposition H_tau, tau neq 0}, which is invariant under $D^H$. Then we can use Proposition \ref{prop expression of horizontal Dirac via representation} to write $D^H$ locally in the form
\begin{equation} \label{eq D^H via right regular representation codim 1}
 D^H \sigma (x) = \sum_{j=1}^d c^H(X_d) \left. \frac{d}{ds} \left(R(\exp sX_j) \sigma\right) (x) \right|_{s=0}.
\end{equation}
To shorten notation in what follows, we set
\begin{equation} \label{eq tilde lambda}
\tilde{\lambda}_j = \tilde{\lambda}_{m+j} := \sqrt{\lambda_j} \ \ \text{for} \ 1 \leq j \leq m \ \ \ \ \text{and} \ \ \ \ \tilde{\lambda}_j = 1 \ \ \text{for} \ 2m+1 \leq j \leq d+1.
\end{equation}
We will keep on working in the polarized coordinates according to the frame $\{\tilde{X}_1, \ldots, \tilde{X}_{d+1}\}$ from \eqref{eq frame tilde X}, such that the group composition on $\G$ is given by \eqref{eq composition rule Heisenberg group polarized coordinates}. By \eqref{eq tilde lambda} we have $X_j = \tilde{\lambda}_j \tilde{X}_j$ for all $1 \leq j \leq d+1$, and the representation of $D^H$ with respect to the frame $\{\tilde{X}_1, \ldots, \tilde{X}_{d+1}\}$ becomes
\begin{equation} \label{eq D^H via right regular representation codim 1 tilde X}
 D^H \sigma (x) = \sum_{j=1}^d c^H(X_d) \left. \frac{d}{ds} \left(R(\exp s \tilde{\lambda}_j \tilde{X}_j) \sigma\right) (x) \right|_{s=0}.
\end{equation}

Since each of the spaces $\mathcal{H}_0^\alpha$ from Lemma \ref{lemma spectral decomposition H_0} and $\mathcal{H}_{\tau,\gamma}^J$ from Lemma \ref{lemma spectral decomposition H_tau, tau neq 0} is invariant under the right regular representation $R$, this decomposition also provides a decomposition of $R$ into its irreducible parts. For $\tau = 0$ these will be the trivial representations, which can be calculated directly, and for $\tau \neq 0$ we can make use of the fact that the representation theory for the Heisenberg group is fully known to detect the corresponding infinite dimensional representations.\smallskip

For $\sigma \in \mathcal{H}_0^\alpha$, which means by \eqref{eq fourier series f_sigma on H_0}
$$ \sigma(x,y,z,t) = a_\alpha e^{2\pi i \langle \alpha, (x,y,z) \rangle} $$
with $a_\alpha \in \C^{2^{[d/2]}}$, we calculate from \eqref{eq D^H via right regular representation codim 1 tilde X} (using the notation $(x_1, \ldots, x_{d+1})$ instead of $(x,y,z,t)$ for the coordinates of $\G$) by \eqref{eq composition rule Heisenberg group polarized coordinates}
\begin{equation} \label{eq D^H on H_0^alpha}
 \begin{split}
D^H \sigma(x_1, \ldots, x_{d+1}) &= \sum_{j=1}^d c^H(X_j) \left. \frac{d}{ds} a_\alpha e^{2\pi i \left\langle \alpha, \left(x_1, \ldots, x_j+\tilde{\lambda}_js, \ldots, x_d\right) \right\rangle} \right|_{s=0}\\
&= \sum_{j=1}^d c^H(X_j) 2\pi i \alpha_j \tilde{\lambda}_j \sigma(x_1, \ldots, x_{d+1}).
 \end{split}
\end{equation}
This shows that the spectrum of $D^H$ restricted to each $\mathcal{H}_0^\alpha$ consists of the eigenvalues of the matrices $\sum_{j=1}^d  2\pi i \alpha_j \tilde{\lambda}_j c^H(X_j)$. To calculate the absolute values of these eigenvalues, we calculate the eigenvalues of $(D^H)^2$ restricted to $\mathcal{H}_0^\alpha$:
\begin{eqnarray*}
 \left(D^H\right)^2 \sigma &=& \left( \sum_{j=1}^d c^H(X_j) 2\pi i \alpha_j \tilde{\lambda}_j \right)^2 \sigma \\
&=& \sum_{j=1}^d \sum_{k=1}^d - 4 \pi^2 \alpha_j \alpha_k \tilde{\lambda}_j \tilde{\lambda}_k c^H(X_j)c^H(X_k) \sigma \\
&=& \sum_{j=1}^d 4 \pi^2 \alpha_j^2 \tilde{\lambda}_j^2 \sigma,
\end{eqnarray*}
since $c^H(X_j) c^H(X_k) + c^H(X_k) c^H(X_j) = 0$ for $j \neq k$ and $(c^H(X_j))^2 = -\id$ by the rules of the Clifford action. But this shows that on each $\mathcal{H}_0^\alpha$ the absolute values of the eigenvalues of $D^H$ are given by
\begin{equation} \label{eq eigenvalues absolute on H_tau^alpha}
 \left| \mu_{0,\alpha} \right| = 2 \pi \sqrt{\sum_{j=1}^d \tilde{\lambda}_j^2 \alpha_j^2} = 2\pi \sqrt{\sum_{j=1}^m \lambda_j \left(\alpha_j^2 + \alpha_{m+j}^2\right) + \sum_{j=2m+1}^d \alpha_j^2},
\end{equation}
each one with multiplicity $2^{[d/2]}$ because $\Sigma^H_\delta M$ is a vector bundle of rank $2^{[d/2]}$. Note that this also shows that the part of the spectrum belonging to $\mathcal{H}_0$ coincides with the spectrum of the Dirac operator $D^{\T^d}$ on the torus $\T^d$ (equipped with a Riemannian metric such that $\{\tilde{X}_1, \ldots, \tilde{X}_d\}$ is orthonormal), from which $D^H$ was constructed. This can also be deduced from Theorem \ref{thm decomposition from center}.\smallskip

To determinate the spectrum on the spaces $\mathcal{H}_{\tau, \gamma}^J \cong L^2(\R^m, \C^{2^{[d/2]}})$, we use results from the representation theory of the Heisenberg group. It is known that $L^2(\R^m)$ is the representation space of irreducible unitary representations of $\Hei^{2m+1}$, and what we have done so far is to decompose the right regular representation $R: \Hei^{2m+1} \times \R^{d-2m} \rightarrow L^2(\Sigma_\delta^H M)$ into its irreducible components. Thereby, the frame $\{\tilde{X}_1, \ldots, \tilde{X}_{d+1}\}$ exactly corresponds to the Carnot group $\Hei^{2m+1} \times \R^{d-2m}$, while the frame $\{X_1, \ldots, X_{d+1}\}$ only corresponds to a Carnot group isomorphic to $\Hei^{2m+1} \times \R^{d-2m}$. \smallskip

Now it is well known by the theorem of Stone and von Neumann how they have to look like on $L^2(\R^m)$ (see e.g. \cite{Fol}, Theorem (1.50)): In the polarized coordinates we are using, the infinite dimensional unitary irreducible (Schrödinger) representations of $\Hei^{2m+1}$ are given by
\begin{equation} \label{eq representations of Heisenberg group}
 \pi_r: \Hei^{2m+1} \rightarrow \mathcal{U}(L^2(M)), \ \ \ \pi_r(x,y,t) f(u) = e^{2\pi i r (t + \langle u,y \rangle)} f(x+u) 
\end{equation}
for any $r \in \R \setminus \{0\}$. In our case, for $\sigma \in \mathcal{H}_{\tau, \gamma}^J \cong L^2(\R^m) \otimes \C^{2^{[d/2]}}$, we have $r = \tau$. From now on, we write $f$ instead of $\sigma$ whenever we consider an element of $L^2(\R^m, \C^{2^{[d/2]}})$. Hence using the expression \eqref{eq D^H via right regular representation codim 1 tilde X} of the horizontal Dirac operator via the right-regular representation $R$ we have, after plugging in \eqref{eq representations of Heisenberg group},
\begin{eqnarray*}
\left. \frac{d}{ds} R(\exp s\tilde{\lambda}_j \tilde{X}_j) f(u) \right|_{s=0} &=& \begin{cases} \left. \frac{d}{ds} f(u + s\tilde{\lambda}_je_j) \right|_{s=0} & j \in \{1, \ldots, m\}\\
      \left. \frac{d}{ds} e^{2\pi i \tau \langle u, s \tilde{\lambda}_{j-m} e_{j-m} \rangle} f(u) \right|_{s=0} & j \in \{m+1, \ldots, 2m\}
     \end{cases}\\
 &=& \begin{cases} \tilde{\lambda}_j \frac{\partial}{\partial u_j} f(u) & j \in \{1, \ldots, m\}\\
      2\pi i \tau u_{j-m} \tilde{\lambda}_{j-m} f(u) & j \in \{m+1, \ldots, 2m\}
     \end{cases}.
\end{eqnarray*}
We still need to express the action of $X_{2m+j} = \tilde{X}_{2m+j}$ for $1 \leq j \leq d-2m$ on $\mathcal{H}_{\tau,\gamma}^J$. But this is just multiplication by $2\pi i \gamma_j$, since $\sigma \in \mathcal{H}_\tau$ is characterized via
\begin{equation*} \label{eq full characterization sigma in H_tau}
 \sigma(x,y,z,t) = e^{2\pi i t \tau} \sum_{\gamma} \sum_{\beta} b_{\gamma,\beta}^\sigma(x) e^{2\pi i \langle \beta,y \rangle} e^{2\pi i \langle\gamma,z \rangle}
\end{equation*}
(see \eqref{eq full characterization sigma in H_tau} in the proof of Lemma \ref{lemma spectral decomposition H_tau, tau neq 0}), and therefore we get
$$\left. \frac{d}{ds} R(\exp s\tilde{X}_{2m+j}) \sigma (x,y,z,t) \right|_{s=0} = 2\pi i \gamma_j \cdot \sigma(x,y,z,t). $$
Plugging everything into the expression \eqref{eq D^H via right regular representation codim 1 tilde X} we see that the action of $D^H$ on every $\mathcal{H}_{\tau,\gamma}^J \cong L^2(\R^m) \otimes \C^{2^{[d/2]}}$ is given by
\begin{equation}\label{eq D^H on H_tau,gamma^J}
\begin{split}
 D^H f(u) &= \sum_{j=1}^m \tilde{\lambda}_j \left( c^H(X_j) \frac{\partial}{\partial u_j} + 2\pi i \tau u_j c^H(X_{m+j}) \right) f(u) + \sum_{k=2m+1}^d 2\pi i \gamma_k c^H(X_k) \cdot f(u)\\
&= \sum_{j=1}^m \sqrt{\lambda_j} \left( c^H(X_j) \frac{\partial}{\partial u_j} + 2\pi i \tau u_j c^H(X_{m+j}) \right) f(u) + \sum_{k=2m+1}^d 2\pi i \gamma_k c^H(X_k) \cdot f(u).
\end{split}
\end{equation}\smallskip

To calculate the absolute values of the eigenvalues of $D^H$ restricted to $\mathcal{H}_{\tau,\gamma}^J$, we once again consider the square of the operator. From \eqref{eq D^H on H_tau,gamma^J} we get
\begin{equation} \label{eq square of D^H on H_tau,gamma^J}
\begin{split}
  \left(D^H\right)^2 f(u) =& \left(\sum_{j=1}^m \sqrt{\lambda_j} \left( c^H(X_j) \frac{\partial}{\partial u_j} + 2\pi i \tau u_j c^H(X_{m+j}) \right) \right)^2 f(u)\\
& + \left(\sum_{k=2m+1}^d 2\pi i \gamma_k c^H(X_k) \right)^2 f(u) \\
& +\sum_{j=1}^m \sqrt{\lambda_j} \left( c^H(X_j) \frac{\partial}{\partial u_j} + 2\pi i \tau u_j c^H(X_{m+j}) \right) \cdot \sum_{k=2m+1}^d 2\pi i \gamma_k c^H(X_k) \\
& + \sum_{k=2m+1}^d 2\pi i \gamma_k c^H(X_k) \cdot \sum_{j=1}^m \sqrt{\lambda_j} \left( c^H(X_j) \frac{\partial}{\partial u_j} + 2\pi i \tau u_j c^H(X_{m+j}) \right),
\end{split}
\end{equation}
and because of the Clifford relation $c^H(X_l)c^H(X_{l'}) + c^H(X_{l'})c^H(X_l) = 0$ for any $l \neq l'$, the sum of the third and the forth term on the right hand side of \eqref{eq square of D^H on H_tau,gamma^J} vanish. For the same reason, and since $c(X_l)^2 = -\id$ for any $l$, we get
$$ \left(\sum_{k=2m+1}^d 2\pi i \gamma_k c^H(X_k) \right)^2 = 4 \pi^2 \sum_{k=2m+1}^d \gamma_k^2 $$
for the second term. Finally, using again the rules of Clifford action and the Leibniz rule for differentiation, we calculate 
\begin{eqnarray*}
 & & \sum_{j=1}^m \sqrt{\lambda_j} c^H(X_j) \frac{\partial}{\partial u_j} \sum_{k=1}^m \sqrt{\lambda_k} \left( c^H(X_k) \frac{\partial}{\partial u_k} f(u) + 2\pi i \tau u_k c^H(X_{m+k}) f(u) \right) \\
&=& \sum_{j=1}^m \sum_{k=1}^m \sqrt{\lambda_j \lambda_k} c^H(X_j) c^H(X_k) \frac{\partial}{\partial u_j} \frac{\partial}{\partial u_k} f(u) + \times \\
& & \times \ \sum_{j=1}^m \sum_{k=1}^m 2 \pi i \tau \sqrt{\lambda_j \lambda_k} c^H(X_j) c^H(X_{m+k}) \frac{\partial}{\partial u_j} \left(u_k f(u) \right) \\
& = & - \sum_{j=1}^m \lambda_j \frac{\partial^2}{\partial u_j^2} f(u) + \times \\
& & \times \  2 \pi i \tau \left( \sum_{j=1}^m \sum_{k=1}^m u_k \sqrt{\lambda_j \lambda_k} c^H(X_j) c^H(X_{m+k}) \frac{\partial}{\partial u_j} f(u) + \sum_{j=1}^m \lambda_j c^H(X_j)c^H(X_{m+j}) f(u) \right)
\end{eqnarray*}
and 
\begin{eqnarray*}
 & & 2\pi i \tau \sum_{j=1}^m u_j \sqrt{\lambda_j} c^H(X_{m+j}) \sum_{k=1}^m \sqrt{\lambda_k} \left( c^H(X_k) \frac{\partial}{\partial u_k} f(u) + 2\pi i \tau u_k c^H(X_{m+k}) f(u) \right) \\
&=& 2 \pi i \tau \sum_{j=1}^m \sum_{k=1}^m u_j \sqrt{\lambda_j \lambda_k} c^H(X_{m+j}) c^H(X_k) \frac{\partial}{\partial u_k} f(u) \times \\
& & \times \ - 4 \pi^2 \tau^2 \sum_{j=1}^m \sum_{k=1}^m u_j u_k \sqrt{\lambda_j \lambda_k} c^H(X_{m+j}) c^H(X_{m+k}) f(u) \\
&=& 2 \pi i \tau \sum_{j=1}^m \sum_{k=1}^m u_j \sqrt{\lambda_j \lambda_k} c^H(X_{m+j}) c^H(X_k) \frac{\partial}{\partial u_k} f(u) + 4 \pi^2 \tau^2 \sum_{j=1}^m \lambda_j u_j^2 f(u).
\end{eqnarray*}
Because of the identity $c^H(X_k) c^H(X_{m+j}) + c^H(X_{m+j}) c^H(X_k) = 0$ for any $k,j \in \{1, \ldots, m\}$, we have
$$\sum_{j=1}^m \sum_{k=1}^m u_k \sqrt{\lambda_j \lambda_k} c^H(X_j) c^H(X_{m+k}) \frac{\partial}{\partial u_j}f + \sum_{j=1}^m \sum_{k=1}^m u_j \sqrt{\lambda_j \lambda_k} c^H(X_{m+j}) c^H(X_k) \frac{\partial}{\partial u_k}f = 0,$$
and thus get for the first term of \eqref{eq square of D^H on H_tau,gamma^J} from these two calculations
\begin{eqnarray*}
 & & \left(\sum_{j=1}^m \sqrt{\lambda_j} \left( c^H(X_j) \frac{\partial}{\partial u_j} + 2\pi i \tau u_j c^H(X_{m+j}) \right) \right)^2 f(u)\\
& = & -\sum_{j=1}^m \lambda_j \frac{\partial^2}{\partial u_j^2} f(u) + 4 \pi^2 \tau^2 \sum_{j=1}^m \lambda_j u_j^2 f(u) + 2\pi i \tau \sum_{j=1}^m \lambda_j c^H(X_j)c^H(X_{m+j}) f(u) \\
& = & \sum_{j=1}^m \lambda_j \left( -\frac{\partial^2}{\partial u_j^2} + 4 \pi^2 \tau^2 u_j^2 + 2\pi i \tau c^H(X_j)c^H(X_{m+j}) \right) f(u).
\end{eqnarray*}
Now after plugging everything into the expression \eqref{eq square of D^H on H_tau,gamma^J} for the restriction of $(D^H)^2$ onto the space $\mathcal{H}_{\tau, \gamma}^J \cong L^2(\R^m, \C^{2^{[d/2]}})$ we see that 
\begin{equation} \label{eq square of D^H on H_tau,gamma^J simplier}
 \left. \left(D^H\right)^2 \right|_{\mathcal{H}_{\tau,\gamma}^J} = \sum_{j=1}^m \lambda_j \left( -\frac{\partial^2}{\partial u_j^2} + 4 \pi^2 \tau^2 u_j^2 + 2\pi i \tau c^H(X_j)c^H(X_{m+j}) \right) + 4 \pi^2 \sum_{k=2m+1}^d \gamma_k^2
\end{equation}
on these spaces.\smallskip

From \eqref{eq square of D^H on H_tau,gamma^J simplier} we observe that we can calculate the eigenvalues of $D^H$ from the eigenvalues of the harmonic oscillator. It is well known that the operator
$$ \sum_{j=1}^m \lambda_j \left(-\frac{\partial^2}{\partial u_j^2} + 4 \pi^2 \tau^2 u_j^2 \right)$$
possesses the eigenvalues $\eta_{\tau, \kappa} = 2\pi |\tau| \sum_{j=1}^m \lambda_j (2\kappa_j + 1)$, where $\kappa = (\kappa_1, \ldots, \kappa_m) \in \N^m$ (see e.g. \cite{Fol}, Section 1.7). Since this operator is acting on the $L^2$-space of $\C^{2^[d/2]}$-valued functions, each of these eigenvalues has the multiplicity $2^{[d/2]}$. In addition, for each $\kappa \in \N^m$ we have to add the eigenvalues of the operator 
$$-\sum_{j=1}^m 2\pi i \tau \lambda_j c^H(X_j)c^H(X_{m+j}).$$
By Proposition \ref{prop eigenvalues sum c(X_j)c(X_k)} from Section 3.2, every eigenvalue $\eta_{\tau,\kappa,l}$, $l \in \{1, \ldots, 2^{[d/2]}\}$ of this $(2^{[d/2]} \times 2^{[d/2]})$-matrix is included in the set 
$$\{-2 \pi |\tau| \sum_{j=1}^m\lambda_j, \ldots, 2 \pi |\tau| \sum_{j=1}^m\lambda_j\}.$$
The third term of \eqref{eq square of D^H on H_tau,gamma^J simplier} is just an additive constant (since $\gamma$ is fixed by choosing a space $\mathcal{H}_{\tau, \gamma}^J$), and hence on each $\mathcal{H}_{\tau,\gamma}^J$ the operator $(D^H)^2$ from \eqref{eq square of D^H on H_tau,gamma^J simplier} has the eigenvalues
\begin{equation} \label{eq exact eigenvalues D^H square on H_tau,gamma^J}
 \tilde{\mu}_{\tau,\gamma,\kappa,l} = 2\pi |\tau| \left( \sum_{j=1}^m \lambda_j (2\kappa_j + 1) + \tilde{\eta}_{\kappa,l} \right) + 4 \pi^2 \sum_{k=2m+1}^d \gamma_k^2
\end{equation}
with $\kappa \in \N^m$ and $- \sum_{j=1}^m\lambda_j \leq \tilde{\eta}_{\kappa,l} \leq \sum_{j=1}^m\lambda_j $ for $ l \in \{1, \ldots, 2^{[d/2]}\}$. But since each number $\tilde{\eta}_{\kappa,l}$ is bounded by constants depending only on the constants $\lambda_1, \ldots, \lambda_m$ which are given by the Levi form on $M$, they do not change the asymptotic behavior of the eigenvalues from \eqref{eq exact eigenvalues D^H square on H_tau,gamma^J}. Therefore we can say that 
\begin{equation} \label{eq asymptotic eigenvalues D^H square on H_tau,gamma^J}
 \left(\mu_{\tau,\gamma,\kappa}\right)^2 \sim 2\pi |\tau| \sum_{j=1}^m \lambda_j (2\kappa_j + 1) + 4 \pi^2 \sum_{k=2m+1}^d \gamma_k^2,
\end{equation}
where each $(\mu_{\tau,\gamma,\kappa})^2$ has the multiplicity $2^{[d/2]}$, for the eigenvalues of $(D^H)^2$ restricted to $\mathcal{H}_{\tau,\gamma}^J$.\smallskip

From \eqref{eq asymptotic eigenvalues D^H square on H_tau,gamma^J} we get for the absolute values of the eigenvalues of $D^H$
\begin{equation} \label{eq eigenvalues absolute on H_tau,gamma^J}
 \left| \mu_{\tau,\gamma,\kappa} \right| \sim \sqrt{2 \pi |\tau| \sum_{j=1}^m \lambda_j (2\kappa_j + 1) + 4 \pi^2 \sum_{k=2m+1}^d \gamma_k^2}
\end{equation}
for $\tau \in \Z \setminus \{0\}$, $\gamma \in \frac{1}{2}\Z^{d-2m}$ such that $e^{2\pi \gamma_j} = \delta_{2m+j}$ and $ \kappa \in \N^m$. Each of these (asymptotic) eigenvalues has the multiplicity $2^{[d/2]} |\tau|^m$, since there are $|\tau|^m$ copies of the spaces $\mathcal{H}_{\tau,\gamma}^J \cong L^2(\R^m, \C^{2^{[d/2]}})$. Thus we have proved statement (a). \smallskip

We still have to show the statements (b) and (c) about the degeneracy of $D^H$. This will be done by showing that for each $\tau \in \Z \setminus \{0\}$, the first term of the operator $(D^H)^2$ from \eqref{eq square of D^H on H_tau,gamma^J simplier} has at least one $0$-eigenvalue on $\mathcal{H}_{\tau,\gamma}^J$: From that it will follow that for every $\gamma \in \Z^{d-2m}$ the number 
$$4 \pi^2 \sum_{k=2m+1}^d \gamma_k^2$$
appears as an eigenvalue on each space $\mathcal{H}_{\tau,\gamma}^J$, and is thus an eigenvalue of infinite multiplicity. Furthermore all the degenerate eigenvalues of $D^H$ are of this form, since every other eigenvalue on $H_{\tau,\gamma}^J$ depends on $\tau$ and is of finite multiplicity on this space. This shows statement (c), and since for $\G \cong \Hei^{2m+1}$ the term 
$$\sum_{k=2m+1}^d 2 \pi i \gamma_k c^H(X_k)$$
does not appear in the expression \eqref{eq D^H on H_tau,gamma^J} of $D^H$ on $H_{\tau,\gamma}^J$ the statement (b) also follows.\smallskip

Indeed the first term 
\begin{equation} \label{eq first term of D^H square on H_tau,gamma}
 \sum_{j=1}^m \lambda_j \left(-\frac{\partial^2}{\partial u_j^2} + 4 \pi^2 \tau^2 u_j^2 \right) - \sum_{j=1}^m 2\pi i \tau \lambda_j c^H(X_j)c^H(X_{m+j})
\end{equation}
of \eqref{eq square of D^H on H_tau,gamma^J simplier} has at least one $0$-eigenvalue. As we already noted, the eigenvalues of the harmonic oscillator are the numbers 
$$\tilde{\mu}_\kappa = 2\pi |\tau| \sum_{j=1}^m \lambda_j (2\kappa_j + 1),$$
where $\kappa = (\kappa_1, \ldots, \kappa_m) \in \N^m$, and the eigenvalues of the matrix 
$$-\sum_{j=1}^m 2\pi i \tau \lambda_j c^H(X_j)c^H(X_{m+j})$$
belong to the set
$$\left\{-2 \pi |\tau| \sum_{j=1}^m \lambda_j, \ldots, 2 \pi |\tau| \sum_{j=1}^m \lambda_j \right\}.$$
Hence, the only chance for the operator \eqref{eq first term of D^H square on H_tau,gamma} to have a $0$-eigenvalue is that $\kappa = 0$ and that the number $-2 \pi |\tau| \sum_{j=1}^m \lambda_j$ is indeed an eigenvalue of $-\sum_{j=1}^m 2\pi i |\tau| \lambda_j c^H(X_j)c^H(X_{m+j})$. But the last statement is true by Proposition \ref{prop eigenvalues sum c(X_j)c(X_k)} from Section 3.2, and hence the above argumentation proves the statements (b) and (c) of the theorem.\eB
\end{thm}\medskip

As mentioned before, the results of Theorem \ref{thm spectral decomposition D^H on Heisenberg group} imply that $D^H$ does not furnish a spectral triple on $M = \Gamma \backslash \G$. But it fits into our definition of a degenerate spectral triple and it detects the Carnot-Carath\'{e}odory metric of our Carnot manifold $M$, which we write down in the following corollary.\medskip

\begin{cor} \label{cor D^H gives degenerate spectral triple}
 Let $M = \Gamma \backslash \G$  be the compact nilmanifold of a Carnot group $\G \cong \Hei^{2m+1} \times \R^{d-2m}$. Let $D^H$ be a horizontal Dirac operator on $M$ arising via pull-back from a spin structure $\Sigma^{\T^d}_\delta$ on the torus $\T^d$. Then the triple $(C(M),L^2(\Sigma_\delta^HM),D^H)$ is a degenerate spectral triple, which detects the Carnot-Carath\'{e}odory metric via Connes' metric formula.\smallskip

\B Since $D^H$ is a horizontal Dirac operator in the sense of Section 3.2, it has been shown in Section 3.3 that it detects the Carnot-Carath\'{e}odory metric on $M$ via Connes metric formula. In particular it has also been shown that 
$$\|[D^H,f]\| = \mathrm{ess} \sup_{x \in M} \left\|\grad^H f \right\|$$
for every $f \in \mathrm{Lip}_{CC}(M)$, and since the number on the right hand side is bounded and $\mathrm{Lip}_{CC}(M)$ is a dense sub-algebra of $C(M)$, condition (i) for a spectral triple (see Definition \ref{def spectral triple}) is fulfilled.\smallskip

Now from Theorem \ref{thm spectral decomposition D^H on Heisenberg group} we know that the spectrum of $D^H$ is discrete. If we exclude the eigenvalues of infinite multiplicity, then for a given number $\Lambda \in \R$ there are only finitely many eigenvalues in the spectrum of $D^H$ which are smaller than $\Lambda$. This shows us that $D^H$ has a compact resolvent if we restrict it to the orthonormal complement of the degenerate eigenspaces, which means we have a degenerate spectral triple according to Definition \ref{def degenerate spectral tripel}. \eB 
\end{cor}\medskip

As a further consequence of Theorem \ref{thm spectral decomposition D^H on Heisenberg group}, we show that besides the Carnot-Carath\'{e}odory metric we can also detect the Hausdorff dimension of the metric space $(M, d_{CC})$ from the degenerate spectral triple $(C(M), L^2(\Sigma^H_\delta M), D^H)$. Therefore we consider the asymptotic behavior of the non-degenerate eigenvalues of $D^H$.\smallskip

To detect the asymptotic behavior in our case, we will make use of the following proposition, which shows the equivalence about the growth function of the number of eigenvalues within a certain radius and their asymptotic growth. We write down the statement in a very general matter, since for its proof it is not important what exactly the elements of the sets whose growth rate we analyze describe. For a proof, we refer to \cite{Shu}, Proposition 13.1. \medskip

\begin{prop} \label{prop asymptotic behaviour versus growth rate}
 Let $\Lambda := \{\lambda_0, \lambda_1, \lambda_2, \ldots\} \subset \R^+$ be a discrete ordered set of positive real numbers, which means we have $\lambda_0 \leq \lambda_1 \leq \ldots$. For $t \in \R^+$, we denote by
$$N_\Lambda(t) := \sum_{\{k \in \N: \lambda_k \leq t\}} 1 = \# \left\{\lambda_k \in \Lambda: \lambda_k \leq t\right\}$$
the number of all elements of $\Lambda$ which are bounded by $t$. In addition, we consider numbers $V_0 \in \R$ and $n \in \N^+$. Then the following statements are equivalent:
\begin{enumerate}[(i)]
 \item $N_\Lambda(t) \sim V_0 t^n$ as $t \rightarrow \infty$.
 \item $\lambda_k \sim V_0^{-1/n} \cdot k^{1/n}$ as $k \rightarrow \infty$. \eB
\end{enumerate}
\end{prop}\medskip

Using this proposition, we can show that the Hausdorff dimension on $(M,d_{CC})$ coincides with the metric dimension of the degenerate spectral triple $(C(M),L^2(\Sigma^H_\delta M), D^H)$.\medskip

\begin{cor} \label{cor Hausdorff dimension from D^H}
 Let $\G \cong \Hei^{2m+1}$ be isomorphic to the $2m+1$-dimensional Heisenberg group, and let $M = \Gamma \backslash \G$  be its compact nilmanifold. Let $D^H$ be a horizontal Dirac operator on $M$ arising via pull-back from a spin structure $\Sigma^{\T^d}_\delta$ on the torus $\T^d$ with $d=2m$.\smallskip

Then the metric dimension of the degenerate spectral triple 
$$\left(C(M),L^2(\Sigma^H_\delta M), D^H\right)$$
is $d+2$, i.e. it coincides with the Hausdorff dimension of the metric space $(M,d_{CC})$. \smallskip

\B By Theorem \ref{thm spectral decomposition D^H on Heisenberg group}, the only degenerate eigenvalue of $D^H$ on $M = \Gamma \backslash \Hei^{2m+1}$ is $0$. Hence let $\Lambda$ denote the set of all eigenvalues of $D^H$ which are not $0$. We have to show by Definition \ref{def metric dimension}:
\begin{equation} \label{eq condition for metric dimension 2m+2}
\sum_{\mu \in \Lambda} \frac{1}{\left| \mu \right|^p} < \infty \ \ \ \Leftrightarrow \ \ \ p > 2m+2.
\end{equation}
Now we can decompose $\Lambda$ into two disjoint sets 
\begin{equation} \label{eq decomposition eigenvalue set}
 \Lambda = \Lambda_1 \dot{\cup} \Lambda_2,
\end{equation}
where $\Lambda_1$ contains all the eigenvalues of $D^H$ listed under (i) and $\Lambda_2$ contains all the eigenvalues of $D^H$ listed under (ii) in Theorem \ref{thm spectral decomposition D^H on Heisenberg group} which are not $0$. The set $\Lambda_1$ contains exactly the eigenvalues of the classical Dirac operator, acting on the spinor bundle $\Sigma_\delta \T^d$ of the torus $\T^d$ (with respect to a Riemannian metric on $\T^d$ such that the vector fields $\tilde{X}_1, \ldots, \tilde{X}_d$ from \eqref{eq frame tilde X} used in the proof of the theorem form an orthonormal frame for $\T^d$), and since this classical Dirac operator is elliptic it is known that they grow proportional to the function $k^d$ by Weyl asymptotic. But this shows
\begin{equation} \label{eq condition for metric dimension on Lambda_1}
 \sum_{\mu \in \Lambda_1} \frac{1}{\left| \mu \right|^p} < \infty \ \ \ \Leftrightarrow \ \ \ p > d = 2m.
\end{equation}
Therefore the crucial point of the asymptotic behavior of our eigenvalues lies in the set $\Lambda_2$.\smallskip

Applying Theorem \ref{thm spectral decomposition D^H on Heisenberg group} to our situation where $\G \cong \Hei^{2m+1}$, we see that the eigenvalues belonging to $\Lambda_2$ are given asymptotically by
\begin{equation} \label{eq eigenvalues Lambda_2}
 \left| \mu_{\tau,\kappa} \right| \sim \sqrt{2 \pi |\tau| \sum_{j=1}^m \lambda_j (2\kappa_j + 1)} \sim \sqrt{|\tau|} \sqrt{\sum_{j=1}^m \kappa_j}
\end{equation}
for $\tau \in \Z \setminus \{0\}$ and $\kappa \in \N^m \setminus \{0\}$. We make a further (disjoint) decomposition
$$\Lambda_2 = \dot{\bigcup}_{\tau \in \Z \setminus \{0\}} \Lambda_{2, \tau}$$
with $\Lambda_{2,\tau}$ containing exactly the eigenvalues $\mu_{\tau,\kappa}$ with $\kappa \in \N^m$. Note that each of these eigenvalues occurs with the multiplicity $2^{[d/2]} |\tau|^m$. For a fixed $\tau \in \Z \setminus \{0\}$ we denote the absolute values of these elements belonging to $\Lambda_{2,\tau}$ by 
\begin{equation} \label{eq denotation elements of Lambda_2,tau}
 \tilde{\mu}_{\tau,0} \leq \tilde{\mu}_{\tau,1} \leq \ldots,
\end{equation}
such that we can work with the number $N_{\Lambda_{2,\tau}}(t)$ from Proposition \ref{prop asymptotic behaviour versus growth rate}. It is well known that for any $C \in \R^+$ we have
$$ \# \left\{ \kappa \in \N^m: \sum_{j=1}^m \kappa_j \leq C \right\} \sim \# \left\{ \kappa \in \N^m: \sqrt{ \sum_{j=1}^m \kappa_j^2} \leq C \right\} \sim C^m$$
(one can consider these numbers as the eigenvalues of the elliptic Dirac operator on the $m$-dimensional torus, and then this statement follows from Weyl asymptotic), which means
$$\# \left\{ \kappa \in \N^m: \sqrt{\sum_{j=1}^m \kappa_j} \leq C \right\} \sim C^{2m}, $$
and hence for every $\tau \in \Z \setminus \{0\}$ it follows from \eqref{eq eigenvalues Lambda_2} that
\begin{equation}
\begin{split} \label{eq asymptotic for N_Lambda_2,tau}
  N_{\Lambda_{2,\tau}}(t) &\sim \# \left\{ \kappa \in \N^m: \sqrt{|\tau|} \sqrt{\sum_{j=1}^m \kappa_j} \leq t \right\}\\
&= \# \left\{ \kappa \in \N^m: \sqrt{\sum_{j=1}^m \kappa_j} \leq \frac{t}{\sqrt{|\tau|}} \right\}\\
&\sim \left( \frac{t}{\sqrt{|\tau|}} \right)^{2m} = \frac{1}{\left|\tau\right|^m} t^{2m}.
\end{split}
\end{equation}
At this point we can apply Proposition \ref{prop asymptotic behaviour versus growth rate}, which tells us that \eqref{eq asymptotic for N_Lambda_2,tau} is equivalent to the fact that
\begin{equation} \label{eq asymptotic behavior lambda_tau,k}
 \tilde{\mu}_{\tau,k} \sim |\tau|^{\frac{1}{2}} k^{\frac{1}{2m}} 
\end{equation}
for the elements of $\Lambda_{2,\tau}$ denoted by \eqref{eq denotation elements of Lambda_2,tau}, where we do not care about the multiplicity $2^{[d/2]}|\tau|^m$ of these numbers as eigenvalues of $D^H$ for a moment.\smallskip

With these results, we can check \eqref{eq condition for metric dimension 2m+2} for the set $\Lambda_2$ which will prove the corollary: Because of our decomposition of $\Lambda_2$ and \eqref{eq asymptotic behavior lambda_tau,k} we have
\begin{eqnarray*}
 \sum_{\mu \in \Lambda_2} \frac{1}{\left|\mu\right|^p}
&=& \sum_{\tau \in \Z \setminus \{0\}} \sum_{\mu \in \Lambda_{2,\tau}} \frac{1}{\left|\mu\right|^p}\\
&\sim& \sum_{\tau \in \Z \setminus \{0\}} \sum_{k \in \N} 2^{\left[\frac{d}{2}\right]} \left| \tau \right|^m \cdot \frac{1}{\left(|\tau|^{\frac{1}{2}} k^{\frac{1}{2m}}\right)^p}\\
&=& \sum_{\tau \in \Z \setminus \{0\}} 2^{\left[\frac{d}{2}\right]} \frac{1}{\left|\tau\right|^{\frac{p}{2}-m}} \cdot \sum_{k \in \N} \frac{1}{k^{\frac{p}{2m}}},\\
\end{eqnarray*}
where we have taken care of the multiplicity $2^{[d/2]} |\tau|^m$ of every eigenvalue $\tilde{\mu}_{\tau,k}$ in the second equation. Now the second geometric series in the last line of this calculation converges if and only if $p>2m$, and the first geometric series in this line converges if and only if $p>2m+2$. But this, together with \eqref{eq condition for metric dimension on Lambda_1} tells us that the series
$$ \sum_{\mu \in \Lambda} \frac{1}{\left|\mu\right|^p} = \sum_{\mu \in \Lambda_1} \frac{1}{\left|\mu\right|^p} + \sum_{\mu \in \Lambda_2} \frac{1}{\left|\mu\right|^p} $$
converges if and only if $p>2m+2$, which shows that \eqref{eq condition for metric dimension 2m+2} is true and hence the corollary is proved. \eB \smallskip

\Bem The same statement should also be true for the case $\G \cong \Hei^{2m+1} \times \R^{d-2m}$, but in this case the argument is a bit more involved. \eBsp 
\end{cor}\medskip

\begin{remark} \normalfont
 From the argumentation of Theorem \ref{thm spectral decomposition D^H on Heisenberg group} it is also possible to calculate the eigenvalues of the horizontal pull-back Dirac operator exactly for given numbers $m \in \N$ and $d \in \N$ (after determining the matrices $c^H(X_j)$ describing the Clifford action). The idea is to use the orthonormal basis of $L^2(\R^m)$ consisting of the Hermite functions
\begin{equation} \label{eq Hermite functions}
 h_k(v) = e^{\left\|v\right\|^2/2} \frac{\partial^{k_1+\ldots+k_m}}{\partial v_1^{k_1} \ldots v_m^{k_m}} e^{-\left\| v \right\|^2},
\end{equation}
which allows us to decompose $\mathcal{H}_{\tau,\gamma}^J$ into finite dimensional subspaces which are invariant under $D^H$. On these subspaces the determination of the spectrum is done by calculating eigenvalues of matrices.\smallskip

In the case $m=1$ where $\G \cong \Hei^3$ is the $3$-dimensional compact Heisenberg nilmanifold, this has been done in \cite{Bae} and \cite{AB} for an ordinary Dirac operator. The argument used there can be transfered to our case. After these calculations, we find that the eigenvalues of $D^H$ are exactly the numbers 
\begin{enumerate}[(i)]
 \item $\mu^{\pm}_{\alpha, \beta} = \pm 2\pi \sqrt{\alpha^2 + \beta^2}$ for $\alpha, \beta \in \Z$ such that $e^{2\pi i \alpha} = \delta_1$, $e^{2 \pi i \beta} = \delta_2$ of multiplicity $1$,
 \item $\mu^{\pm}_{\tau, \kappa} = \pm 2 \sqrt{\kappa\pi |\tau|}$ for $\kappa \in \Z^+$, $\tau \in \Z \setminus \{0\}$ of multiplicity $|\tau|$,
 \item $\mu_{\tau,0} = 0$ for $\tau \in \Z \setminus \{0\}$ of multiplicity $|\tau|$.
\end{enumerate}
Considering the inverses of the absolute values of these eigenvalues which are not zero, we find (for the set $\Lambda_2$ from \eqref{eq decomposition eigenvalue set} in the proof of Theorem \ref{thm spectral decomposition D^H on Heisenberg group}) that
$$\sum_{\tau \in \Z \setminus \{0\}} \sum_{\kappa \in \N} \frac{1}{\left|\mu_{\tau,\kappa}^{\pm} \right|^p} \leq \infty \ \ \ \Leftrightarrow \ \ \ p > 4,$$
which shows the statement of Corollary \ref{cor Hausdorff dimension from D^H} for this example.\smallskip

Using the Hermite functions from \eqref{eq Hermite functions}, it is possible to describe the kernel of $D^H$ on the sub-spaces $L^2(\R^m, \C^{2^{[d/2]}})$ of the decomposition \eqref{eq decomposition H_tau codim 1} of each space $\mathcal{H}_\tau$ from Lemma \ref{lemma spectral decomposition H_tau, tau neq 0}. In the case $m=1$ and $\G \cong \Hei^3$, one finds that on each space $L^2(\R,\C^2)$ belonging to $\mathcal{H}_\tau$ the kernel of $D^H$ is spanned by the $\C^2$-valued functions
$$ f_0(t) = h_0(\sqrt{2\pi |\tau|} t) \begin{pmatrix} 1 \\ 0 \end{pmatrix} = e^{-\pi |\tau| t^2} \begin{pmatrix} 1 \\ 0 \end{pmatrix}$$
for $\tau > 0$ and
$$f_0(t) = h_0(\sqrt{2\pi |\tau|} t) \begin{pmatrix} 0 \\ 1 \end{pmatrix} = e^{-\pi |\tau| t^2} \begin{pmatrix} 0 \\ 1 \end{pmatrix}$$
for $\tau < 0$. In particular we observe that these functions are smooth.\smallskip

This observation leads to the idea that $D^H$ may be an example for a so-called \emph{weakly hypoelliptic operator}, in the sense that $D^H \varphi = 0$ implies that $\varphi$ is smooth. Recently, these weakly hypoelliptic operators have been considered by Christian Bär in \cite{Bae2}. \eBsp
\end{remark}\medskip

To close this section, let us once again summarize the results we achieved so far: We have seen that for the case $\G \cong \Hei^{2m+1} \times \R^{d-2m}$, the horizontal pull-back Dirac operator on $\Gamma \backslash \G$ does not have a compact resolvent. Hence it only furnishes a degenerate spectral triple, but from this degenerate spectral triple the most important ingredients of the Carnot-Carath\'{e}odory geometry on $M$, which are the Carnot-Carath\'{e}odory metric and the Hausdorff dimension, can be detected.\smallskip

In the next section, we will deduce from this co-dimension $1$ case that the absence of a compact resolvent of $D^H$ occurs on any compact Carnot nilmanifold. \bigskip

\section{Degeneracy of $D^H$ in the General Case}

This section will provide the most general result of this chapter concerning compact Carnot nilmanifolds: On an arbitrary compact Carnot nilmanifold the horizontal pull-back Dirac operator has (at least one) infinite dimensional eigenspace and does therefore not have a compact resolvent. Later we will prove that any horizontal Dirac operator on an arbitrary Carnot manifold fails to be hypoelliptic, which can be seen as a generalization of this statement.\smallskip

The starting point is Theorem \ref{thm decomposition from center}, which allows us to reduce the problem from an arbitrary Carnot group to a $2$-step nilpotent one. In the last section we have made a detailed treatment of compact Heisenberg nilmanifolds, which are exactly those $2$-step nilmanifolds whose horizontal distribution has co-dimension $1$, and seen that we have a degenerate eigenspace in this case. So the only step missing is to get from an arbitrary space of step $2$ to a compact Heisenberg nilmanifold. \smallskip

To do this, we consider a submersion of the type introduced in Section 2.4. Let $M_2 = \Gamma_2 \backslash \G_2$ be the compact nilmanifold of a Carnot group $\G_2$ of rank $2$. We assume that the Lie algebra $\mathfrak{g}_2$ of $\G$ has the grading $\mathfrak{g}_2 = V_1 \oplus V_2$, where $\{X_{1,1}, \ldots, X_{1,d_1}\}$ is an orthonormal frame for $V_1$ and $\{X_{2,1}, \ldots, X_{2,d_2}\}$ is an orthonormal frame for $V_2$. For a $\nu \in \{1, \ldots, d_2\}$ we consider the orthonormal projection
\begin{equation} \label{eq projection algebra level for homogeneous Carnot spaces}
\begin{split}
 \mathrm{pr}_\nu: & \ \ \mathfrak{g}_2 \rightarrow \mathfrak{g}_{2,\nu} \simeq V_1 \oplus \spa \{X_{2,\nu}\}, \\
 & \ \ v \mapsto v \mod \spa \left( \{X_{2,1}, \ldots , X_{2,d_2} \} \setminus \{X_{2,\nu}\} \right).
\end{split}
\end{equation}
We have seen in Section 2.4 that the vector space $\mathfrak{g}_{2,\nu}$ can be canonically equipped with a Lie bracket, such that it is a graded nilpotent Lie algebra of rank $2$, where $V_1$ is bracket generating of step $2$ and co-dimension $1$ for $\mathfrak{g}_{2,\nu}$. We denote the Carnot group arising from $\mathfrak{g}_{2,\nu}$ by $\G_{2,\nu}$, and the Lie group homomorphism arising from $\mathrm{pr}_\nu$ by
\begin{equation} \label{eq projection group level for homogeneous Carnot spaces}
 \psi_\nu := \exp_{\G_{2,\nu}} \circ \mathrm{pr}_\nu \circ \exp_{\G_2}^{-1}: \G_2 \rightarrow \G_{2,\nu}.
\end{equation}
For $\Gamma_{2,\nu} := \psi_\nu(\Gamma_2)$, we further define the compact Carnot nilmanifold $M_{2,\nu} := \Gamma_{2,\nu} \backslash \G_{2,\nu}$ over $\G_{2,\nu}$. \smallskip

Assume that a horizontal Clifford bundle $\Sigma_\delta^H M_2$ and a horizontal pull-back Dirac operator $D^H_{M_2}$ on $\Gamma_2 \backslash \G_2$, arising from a spin structure $\Sigma^{\T^d}_\delta$ on $\T^d$, are given. Now $\Sigma^{\T^d}_\delta$ also defines a horizontal Clifford bundle $\Sigma_\delta^H M_{2,\nu}$ on $M_{2,\nu}$, which is a vector bundle of the same rank as $\Sigma^H_\delta M_2$; and since the horizontal distributions of $M_2$ and $M_{2,\nu}$ can both be identified with $T\T^d$ we have a horizontal Clifford action on $M_{2,\nu}$ which coincides with the horizontal Clifford action on $M_2$, meaning $c^H_{M_2}(X_{1,j}) = c^H_{M_{2,\nu}}(\pr_\nu(X_{1,j}))$ for all $1 \leq 1 \leq d$ as endomorphisms on the vector bundles of rank $2^{[d/2]}$. Hence we have a horizontal pull-back Dirac operator
\begin{equation} \label{eq horizontal Dirac operator on M_2,nu}
 D^H_{M_{2,\nu}} := \sum_{j=1}^{d_1} c^H_{M_{2,\nu}}(\pr_\nu(X_{1,j})) \partial_{\pr_\nu (X_{1,j})}
\end{equation}
on $M_{2,\nu}$. But this operator allows us to prove the following lemma. \medskip

\begin{lemma} \label{lemma reduction to codimension 1 case on homogeneous spaces}
Let $D^H_{M_2}$ be the horizontal pull-back Dirac operator on the compact Carnot nilmanifold $M_2 = \Gamma_2 \backslash \G_2$, where $\G_2$ is a Carnot group of nilpotency step $2$, and for a $\nu \in \{1, \ldots, d_2\}$ let $D^H_{M_{2,\nu}}$ be the horizontal pull-back Dirac operator from \eqref{eq horizontal Dirac operator on M_2,nu} on the compact Heisenberg nilmanifold $M_{2,\nu} = \Gamma_{2,\nu} \backslash \G_{2,\nu}$ constructed above.\smallskip

Then if there is a $\nu \in \{1, \ldots, d_2\}$ such that the section $\sigma_\nu \in \Gamma^\infty(\Sigma_\delta^H M_{2,\nu})$ lies in the kernel of $D^H_{M_{2,\nu}}$, the section
$$\sigma := \sigma_\nu \circ \psi_\nu$$
is an element of the kernel of $D^H_{M_2}$.\smallskip

\B The idea is to use the expression from Theorem \ref{prop expression of horizontal Dirac via representation} of the horizontal Dirac operators involving the right-regular representation $R$ of the Carnot groups $\G_2$ and $\G_{2,\nu}$ on the spaces $L^2(\Sigma^H_\delta M_2)$ and $L^2(\Sigma^H_\delta M_{2,\nu})$. This means in the current situation that the operators $D_{M_2}^H$ and $D_{M_{2,\nu}}^H$, applied to $\sigma$ and $\sigma_\nu$, are given by
\begin{equation} \label{eq D^H_M2 via R}
 D^H_{M_2} \sigma = \sum_{j=1}^{d_1} c^H_{M_2}(X_{1,j}) \left. \frac{d}{dt} R(\exp_{\G_2} tX_{1,j}) \sigma \right|_{t=0}
\end{equation}
and
\begin{equation} \label{eq D^H_m2nu via R}
 D^H_{M_{2,\nu}} \sigma_\nu = \sum_{j=1}^{d_1} c^H_{M_{2,\nu}}\left(\pr_\nu(X_{1,j})\right) \left. \frac{d}{dt} R\left(\exp_{\G_{2,\nu}} t \pr_\nu(X_{1,j})\right) \sigma_\nu \right|_{t=0}.
\end{equation}
We want to show that for any $\sigma_\nu$ such that $D^H_{M_{2,\nu}} \sigma_\nu = 0$ we have $D^H_M \sigma = 0$, where 
$$\sigma = \sigma_\nu \circ \psi_\nu.$$
Following the discussion preceding this lemma, we have $c^H_{M_2}(X_{1,j}) = c^H_{M_{2,\nu}}(\pr_\nu(X_{1,j}))$ as endomorphisms on the vector bundles of rank $2^{[d/2]}$. Hence, from the expressions \eqref{eq D^H_M2 via R} and \eqref{eq D^H_m2nu via R} the desired statement will follow if we have 
\begin{equation} \label{eq desired identity for R}
 \left(R(\exp_{\G_2} tX_{1,j}) \sigma\right) (x) = \left(R\left(\exp_{\G_{2,\nu}} t \pr_\nu(X_{1,j})\right) \sigma_\nu\right) (\psi_\nu(x))
\end{equation}
for any $x \in M_2$. Will check this via a small calculation.\smallskip

Let $t \in \R$. Then for the left hand side of \eqref{eq desired identity for R} we get for every $\sigma \in \Gamma^\infty(\Sigma_\delta^H M_2)$, using exponential coordinates of $\G_2$ and the Baker-Campbell-Hausdorff formula (see Equation \eqref{eq Baker-Campbell-Hausdorff formula} in Section 2.2),
\begin{eqnarray*}
& & R(\exp_{\G_2} tX_{1,j}) \sigma \left(x^{(1)}, x^{(2)}\right)\\
&=& \sigma \left( \left(x^{(1)},x^{(2)}\right) . \exp_{\G_2} tX_{1,j}\right)\\
&=& \sigma \left( \exp_{\G_2}\left( \sum_{k=1}^{d_1} x_{1,k}X_{1,k} + \sum_{\mu=1}^{d_2}x_{2,\mu}X_{2,\mu}\right) . \exp_{\G_2} tX_{1,j} \right)\\
&=& \sigma \left( \exp_{\G_2}\left( \sum_{k=1}^{d_1} x_{1,k}X_{1,k} + \sum_{\mu=1}^{d_2}x_{2,\mu}X_{2,\mu} + tX_{1,j} + \frac{1}{2} \left[\sum_{k=1}^{d_1} x_{1,k}X_{1,k}, tX_{1,j}\right] \right) \right).  
\end{eqnarray*}
Now we can calculate the commutators occurring in the last line using for $1 \leq \mu \leq d_2$ the $\mu$-Levi forms (see Definition \ref{def generalized Levi form}) $L^{(\mu)}$ of $\G_2$, which gives us the identity
$$[X_{1,k},X_{1,j}] = \sum_{\mu=1}^{d_2} L_{jk}^{(\mu)} X_{2,\mu}$$
for all $1 \leq k \leq d$. Plugging this into the above calculation we get
\begin{equation} \label{eq R applied to sigma} 
\begin{split}
   & R(\exp_{\G_2} tX_{1,j}) \sigma \left(x^{(1)}, x^{(2)}\right)\\
   & = \sigma \left( \exp_{\G_2}\left( \sum_{k=1}^{d_1} x_{1,k}X_{1,k} + \sum_{\mu=1}^{d_2}x_{2,\mu}X_{2,\mu} + tX_{1,j} + \frac{1}{2} t \sum_{\mu=1}^{d_2} \sum_{k=1}^{d_2} x_{1,k} L_{kj}^{(\mu)} X_{2,\mu}  \right) \right).
\end{split}
\end{equation}
Now we set 
$$\sigma = \sigma_{\nu} \circ \psi_{\nu} = \sigma_{\nu} \circ \left(\exp_{\G_{2,\nu}} \circ \pr_{\nu} \circ \exp_{\G_2}^{-1} \right)$$
in \eqref{eq R applied to sigma} and get from the definition \eqref{eq projection algebra level for homogeneous Carnot spaces} of $\pr_\nu$
\begin{equation} \label{eq desired identity for R right hand side}
\begin{split}
 &  R(\exp_{\G_2} tX_{1,j}) \sigma \left(x^{(1)}, x^{(2)}\right)\\
 & = \sigma_\nu \left( \exp_{\G_{2,\nu}}\left( \sum_{k=1}^{d_1} x_{1,k}\pr_\nu(X_{1,k}) + x_{2,\nu} \pr_\nu(X_{2,\nu}) + t\pr_\nu(X_{1,j}) + \frac{1}{2} t \sum_{k=1}^{d_1} x_{1,k} L_{kj}^{(\nu)} \pr_\nu(X_{2,\nu})  \right) \right)
\end{split}
\end{equation}
since $\pr(X_{2,\mu}) = 0$ for all $\mu \neq 0$. Since $\pr_\nu$ is a Lie algebra homomorphism, we have
$$[\pr_\nu(X_{1,j}), \pr_\nu(X_{1,k})] = \pr_\nu( [X_{1,j}, X_{1,k}])$$
for all $1 \leq j,k \leq d_1$. Thus we see immediately (after using the Baker-Campbell-Hausdorff formula on $\G_{2,\nu}$ in the same way we did above on $\G_2$) that \eqref{eq desired identity for R right hand side} is exactly the right hand side of \eqref{eq desired identity for R}, and therefore we have
$$ D^H_{M_2} \sigma (x) = D^H_{M_{2,\nu}} \sigma_\nu (\pi_\nu(x)) = 0.$$
This shows the statement of the lemma. \eB
\end{lemma}\medskip

Now we can put the things together to prove that the horizontal pull-back Dirac operator we constructed has an infinite dimensional eigenspace on any Carnot group one can choose.\medskip

\begin{thm} \label{thm D^H has infinite dimensional kernel on any homogeneous Carnot space}
Let $M = \Gamma \backslash \G$ be the compact nilmanifold of a Carnot group $\G$ of rank $d$ and step $R$. Let $D^H$ be the horizontal pull-back Dirac operator acting on the horizontal Clifford bundle $\Sigma_\delta^H M$ which is arising from a spin structure $\Sigma^{\T^d}_\delta$ of the torus $\T^d$.\smallskip

Then the kernel of $D^H$ is infinite-dimensional. This means in particular that $D^H$ does not have a compact resolvent. \smallskip

 \B Let $\mathfrak{g} = \bigoplus_{S=1}^R V_S$ be the grading of the Lie algebra of $\G$. We consider the decomposition
$$L^2(\Sigma_\delta^H M) = \mathcal{H}_0 \oplus \bigoplus_{\tau \in \Z^{\dim V_R} \setminus \{0\}} \mathcal{H}_\tau$$
of $L^2(\Sigma_\delta^H M)$ from Theorem \ref{thm decomposition from center}, where all the spaces $\mathcal{H}_\tau$ are invariant under $D^H$. By the second statement of Theorem \ref{thm decomposition from center}, we have 
$$\mathcal{H}_0 \cong L^2 (\Sigma_\delta^H M_{R-1}),$$
where $M_{R-1}$ is the compact nilmanifold of the Carnot group $\G_{R-1} \cong \G / Z_R(\G)$ and the restriction of $D^H$ to $\mathcal{H}_0$ can be identified with a horizontal pull-back Dirac operator $\tilde{D}^H$ acting on $\Sigma^H_\delta M_{R-1}$. But this means that we can apply the same decomposition to the space $\mathcal{H}_0$, which contains a Hilbert space isomorphic to $L^2(\Sigma_\delta^H M_{R-2})$, with $M_{R-2}$ the compact nilmanifold of the step $R-2$ Carnot group $\G_{R-2} \cong \G_{R-1} \left/ Z_{R-1}(\G_{R-1}) \right.$, and so on.\smallskip

Inductively, we find a Hilbert space $\tilde{\mathcal{H}} \cong L^2(\Sigma_\delta^H M_2)$, where $M_2$ is the compact nilmanifold of a Carnot group $\G_2$ of step $2$, which is invariant under $D^H$ and on which $D^H$ can be identified with a horizontal pull-back Dirac operator $D^H_{M_2}$, acting on $\Sigma^H_\delta M_2$. But for this operator, we find an infinite dimensional kernel by Lemma \ref{lemma reduction to codimension 1 case on homogeneous spaces}: Since $\G$ is not abelian, we find a $\nu \in \{1, \ldots, d_2\}$, where $d_2$ is the dimension of the space $V_2$ from the grading $\mathfrak{g}_2 = V_1 \oplus V_2$ of the Lie algebra of $\G_2$, such that 
$$\G_{2,\nu} = \psi_\nu(\G_2) \cong \Hei^{2m+1} \times \R^{d-2m}$$
for some $m \geq 1$. We can define a horizontal Clifford bundle $\Sigma^H_\delta M_{2,\nu}$ and a horizontal pull-back operator $D^H_{M_{2,\nu}}$ from the corresponding objects on $M_2$ like we did in the discussion preceding Lemma \ref{lemma reduction to codimension 1 case on homogeneous spaces}.\smallskip

Now we know by Theorem \ref{thm spectral decomposition D^H on Heisenberg group} that $D^H_{M_{2,\nu}}$ has an infinite dimensional kernel, and that we can choose a basis $\{\tilde{\sigma}_1, \tilde{\sigma}_2, \ldots \}$ of this kernel. Using Lemma \ref{lemma reduction to codimension 1 case on homogeneous spaces}, we can lift this basis to an orthonormal system of infinitely many independent sections of $L^2(\Sigma_\delta^H M_2)$, given by $\sigma_j := \tilde{\sigma}_j \circ \psi_\nu$, which all lie in the kernel of $D^H$. (Note that these sections are indeed linear independent since this is the case for the $\tilde{\sigma}_j$'s and $\psi_\nu$ is a submersion.) Hence we have shown that $D^H$ has an infinite-dimensional kernel on $L^2(\Sigma_\delta^H M_2)$, and from the above argumentation this is also the case on $L^2(\Sigma_\delta^H M)$.\smallskip

The statement that $D^H$ cannot have a compact resolvent follows trivially from the fact that $\ker D^H$ is infinite-dimensional. \eB 
\end{thm}\medskip

We have shown by Theorem \ref{thm D^H has infinite dimensional kernel on any homogeneous Carnot space} that the horizontal pull-back Dirac operator $D^H$ does not furnish a spectral triple on arbitrary compact Carnot nilmanifolds. Theoretically it is possible to do spectral decompositions like in Section 4.3 for any given Carnot group and thus get statements about the asymptotic behavior of the non-degenerate eigenvalues: We have to know about the representation theory of $\G$.\smallskip

Now there is an algorithm to determine all the irreducible unitary representation of a Carnot group $\G$ (up to equivalence) developed by Alexander Kirillov (\cite{Kir1}), which is also referred to as the orbit method (see e.g. \cite{CG} or \cite{Kir2}). But since this algorithm makes use of the concrete structure of $\G$ it is hard to get general results in our context. And even for a given Carnot group which is not isomorphic to $\Hei^{2m+1} \times \R^{d-2m}$ we expect the calculations to be very long and complicated.\smallskip

But anyway we will show in the following chapters that the phenomenon of the degeneracy of the horizontal Dirac operator occurs in general, such that the description of the Carnot-Carath\'{e}odory geometry via spectral triples does not work as one would expect.\bigskip

\chapter{Calculus on Heisenberg Manifolds}

In the previous chapter we presented an explicit construction for horizontal Dirac operators on compact Carnot nilmanifolds and we saw that these operators do not have a compact resolvent. To put these observations into a greater generality, we want to adopt tools from pseudodifferential calculus. In the classical case, the Dirac operator $D$ on a compact manifold is an elliptic operator of order $1$, and therefore it follows from pseudodifferential theory that it admits a parametrix of order $-1$, which is compact because of the Sobolev embedding theorem. The existence of a parametrix leads to the possibility to construct complex powers within the calculus, which shows that the resolvent $(D^2+I)^{-1/2}$ of $D$ is compact. Now we intend to present something analogous for operators on Carnot manifolds, respecting the grading of a graded nilpotent Lie algebra.\smallskip

Indeed there is a pseudodifferential symbol calculus for Heisenberg manifolds: It has been developed simultaneously by Richard Beals and Peter Greiner (see \cite{BG}) and Michael Taylor (see \cite{Tay}) in the 1980s. In the last decade, some properties which are important for our work have been presented by Rapha\"{e}l Ponge (see \cite{Pon1}). We will see that in this calculus hypoellipticity takes the place of ellipticity, since hypoellipticity implies the existence of complex powers. In addition, we will see that on a compact manifold operators of negative order are compact, such that one can get a compact resolvent for a given operator of positive Heisenberg order. For hypoelliptic self-adjoint horizontal Laplacians on compact Heisenberg manifolds it is also known that their eigenvalues grow polynomial with a rate which gives back the graded dimension of the manifold. This can be seen as an analogy to the Weyl asymptotics in the elliptic case. \smallskip

In this chapter we give an overview over the Heisenberg calculus developed by Richard Beals and Peter Greiner, explain the composition of symbols and explore the role of hypoellipticity. For details we refer to the books by Beals and Greiner (\cite{BG}) and by Rapha\"{e}l Ponge (\cite{Pon1}). Afterwards we present results concerning the asymptotic growth of the eigenvalues and the existence of complex powers of hypoelliptic operators, which can be derived by an expansion of the heat kernel of the operator similarly to the classical case. Here we refer to \cite{BGS} and also to \cite{Pon1}. It is expected that there are analogous results for arbitrary Carnot manifolds, see \cite{Pon1}, but we did not find this generalization worked out in the literature and hence we restrict ourselves to Heisenberg manifolds. \bigskip

\section{The Heisenberg Calculus}

Let $M$ be a Heisenberg manifold of dimension $n = d+1$, which means by Section 2.3 that we have a grading $TM = HM \oplus VM$ of the tangent bundle such that $HM$ is a bracket generating horizontal distribution of rank $d$ and $VM = [HM,HM]$ is of rank $1$. An orthonormal frame for $HM$ shall be given by the vector fields $\{X_1, \ldots, X_d\}$, while $\{X_{d+1}\}$ shall span $VM$. As we know by Section 2.2, the graded (and therefore the Hausdorff) dimension of the metric space $(M,d_{CC})$, where $d_{CC}$ is the Carnot-Carath\'{e}odory metric, is equal to $d+2$.\smallskip

In Section 2.3 we have seen how one can identify the tangent space of $M$ with a bundle $\mathfrak{g} M$ of graded nilpotent Lie algebras, and hence we also have this structure on the cotangent bundle $T^\ast M$ which we denote by $\mathfrak{g}^\ast M$. At a point $x_0 \in M$, we have $\mathfrak{g}_{x_0}M \cong \R^n$ (as a vector space), and we have the dilations
\begin{equation} \label{eq dilation on cotangent Lie algebra bundle}
 \lambda . (\xi + \xi_{d+1}) = \lambda \xi + \lambda^2 \xi_{d+1}.
\end{equation}
for coordinates $(\xi, \xi_{d+1}) \in \R^n$ (with $\xi = (\xi_1, \ldots, \xi_d)$). We will further use the Koranyi gauge
\begin{equation} \label{eq Koranyi gauge Heisenberg}
 \left\| \xi \right\|_\Hei = \left( \sum_{j=1}^d |\xi_j|^4 + |\xi_{d+1}|^2 \right)^\frac{1}{4},
\end{equation}
see Definition \ref{def Koranyi gauge}.\smallskip

First we will consider $M = U \subset \R^n$ to be an open subset of $\R^n$; the generalization to vector bundles and manifolds will be the content of a theorem we mention at a later point of this section. Let $\sigma_j(x,\xi) = \sigma(-iX_j)$ denote the (classical) symbol of the vector fields $-iX_j$, and set $\sigma(x, \xi) := (\sigma_1(x, \xi), \ldots, \sigma_n(x,\xi))$. Note that by Proposition \ref{prop Carnot Lie algebra as vector fields} we can detect $\xi$ from the symbol $\sigma$ of the homogeneous operator. Using the Koranyi gauge $\| \cdot \|_\Hei$ and the notation
\begin{equation} \label{eq weighted multiindices}
 \langle \alpha \rangle = \sum_{j=1}^d \alpha_j + 2 \alpha_{d+1}
\end{equation}
for a multi-index $\alpha \in \N^{d+1}$ we take care of the homogeneity of the $X_j$ considered as differential operators. We can define the following symbol classes on which the Heisenberg calculus will be based (see \cite{BG}, (10.5)-(10.18), and \cite{Pon1}, Section 3.1.2).\medskip

\begin{defin} \label{def symbol classes Heisenberg calculus} 
 Let $U \subset \R^n$ be an open subset. Then:
\begin{enumerate}[(i)]
  \item For $m \in \Z$, we set
$$F_{\Hei,m}(U) := \left\{f \in C^\infty(U \times \R^n \setminus \{0\}): f(x, \lambda . \sigma) = \lambda^m f(x,\sigma) \ \forall \lambda > 0\right\}.$$
We will further denote the class of these homogeneous functions, which do not depend on $x$, by $F_{\Hei,m}$.
\item For $m \in \Z$, we set 
$$ F^m_\Hei (U) := \left\{ f \in C^\infty (U \times \R^n): f \sim \sum_{j=0}^\infty f_{m-j}, \ f_k \in F_{\Hei,k}(U) \right\},$$
where the asymptotic expansion $f \sim \sum f_{m-j}$ is meant in the sense that for all multi-indices $\alpha, \beta \in \N^n$ and all $N > 0$, we have
\begin{equation} \label{eq asymptotic expansion F^m}
\left| D_x^\alpha D_{\xi}^\beta \left(f(x,\sigma) - \sum_{j<N} f_{m-j}(x,\sigma) \right) \right| \leq C_{\alpha \beta N}(x) \| \sigma \|_{\Hei}^{m-N-\langle \beta \rangle}, 
\end{equation}
for a locally bounded function $C_{\alpha \beta N}$ on $U$. \smallskip

 We will further denote the class of these functions, which do not depend on $x$, by $F^m_\Hei$:
$$f \in F^m_\Hei \ \Leftrightarrow \ f \sim \sum_{j=0}^\infty f_{m-j}, \ f_k \in F_{\Hei,k}.$$

\item For $m \in \Z$, we set
$$S_{\Hei,m}(U) := \left\{q \in C^\infty(U \times \R^n \setminus \{0\}): \exists f \in F_{\Hei,m}(U) \ \text{with} \ q(x, \xi) = f(x, \sigma(x,\xi)) \right\}.$$

\item For $m \in \Z$, we set
$$S_{\Hei}^m(U) := \left\{q \in C^\infty(U \times \R^n): \exists f \in F_{\Hei}^m(U) \ \text{with} \ q(x, \xi) = f(x, \sigma(x,\xi)) \right\}.$$
For $f \sim \sum_{j=0}^\infty f_{m-j}$, the asymptotic expansion of $q(x,\xi) = f(x, \sigma(x,\xi))$ is given by
\begin{equation} \label{eq asymptotic expansion S^m}
  q \sim \sum_{j=0}^\infty q_{m-j} \ \ \ \text{with} \ q_k(x,\xi) = f_k(x, \sigma(x,\xi)).
\end{equation}
We call elements belonging to the class $S_{\Hei}^m(U)$ \emph{Heisenberg symbols} of order $m$.

\item The symbol class
$$S_{\Hei}^\infty(U) := \bigcup_{m \in \Z} S_{\Hei}^m(U) $$
induces the class of \emph{Heisenberg pseudodifferential operators} on $U$, which we denote by $\Psi_{\Hei}(U)$. \smallskip

 In detail, for $q \in S^m(U)$, the corresponding operator $\mathrm{Op}(q) \in \Psi^m_{\Hei}(U)$ is given by
\begin{equation} \label{eq Heisenberg Op(q)}
  \mathrm{Op}(q) u(x) = \int_{\R^n} e^{i \langle x,\xi \rangle} q(x,\xi) \hat{u}(\xi) \dbar \xi
\end{equation}
for any function $u \in C^\infty_c(U)$. On the other hand, if an operator $Q$ can be written in the form \eqref{eq Heisenberg Op(q)} for a function $q \in C^\infty(U \times \R^n)$, we call $\sigma_\Hei(Q) := q$ the Heisenberg symbol of $Q$.

\item The class
$$S_\Hei^{-\infty}(U) := \bigcap_{m \in \Z} S_\Hei^m(U)$$
is called the class of \emph{smoothing operators} in the Heisenberg calculus.
 \end{enumerate}
\eBsp
\end{defin}\medskip

Before we go on with the theory, we will state an example for Heisenberg pseudodifferential operators we already know.\medskip

\begin{example} \label{ex graded differential operators as Heisenberg pseudos} \normalfont
 Let $M = \R^n$ equipped with a frame $\{X_1, \ldots, X_{d+1}\}$ such that $\{X_1, \ldots, X_d\}$ forms a bracket generating horizontal distribution. Then any graded differential operator
$$D = p(X_1, \ldots, X_d, X_{d+1})$$
of order $\mu \in \N$ with respect to this frame (see Definition \ref{def graded differential operator}) is a Heisenberg pseudodifferential operator of Heisenberg order $\mu$. Its Heisenberg symbol is a polynomial 
$p(\sigma_1, \ldots, \sigma_{d+1}) \in F^\mu_{\Hei}(\R^n)$
of (homogeneous) degree $\mu$, and its asymptotic expansion is given by the homogeneous terms of this polynomial. In case the coefficients of $D$ are constant, we have a Heisenberg symbol belonging to the class $F^\mu_\Hei$.\smallskip

In particular, a horizontal Laplacian is a Heisenberg pseudodifferential operator of Heisenberg order $2$. \eBsp
\end{example}\medskip

We do not know yet if the expression \eqref{eq Heisenberg Op(q)} from the last item of Definition \ref{def symbol classes Heisenberg calculus} makes sense. But this will be the case, since a symbol of $S_{\Hei}^m(U)$ belongs to a Hörmander symbol class
\begin{equation} \label{eq Hormander rho delta}
  S^m_{\rho,\delta}(U) = \left\{q \in C^\infty(U \times \R^n): \left| D^\alpha_x D^\beta_\xi q(x,\xi) \right| \leq C_{\alpha \beta}(x) \left(1 + |\xi| \right)^{m + \delta|\alpha| - \rho|\beta|} \right\},
\end{equation}
where $C_{\alpha \beta}$ is again a locally bounded function on $U$. These symbol classes were established by Lars Hörmander in \cite{Hor1}, where it was also shown that the corresponding operators can be extended to bounded operators on certain Sobolev spaces. For the following theorem, we refer to \cite{BG}, Proposition (10.22).\medskip

\begin{thm} \label{thm embedding of S_Hei^m}
 For every $m \in \Z$, we have
$$S_{\Hei}^m(U) \subset \begin{cases} S^m_{\frac{1}{2},\frac{1}{2}}(U) & \text{for} \ m \geq 0 \\ S^{\frac{1}{2}m}_{\frac{1}{2}, \frac{1}{2}}(U) & \text{for} \ m < 0 \end{cases}.$$
\eB
\end{thm}\medskip

From this embedding, we see that the operator from equation \eqref{eq Heisenberg Op(q)} is well defined for every test function $u \in C^\infty_c(U)$. Moreover, we immediately get some regularity properties for Heisenberg pseudodifferential operators, which follow immediately from the corresponding regularity properties of operators belonging to the class $S^m_{1/2, \, 1/2}(U)$.\medskip

\begin{cor} \label{cor Heisenberg pseuDOs regularity}
 Let $m \in \Z$, $q \in S_{\Hei}^m(U)$. We can define an operator
$$ Qu(x) := \int_{\R^n} e^{i \langle x, \xi \rangle} q(x,\xi) \hat{u}(\xi) \dbar \xi $$
which has the following properties:
\begin{enumerate}[(i)]
 \item $Q: C^\infty_c(U) \rightarrow C^\infty(U)$ is a continuous linear operator.
 \item For every $s \in \R$, $Q$ can be extended to a bounded linear operator 
$$ Q: H^{s+m}(U) \rightarrow H^{s}(U),$$
 where $H^{s+m}(U)$ and $H^{s}(U)$ denote the $L^2$-Sobolev spaces.\smallskip

In particular, this means $Q \in \mathcal{B}(L^2(U))$ if $m \leq 0$. \eB
\end{enumerate}
\end{cor}\medskip

In \cite{BG}, it is also shown that the class $S^m_\Hei(U)$ does not depend on the choice of the frame $\{X_1, \ldots, X_d, X_{d+1}\}$ of $TM$, as long as it respects the grading structure (see \cite{BG}, Proposition (10.46)).\smallskip

After defining our symbol classes, we show that they induce a meaningful calculus. This is in general not the case for symbols of the class $S_{1/2, \, 1/2}(U)$, but it turns out that the composition of two Heisenberg pseudodifferential operators gives  a $\Psi_\Hei DO$ again. We briefly sketch how this composition is defined, referring to \cite{BG}, Chapters 12-14, and to \cite{Pon1}. Section 3.1.3, for the details. Assume first that the symbols of the operators which shall be composed are given by $p_1 \in S_{\Hei,m_1}(U)$ and $p_2 \in S_{\Hei,m_2}(U)$, where $p_k(x,\xi) = f_k(x,\sigma(x,\xi))$ for homogeneous functions $f_k \in F_{\Hei,m_k}$.

\begin{itemize}
 \item We fix an $x \in U$ and choose coordinates on $U$ which are centered at $x$. This provides us with symbols 
$$f_k^{(x)} \left(\sigma\right) := f_k\left(x, \sigma\right)$$
with $f_k^{(x)} \in F_{\Hei,m_1}$ for $k \in \{1,2\}$, defined on the dual $\mathfrak{g}^\ast_x U$ of the tangent graded Lie algebra $\mathfrak{g}_xU$ at $x$, which corresponds to a Carnot group $\G^x$ via the exponential mapping. Then $\G^x$ can be identified with the tangent Carnot group $\G_xU$ of $U$.

 \item After these identifications, we can use the Lie group composition $.$ of the tangent Carnot group $\G_xM$ to define a convolution product
\begin{equation} \label{eq convolution product Heisenberg symbols}
 \ast^{(x)}: F_{\Hei,m_1} \times F_{\Hei,m_2} \rightarrow F_{\Hei, m_1+m_2}
\end{equation}
 in the following way: For $f_k^{(x)}$, $k=1,2$, the operator associated to the symbol is simply the convolution operator
$$\mathrm{Op}\left(f_k^{(x)}\right) u(z) = \int_{\G_xM} \check{f}_k^{(x)}(y) u(z.y^{-1}) \dbar y,$$
acting on the Carnot group $\G_xU$. The composition of two convolution operators of this type gives the bilinear mapping \eqref{eq convolution product Heisenberg symbols} (see \cite{BG}, Proposition (12.14)).

\item Finally, one finds that the product $\ast^{(x)}$ from \eqref{eq convolution product Heisenberg symbols} depends smoothly on $x$ (see \cite{BG}, Proposition (13.3)), which gives us a continuous bilinear product
\begin{equation} \label{eq convolution product on open subset Heisenberg symbols}
\begin{split}
 \ast: & S_{\Hei, m_1}(U) \times S_{\Hei,m_2}(U) \rightarrow S_{\Hei, m_1+m_2}(U) \\
& \left(p_1 \ast p_2 \right)(x,\xi) = \left(f_1^{(x)} \ast^{(x)} f_2^{(x)} \right)(\sigma(x,\xi)).
\end{split}
\end{equation}
\end{itemize}

The above construction shows how homogeneous symbols can be composed. Now every Heisenberg symbol has an asymptotic expansion \eqref{eq asymptotic expansion S^m} into homogeneous symbols, and therefore one can find an asymptotic expansion for the composition of two arbitrary Heisenberg symbols, which is again a Heisenberg symbol and furnishes a $\Psi_\Hei DO$. The details for this expansion are formulated within the next theorem (see \cite{BG}, Theorems (14.1) and (14.7), and \cite{Pon1}, Proposition 3.1.9).\medskip

\begin{thm} \label{thm composition of Heisenberg pseudos}
 For $j = 1,2$ let $P_j \in \Psi_\Hei^{m_j}(U)$ have the symbol
$$p_j \sim \sum_{k \geq 0} p_{j,m_j-k}$$
in the sense of \eqref{eq asymptotic expansion S^m} and assume that one of these operators is properly supported. Then for the operator $P = P_1P_2$ we have $P \in \Psi_\Hei^{m_1+m_2}(U)$, and its symbol $p$ is given by
$$p \sim \sum_{k \geq 0} p_{m_1+m_2-k}$$
in the sense of \eqref{eq asymptotic expansion S^m}, where
\begin{equation} \label{eq composition Heisenberg pseudos homogeneous terms}
 p_{m_1+m_2-k}(x,\xi) = \sum_{k_1+k_2 \leq k} \sum_{\alpha,\beta,\gamma,\delta}^{(k-k_1-k_2)} h_{\alpha,\beta,\gamma,\delta}(x) \cdot \left(D^\delta_\xi p_{1,m_1-k_1}(x,\xi)\right) \ast \left(\xi^{\gamma} \partial_x^{\alpha} \partial_\xi^{\beta} p_{2,m_2-k_2}(x,\xi) \right).
\end{equation}
 In \eqref{eq composition Heisenberg pseudos homogeneous terms}, $\sum_{\alpha,\beta,\gamma,\delta}^{(l)}$ denotes the sum over all the indices such that
$$|\alpha| + |\beta| \leq \langle \beta \rangle - \langle \gamma \rangle + \langle \delta \rangle = l \ \ \ \text{and} \ \ \ |\beta| = |\gamma|,$$
and the functions $h_{\alpha,\beta,\gamma,\delta}$ are polynomials in the derivatives of the coefficients of the vector fields $X_1, \ldots, X_d, X_{d+1}$. \eB
\end{thm}\medskip

So far, we have developed the Heisenberg calculus on open subsets of $\R^n$. To define this calculus for operators acting on vector bundles over arbitrary Heisenberg manifolds, one has to show that the class of $\Psi_\Hei DO$s is invariant under the change of charts respecting the Heisenberg structure. In detail, we have the following theorem which shows us that we are able to extend the theory to the manifold case (see \cite{BG}, Theorem (10.67), and \cite{Pon1}, Proposition 3.1.18).\medskip

\begin{thm} \label{thm invariance theorem for Heisenberg calculus}
 Let $U_1$ and $U_2$ be open subsets on $\R^{d+1}$ together with hyperplane bundles $HU_1 \subset TU_1$ and $HU_2 \subset TU_2$, and let $\phi: (U_1,HU_1) \rightarrow (U_2,HU_2)$ be a Heisenberg diffeomorphism, which means we have $D\phi HU_1 = HU_2$ for the differential of $\phi$.\smallskip

 Then if $P_2 \in \Psi^m_\Hei$ is a Heisenberg pseudodifferential operator of order $m$ on $U_2$, the pullback $P_1 := \phi^\ast P_2$ of this operator to $U_1$ is a Heisenberg pseudodifferential operator of order $m$ on $U_1$. \eB
\end{thm}\medskip

Because of this theorem, we can consider Heisenberg pseudodifferential operators on Heisenberg manifolds $M$, acting on vector bundles $E$ from now on, and write $\Psi_\Hei(M,E)$ for this class of operators. In particular, we can derive the following consequence from Corollary \ref{cor Heisenberg pseuDOs regularity} which states that $\Psi_\Hei DO$s of negative order defined on a compact manifold are compact. It can be used if we want to show that a certain Heisenberg pseudodifferential operator has a compact resolvent. \medskip

\begin{cor} \label{cor Heisenberg pseuDOs of negative order are compact}
 For $m < 0$, let $Q \in \Psi^m(M,E)$ be a Heisenberg pseudodifferential operator acting on a vector bundle $E$ over a Heisenberg manifold $M$. Then $Q$ is a compact operator on the Hilbert space $L^2(M,E)$.\smallskip

\B By Corollary \ref{cor Heisenberg pseuDOs regularity} and Theorem \ref{thm invariance theorem for Heisenberg calculus}, for every $s \in \R$ the operator 
$$Q: H^{s+m}(M,E) \rightarrow H^s(M,E)$$
is bounded on the $L^2$-Sobolev spaces $H^{s+m}(M,E)$. We can now choose $s=-m$ and use the fact that the embedding $H^{-m}(M,E) \hookrightarrow H^0(M,E) = L^2(M,E)$ is compact by the Sobolev embedding theorem to derive that that operator
$$Q: L^2(M,E) \rightarrow H^{-m}(M,E) \hookrightarrow L^2(M,E)$$
is compact as the composition of a bounded and a compact operator. \eB
\end{cor}\medskip

After having established a composition rule inside the class of $\psi_\Hei DO$s in Theorem \ref{thm composition of Heisenberg pseudos}, the next thing to examine is the existence and the regularity of parametrices inside this class. First of all let us recall the definition of a parametrix. \medskip

\begin{defin} \label{def parametrix}
 Let $E$ be a vector bundle over a Heisenberg manifold $M$, and let $P: C_c^\infty(M,E) \rightarrow C^\infty(M,E)$ be a $\Psi_\Hei DO$. Then an operator $Q: C^\infty_c(M,E) \rightarrow C^\infty(M,E)$ with $Q \in \Psi_\Hei (M,E)$ is called a \emph{Heisenberg parametrix} or a \emph{Heisenberg pseudodifferential inverse} of $P$, if we have
$$QP = PQ = I \mod \Psi^{-\infty}_\Hei(M,E),$$
which means that the operators $PQ-I$ and $QP-I$ are smoothing. \eBsp
\end{defin}\medskip

In the classical calculus for symbols of the type $S_{1,0}$, a parametrix for a pseudodifferential operator $P$ exists if the operator is elliptic.  We will see that a necessary condition for the existence of a parametrix in the Heisenberg calculus is the hypoellipticity of the operator. The following classical definition of hypoellipticity is due to Lars Hörmander and can also be found in \cite{Pon1} and \cite{BG}.\medskip

\begin{defin} \label{def hypoelliptic}
 Let $P$ be a $\Psi_\Hei DO$ of order $m \in \Z$, acting on a vector bundle $E$ over a Heisenberg manifold $M$. Then $P$ is called \emph{hypoelliptic}, if for any distribution $u \in \mathcal{D}'(M,E)$ we have
$$Pu \in C^\infty(M,E) \ \ \ \Rightarrow \ \ \ u \in C^\infty(M,E).$$
In more detail, we call $P$ \emph{hypoelliptic with the loss of $k$ derivatives}, if we have for any $s \in \R$:
$$Pu \in H^s(M,E) \ \ \ \Rightarrow \ \ \ u \in H^{s+k}(M,E).$$
\eBsp
\end{defin}\medskip

In the classical case, a pseudodifferential operator is elliptic if its principal symbol is invertible. This generalizes to the case of Heisenberg pseudodifferential operators in terms of hypoellipticity: The main result will be that a $\Psi_\Hei DO$ is hypoelliptic if the $\Psi_\Hei DO$ associated to the principal part in the asymptotic (Heisenberg) expansion of its symbol is invertible in the Heisenberg calculus. We first introduce the notion of the principal symbol and the model operator, as it is done in \cite{BG} and \cite{Pon1}. \medskip

\begin{defin} \label{def principal symbol and model operator} 
Let $P$ be a $\Psi_\Hei DO$ of order $m \in \Z$ with Heisenberg symbol $p \in S^m_\Hei(M,E)$, acting on a vector bundle $E$ over a Heisenberg manifold $M$. For a point $a \in M$, we choose an open subset $U \ni a$ of $M$ and local coordinates of $U$ centered at $a$ to consider the asymptotic expansion
\begin{equation} \label{eq asymptotic expansion in local coordinates}
 p \sim \sum_{j=0}^\infty \tilde{p}_{m-j} \ \ \ \text{with} \ \tilde{p}_{m-j} \in S_{\Hei,m}(U,E)
\end{equation}
in the sense of Definition \ref{def symbol classes Heisenberg calculus}.\smallskip

For $\sigma(x,\xi) = (\sigma_1(x,\xi), \ldots, \sigma_{d+1}(x,\xi))$ with $\sigma_j(x,\xi) = \sigma(-iX_j)$ the classical symbol of $-iX_j$, we find a $p_m \in F^m_{\Hei}(U)$ such that in \eqref{eq asymptotic expansion in local coordinates} $\tilde{p}_m(x,\xi) = p_m(x,\sigma(x,\xi))$. Then the symbol
$$p^{(a)}_m(\sigma) := p_{m}(a,\sigma(a,\xi)) \ \in F_{\Hei m}(U).$$\smallskip
is called the \emph{(homogeneous) principal symbol} of $P$ at $a$. The corresponding operator
\begin{equation} \label{eq model operator}
 P^{(a)} := \mathrm{Op}\left(p^{(a)}_m\right) : C^\infty_c(\G_aM, E_a) \rightarrow C^\infty(\G_aM, E_a),
\end{equation}
where $\G_aM$ is the tangent Carnot group at $a$, is called the \emph{(homogeneous) model operator} of $P$ at $a$. \eBsp\smallskip

\Bem It is also possible to define a global principal symbol $\sigma_m(P) \in F_{\Hei,m}(M,E)$ on $M$ using the kernel representation of Heisenberg pseudodifferential operators, for which we refer to \cite{Pon1}, Theorem 3.2.2. \eBsp
\end{defin}\medskip

\begin{example} \normalfont \label{ex principal symbol graded differential operator}
 For a graded differential operator
$$D = p(X_1, \ldots, X_{d+1})$$
from Example \ref{ex graded differential operators as Heisenberg pseudos}, the principal symbol of $D$ at a point $a \in M$ is given by the leading homogeneous term of the polynomial $p$ after freezing the coefficients of $p$ in $a$. Hence the model operator of $D$ is given by the homogeneous graded differential operator with constant coefficients belonging to this homogeneous term.\smallskip

In particular, for the case of a horizontal Laplacian of the form
$$ \Delta^\mathrm{hor} = -\sum_{j=1}^d X_j^2 + B(x)X_{d+1} + \sum_{j=1}^d a_j(x)X_j + a(x), $$
the model operator of $\Delta^\mathrm{hor}$ at $a \in M$ is the operator
$$\Delta^\mathrm{hor}_\mathrm{mod} = -\sum_{j=1}^d X_j^2 + B(a) X_{d+1}.$$
\eBsp
\end{example}\medskip

One can show that the convolution of two principal symbols gives the principal symbol of the composition of the corresponding operators, and hence the composition of two model operators gives the model operator of the composition of the original operators. See for example \cite{Pon1}, Proposition 3.2.9.\smallskip

In particular, the existence of a parametrix for a Heisenberg pseudodifferential operator implies the existence of a parametrix of its model operator at each point. We even have equivalence for these two statements, and both statements imply the hypoellipticity of each model operator, which is formulated in the next theorem (see \cite{Pon1}, Proposition 3.3.1, Theorem 3.3.18 and Proposition 3.3.20). \medskip

\begin{thm} \label{thm hypoellipticity and existence of parametricies}
 Let $E$ be a vector bundle over a Heisenberg manifold $M$, and let 
$$P: C^\infty_c(M,E) \rightarrow C^\infty(M,E)$$
 be a $\Psi_\Hei DO$ of order $m \in \N$. Then the following statements are equivalent.
\begin{enumerate}[(i)]
 \item $P$ has a Heisenberg parametrix $Q$ with symbol $q \in S_\Hei^{-m}(M,E)$.
 \item At each point $a \in M$, the model operator $P^{(a)}$ from \eqref{eq model operator} of $P$ has a Heisenberg parametrix $Q^{(a)}$ with symbol $q^{(a)} \in F_{\Hei,-m}$.
 \item The global principal symbol $\sigma_m(P)$ of $P$ is invertible with respect to the convolution product for homogeneous symbols.
\end{enumerate}
If any of these conditions is fulfilled, $P$ is hypoelliptic with loss of $\frac{m}{2}$ derivatives. \eB
\end{thm}\medskip

We note that in \cite{Pon1} it is only shown that hypoellipticity is implied by the conditions (i) - (iii) of Theorem \ref{thm hypoellipticity and existence of parametricies}. But in the case of a graded differential operator $D$ with constant coefficients, the hypoellipticity of $D$ is equivalent to the invertibility of $D$ in the Heisenberg calculus. This can be seen by the so-called \emph{Rockland condition}, which we will introduce in Chapter 6: It is shown in \cite{Pon1}, that this Rockland condition (in a more general version compared to the one we will introduce) for a Heisenberg pseudodifferential operator $P$ is equivalent to the statements (i) - (iii) of Theorem \ref{thm hypoellipticity and existence of parametricies} (see \cite{Pon1}, Theorem 3.3.18). But for graded differential operators with constant coefficients, the Rockland condition is equivalent to the hypoellipticity of $D$, as we will see in Section 6.1.\smallskip

In case of a horizontal Laplacian of the form 
\begin{equation} \label{eq horizontal Laplacian Section 5.2}
 \Delta^\mathrm{hor} = -\sum_{j=1}^d X_j^2 + B(x)X_{d+1} + \sum_{j=1}^d a_j(x)X_j + a(x)
\end{equation}
it was proved in \cite{BG} that the hypoellipticity of this operator is equivalent to its invertibility in the Heisenberg calculus. Since we will work in particular with horizontal Laplacians, we state this theorem here (see \cite{BG}, Theorem (18.4)). \medskip

\begin{thm} \label{thm hypoellipticity and existence of parametricies horizontal Laplacian}
 For the horizontal Laplacian $\Delta^\mathrm{hor}$ from \eqref{eq horizontal Laplacian Section 5.2}, the following statements are equivalent:
\begin{enumerate}[(i)]
 \item $\Delta^\mathrm{hor}$ has a Heisenberg parametrix $Q \in \Psi_\Hei^{-2}(U)$
 \item At each $a \in U$, the model operator of $\Delta^\mathrm{hor}$ has a Heisenberg pseudodifferential inverse.
 \item $\Delta^\mathrm{hor}$ is hypoelliptic with loss of one derivative. \eB
\end{enumerate} 
\end{thm}\medskip

We will see in the next section that from the point of view of constructing and analyzing spectral triples using the techniques of Heisenberg calculus indeed hypoellipticity is the central property. Right now, we close this section with a remark concerning a possible generalization to arbitrary Carnot manifolds. \medskip

\begin{remark} \label{rem calculus on general Carnot manifold} \normalfont
 Let $M$ be an arbitrary Carnot manifold of step $R$ with grading $TM \cong V_1M \oplus \ldots \oplus V_RM$ such that for any $1 \leq S \leq R$ $X_{S,1}, \ldots, X_{S,d_S}$ forms an orthonormal frame for $V_SM$ (with $d_S = \rank V_SM$). Then we assume that we can generalize Definition \ref{def symbol classes Heisenberg calculus} as follows: For the Carnot group $\G$ corresponding to $M$ we consider the exponential coordinates $\xi \in \mathfrak{g} = \exp^{-1} \G$ with $\xi = (\xi^{(1)}, \ldots, \xi^{(R)})$, where $\xi^{(S)} \in \R^{d_S}$, the dilations
$$\lambda.\xi = \left( \lambda \xi^{(1)}, \lambda^2 \xi^{(2)}, \ldots, \lambda^R \xi^{(R)} \right)$$
and the Koranyi gauge
$$  \left\| x \right\|_\G := \left(\sum_{S=1}^R \sum_{j=1}^{d_S} \left| x_{S,j} \right|^\frac{2R!}{S}\right)^{\frac{1}{2R!}}, $$
see Definition \ref{def Koranyi gauge}. In addition, for multi-indices $\alpha \in \N^{d_1 + \ldots + d_R}$ we set
$$ \langle \alpha \rangle_\G := \sum_{S=1}^R S \cdot \sum_{j=1}^{d_S} \alpha_{S,j}.$$
Then we can define symbol classes in analogy to Definition \ref{def symbol classes Heisenberg calculus}, and it seems natural that it is possible to generalize the further definitions and theorems of this section to the general Carnot case. The proofs should work more or less analogous to the presentation in \cite{BG} and \cite{Pon1}, but it would by very laborious to write everything down in detail.\eBsp
\end{remark}\bigskip

\section{Complex Powers and Eigenvalue Asymptotics}

In this short section we present some results in Heisenberg calculus which open the door to constructing meaningful spectral triples from this calculus. As we will see, the crucial assumption in all these results is the hypoellipticity of the operator.\smallskip

First of all, we provide a theorem from which it can be shown that the operator $(D^2 + I)^{-1/2}$ is compact for a hypoelliptic self-adjoint operator $D$ of Heisenberg order $1$. We formulate a combinations of the Theorems 5.3.1 and 5.4.10 from \cite{Pon1}. \medskip

\begin{thm} \label{thm complex powers of Heisenberg pseudos}
 Let $M$ be a Heisenberg manifold. Suppose that $P$ is a hypoelliptic self-adjoint Heisenberg pseudodifferential operator of order $\nu \in \Z$ which is bounded from below and which is satisfying $\ker P = \{0\}$. Then, for any $s \in \C$, the operator $P^s$ defined via functional calculus is a $\Psi_\Hei DO$ of order $\nu s$. \eB \smallskip

\Bem Note that we have not introduced $\Psi_\Hei DO$s of non-integer order properly in this thesis. For this, we refer to \cite{Pon1}. Since in our case every operator appearing in this context will be of integer order, the theorems formulated in Section 5.1 are sufficient.\smallskip

Note also that the formulation in \cite{Pon1} of this theorem is more general in the sense that he does not assume that the kernel of $P$ is only assumed to be finite dimensional, which is the case because of the hypoellipticity of $P$. In this general case, one can also construct complex powers using projections onto the orthonormal complement of the kernel. \eBsp
\end{thm} \medskip

The next thing is to discover the asymptotic behavior of the eigenvalues of a self-adjoint Heisenberg pseudodifferential operator, from which one can detect the metric dimension of a spectral triple. In classical pseudodifferential calculus, one can consider the heat kernel of a positive self-adjoint elliptic operator $P$ to get asymptotics for the growth of its eigenvalues. This is done by expanding the trace of the operator $e^{-tP}$. Now something similar works for hypoelliptic self-adjoint horizontal Laplacians which are bounded from below: Mostly, this is the content of the paper \cite{BGS} by Richard Beals, Peter Greiner and Nancy Stanton. Some further considerations have been carried out by Rapha\"{e}l Ponge (see e.g. \cite{Pon1}). Without going into the details, we just state the results here. \smallskip

The following theorem shows how the trace expansion of the heat kernel of such a horizontal Laplacian looks like (see \cite{BGS}, Theorem (5.6), or \cite{Pon1}, Proposition 6.1.1).\medskip

\begin{thm} \label{thm heat trace expansion horizontal Laplacian}
 Let $\Delta^{\mathrm{hor}}$ be a hypoelliptic and self-adjoint horizontal Laplacian which is bounded from below, acting on a vector bundle $E$ over a Heisenberg manifold $M$ of dimension $d+1$ (which means that the horizontal distribution of $M$ has rank $d$). Then for $t \rightarrow 0^+$ we have the expansion
\begin{equation} \label{eq heat kernel expansion horizontal Laplacian}
 \mathrm{Tr} \; e^{-t\Delta^{\mathrm{hor}}} \sim t^{-\frac{d+2}{2}} \sum_{j=0}^\infty t^\frac{2j}{m} A_j(\Delta^\mathrm{hor})
\end{equation}
with $A_j(\Delta^\mathrm{hor}) = \int_M \tr_E a_j(\Delta^\mathrm{hor})(x) dx$, where $a_j$ can be computed from the term of degree $-2-2j$ in the asymptotic expansion of the symbol of the parametrix of $\Delta^\mathrm{hor}$ in local coordinates. \eB
\end{thm}\medskip

Now we denote by 
$$\lambda_0(\Delta^{\mathrm{hor}}) \leq \lambda_1(\Delta^{\mathrm{hor}}) \leq \ldots$$
the eigenvalues of $\Delta^{\mathrm{hor}}$, counted with multiplicity. A consequence of \eqref{eq heat kernel expansion horizontal Laplacian}, in connection with Karamata`s Tauberian Theorem, is that the following asymptotic behavior of these eigenvalues holds (see \cite{Pon1}, Proposition 6.1.2).\medskip

\begin{thm} \label{thm asymptotic behaviour of eigenvalues horizontal Laplacian}
 Let $\Delta^{\mathrm{hor}}$ be a hypoelliptic and self-adjoint horizontal Laplacian which is bounded from below, acting on a vector bundle $E$ over a Heisenberg manifold $M$ of dimension $d+1$. Then for $j \rightarrow \infty$ we have
\begin{equation} \label{eq asymptotic behaviour of eigenvalues horizontal Laplacian}
 \lambda_j(\Delta^{\mathrm{hor}}) \sim \left( \frac{j}{\nu_0\left(\Delta^{\mathrm{hor}}\right)} \right)^\frac{2}{d+2}
\end{equation}
for the eigenvalues $\lambda_j$ of $\Delta^{\mathrm{hor}}$, where $\nu_0(\Delta^{\mathrm{hor}})$ is a constant depending on the dimension $d$ of the horizontal distribution and the term $A_0(\Delta^{\mathrm{hor}})$ in the heat trace expansion \eqref{eq heat kernel expansion horizontal Laplacian} of $\Delta^{\mathrm{hor}}$. \eB
\end{thm}\medskip

Note that the above theorem just gives a qualitative statement about the growth of the eigenvalues, which suffices to detect the metric dimension of a spectral triple constructed from a first order hypoelliptic and self-adjoint operator on $M$. We will carry out this construction in Section 7.1. But first of all, we want to return to the horizontal Dirac operator and show that this theory cannot be applied to $D^H$, since $D^H$ is not hypoelliptic.\bigskip

\chapter{Hypoellipticity of Graded Differential Operators}

As we saw in the last chapter, the condition of hypoellipticity allows us to define complex powers of a self-adjoint Heisenberg pseudodifferential operator within the Heisenberg calculus. Since operators of negative Heisenberg order on a compact Carnot manifold are compact, one can argue that a hypoelliptic graded differential operator has a compact resolvent. This means that if a horizontal Dirac operator $D^H$ acting on a compact Carnot manifold $M$ is hypoelliptic, condition (ii) for a spectral triple will be fulfilled on $(C(M),L^2(M),D^H)$. But in Chapter 4 we already constructed a horizontal Dirac operator on an arbitrary compact Carnot nilmanifold which does not have a compact resolvent.\smallskip

The aim of this chapter is to show a generalization of the results of Chapter 4 in the setting of pseudodifferential calculus: We show that any horizontal Dirac operator on a Carnot manifold is not hypoelliptic. From this we will draw the conclusion that a horizontal Dirac operator does not have a compact resolvent, and hence does not furnish a spectral triple. \smallskip

We start this chapter by reviewing some well-known hypoellipticity criteria, starting with Hörmanders \emph{Sum-of-Squares Theorem} and leading to the \emph{Rockland Condition}, which states an equivalence between hypoellipticity of a graded differential operator and the non-degeneracy in the irreducible representations of its associated Lie algebra. Then we will draw special attention to the situation where this Lie algebra is a Heisenberg algebra, since in this situation one has a good classification for the hypoellipticity of horizontal Laplacians. This will pay off, because in the second section we develop a criterion to exclude hypoellipticity of a graded differential operator by reducing the case to the co-dimension $1$ case. Similar to the argument given in Chapter 4 for a specific example on compact Carnot nilmanifolds, we will make use of the submersions between graded nilpotent Lie algebras introduced in Section 2.4.\smallskip

Finally, in Section 6.3 we prove that any horizontal Dirac operator is not hypoelliptic. This follows quickly from the previous work by considering its square.\bigskip

\section{Some Classical Theorems}

The development of hypoellipticity criteria for certain differential operators was a great matter in the 1970s and 1980s, and there are some celebrated results. The origin of all these criteria is the famous \emph{sum-of-squares theorem} by Lars Hörmander (see \cite{Hor2}). \medskip

\begin{thm} \label{thm Hormander sum of squares}
 Let $X_0, X_1, \ldots, X_d$ be homogeneous vector fields on $\R^n$ with real $C^\infty$-coefficients on an open set $\Omega \subset \R^n$ and $c \in C^\infty(\Omega)$ real valued. Then the operator
$$P = \sum_{j=1}^d X_j^2 + X_0 + c$$
is hypoelliptic, if among the operators $X_j$ and all their commutators there exist $n$ which are linearly independent at any given point in $\Omega$. \eB
\end{thm}\medskip

The problem is that this theorem only works for vector fields with real coefficients, so it will not apply to our case of the square of a horizontal Dirac operator where the Clifford action causes complex coefficients. During the following years there were several generalizations of Hörmander's theorem, for example by Kohn (\cite{Koh1} and \cite{Koh2}) or Rothschild and Stein (\cite{RS}). Rothschild and Stein developed a close-to-complete characterization for the hypoellipticity of horizontal Laplacians of the form
\begin{equation} \label{eq hoizontal laplacian for rothschild stein}
 \Delta^{\mathrm{hor}}_{RS} = -\sum_{j=1}^d X_j^2 - \frac{i}{2} \sum_{j,k=1}^d b_{jk}[X_j,X_k],
\end{equation}
where the $X_j$ are homogeneous vector fields on $\R^n$ and $b = \left(b_{jk}\right) \in \mathrm{Skew}_{d \times d}(\R)$ is assumed to be a real skew-symmetric matrix. Then Rothschild and Stein proved the following theorem in a slightly more general version (see \cite{RS}, Theorem 1' and Theorem 2): \medskip

\begin{thm} \label{thm hypoellipticity Rothschild Stein}
 Consider the space
\begin{equation} \label{eq space R from rothschild stein}
 \mathcal{R} := \left\{ r = \left( r_{jk}\right) \in \mathrm{Skew}_{d \times d}(\R): \sum_{j,k = 1}^d r_{jk}[X_j,X_k] = 0 \right\},
\end{equation}
and its orthonormal complement $\mathcal{R}^\perp$ with respect to the inner product $(s_1,s_2) = - \mathrm{tr} (s_1s_2)$ on $\mathrm{Skew}_{d \times d}(\R)$.

\begin{enumerate} [(i)]
 \item If
$$\sup_{\rho \in \mathcal{R}^\perp, \ \| \rho \|_1 \leq 1} \left| \tr\left(b \rho\right) \right| < 1,$$
 then $\Delta^{\mathrm{hor}}_{RS}$ from \eqref{eq hoizontal laplacian for rothschild stein} is hypoelliptic.
 \item Assume that the Lie algebra $\mathfrak{g}$ spanned by $\{X_1, \ldots, X_d\}$ is graded, which means $\mathfrak{g} = \bigoplus_{S=1}^R V_R$ such that $[V_S,V_T] \subset V_{S+T}$ if $S+T \leq R$ and $[V_S,V_T] = 0$ if $S+T>R$ (see also Section 2.2). Then if 
 $$\sup_{\rho \in \mathcal{R}^\perp, \ \| \rho \|_1 \leq 1} \left| \tr\left(b \rho\right) \right| \geq 1$$
and if the algebra $\mathfrak{g}_2 := \mathfrak{g} \left/ \bigoplus_{S=3}^R V_S \right.$ is not the Lie algebra of a Heisenberg group $\Hei^{2m+1}$, $\Delta^{\mathrm{hor}}_{RS}$ from \eqref{eq hoizontal laplacian for rothschild stein} is not hypoelliptic.
\end{enumerate}
\eB \smallskip

\Bem For $\exp \mathfrak{g}_2 \cong \Hei^{2m+1}$, the situation is more involved: It can be shown that there are situations where $\mathfrak{g}_2$ is the Lie algebra of a Heisenberg group, in which the operator $\Delta^{\mathrm{hor}}_{RS}$ is hypoelliptic, even though
$$\sup_{\rho \in \mathcal{R}^\perp, \ \| \rho \|_1 \leq 1} \left| \tr\left(b \rho\right) \right| = 1.$$
For a more detailed treatment of the co-dimension $1$ case, we refer to Theorem \ref{thm hypoellipticity of laplacian codimension 1} below. \eBsp
\end{thm}\medskip

The idea for the proof of the second statement of this theorem is to reduce the situation to the case $R=2$, and to describe the operator $\Delta^\mathrm{hor}$ on $\mathfrak{g}_2$ via irreducible unitary representation of $\mathfrak{g}_2$. Doing this, it is possible to write down a function which is not $C^\infty$ but belongs to the kernel of $\Delta^\mathrm{hor}$, which is a contradiction to the hypoellipticity of the operator.\smallskip

Remember that we already used techniques from representation theory for the case of horizontal Dirac operators on compact Carnot nilmanifolds in Chapter 4, and it turned out that there is indeed a close connection between the hypoellipticity of graded differential operators and the representation theory of the underlying Carnot group. We will describe this connection now.\smallskip

As we noted in Section 2.2, a graded differential operator is an operator of the form
$$D = p\left(X_1, \ldots, X_d, X_{2,1}, \ldots, X_{R,d_R} \right),$$
where $p$ is a polynomial with matrix-valued $C^\infty$ coefficients and where $\{X_{S,1}, \ldots, X_{S,d_S}\}$ is a frame for the vector space $V_S$ appearing in the grading $\mathfrak{g} = V_1 \oplus \ldots \oplus V_R$ of a graded Lie algebra $\mathfrak{g}$. If the coefficients of $D$ are constant, This suggests to consider $D$ as an element of the universal enveloping algebra $\mathcal{U}(\mathfrak{g})$ of $\mathfrak{g}$. Since $\mathfrak{g}$ is nilpotent, the exponential mapping is an isomorphism from $\mathfrak{g}$ onto its Lie group $\G$ which is a Carnot group. Now let $\pi$ be an irreducible unitary representation of $\G$ on a Hilbert space $\mathcal{H}$, consisting of $L^2$-functions. Since $\pi$ gives rise to an irreducible, unitary representation of $\mathfrak{g}$ on $\mathcal{H}$ via
$$(d \pi(X) \varphi)(x) = \left. \frac{d}{dt} \pi( \exp t_X) \varphi (x) \right|_{t=0},$$
we obtain a representation $d \pi(D)$ of the graded differential operator $D$.\smallskip

Let $\hat{\G}$ denote the unitary dual of $\G$, which is the space of all irreducible unitary representations of $\G$. Using the above concept, there is a representation theoretic criterion which characterizes the hypoellipticity of graded differential operators completely. We will now state this criterion, which is called the \emph{Rockland condition} since it was developed by Charles Rockland \cite{Roc1} for Heisenberg algebras. B. Helffer and J. Nourrigat extended the Rockland condition to the case of arbitrary graded nilpotent Lie algebras (see \cite{HN} or also \cite{Rot}).\medskip

\begin{thm} \label{thm Rockland condition}
 Let $\mathfrak{g}$ be a graded nilpotent Lie algebra, $\G = \exp(\mathfrak{g})$, and let $L \in \mathcal{U}(\mathfrak{g})$ be a graded differential operator which is homogeneous of degree $m$. Then $L$ is hypoelliptic if and only if $\pi(L)$ is injective for all nontrivial $\pi \in \hat{\G}$. \eB
\end{thm}\medskip

Although there is a way to determine the irreducible unitary representations of an arbitrary nilpotent Lie group (see e.g. \cite{CG}, we also mentioned this at the end of Section 4.4), this is a quite difficult task for specific examples. Hence we follow the same approach as in Chapter 4: For the co-dimension $1$ case, it is not too difficult to formulate hypoellipticity criteria deduced from the Rockland condition. Then in the next section, we will see how this case enables us to make more general statements.\smallskip

For the rest of this section, we assume $n=d+1$ and $\mathfrak{g} \cong V_1 \oplus V_2$ with $d = \dim V_1$. Let $\{X_1, \ldots, X_{d+1}\}$ be a frame for $T\R^n$ such that $V_1 = \spa \{X_1, \ldots, X_d\}$ and $V_2 = \spa \{X_{d+1}\}$ for the representation of this Lie algebra as vector fields on $\R^n$ from Proposition \ref{prop Carnot Lie algebra as vector fields}. We want to check horizontal Laplacians of the form
\begin{equation} \label{eq horizontal Laplacian codimension 1}
 \Delta^{\mathrm{hor}} := -\sum_{j=1}^d X_j^2 -i A X_{d+1} + O_H(1),
\end{equation}
acting on a vector bundle $E$ of rank $p \in \N$ over an open subset $\Omega \subset \R^{d+1}$ with $A \in \mathrm{Mat}_{p \times p} (\C)$ for hypoellipticity. Here, the term $O_H(1)$ denotes a graded differential operator of order smaller or equal to $1$ (which means it is a differential operator of order $1$ only depending on the horizontal vector fields $X_1, \ldots, X_d$). To formulate the criterion we recall the notion of the Levi form from Chapter 2.3 (see Definition \ref{defin Levi form}), which is the bilinear form  
$$\mathcal{L}: V_1 \times V_1 \rightarrow V_2, \ \ \ (Y_1,Y_2) \mapsto [Y_1,Y_2] \mod V_1.$$
For $\mathcal{L}(X_j,X_k) = L_{ik} X_{d+1}$ with $L_{ik} \in \R$, we denote by $L = \left(L_{jk}\right)$ the antisymmetric matrix describing $\mathcal{L}$.\smallskip

The following result states that the hypoellipticity of $\Delta^{\mathrm{hor}}$ from \eqref{eq horizontal Laplacian codimension 1} only depends on how the eigenvalues of $A$ behave in comparison with the eigenvalues of $\mathcal{L}$. It is well known and can be found at various places in the literature (see e.g. \cite{Pon1} or \cite{BG}), but because of its importance for our future arguments we give a proof, orientated towards the one given in \cite{Pon1}. \medskip

\begin{thm} \label{thm hypoellipticity of laplacian codimension 1}
 Let $\Delta^{\mathrm{hor}}$ be the horizontal Laplacian given by \eqref{eq horizontal Laplacian codimension 1} with corresponding graded nilpotent Lie algebra $\mathfrak{g} = \spa \{X_1, \ldots, X_d\} \oplus \spa \{X_{d+1}\}$, Carnot group $\G = \exp \mathfrak{g}$ and Levi form $\mathcal{L}$, which is described by the Levi matrix $L \in \mathrm{Skew}_{d \times d}(\R)$ corresponding to this basis. The non-zero eigenvalues of $\mathcal{L}$ are denoted by $\pm i\lambda_1, \ldots, \pm i\lambda_m$ (including multiplicity) with $\lambda_j > 0$ for all $j \in \{1, \ldots, m\}$, where $2m \leq d$ is the rank of $\mathcal{L}$..\smallskip

Then the hypoellipticity of $\Delta^{\mathrm{hor}}$ can be characterized as follows:
\begin{enumerate}[(i)]
 \item If $\G \cong \mathbb{H}^{2m+1}$ with $\mathbb{H}^{2m+1}$ the $(2m+1)$-dimensional Heisenberg group (with $d = 2m$), then  $\Delta^{\mathrm{hor}}$ is hypoelliptic if and only if no eigenvalue of $A$ belongs to the set
\begin{equation} \label{eq singular set Heisenberg}
 \Lambda := \left\{ \pm \left( \frac{1}{2} \left\|L\right\|_1 + 2 \sum_{1 \leq j \leq m} \alpha_j \left|\lambda_j\right| \right): \alpha_j \in \N^{m} \right\}.
\end{equation}
 \item If $\G \cong \Hei^{2m+1} \times \R^{d-2m}$ with $2m<d$, then $\Delta^{\mathrm{hor}}$ is hypoelliptic if and only if no eigenvalue of $A$ belongs to the set
\begin{equation} \label{eq singular set non-Heisenberg}
 \Lambda := \left( -\infty, -\frac{1}{2} \left\|L\right\|_1 \right] \cup \left[ \frac{1}{2} \left\|L\right\|_1, \infty \right).
\end{equation}
\end{enumerate}
In both cases, $\left\|L\right\|_1 = \mathrm{tr}(|L|) = 2\sum_{j=1}^m |\lambda_j|$ denotes the trace norm of $L$.\smallskip

\B Since $\Delta^\mathrm{hor}$ is a horizontal Laplacian, by Theorem \ref{thm hypoellipticity and existence of parametricies horizontal Laplacian} its hypoellipticity is equivalent to the hypoellipticity of its model operator, which is in this case the homogeneous horizontal Laplacian
\begin{equation} \label{eq horizontal Laplacian codimension 1 homogeneous}
 \Delta^\mathrm{hor}_{\mathrm{mod}} = -\sum_{j=1}^d X_j^2 - i A X_{d+1}.
\end{equation}
Hence we only have to check the operator \eqref{eq horizontal Laplacian codimension 1 homogeneous} for hypoellipticity to prove the theorem, which can be done by the Rockland condition (see Theorem \ref{thm Rockland condition}). \smallskip

First of all, for every Carnot group $\G$ of co-dimension $1$ we have $\G \cong \mathbb{H}^{2m+1} \times \R^{d-2m}$, where $2m = \rank L \leq d$, by Proposition \ref{prop tangent Lie group bundle and rank of Levi form}. By the same proposition, there exists an orthonormal basis transformation of $V_1 = \spa \{X_1, \ldots, X_d\} \cong \R^d$ such that after this transformation we have
\begin{equation} \label{eq Levi form normal form}
 L = \begin{pmatrix} 0 & D & 0 \\ -D & 0 & 0 \\ 0 & 0 & 0 \end{pmatrix}
\end{equation}
for the Levi matrix of $\mathcal{L}$, where $D$ is a diagonal matrix with diagonal entries $\lambda_1, \ldots, \lambda_m > 0$. Since an orthonormal change of the frame of $V_1$ does not change the form of the horizontal Laplacian \eqref{eq horizontal Laplacian codimension 1 homogeneous}, as one can see after a small calculation, we can assume \eqref{eq Levi form normal form} to be the matrix of the Levi form $\mathcal{L}$ of $\G$. But this gives the commutator relations
$$[X_j, X_{j+m}] = \lambda_j X_{d+1} \ \ \ \text{and} \ \ \ [X_{j+m},X_j] = -\lambda_j X_{d+1}$$
if $1 \leq j \leq m$, while all the other commutators are zero. For $\lambda_1 = \ldots = \lambda_m = 1$ these are exactly the commutator relations for the Lie algebra of $\mathbb{H}^{2m+1} \times \R^{d-2m}$. Thus, the isomorphism $\phi: \mathbb{H}^{2m+1} \times \R^{d-2m} \rightarrow \G$ is coordinate-wise defined by
\begin{equation} \label{eq isomorphism codimension 1 group to Heisenberg group}
\phi(x_j) = \begin{cases}
             \sqrt{\lambda_j} \; x_j & \text{for } 1 \leq j \leq m\\
             \sqrt{\lambda_{j-m}} \; x_{j} & \text{for } m+1 \leq j \leq 2m\\
             x_j & \text{for } j>2m
            \end{cases},
\end{equation}
where we work on exponential coordinates of the Carnot groups. \smallskip

The next step is to show that the general case of $\Delta^{\mathrm{hor}}_{\mathrm{mod}}$ acting on a vector bundle $E$ of rank $p$ can be restricted to the scalar case: We can choose (point-wise) a basis of $E$ such that in this basis the matrix $A$ in \eqref{eq horizontal Laplacian codimension 1 homogeneous} is given by
$$A = \begin{pmatrix} \mu_1 & \ast & \ast \\ 0 & \ddots & \ast \\ 0 & 0 & \mu_p \end{pmatrix}, $$
where $\mu_1, \ldots, \mu_p$ are the eigenvalues of $A$. With respect to this basis the operator $\Delta^{\mathrm{hor}}_{\mathrm{mod}}$ takes the form
\begin{equation} \label{eq horizontal Laplacian triangular}
 \Delta^\mathrm{hor}_{\mathrm{mod}} = \begin{pmatrix} \Delta^\mathrm{hor}_1 & \ast & \ast \\ 0 & \ddots & \ast \\ 0 & 0 & \Delta^\mathrm{hor}_p \end{pmatrix},
\end{equation}
where for each $j$ with $1 \leq j \leq p$ we have a horizontal Laplacian
$$\Delta^\mathrm{hor}_j = - \sum_{j=1}^d X_j^2 - i \mu_j X_{d+1}$$
acting on scalar valued functions. Now it is obvious that the operator \eqref{eq horizontal Laplacian triangular} fulfills the Rockland condition from Theorem \ref{thm Rockland condition} if and only if each of the scalar operators $\Delta^\mathrm{hor}_j$ does (because otherwise the matrix would not be invertible which would be a contradiction to the injectivity of the corresponding irreducible representations). Therefore we can restrict our considerations to the scalar case, working with the eigenvalues $\mu_1, \ldots, \mu_p$ of $A$. From now on, we will assume $A$ is scalar.\smallskip

After the above simplifications, our task is to check the Rockland condition for the scalar case with Levi form of $\G$ given by \eqref{eq Levi form normal form}. We start with the situation where $2m=d$ (which means $\G \cong \mathbb{H}^{2m+1})$). For the Lie algebra of the Heisenberg group the representation theory is well known: Up to equivalence, the nontrivial irreducible representations of the operators $X_1, \ldots, X_{d+1}$, considered as basis vectors of the Heisenberg algebra $\mathfrak{h}_{2m+1}$ are of two types (see e.g. \cite{Pon1}, (3.3.5) - (3.3.7), or \cite{Fol}; compare also to the proof of Theorem \ref{thm spectral decomposition D^H on Heisenberg group}): 
\begin{enumerate}[(i)]
 \item The infinite dimensional representations $\pi_{t}: \mathbb{H}^{2m+1} \rightarrow \mathcal{U}(L^2(\R^m))$ of $\Hei^{2m+1}$, parametrized by $t \in \R \setminus \{0\}$, give rise to the representations 
\begin{eqnarray*}
 d\pi_{t}(X_j) f(\xi) & = & |t| \partial_{\xi_j} f(\xi) \ \ \ \text{for } 1 \leq j \leq m\\
 d\pi_{t}(X_{m+j}) f(\xi) & =  & it\xi_{j} f(\xi) \ \ \ \ \ \text{for } 1 \leq j \leq m\\
 d\pi_{t}(X_{2m+1}) f(\xi) & = & it|t| f(\xi).
\end{eqnarray*}
 \item The one dimensional representations $\pi_\zeta: \mathbb{H}^{2m+1} \rightarrow \C$ of $\Hei^{2m+1}$, parametrized by $\zeta \in \R^{2m} \setminus \{0\}$, give rise to the representations
\begin{eqnarray*}
 d\pi_\zeta(X_j) &=& i \zeta_j \ \ \ \text{for } 1 \leq j \leq 2m\\
 d\pi_\zeta(X_{2m+1}) &=& 0.
\end{eqnarray*}
\end{enumerate}
Under the Lie group isomorphism $\phi: \mathbb{H}^{2m+1} \rightarrow \G$ given by \eqref{eq isomorphism codimension 1 group to Heisenberg group}, the representations of the basis vectors of the Lie algebra $\mathfrak{g}$ we are interested in are given by the above representations applied to $\sqrt{\lambda_j}X_j$ and $\sqrt{\lambda_j}X_{m+j}$ for $1 \leq j \leq m$. Hence all the irreducible representation for $\Delta^\mathrm{hor}$ are
\begin{equation} \label{eq infinite dimensional representation horizontal Laplace}
  d\pi_{t,\xi}\left(\Delta^\mathrm{hor}_{\mathrm{mod}}\right) = -\sum_{j=1}^m \lambda_j \left(\left|t\right|^2 \partial_{\xi_j}^2 - t^2 \xi_j^2\right) + A t|t| \ \ \ \text{for} \ t \in \R \setminus \{0\}
\end{equation} 
and
\begin{equation} \label{eq finite dimensional representation horizontal Laplace}
 d\pi_\zeta\left(\Delta^\mathrm{hor}_{\mathrm{mod}}\right) = -\sum_{j=1}^{2m} \zeta_j^2 \ \ \ \text{for} \ \xi \in \R^{2m} \setminus \{0\}.
\end{equation}
Obviously \eqref{eq finite dimensional representation horizontal Laplace} is injective for any $\zeta \in \R^{2m} \setminus \{0\}$, so we only have to check the operators given by \eqref{eq infinite dimensional representation horizontal Laplace} for the Rockland condition. For these operators, their injectivity is equivalent to the injectivity of the operators
\begin{equation} \label{eq infinite dimensional representation horizontal Laplace harmonic oszillator}
 -\sum_{j=1}^m \lambda_j \left(\partial_{\xi_j}^2 - \xi_j^2 \right) \pm A.
\end{equation}
It is known that the spectrum of the harmonic oscillator $\sum_{j=1}^m \lambda_j \left(\partial_{\xi_j}^2 - \xi_j^2 \right)$ is exactly the set $\sum_{j=1}^m \lambda_j \left(1 + 2 \N\right)$ (see e.g. \cite{Pon1}, Section 3.4, or \cite{Fol}; compare also to the proof of Theorem \ref{thm spectral decomposition D^H on Heisenberg group}), and therefore the invertibility of \eqref{eq infinite dimensional representation horizontal Laplace harmonic oszillator} is equivalent to the condition
$$A \notin \left\{ \pm \sum_{j=1}^m \lambda_j \left(1 + 2\N\right) \right\} = \left\{ \pm \left(\frac{1}{2} \left\|L\right\|_1 + 2 \sum_{j=1}^m \alpha_j \left|\lambda_j\right| \right): \alpha_j \in \N^m \right\}.$$
But this shows statement (i) of the theorem for a scalar $A \in \C$, and because of the above simplifications statement (i) is also proved for $A \in \mathrm{Mat}_{p \times p}(\C)$. \smallskip

The second statement where $\G \cong \Hei^{2m+1} \times \R^{d-2m}$ with $2m < d$ can be proved via the Rockland condition in a similar way: One can use that the irreducible representations of the abelian group $\R^{d-2m}$ are simply the trivial ones. Note that for the case $A \in \R$ the second statement can also be proved using Theorem \ref{thm hypoellipticity Rothschild Stein}, if we set
$$ b_{j,k} = \begin{cases}
	      \frac{|\lambda_j| A}{\sum_{j=1}^m |\lambda_j|} & \ \text{for} \ 1 \leq j \leq m, \ k=m+j\\
	      -\frac{|\lambda_k| A}{\sum_{k=1}^m |\lambda_k|} & \ \text{for} \ 1 \leq k \leq m, \ j=k+m\\
	      0, & \ \text{otherwise}
             \end{cases}
$$
in \eqref{eq hoizontal laplacian for rothschild stein}. This works for both directions of the equivalence we want to prove since in this situation $\G$ is not a Heisenberg group. By a more general version of the Theorem by Rothschild and Stein (see \cite{RS}), the case $A = B+iC \in \C$ is also covered. \eB \smallskip

\Bem If we look at the proof of this theorem, we note that the strategy is very similar to the strategy of proving Theorem \ref{thm spectral decomposition D^H on Heisenberg group}: In both cases, we make use of the fact that we can describe our operator using the irreducible representations of $\G$. While in Theorem \ref{thm spectral decomposition D^H on Heisenberg group} we deduced that the horizontal pull-back Dirac operator on a compact Carnot nilmanifold has infinite dimensional eigenspaces, we will deduce from this theorem that the square of any horizontal Dirac operator (and hence the operator $D^H$ itself) on a Heisenberg manifold is not hypoelliptic in Section 6.3. \eBsp
\end{thm} \medskip

Before we close this section, we remark that of course one can combine Theorem \ref{thm hypoellipticity of laplacian codimension 1} with Theorem \ref{thm hypoellipticity Rothschild Stein} to deduce a better characterization of hypoellipticity for horizontal Laplacians of the form \eqref{eq hoizontal laplacian for rothschild stein}, acting on scalar-valued functions. Since this is straight-forward, we will not write it down here. \bigskip

\section{A Reduction Criterion for Non-hypoellipticity}

In this section, we want to develop a criterion which we can use to show that any horizontal Dirac operator on a Carnot manifold is not hypoelliptic. Note that Theorem \ref{thm hypoellipticity Rothschild Stein} is not very practical for this situation: First of all, it has to be extended to the case of vector bundles, and even after doing this is seems to be a quite complicated task to work with the spaces $\mathcal{R}$ from \eqref{eq space R from rothschild stein}, which depends heavily on the concrete structure of the Carnot group, in a general setting. But on the other hand Theorem \ref{thm hypoellipticity of laplacian codimension 1} can be applied easily to the square of such a horizontal Dirac operator, as long as we are in the case of Heisenberg manifolds.\smallskip

Thus the idea we want to follow is to reduce the general problem of excluding hypoellipticity on an arbitrary Carnot manifold to the co-dimension $1$ case. We have already seen that this approach works for the case of horizontal pull-back Dirac operators on compact Carnot nilmanifolds in Chapter 4: It was possible to lift the infinite dimensional eigenspaces of $D^H$ to a higher dimensional nilmanifold. In Chapter 4 we used the submersions between the corresponding Carnot group introduced in Section 2.4, and for our current (and more general) situation this strategy will also work.\smallskip

Like in Section 6.1, we consider a graded differential operator with constant coefficients
\begin{equation*} \label{eq graded differential operator}
 D = D(X_1, \ldots, X_n) \in \mathcal{U}(\mathfrak{g}) \otimes \C^p
\end{equation*}
as an element of the universal enveloping algebra of a graded Lie algebra $\mathfrak{g} = \bigoplus_{S=1}^R V_S$ belonging to a Carnot group $\G$ of step $R$. The (vector-valued) functions $D^H$ is acting on are supposed to be defined on a Carnot group $\G$, realized as a non-abelian group structure on $\R^n$. We remember the orthogonal projection
\begin{equation} \label{eq projection Lie algebra 2}
 \mathrm{pr}: \mathfrak{g} \rightarrow \tilde{V}, \ \ \ v \mapsto v \mod \tilde{V}^\perp,
\end{equation}
from Section 2.4, where
\begin{equation*} \label{eq projection Lie algebra subspace 2}
 \tilde{V} := \bigoplus_{S=1}^{M-1} V_S \oplus \tilde{V}_{M} 
\end{equation*}
for a linear subspace $\tilde{V}_M \subset V_M$ for some $1 \leq M \leq R$. By Proposition \ref{prop properties projection lie algebra}, $\tilde{\mathfrak{g}} = \pr(\mathfrak{g})$ has the structure of a nilpotent graded Lie algebra which is induced by $\mathfrak{g}$. After applying the map $\mathrm{pr}$ to the elements of $\mathfrak{g}$, we get a new differential operator
\begin{equation*} \label{eq projected graded differential operator}
 \pr(D) := D \left(\pr(X_1), \ldots, \pr(X_n)) \in \mathcal{U}(\tilde{\mathfrak{g}} \right) \otimes \C^q
\end{equation*}
considered as an element of the universal enveloping algebra of $\tilde{\mathfrak{g}}$. Note that of course the operator $\pr(D)$ is supposed to act on the (compared to $\G$ lower dimensional) Carnot group $\tilde{G} = \exp( \tilde{\mathfrak{g}})$, realized on the Euclidean space $\R^{\dim \tilde{\mathfrak{g}}}$. \smallskip

By Proposition \ref{prop properties projection lie algebra} the projection $\pr$ gives rise to a submersive Lie group homomorphism
\begin{equation*} \label{eq projection Lie groups 2}
 \psi := \exp_{\tilde{\G}} \circ \ \mathrm{pr} \ \circ \exp_{\G^{-1}}: \G \rightarrow \tilde{\G}
\end{equation*}
between the Carnot groups $\G = \exp \mathfrak{g}$ and $\tilde{\G} = \exp \tilde{\mathfrak{g}}$. The operator $\pr(D)$ is supposed to act on vector-valued functions living on the Carnot group $\tilde{\G}$. Since $\psi$ is a submersion, a function or distribution on $\tilde{\G}$ can be extended to a function or distribution on $\G$ via pullback along $\psi$. This observation leads us to the idea to deduce that $D$ is not hypoelliptic if $\pr(D)$ is not hypoelliptic for some projection of the type \eqref{eq projection Lie algebra 2}, as we will do via the following theorem.\medskip

\begin{thm} \label{thm non-hypoellipticity on vector bundles over carnot groups}
 Let $\G$ be a Carnot group with grading $\mathfrak{g} = \bigoplus_{S=1}^R V_S$ of its Lie algebra, and let $D \in \mathcal{U}(\mathfrak{g}) \otimes \C^p$ be a graded differential operator with constant coefficients.\smallskip

Assume that for some $2 \leq M \leq R$ there is a linear space $\tilde{V}_M \subset V_M$ such that the corresponding orthonormal projection
$$\mathrm{pr}: \mathfrak{g} \rightarrow \tilde{\mathfrak{g}} := \tilde{V} = \bigoplus_{S=1}^{M-1} V_S \oplus \tilde{V}_M, \ \ \ X \mapsto \tilde{X} = \pr (X)$$ 
defined by \eqref{eq projection Lie algebra 2} furnishes a graded differential operator 
$$\pr(D) = D(\tilde{X}_1, \ldots, \tilde{X}_n) \in \mathcal{U}(\tilde{\mathfrak{g}}) \otimes \C^q,$$
which is not hypoelliptic. Then $D$ is not hypoelliptic.\smallskip

\B We write down the argument for the case $p=1$ where $D$ is acting on complex-valued functions functions; the vector-valued case works analogously. Assume the operator $\tilde{D}$ constructed above is not hypoelliptic, which means that there exists a distribution $\tilde{\varphi} \in \mathcal{D}'(\tilde{\G})$ such that $\tilde{\varphi} \notin C^\infty(\tilde{\G})$, but $\tilde{D}\tilde{\varphi} \in C^\infty(\tilde{\G})$ for the Carnot group $\tilde{\G} = \exp(\tilde{\mathfrak{g}})$. We have to find a $\varphi \in \mathcal{D}'(\G)$ which is not a $C^\infty$-function, such that $D\varphi \in C^\infty(\G)$.\smallskip

We consider the submersive Lie group homomorphism
\begin{equation}
 \psi := \exp_{\tilde{\G}} \circ \ \mathrm{pr} \ \circ \exp_\G^{-1}: \G \rightarrow \tilde{\G}
\end{equation}
arising from $\mathrm{pr}$, see \eqref{eq projection Lie groups 2}. In addition we use the projection $\pr^\perp: \mathfrak{g} \rightarrow \ker \, (\pr)$, which is given by $\pr^\perp(v) = v \mod \tilde{V}$, of $\mathfrak{g}$ onto the kernel of $\pr$. Then from Proposition \ref{prop properties projection lie algebra} we get an isomorphism
\begin{equation} \label{eq isomorphism G to tildeG times N}
\alpha: \G \ \tilde{\rightarrow} \ \tilde{\G} \times N, \ \ \ x \mapsto (\tilde{x},n) = (\psi(x), \nu(x))
\end{equation}
from $\psi$, where $N = \ker (\psi)$ and 
$$\nu := \exp_N \circ \ \pr^\perp \ \circ \exp_{\G}^{-1}: \G \rightarrow N$$
is the projection onto the kernel of $\psi$, arising from $\pr^\perp$. From now on, we will use the coordinates $(\tilde{x},n)$ on $\G \cong \tilde{\G} \times N$ which are induced by the isomorphism $\alpha$ from \eqref{eq isomorphism G to tildeG times N}.\smallskip

If $x = \exp_{\G}X$ are exponential coordinates (corresponding to a vector field $X \in \mathfrak{g}$) on $\G$, then we have 
\begin{equation} \label{eq isomorphism G to tildeG times N exponential coordinates}
 \alpha\left( \exp_{\G} X \right) = \left(\exp_{\tilde{\G}} \pr(X), \exp_N \pr^\perp (X) \right),
\end{equation}
as one sees immediately from the definition of $\alpha$ via $\pr$. But from these exponential coordinates we see that for the differential of $\alpha$ we have
$$D\alpha(X) = \left( \pr(X), \pr^\perp(X) \right),$$
and hence we see how the application of a vector field $X \in \mathfrak{g}$ to a function $f \in C^\infty(\G)$ carries over to the push-forward $\alpha_\ast f$ of $f$ on $\tilde{\G} \times N$: We have
\begin{equation} \label{eq application vector field to tildeG times N}
 Xf(x) = \left(\pr(X), \pr^\perp(X) \right) \alpha_\ast f(\tilde{x},n),
\end{equation}
where $(\tilde{x},n) = \alpha(x)$. From now on, we will consider functions and distributions on $\G$ as functions and distributions on $\tilde{\G} \times N$, where the identification is given via the isomorphism $\alpha$. \smallskip

After these preparations, we are ready to prove the theorem. Let $\tilde{\varphi} \in \mathcal{D}'(\tilde{\G})$ be given such that $\tilde{\varphi} \notin C^\infty(\tilde{\G})$, but $\pr(D)\tilde{\varphi} \in C^\infty(\tilde{\G})$ for our graded differential operator $D$. For a test function $f \in C^\infty_c(\tilde{\G} \times N)$ we can define the push-forward $\psi_\ast f$ of $f$ along $\psi$ via
\begin{equation} \label{eq f_ast for test function f}
 \psi_\ast f(\tilde{x}) := \int_N f(\tilde{x},n) d\mu(n),
\end{equation}
where $\mu(n)$ is the Haar measure on the nilpotent Lie group $N$. Note that the expression \eqref{eq f_ast for test function f} is well-defined since $N$ is a normal subgroup of $\G$ as the kernel of the Lie group homomorphism $\psi$, and that $\psi_\ast f$ is a $C^\infty$-function with compact support on $\tilde{\G}$ since this is the case for $f$ on $\tilde{\G} \times N$. But this means we get a distribution $\varphi \in \mathcal{D}'(\tilde{\G} \times N)$ via
\begin{equation} \label{eq varphi^ast for distribution varphi}
 \langle \varphi,f \rangle := \langle \tilde{\varphi}, \psi_\ast f \rangle.
\end{equation}
Obviously, since we assumed $\tilde{\varphi} \notin C^\infty(\tilde{\G})$, we have $\varphi \notin C^\infty(\tilde{\G} \times N)$. We will show now that $D \varphi \in C^\infty(\G)$, from which the theorem will be proved. \smallskip

If we apply a vector field $X \in \mathfrak{g}$ to the distribution, which means by \eqref{eq application vector field to tildeG times N} to apply $(\pr(X), \pr^\perp(X)) \in \tilde{\mathfrak{g}} \times \mathfrak{n}$ to $\varphi \in \mathcal{D}'(\tilde{\G} \times N)$, we get from \eqref{eq varphi^ast for distribution varphi}
\begin{equation} \label{eq vector field on distribution first part}
\begin{split}
 \left\langle \left(\pr(X), \pr^\perp(X)\right) \varphi, f \right\rangle & = \left\langle \varphi, \left(\pr(X), \pr^\perp(X)\right)f \right\rangle\\
 & = \left\langle \tilde{\varphi}, \psi_\ast \left[\left(\pr(X), \pr^\perp(X)\right)f\right] \right\rangle 
\end{split}
\end{equation}
for any test function $f \in C^\infty_c(\tilde{\G} \times N)$. Now for the push-forward of $(\pr(X), \pr^\perp(X))f$ to $\tilde{\G}$ via $\psi$ we calculate, using \eqref{eq f_ast for test function f},
\begin{eqnarray*}
 \psi_\ast \left[\left(\pr(X), \pr^\perp(X)\right)f\right](\tilde{x}) &=& \int_N \left(\pr(X), \pr^\perp(X)\right)f(\tilde{x},n) d\mu(n)\\
&=& \int_N \left. \frac{d}{dt} f\left(\tilde{x} \, ._{\tilde{\G}} \, \exp_{\tilde{\G}} t \pr(X), \ n \, ._N \, \exp_N t \pr^\perp(X) \right) \right|_{t=0} d\mu(n) \\
&=& \left. \frac{d}{dt} \int_N f\left(\tilde{x} \, ._{\tilde{\G}} \, \exp_{\tilde{\G}} t \pr(X), \ n \, ._N \, \exp_N t \pr^\perp(X) \right) d\mu(n) \right|_{t=0} \\
&=& \left. \frac{d}{dt} \int_N f\left(\tilde{x} \, ._{\tilde{\G}} \, \exp_{\tilde{\G}} t \pr(X), \ n \right) d\mu(n) \right|_{t=0} \\
&=& \left.\frac{d}{dt} \psi_\ast f \left(\tilde{x} \, ._{\tilde{\G}} \, \exp_{\tilde{\G}} t \pr(X)\right) \right|_{t=0}\\
&=& \pr(X) \psi_\ast f(\tilde{x}). 
\end{eqnarray*}
Here we are allowed to interchange differentiation and integration since $f$ is a test function, and the forth equation is true since the Haar measure is translation invariant. The same calculation works if we apply another vector field $Y$ to $X\varphi$, which provides
\begin{equation} \label{eq derivatives of pushed forward functions}
 \psi_\ast \left[\left(\pr(Y), \pr^\perp(Y)\right) \left(\pr(X), \pr^\perp(X)\right)f\right](\tilde{x}) = \pr(Y) \pr(X) \psi_\ast f(\tilde{x}),
\end{equation}
and we can go on inductively for higher order applications of vector fields. Now we can apply these identities to \eqref{eq vector field on distribution first part} to find that
\begin{equation*} \label{eq vector fields on distributions}
 \left\langle \left(\pr(X), \pr^\perp(X)\right) \varphi, f \right\rangle = \left\langle \tilde{\varphi}, \pr(X) \psi_\ast f \right\rangle = \langle \pr(X) \tilde{\varphi}, \psi_\ast f \rangle,
\end{equation*}
and because of \eqref{eq derivatives of pushed forward functions} the same identity is true if we apply further vector fields to $X\varphi$. But this means that for any graded differential operator $D$ we get
\begin{equation} \label{eq graded differential operator lowers on distributions}
 \langle D\varphi,f \rangle = \langle \pr(D) \tilde{\varphi}, \psi_\ast f \rangle.
\end{equation} \smallskip

Now we assumed $\pr(D) \tilde{\varphi} \in C^\infty(\tilde{\G})$ to be a smooth function. But from this it follows immediately from \eqref{eq graded differential operator lowers on distributions} that $D\varphi$ is a smooth function on $\tilde{\G} \times N$: We have by this equation and by the definition of $\psi_\ast f$
\begin{eqnarray*}
 \langle D\varphi,f \rangle &=& \langle \pr(D) \tilde{\varphi}, \psi_\ast f \rangle\\
&=& \int_{\tilde{\G}} \pr(D) \tilde{\varphi}(\tilde{x}) \int_N f(\tilde{x},n) \ d \mu(n) d\mu(\tilde{x}) \\
&=& \int_{\tilde{\G} \times N} \pr(D) \tilde{\varphi} \left( \tilde{x} \right) f(\tilde{x},n) \ d\mu(n) d\mu(\tilde{x}),
\end{eqnarray*}
where $d\mu(n)$ and $d\mu(\tilde{x})$ are the corresponding Haar measures on $N$ and $\tilde{\G}$. But this shows that we must have
$$D\varphi(\tilde{x},n) = \pr(D) \tilde{\varphi} \left(\psi(\tilde{x},n)\right),$$
and since $\pr(D) \tilde{\varphi}$ is smooth this must also be the case for $D\varphi$.\smallskip

Altogether we have found a distribution $\varphi$ on $\G \cong \tilde{\G} \times N$, for which we have $D\varphi \in C^\infty(\G)$, but $\varphi \notin C^\infty(\G)$, whenever there is a projection $\pr$ of the type \eqref{eq projection Lie algebra 2} such that $\pr(D)$ is not hypoelliptic. But this shows that the graded differential operator $D$ is not hypoelliptic in this situation, such that the theorem is proved. \eB 
\end{thm}\medskip

\begin{remark} \normalfont
 We did the proof of Theorem \ref{thm non-hypoellipticity on vector bundles over carnot groups} in the very general setting of distributions. If in a more specific sense the operator $\pr(D)$ is assumed to be hypoelliptic in the sense that there is a function $\tilde{\varphi}$ on $\tilde{\G}$, which lies in the domain of $D$ but which is not $C^\infty$, such that $\pr(D) \tilde{\varphi} \in C^\infty(\tilde{\G})$, then the hypoellipticity of $D$ is neglected by the function $\varphi(x) = \tilde{\varphi}(\psi(x))$ on $\G$. The calculation is straight forward, making use of the isomorphism $\G \cong \tilde{\G} \times N$ like in the proof of the above theorem.\smallskip

 Via the same calculation one can show that $\pr(D) \tilde{\varphi} = 0$ implies $D\varphi = 0$. But from this, we immediately get Theorem \ref{thm D^H has infinite dimensional kernel on any homogeneous Carnot space}, which states the degeneracy of the horizontal pull-back operator on any compact Carnot nilmanifold (if we use the degeneracy in the Heisenberg case). This is no surprise, since the techniques of lifting an operator from a low-dimensional Carnot manifold to a higher dimensional one work very similar in both cases. In this way, the results from this chapter are the more general ones, while we get some extra features (like the detection of the metric dimension) from the specific examples in Chapter 4 which we were not able to derive in the general case. \eBsp
\end{remark}\medskip

In the next section, we want to use Theorem \ref{thm non-hypoellipticity on vector bundles over carnot groups} to show that any horizontal Dirac operator constructed in Section 3.2 is not hypoelliptic. The idea is to use Theorem \ref{thm non-hypoellipticity on vector bundles over carnot groups} after reducing the operator to an operator acting on a Carnot group of co-dimension $1$, which will be assumed to be of the type $\Hei^{2m+1} \times \R^{d-2m}$. For this operator, the hypoellipticity can be neglected by Theorem \ref{thm hypoellipticity of laplacian codimension 1}.\smallskip

Like before, let $d_S := \dim V_S$ be the dimension of $V_S$ and let $\left\{X_{S,1}, \ldots, X_{S,d_S}\right\}$ be an orthonormal basis of $V_S$ for the grading $\mathfrak{g} = \bigoplus_{S=1}^R V_S$ of $\mathfrak{g}$. Then a projection of $\mathfrak{g}$ onto a graded nilpotent Lie algebra of co-dimension $1$ can be constructed as follows: For any $\nu \in \{1, \ldots, d_2\}$, consider the $1$-dimensional linear subspace $\tilde{V}_{2,\nu} := \spa \{X_{2,\nu}\}$ of $V_2$. We then define $\tilde{V}_\nu := V_1 \oplus \tilde{V}_{2,\nu} \subset \mathfrak{g}$, such that our orthogonal projection \eqref{eq projection Lie algebra 2} becomes
\begin{equation} \label{eq projection onto codimension 1}
 \pr_\nu: \mathfrak{g} \rightarrow \tilde{V}_\nu, \ \ \ v \mapsto v \mod \left(\tilde{V}_\nu\right)^\perp.
\end{equation}
We can apply Proposition \ref{prop properties projection lie algebra} to $\pr_\nu$, which gives us graded nilpotent Lie algebras $\mathfrak{g}_{2,\nu} := \pr_\nu(\mathfrak{g})$ and $\mathfrak{n}_\nu := \ker (\pr_\nu)$ with corresponding Carnot groups $\G_{2,\nu}$ and $N_\nu$.\smallskip

Now let 
\begin{equation} \label{eq horizontal laplacian vector bundles}
 \Delta^\mathrm{hor} = -\sum_{j=1}^{d_1} X_j^2 - i \sum_{j<k} A_{j,k} [X_j,X_k] + O_H(1) \ \ \ \in \mathcal{U}(\mathfrak{g}) \otimes \C^q
\end{equation}
be a horizontal Laplacian acting on a vector bundle $E$ of rank $p$, which is considered as an element of the universal enveloping algebra of $\mathfrak{g}$ tensored with $\C^p$. The term $O_H(1)$ denotes a graded differential operator of order smaller or equal to $1$ (which is a first order differential operator depending only on the vector fields $X_1, \ldots, X_d$). The $A_{j,k}$ are $(p \times p)$-matrices with complex valued entries.\smallskip

Using the Lie algebra homomorphism $\pr_\nu$ from \eqref{eq projection onto codimension 1} we can define a horizontal Laplacian $\tilde{\Delta}^\mathrm{hor}_\nu \in \mathcal{U}(\mathfrak{g}_{2,\nu}) \otimes \C^p$, acting on a Carnot group with a horizontal distribution $\{\tilde{X_1}, \ldots, \tilde{X}_d\}$ of co-dimension $1$: For $\tilde{X}_j := \pr_\nu(X_j)$ we get
\begin{equation*}
  \tilde{\Delta}_\nu^\mathrm{hor} := \pr(\Delta^\mathrm{hor}) = -\sum_{j=1}^{d_1} \tilde{X}_j^2 - i \sum_{j<k} A_{j,k} \left( [\tilde{X}_j,\tilde{X}_k] \right) + O_H(1) \ \ \ \in \mathcal{U}(\mathfrak{g}_{2,\nu}) \otimes \C^p.
\end{equation*}
After calculating the commutators we have an $A_\nu \in \mathrm{Mat}_{p \times p}(\C)$ such that (for $d=d_1$)
\begin{equation} \label{eq reduced horizontal laplacian codimension 1}
 \tilde{\Delta}_\nu^\mathrm{hor} = -\sum_{j=1}^{d} \tilde{X}_j^2 - i A_\nu \tilde{X}_{d+\nu} + O_H(1),
\end{equation}
with $\tilde{X}_{d+\nu} = \pr_\nu(X_{2,\nu})$, which is a horizontal Laplacian for which its hypoellipticity can be determined by Theorem \ref{thm hypoellipticity of laplacian codimension 1}. \smallskip

This argument enables us to formulate non-hypoellipticity criteria involving the $\nu$-Levi form introduced in Section 2.3: In the above situation, for $\nu \in \{1, \ldots, d_2\}$ this is the bilinear form
\begin{equation} \label{eq nu-Levi form for non-hypoellipticity}
 \mathcal{L}_\nu: V_1 \times V_1 \rightarrow \spa\{X_{2,\nu}\}, \ \ \ (Y_1,Y_2) \mapsto [Y_1,Y_2] \mod \left( \spa \{X_{2,\nu}\} \right)^\perp.
\end{equation}
If we fix a basis for $V_1$, $\mathcal{L}_\nu$ is described by an antisymmetric matrix $L_\nu = \left(L_{jk}^{(\nu)}\right)$, such that $\mathcal{L}_\nu(X_j,X_k) = L_{jk} X_{2,\nu}$. Now we can use Theorem \ref{thm non-hypoellipticity on vector bundles over carnot groups} together with Theorem \ref{thm hypoellipticity of laplacian codimension 1} to formulate the following criterion.\medskip

\begin{cor} \label{cor non-hypoellipticity from reduced group of codimension 1}
 Let $\G$ be a Carnot group with corresponding Lie algebra $\mathfrak{g} = \bigoplus_{S=1}^R V_S$. For any number $\nu \in \{1, \ldots, \dim V_2\}$, we consider the projection $\mathrm{pr}_\nu$ from \eqref{eq projection onto codimension 1} of $\mathfrak{g}$ onto the Lie algebra $\mathfrak{g}_{2,\nu} = \mathrm{pr}_\nu(\mathfrak{g})$ together with its corresponding Carnot group $\G_{2,\nu}$. Let $\Delta^\mathrm{hor} \in \mathcal{U}(\mathfrak{g}) \otimes \C^p$ be a horizontal Laplacian of the type \eqref{eq horizontal laplacian vector bundles}, such that the horizontal Laplacian $\tilde{\Delta}^\mathrm{hor}_\nu = \mathrm{pr}_\nu(\Delta^\mathrm{hor})$ on $\G_{2,\nu}$ is given via
\begin{equation*}
 \tilde{\Delta}_\nu^\mathrm{hor} = -\sum_{j=1}^d \tilde{X}_j^2 - i A_\nu \tilde{X}_{d+\nu} +O_H(1),
\end{equation*}
like in \eqref{eq reduced horizontal laplacian codimension 1}, where $A_\nu \in \mathrm{Mat}_{p \times p}(\C)$ and $O_H(1)$ denotes a graded differential operator of order smaller or equal to $1$.\smallskip

Assume there is a $\nu \in \{1, \ldots, \dim V_2\}$ such that there is an eigenvector of $A_\nu$ which is contained in the singular set $\Lambda_\nu$ defined as follows:
\begin{enumerate}[(i)]
 \item If $\G_{2,\nu} \cong \mathbb{H}^{2m+1}$ is isomorphic to the $(2m+1)$-dimensional Heisenberg group (with $2m=d$), we define
\begin{equation*}
 \Lambda_\nu := \left\{ \pm \left( \frac{1}{2} \left\|L_\nu\right\|_1 + 2 \sum_{1 \leq j \leq m} \alpha_j \left|\lambda_j\right| \right): \alpha_j \in \N^{m} \right\},
\end{equation*}
where $\pm i \lambda_1, \ldots, \pm i \lambda_m$ denote the non-zero eigenvalues and $\left\| L_\nu \right\|_1 = 2 \sum_{j=1}^m |\lambda_j|$ denotes the trace norm of the $\nu$-Levi form $L_\nu$.
 \item If $\G \cong \Hei^{2m+1} \times \R^{d-2m}$ with $2m<d$ is not isomorphic to a Heisenberg group, then we define
\begin{equation*}
 \Lambda_\nu := \left( -\infty, -\frac{1}{2} \left\|L_\nu\right\|_1 \right] \cup \left[ \frac{1}{2} \left\|L_\nu\right\|_1, \infty \right),
\end{equation*}
 where $\left\| L_\nu \right\|_1$ denotes the trace norm of the $\nu$-Levi form $L_\nu$.
\end{enumerate}
Then $\Delta^\mathrm{hor}$ is not hypoelliptic.\smallskip

\Bem Note that for an arbitrary $\nu \in \{1, \ldots \dim V_2\}$, the cases $\G_{2,\nu} \cong \Hei^{2m+1}$ and $\G_{2,\nu} \cong \Hei^{2m+1} \times \R^{d-2m}$, with $m \geq 1$, are indeed the only possible cases since $\G_{2,\nu}$ is a Carnot group of step $2$ and horizontal co-dimension $1$.\smallskip

\B By Theorem \ref{thm non-hypoellipticity on vector bundles over carnot groups}, $\Delta^\mathrm{hor}$ is not hypoelliptic if $\tilde{\Delta}^\mathrm{hor}_\nu$ is not hypoelliptic for some $\nu \in \{1, \ldots, \dim V_2\}$. But to the operator $\tilde{\Delta}^\mathrm{hor}_\nu$ we can apply Theorem \ref{thm hypoellipticity of laplacian codimension 1}, which states that the operator $\tilde{\Delta}^\mathrm{hor}_\nu$ is not hypoelliptic if and only if there is an eigenvalue of $A_\nu$ belonging to one of the singular sets $\Lambda_\nu$ from (i) and (ii), depending on how the structure of $\G_{2,\nu}$ concretely looks like. \eB
\end{cor}\medskip

From Corollary \ref{cor non-hypoellipticity from reduced group of codimension 1}, we can immediately derive a simple criterion which ensures us that a horizontal Laplacian is not hypoelliptic if one of the matrices from \eqref{eq reduced horizontal laplacian codimension 1} has a certain eigenvalue. In detail we can use the matrices $L^{(\nu)}$ of the $\nu$-Levi forms $\mathcal{L}_\nu$ from \eqref{eq nu-Levi form for non-hypoellipticity} according to the frame $\{X_1, \ldots, X_d\}$ of $V_1$, such that for all $j,k \in \{1, \ldots, d\}$ the commutators have the form
$$[X_j,X_k] = \sum_{\nu_1}^{d_2} L_{jk}^{(\nu)} X_{2,\nu}.$$
Thus the horizontal Laplacian \eqref{eq horizontal laplacian vector bundles} can be rewritten in the form
\begin{equation} \label{eq horizontal laplacian vector bundles with Levi}
\begin{split}
  \Delta^\mathrm{hor} &= -\sum_{j=1}^{d_1} X_j^2 - i \sum_{\nu=1}^{d_2} \sum_{j<k} A_{j,k} L_{jk}^{(\nu)} X_{2,\nu} + O_H(1)\\
&= -\sum_{j=1}^{d_1} X_j^2 - i \sum_{\nu=1}^{d_2} A_\nu X_{2,\nu} + O_H(1),
\end{split}
\end{equation}
such that the matrices $A_\nu$ are given by
\begin{equation} \label{eq A_nu}
 A_\nu = \sum_{j<k} A_{j,k} L_{jk}^{(\nu)}.
\end{equation}
But from these matrices one can check immediately that a given horizontal Laplacian is not hypoelliptic: \medskip

\begin{cor} \label{cor non-hypoellipticity by eigenvalue of matrix}
 If there is a $\nu \in \{1, \ldots, \dim V_2\}$ such that there is an eigenvalue $\mu$ of the $(p \times p)$-matrix $A_\nu$ from \eqref{eq horizontal laplacian vector bundles with Levi} with 
 $$ \mu  = \pm \frac{1}{2} \left\| L_\nu \right\|_1 = \pm \sum_{j=1}^m |\lambda_j|,$$
then the operator $\Delta^\mathrm{hor}$ is not hypoelliptic. Here, for $j \in \{1, \ldots, m\}$ the numbers $\pm i \lambda_j$ are supposed to be the non-zero eigenvalues of the $\nu$-Levi form $\mathcal{L}_\nu$ of $\G$.\smallskip

\B Using the projection $\pr_\nu$ for any $\nu \in \{1, \ldots, d_2\}$ from \eqref{eq projection onto codimension 1}, we get from \eqref{eq horizontal laplacian vector bundles with Levi}
$$\pr_\nu\left(\Delta^\mathrm{hor}\right) = -\sum_{j=1}^{d_1} \pr_\nu(X_j)^2 - i A_\nu \pr_\nu(X_{2,\nu}) + O_H(1). $$
Now we can apply Corollary \ref{cor non-hypoellipticity from reduced group of codimension 1} to this operator and see that it is not hypoelliptic if one of the eigenvalues of $A_\nu$ belongs to one the sets $\Lambda_\nu$ from Corollary \ref{cor non-hypoellipticity from reduced group of codimension 1}. Since the numbers $\pm \frac{1}{2} \left\| L_\nu \right\|_1$ are included in both sets, the operator $\Delta^\mathrm{hor}$ fails to be hypoelliptic if one of them is an eigenvalue of one of the matrices $A_\nu$, no matter how the Carnot groups $\G_{2,\nu} = \psi_\nu(\G)$ looks like.  \eB
\end{cor}\bigskip

\section{Non-hypoellipticity of $D^H$}

Using the preparing work done in the last section we are finally ready to show that the horizontal Dirac operators we constructed in Chapter 3 cannot be hypoelliptic.\smallskip

We have constructed our horizontal Dirac operators $D^H$ more or less analogously to classical Dirac operators, acting on a horizontal Clifford bundle which is arising from the natural horizontal connection on a Carnot manifold $M$. We have seen that these operators detect the Carnot-Carath\'{e}odory metric and that they are therefore also a generalization of the classical case from the Connes metric point of view. If $D^H$ would be hypoelliptic, it would follow that it has a compact resolvent by Heisenberg calculus (which is developed for the case of Heisenberg manifolds in Chapter 5), and hence it would give a spectral triple. But this is not the case: We will prove that $D^H$ is not hypoelliptic, and therefore the machinery of graded pseudodifferential calculus does not work. We have already seen this in detail for a specific example in Chapter 4, and from the following theorems it will follow that these results fit into the general situation.\smallskip

To prove the non-hypoellipticity of $D^H$, we will consider its square calculated locally in Proposition \ref{prop D^H squared} from Chapter 3. This strategy is justified by the following simple observation. \medskip

\begin{prop} \label{prop hypoellipticity D and D square}
 Let $D$ be a (pseudo-)differential operator acting on a vector bundle $E$ over a manifold $M$. Then $D$ is hypoelliptic of and only if $D^2$ is hypoelliptic.\smallskip

\B Assume $D$ is hypoelliptic, which means for any distribution $\varphi$ we have that $D\varphi \in C^\infty(M,E)$ implies $\varphi \in C^\infty(M,E)$. This gives the implication
$$D^2 \varphi \in C^\infty(M,E) \ \Rightarrow \ D \varphi \in C^\infty(M,E) \ \Rightarrow \ \varphi \in C^\infty(M,E),$$
and hence $D^2$ is hypoelliptic.\smallskip

On the other hand, let $D^2 \varphi \in C^\infty(M,E)$ imply that $\varphi \in C^\infty(M,E)$. We have to show that $D$ is hypoelliptic. But this follows immediately since for $D\varphi \in C^\infty(M,E)$, we also have $D^2\varphi \in C^\infty(M,E)$, which shows $\varphi \in C^\infty(M,E)$ by assumption. Altogether the proposition is proved. \eB
\end{prop}\medskip

Since $(D^H)^2$ is a horizontal Laplacian, we can use Theorem \ref{thm hypoellipticity of laplacian codimension 1} to decide about its hypoellipticity for the case it is acting on a Heisenberg manifold. For the case of a general Carnot manifold, we can use the criterion developed in Section 6.2 to get a corresponding horizontal Laplacian on a Heisenberg manifold to which we can apply Theorem \ref{thm hypoellipticity of laplacian codimension 1}. We have already fulfilled this reduction step in Section 6.2, such that we can simply use the Corollaries \ref{cor non-hypoellipticity from reduced group of codimension 1} and \ref{cor non-hypoellipticity by eigenvalue of matrix} to prove the following theorem.\medskip

\begin{thm} \label{thm horizontal Dirac is not hypoelliptic}
 Let $M$ be a Carnot manifold with horizontal distribution $HM$, equipped with a horizontal Clifford bundle $S^H$ arising from the horizontal connection on $M$, and let $D^H$ be any horizontal Dirac operator acting on $S^HM$. Then $D^H$ is not hypoelliptic.\smallskip

\B Assume the grading of $TM$ is given by 
$$TM \cong HM \oplus V_2M \oplus \ldots \oplus V_RM,$$
and that we have a Riemannian metric on $M$ such that these sub-bundles are point-wise orthogonal to each other. Further we assume $\{X_1, \ldots, X_d\}$ (with $d = \rank HM$) to be an orthonormal frame for $HM$ and $\{X_{2,1}, \ldots, X_{2,d_2}\}$ to be an orthonormal frame for $V_2M$ (with $d_2 = \rank V_2M$). After fixing a $\nu_0 \in \{1, \ldots, d_2\}$ we can assume the frame of $HM$ to have the additional property that the matrix $L^{(\nu_0)} \in \mathrm{Skew}_{d \times d}(\R)$ describing the $\nu_0$-Levi form $\mathcal{L}_{\nu_0}$ with respect to this frame is given by
\begin{equation} \label{eq nu_0 Levi form for proof of non-hypoellipticity D^H}
 L^{(\nu_0)} = \begin{pmatrix} 0&D&0 \\ -D&0&0 \\ 0&0&0 \end{pmatrix},
\end{equation}
where $D$ is a diagonal matrix carrying the absolute values $\lambda_1, \ldots, \lambda_m$ of the non-zero eigenvalues of $\mathcal{L}_{\nu_0}$ on its diagonal. This can always be achieved by on orthonormal transformation of the horizontal frame because any $\nu$-Levi matrix is skew symmetric. For the proof of this theorem, we work with the expression of $D^H$ according to this horizontal frame $\{X_1, \ldots, X_d\}$, see e.g. Equation \eqref{eq D^H} from Theorem \ref{thm D^H selfadjoint}. \smallskip

We show that $(D^H)^2$ is not hypoelliptic, then the non-hypoellipticity of $D^H$ follows from Proposition \ref{prop hypoellipticity D and D square}. Note that it suffices to show the non-hypoellipticity locally in an environment of any point $x \in M$. By Proposition \ref{prop D^H squared} we have locally
\begin{equation} \label{eq D^H squared for proof non-hypoellipticity}
 \left(D^H\right)^2 = -\sum_{j=1}^d X_j^2 + \sum_{j<k} c^H(X_j)c^H(X_k) \left[X_j,X_k\right] + O_H(1) 
\end{equation}
with horizontal Clifford action $c^H: HM \rightarrow \, \en S^HM$ on the horizontal Clifford bundle $S^HM$. As before $X_j$ is to be understood as a component wise directional derivative in a local chart and $O_H(1)$ denotes a graded differential operator of order smaller than or equal to $1$.\smallskip

For any $\nu \in \{1, \ldots, \rank V_2M\}$ let $L^{(\nu)}$ denote the matrix of the $\nu$-Levi form corresponding to the frame $\{X_1, \ldots, X_d\}$ of $V_1M$. Hence for any pair $j,k \in \{1, \ldots, d\}$ we have
$$[X_j,X_k] = \sum_{\nu=1}^{d_2} L_{j,k}^{(\nu)} X_{2,\nu},$$
and plugging this into \eqref{eq D^H squared for proof non-hypoellipticity} we find that
\begin{eqnarray*}
 \left(D^H\right)^2 &=& -\sum_{j=1}^d X_j^2 + \sum_{\nu=1}^{d_2} \sum_{j<k} c^H(X_j)c^H(X_k) L_{j,k}^{(\nu)} X_{2,\nu} + O_H(1)\\
&=& -\sum_{j=1}^d X_j^2 - i \sum_{\nu=1}^{d_2} \sum_{j<k} i c^H(X_j)c^H(X_k) L_{j,k}^{(\nu)} X_{2,\nu} + O_H(1).
\end{eqnarray*}
Now we can use Corollary \ref{cor non-hypoellipticity by eigenvalue of matrix} applied to the matrix
$$A_{\nu_0} = \sum_{j<k} i c^H(X_j)c^H(X_k) L_{j,k}^{(\nu_0)},$$
where $\nu_0 \in \{1, \ldots, d_2\}$ is the number we chose in the beginning of the proof such that the $\nu_0$-Levi matrix $L^{(\nu_0)}$ has the form \eqref{eq nu_0 Levi form for proof of non-hypoellipticity D^H}. But this means that we have
\begin{equation} \label{eq A_nu_0 for proof non-hypoellipticity}
 A_{\nu_0} = \sum_{j=1}^m i c^H(X_j)c^H(X_{j+m}) \lambda_j,
\end{equation}
where $\lambda_1, \ldots, \lambda_m$ denote the absolute values of the non-zero eigenvalues of $\mathcal{L}_\nu$.\smallskip

According to Corollary \ref{cor non-hypoellipticity by eigenvalue of matrix}, we have to show that there is an eigenvalue of the matrix $A_\nu$ from \eqref{eq A_nu_0 for proof non-hypoellipticity} with absolute value $\sum_{j=1}^m \lambda_j$. But this has already been done by Proposition \ref{prop eigenvalues sum c(X_j)c(X_k)} from Section 3.2: This proposition states that $i \sum_{j=1}^m \lambda_j$ and $-i \sum_{j=1}^m \lambda_j$ are eigenvalues of the matrix $\sum_{j=1}^m \lambda_j c(X_j)c(X_{m+j})$. Hence the Theorem is proved because of Corollary \ref{cor non-hypoellipticity by eigenvalue of matrix}. \eB \smallskip

\Bem The property of being hypoelliptic only depends on the leading term of the local expression \eqref{eq D^H squared for proof non-hypoellipticity} for the square of $D^H$. As we already noted in the remark after Proposition \ref{prop D^H squared}, this term does not change if we modify $D^H$ by adding a section of the endomorphism bundle of $S^HM$, meaning an order zero term in the language of differential operators. Hence there is no chance of getting a hypoelliptic first order horizontal differential operator on $S^H M$ via a modification of the connection on $S^HM$ by an endomorphism. \eBsp
\end{thm}\medskip

In Chapter 5, Theorem \ref{thm complex powers of Heisenberg pseudos} we have seen that on a compact Heisenberg manifold $M$ we need the condition of being hypoelliptic to get a resolvent for a Heisenberg pseudodifferential operator of positive order which is compact. On the other hand, it is clear that from the existence of a compact resolvent of a differential operator $D$ of order $m$ one can expect this operator to be hypoelliptic: For every $s \in \R$, $D$ can be extended to a bounded operator
$$D: H^s(M) \rightarrow H^{s-m}(M)$$
between the $L^2$-Sobolev spaces $H^s(M)$ and $H^{s-m}(M)$. Hence we expect a (compact) resolvent of $D$ to be a mapping from $H^{s-m}(M)$ to $H^s(M)$ for any $s \in \R$. Now assume $D$ is not hypoelliptic, i.e. that there is an element $\varphi \notin C^\infty(M)$ such that $D\varphi \in C^\infty(M)$. But this means that the resolvent of $D$ maps the $C^\infty$-function $D\varphi$ (which is an element of $H^s(M)$ for every $s \in \R$) to a distribution which does not belong to $H^{s-m}(M)$ for one $s \in \R$. This is a contradiction.\smallskip

From the above discussion, we can formulate the following corollary.\medskip

\begin{cor} \label{cor horizontal Dirac has no compact resolvent}
Let $M$ be a Carnot manifold with horizontal distribution $HM$, equipped with a horizontal Clifford bundle $S^H$ arising from the horizontal connection on $M$, and let $D^H$ be any horizontal Dirac operator acting on $S^HM$. Then $D^H$ does not have a compact resolvent. \eB
\end{cor}\medskip

We emphasize once again that Theorem \ref{thm horizontal Dirac is not hypoelliptic} is true for any example of a horizontal Dirac operator according to the horizontal Levi-Civita connection on $HM$ one can imagine. It is true for the horizontal pull-back Dirac operators we discussed in Chapter 4, and it is also true for horizontal Dirac operators defined on classical Clifford or spinor bundles (by skipping the non-horizontal derivatives), see Proposition \ref{prop horizontal Clifford bundle from spinor bundle}. In particular, it is true for the operator $d^H + (d^H)^\ast$, where $d^H$ is the horizontal exterior derivative (see Example \ref{ex exterior bundle}). Note that all these constructions are based on the horizontal connection on $M$, which is induced by the Levi-Civita connection.\smallskip

This shows that the degeneracy we detected in Chapter 4 for the horizontal pull-back Dirac operator on compact Carnot nilmanifolds is not because of a bad choice for $D^H$. Rather it reflects a general phenomenon: The natural differential operator, which detects the Carnot-Carath\'{e}odory metric on a Carnot manifold via Connes' formula, does not give a spectral triple, and the classical construction of a spectral triple on a compact spin manifold cannot be transported to the Carnot case.\bigskip

\chapter{Spectral Triples from Horizontal Laplacians}

In the last chapter we saw that the canonical candidate for a spectral triple over a Carnot manifold detecting the horizontal geometry, the horizontal Dirac operator $D^H$, reproduces the Carnot-Carath\'{e}odory metric but is not hypoelliptic, and therefore does not have a compact resolvent in the Heisenberg calculus. But on the other hand the classical hypoellipticity criteria imply that there are a lot of horizontal Laplacians which are hypoelliptic, and in addition (at least in the Heisenberg case) give back the Hausdorff dimension of $(M,d_{CC})$ via their eigenvalue asymptotics.\smallskip

Now this chapter is devoted to studying horizontal Laplacians and to discussing how they furnish the geometry in the sense of Connes. In the first section we show explicitly how one can construct a spectral triple from a positive hypoelliptic horizontal Laplacian using the Heisenberg calculus. We present some operators which are induced by a small perturbation of the square of a horizontal Dirac operator. Then in Section 7.2 we will show that any horizontal Laplacian detects the Carnot-Carath\'{e}odory metric of a Carnot manifold via a formula similar to Connes' one. We also discuss what this means for a first order operator whose square is a horizontal Laplacian.\smallskip

The last section of this chapter is rather speculative: We intend to give some ideas how one can find estimates for the Connes metric of spectral triples arising from horizontal Laplacians towards the Carnot-Carath\'{e}odory metric (which is exactly detected by a horizontal Dirac operator $D^H$). Therefore we use some of the observations we made in Section 1.2 concerning the convergence of the metrics belonging to compact quantum order unit spaces. But sadly we have not been able to prove such an estimate, so there are a few open problems formulated in that section. \bigskip

\section{Spectral Triples via Heisenberg Calculus}

We now point out how we can construct spectral triples on compact Heisenberg manifolds using the Heisenberg calculus we introduced in Chapter 5. The idea is to start with a hypoelliptic, positive horizontal Laplacian $\Delta^\mathrm{hor}$ and to consider the operator $D_\mathrm{hor} = \sqrt{\Delta^\mathrm{hor}}$, defined via functional calculus.\smallskip

To show that such a construction indeed furnishes a spectral triple, we check that the commutator of any Heisenberg pseudodifferential operator of order $1$ with a smooth function $f$ is bounded.\medskip

\begin{prop} \label{prop boundedness of commutator}
 Let $M$ be a Heisenberg manifold and $P \in \Psi_\Hei^1(M,E)$ a Heisenberg pseudodifferential operator of order $1$ acting on a vector bundle $E$ over $M$. For any function $f \in C^\infty(M)$ we denote by $M_f$ the operator of multiplication by $f$ on the Hilbert space $L^2(M,E)$. Then the commutator $[P,M_f]$ is bounded. \smallskip

\B We assume that $E$ is the trivial line bundle over $M$, the general case works analogously. Our strategy is to show that the symbol of the commutator $[P,M_f]$ lies in $S_\Hei^0(M)$, and hence provides a Heisenberg pseudodifferential operator of order $0$. Since $S_\Hei^0(M) \subset S^0_{1/2, \, 1/2}$ for the Hörmander class $S^0_{1/2,1/2}$ by Theorem \ref{thm embedding of S_Hei^m}, this implies by classical pseudodifferential calculus that $[P,M_f]$ is bounded (see also Corollary \ref{cor Heisenberg pseuDOs regularity}).\smallskip

First of all we note that the multiplication operator $M_f$ is a Heisenberg pseudodifferential operator of order $0$, and its symbol is exactly the function $f$ by the definition of the symbol classes (see Definition \ref{def symbol classes Heisenberg calculus}). Let $p \in S^1_\Hei(M)$ denote the symbol of $P$. Then we can calculate the symbol of the commutator $[P,M_f]$ using Theorem \ref{thm composition of Heisenberg pseudos}: For general Heisenberg symbols $q_1$ of order $m_1$ and $q_2$ of order $m_2$, the asymptotic expansion of the symbol 
$$q = q_1 \# q_2 \sim \sum_{k \geq 0} q_{m_1+m_2-k},$$
 denoting the symbol of the composition of $\mathrm{Op}(q_1) \mathrm{Op}(q_2)$, is given by the terms
\begin{equation} \label{eq composition Heisenberg pseudos homogeneous terms section 7.1}
 q_{m_1+m_2-k}(x,\xi) = \sum_{k_1+k_2 \leq k} \sum_{\alpha,\beta,\gamma,\delta}^{(k-k_1-k_2)} h_{\alpha,\beta,\gamma,\delta}(x) \cdot \left(D^\delta_\xi q_{1,m_1-k_1}(x,\xi)\right) \ast \left(\xi^{\gamma} \partial_x^{\alpha} \partial_\xi^{\beta} q_{2,m_2-k_2}(x,\xi) \right),
\end{equation}
 where $\sum_{\alpha,\beta,\gamma,\delta}^{(l)}$ denotes the sum over all the indices such that
\begin{equation} \label{eq relation alpha beta gamma delta}
 |\alpha| + |\beta| \leq \langle \beta \rangle - \langle \gamma \rangle + \langle \delta \rangle = l \ \ \ \text{and} \ \ \ |\beta| = |\gamma|,
\end{equation}
and the functions $h_{\alpha,\beta,\gamma,\delta}$ are polynomials in the derivatives of the coefficients of the vector fields $X_1, \ldots, X_d, X_{d+1}$ forming a graded frame for $TM$. The operation $\ast$ denotes the operation of point-wise convolution, varying smoothly over $x$ (see Equations \eqref{eq convolution product Heisenberg symbols} and \eqref{eq convolution product on open subset Heisenberg symbols} from Section 5.1). See the discussion before Theorem \ref{thm composition of Heisenberg pseudos} for a more detailed explanation.\smallskip

To get the boundedness of the commutator we only have to show that the order $1$ part of this asymptotic expansion vanishes, since all the lower order terms have maximal Heisenberg order $0$ which gives rise to a bounded operator by Corollary \ref{cor Heisenberg pseuDOs regularity}. Since $m_1+m_2 = 1$ in our situation, these are exactly the terms in \eqref{eq composition Heisenberg pseudos homogeneous terms section 7.1} such that $k=0$. But this means $k_1=k_2 = 0$ and therefore also $\alpha, \beta, \gamma, \delta = 0$ by \eqref{eq relation alpha beta gamma delta}. Now the Heisenberg symbol of the commutator $[P,M_f] = Pf - fP$ is given by
$$ \sigma_\Hei([P,M_f]) = p \# f - f \# p,$$
and for its leading (Heisenberg order $1$) symbol we have because of the above argumentation
\begin{equation} \label{eq leading term for symbol commutator}
 \left(p\# f - f\# p\right)_1(x,\xi) = h_{0,0,0,0}(x) \cdot \left( p_1(x,\xi) \ast f(x) - f(x) \ast p_1(x,\xi) \right), 
\end{equation}
where $p_1$ denotes the leading term of the Heisenberg symbol of $P$. Note that the polynomial $h_{0,0,0,0}(x)$ only depends on the coefficients of the vector fields $X_1, \ldots, X_{d+1}$, but not on the symbols of the operators to be composed.\smallskip

Hence all we have to show is that the convolution of $p_1$ and $f$ commutes. But this is clear by the definition of the convolution: Whenever we fix an $x \in M$, the terms commute point-wise in $x$ since in this case $f(x)$ is a constant. Now since $\ast$ varies smoothly in $x$, this shows that 
$$p_1(x,\xi) \ast f(x) - f(x) \ast p_1(x,\xi) = 0,$$
and from \eqref{eq leading term for symbol commutator} we get $\left(p\# f - f\# p\right)_1 = 0$.\smallskip

We have seen that the leading term in the asymptotic expansion of the symbol belonging to the commutator $[P,M_f]$ has mostly Heisenberg order zero, and we can conclude that this commutator is bounded by Corollary \ref{cor Heisenberg pseuDOs regularity}. Therefore the proposition is proved. \eB
\end{prop}\medskip

To check the compactness of the resolvent, we refer to Theorem \ref{thm complex powers of Heisenberg pseudos}, which states that for a hypoelliptic self-adjoint $\Psi_\Hei DO$ $P$ of order $\nu$ the operator $P^s$ is a $\Psi_\Hei DO$ of order $\nu s$. Then it is clear that $\Psi_\Hei DO$s of order one can provide spectral triples. \medskip

\begin{thm} \label{thm spectral triple from operator Heisenberg order 1}
 Let $M$ be a compact Heisenberg manifold, and let $P$ be a hypoelliptic self-adjoint $\Psi_\Hei DO$ of order $1$ acting on a vector bundle $E$ over $M$, which is bounded from below.\smallskip

Then the triple
$$\left( \mathcal{A}, \mathcal{H}, D \right) = \left( C(M), L^2(M,E), P \right)$$
is a spectral triple, where the representation $\pi: C(M) \rightarrow \mathcal{B}(L^2(M,E))$ is given by left multiplication of a function $f \in C(M)$. \smallskip

\B We have to show that 
\begin{enumerate}[(i)]
 \item There is a dense sub-algebra $\mathcal{A}' \subset C(M)$ such that $[P,\pi(f)]$ is bounded for any $f \in \mathcal{A}'$.
 \item The operator $P$ has a compact resolvent.
\end{enumerate}
Statement (i) is immediately clear from Proposition \ref{prop boundedness of commutator}, since $C^\infty(M)$ is a dense sub-algebra of $C(M)$. And statement (ii) follows immediately from Theorem \ref{thm complex powers of Heisenberg pseudos}: By this theorem, the operator $(P^2 + I)^{-1/2}$ exists in the Heisenberg calculus and is a $\Psi_\Hei DO$ of order $-1$. Since $M$ is compact, this is a compact operator by Corollary \ref{cor Heisenberg pseuDOs of negative order are compact}, and therefore $P$ possesses a compact resolvent. \eB \end{thm}\medskip

From the above considerations we find that it is in general possible to construct spectral triples via the Heisenberg calculus. We will see that concrete examples can be comfortably derived from horizontal Laplacians, and by the results stated in Section 5.2 concerning the eigenvalue asymptotics we will also see that the metric dimension of these spectral triples provides the Hausdorff dimension of $(M,d_{CC})$. \medskip

\begin{thm} \label{thm spectral triple from horizontal Laplacian}
 Let $M$ be a compact Heisenberg manifold equipped with a horizontal distribution of rank $d$, and let $\Delta^\mathrm{hor}$ be a hypoelliptic self-adjoint horizontal Laplacian which is bounded from below, acting on a vector bundle $E$ over $M$. Then the triple
\begin{equation} \label{eq spectral triple from horizontal Laplacian}
 \left( C(M), L^2(M,E), \left( \Delta^\mathrm{hor} \right)^{\frac{1}{2}} \right)
\end{equation}
is a spectral triple of metric dimension $d+2$.\smallskip

In particular, the metric dimension of this spectral triple coincides with the Hausdorff dimension of the Carnot manifold $(M,d_{CC})$. \smallskip

\B Without loss of generality we assume $\Delta^\mathrm{hor}$ to be invertible (since it is bounded from below, this can always be achieved by adding a constant). Thus, by Theorem \ref{thm complex powers of Heisenberg pseudos}, $\left( \Delta^\mathrm{hor} \right)^{1/2}$ exists and is a $\Psi_\Hei DO$ of order one, which is hypoelliptic, self-adjoint and positive. Then the statement that \eqref{eq spectral triple from horizontal Laplacian} is a spectral triple follows from Theorem \ref{thm spectral triple from operator Heisenberg order 1}. \smallskip

To get the additional statement about the metric dimension of this spectral triple, we can apply Theorem \ref{thm asymptotic behaviour of eigenvalues horizontal Laplacian}, which says that the eigenvalues of $\Delta^\mathrm{hor}$ have the asymptotic behavior
$$ \lambda_k\left(\Delta^{\mathrm{hor}}\right) \sim \left( \frac{k}{\nu_0\left(\Delta^{\mathrm{hor}}\right)} \right)^\frac{2}{d+2},$$
which gives us
$$ \lambda_k\left( \left(\Delta^{\mathrm{hor}}\right)^\frac{1}{2} \right) \sim \left( \frac{k}{\nu_0\left(\Delta^{\mathrm{hor}}\right)} \right)^\frac{1}{d+2}.$$
This asymptotic behavior shows that $\left(\Delta^\mathrm{hor}\right)^{-p/2}$ is trace class if and only if $p>d+2$. Thus the metric dimension of the spectral triple is $d+2$, which is also the Hausdorff dimension of $(M,d_{CC})$ by Theorem \ref{thm Mitchells measure theorem}. \eB
\end{thm}

From now on, we consider once again the orthonormal frame $\{X_1, \ldots, X_d, X_{d+1}\}$ of the tangent bundle of our Heisenberg manifold $M$, such that $\{X_1, \ldots, X_d\}$ span the horizontal distribution $HM$. We assume that the Levi form of $M$ according to this frame is given by the matrix
\begin{equation} \label{eq Levi form local frame section 7.1}
L = \begin{pmatrix} 0&D_m&0 \\ -D_m&0&0 \\ 0&0&0 \end{pmatrix} \ \ \ \text{with} \ D_m = \begin{pmatrix} \lambda_1& & \\ & \ddots & \\ & & \lambda_m \end{pmatrix}, \ \lambda_1, \ldots, \lambda_m > 0
\end{equation}
for an $m \leq d/2$. Note that such a Levi-form can always be achieved by an orthonormal change of the frame (see Proposition \ref{prop tangent Lie group bundle and rank of Levi form}). But this means we have $\G \cong \Hei^{2m+1} \times \R^{d-2m}$ for the underlying Carnot group, and the commutator relation of the vector fields belonging to the frame of $M$ become
$$[X_j,X_k] = \begin{cases}
               \lambda_j, & \text{for} \, 1 \leq j \leq m, \, k = m+j \\
	       -\lambda_k & \text{for} \, 1 \leq k \leq m, \, j = k+m\\
	       0 & \text{otherwise}
              \end{cases}. $$

Now an obvious candidate for a horizontal Laplacian which provides a spectral triple is the sum-of-squares operator
\begin{equation} \label{eq sum of squares}
 \Delta^H = \nabla^\ast_{X_1} \nabla_{X_1} + \ldots + \nabla^\ast_{X_d} \nabla_{X_d},
\end{equation}
where $\nabla$ is any connection acting on a vector bundle $E$ over $M$. This operator is obviously self-adjoint, and after applying Hörmander's sum-of-squares theorem (see Theorem \ref{thm Hormander sum of squares}) in any local chart we see that is is  hypoelliptic. It is also known that this operator is bounded from below, as it is the case for any self-adjoint horizontal Laplacian, whose principal symbol is invertible in the Heisenberg calculus. See \cite{Pon1}, Remark 5.2.10 for that. \smallskip

As we have seen in the previous chapters, the horizontal Dirac operator $D^H$ detecting the Carnot-Carath\'{e}odory metric on $M$ is not hypoelliptic and therefore does not provide a spectral triple. To see why the hypoellipticity fails, we take a look at the proof of Theorem \ref{thm horizontal Dirac is not hypoelliptic} adapted to this situation: The square of $D^H$ is (locally) given by
\begin{equation} \label{eq D^H squared locally Section 7.1}
 \left(D^H\right)^2 = -\sum_{j=1}^d X_j^2 + \sum_{j=1}^d \lambda_j c^H(X_j)c^H(X_{m+j}) X_{d+1} + O_H(1),
\end{equation}
and this operator fails to be hypoelliptic since $i\sum_{j=1}^m \lambda_j$ and $-i\sum_{j=1}^m \lambda_j$ are eigenvalues of the matrix $\sum_{j=1}^d \lambda_j c(X_j)c(X_{m+j})$. This is the content of Theorem \ref{thm hypoellipticity of laplacian codimension 1}, but the same theorem also states that the operator from \eqref{eq D^H squared locally Section 7.1} would be hypoelliptic if the absolute value of every eigenvalue of the matrix coefficient of $X_{d+1}$ would be smaller than $\sum_{j=1}^m \lambda_j$. Since this is fulfilled for any other eigenvalue of $\sum_{j=1}^d \lambda_j c(X_j)c(X_{m+j})$ (see Proposition \ref{prop eigenvalues sum c(X_j)c(X_k)}), one can disturb this operator by a small number $\theta > 0$ and conclude that the operator
$$\left(D^H_\theta\right)^2 := -\sum_{j=1}^d X_j^2 + (1-\theta) \sum_{j=1}^d \lambda_j c^H(X_j)c^H(X_{m+j}) X_{d+1} + O_H(1)$$
is hypoelliptic for any $0 < \theta \leq 1$.\smallskip

We thus construct a spectral triple which shall be a small perturbation of the horizontal Dirac operator, using the sum-of-squares operator from \eqref{eq sum of squares}. The connection we are using is the horizontal Clifford connection $\nabla^{S^H}$ from Section 3.2, acting on a Clifford bundle $S^HM$ over $M$ on which the horizontal Dirac operator $D^H$ is defined. \medskip

\begin{cor} \label{cor disturbed spectral triple on Heisenberg manifold}
 Let $M$ be a closed Heisenberg manifold, and let $S^HM$ be a horizontal Clifford bundle over $M$ on which a horizontal Dirac operator $D^H$ is acting. We further define the horizontal Laplacian
$$ \Delta^H = \left(\nabla^{S^H}_{X_1}\right)^\ast \nabla^{S^H}_{X_1} + \ldots + \left(\nabla^{S^H}_{X_d}\right)^\ast \nabla^{S^H}_{X_d}. $$
Then for any $0 < \theta \leq 1$ the operator
\begin{equation} \label{eq D^H_theta on Heisenberg manifold}
 D^H_\theta := \left( (1-\theta)\left(D^H\right)^2 + \theta \Delta^H \right)^\frac{1}{2}
\end{equation}
provides a spectral triple $(C(M), L^2(S^HM), D^H_\theta)$ of metric dimension $d+2$.\smallskip

\B We have to show that for any $\theta \in (0,1]$ the operator 
\begin{equation} \label{eq Delta^H_theta}
 \Delta^H_\theta := (1-\theta)\left(D^H\right)^2 + \theta \Delta^H \ \in \Psi_\Hei^2(S^HM)
\end{equation}
is positive (which also means that \eqref{eq D^H_theta on Heisenberg manifold} is well-defined) and hypoelliptic. In this case we can apply Theorem \ref{thm spectral triple from horizontal Laplacian} to get that $D^H_\theta$ provides a spectral triple of metric dimension $d+2$ and the corollary is proved. \smallskip

We assume that $\Delta^H \geq 0$, which can always be achieved by adding a constant because it is bounded from below. Since the horizontal Dirac operator $D^H$ is self-adjoint by Theorem \ref{thm D^H selfadjoint}, its square $\left(D^H\right)^2$ is positive, and therefore $\Delta^H_\theta$ from \eqref{eq Delta^H_theta} is positive as a convex combination of positive operators. \smallskip

To show that $\Delta^H_\theta$ is hypoelliptic, we note that
\begin{equation} \label{eq Delta^H_theta via difference}
 \Delta^H_\theta = \theta \Delta^H + \left(1-\theta\right) \left(D^H\right)^2 = \Delta^H + \left(1-\theta\right) \left(\left(D^H\right)^2 - \Delta^H\right).
\end{equation}
We consider this operator in the local frame $\{X_1, \ldots, X_{d+1}\}$ introduced above, such that $HM = \spa \{X_1, \ldots, X_d\}$ and the matrix of the Levi form corresponding to this basis is given by \eqref{eq Levi form local frame section 7.1}. But in these local coordinates, the argument we sketched in the discussion preceding this corollary applies: We have locally $\Delta^H = -\sum_{j=1}^d X_j^2$ and therefore 
$$\left( D^H \right)^2 - \Delta^H = \sum_{j=1}^m c^H(X_j)c^H(X_{m+1}) X_{d+1} + O_H(1) $$
by the local expression of $(D^H)^2$ from Proposition \ref{prop D^H squared}. Hence \eqref{eq Delta^H_theta via difference} has locally the form
\begin{equation} \label{eq Delta^H_theta local coordinates}
 \Delta^H_\theta = -\sum_{j=1}^d X_j^2 + \left(1-\theta\right) \sum_{j=1}^m c^H(X_j)c^H(X_{m+1}) X_{d+1} + O_H(1),
\end{equation}
where like in the previous chapters $O_H(1)$ denotes a graded differential operator of order smaller than or equal to $1$. Now by Proposition \ref{prop eigenvalues sum c(X_j)c(X_k)}, the eigenvalues of the matrix 
$$\left(1-\theta\right) \sum_{j=1}^m \lambda_j c^H(X_j)c^H(X_{m+j})$$
 are included in the interval 
$$\left[-(1-\theta) i\sum_{j=1}^d \lambda_j, (1-\theta) i\sum_{j=1}^d \lambda_j\right] \ \subset \R i$$
on the imaginary line, and since $1 \geq \theta > 0$ the absolute value of any of these eigenvalues is smaller than $\sum_{j=1}^d \lambda_j = \frac{1}{2} \tr | L|$. But this means by Theorem \ref{thm hypoellipticity of laplacian codimension 1} that $\Delta^H_\theta$ is (locally) hypoelliptic.\smallskip

The global hypoellipticity then follows by Theorem \ref{thm hypoellipticity and existence of parametricies horizontal Laplacian}, and the theorem is proved. \eB
\end{cor}\bigskip

\section{Detection of the Metric via Horizontal Laplacians}

After we have constructed spectral triples from horizontal Laplacians, we ask ourselves whether one can detect the Carnot-Carath\'{e}odory metric from these spectral triples. While we postpone the discussion of estimates for the Connes metric formula of these operators to Section 7.3, we want to show now that there is a formula to detect the Carnot-Carath\'{e}odory metric directly from any horizontal Laplacian. In the case where we have an ordinary Laplacian $\Delta = -\sum_{j=1}^n X_j^2$, it is known that for $f \in C^\infty(M)$ we have
$$ \left[ [\Delta,f],f \right] = -2 \left\|df\right\|^2,$$
see for example \cite{BGV}, Proposition 2.3. Now an analogous result is true for horizontal Laplacians: The key observation is the following lemma.\medskip

\begin{lemma} \label{lemma double commutator from horizontal Laplacian}
 Let $\Delta^\mathrm{hor}$ be a horizontal Laplacian, acting on a vector bundle $E$ over a Carnot manifold $M$. Then we have for any function $f \in C^\infty(M)$
\begin{equation} \label{eq double commutator horizontal Laplacian}
 \frac{1}{2} \left[\left[\Delta^\mathrm{hor},f\right],f\right] \sigma = - \left\| \grad^H f \right\|^2 \sigma
\end{equation} 
for $\sigma \in \Gamma^\infty(M,E)$. Here, $\| \cdot \|$ denotes the (point-wise) norm of a vector in $E_xM$ induced by the Riemannian metric $g$ on $M$.\smallskip

\B Let $TM = V_1M \oplus \ldots \oplus V_RM$ be the grading of $M$, such that $\{X_1, \ldots, X_d\}$ is an orthonormal frame of $V_1M$ and $\{X_{d+1}, \ldots, x_{d+d_2}\}$ is an orthonormal frame of $V_2M$. We prove the statement locally, meaning that $\Delta^\mathrm{hor}$ is given in the form
\begin{equation} \label{eq horizontal Laplacian locally for double commutator}
 \Delta^\mathrm{hor} = -\sum_{j=1}^d \partial_{X_j}^2 + B(\partial_{X_1}, \ldots, \partial_{X_{d+d_2}}) + b,
\end{equation}
where $\partial_{X_k}$ denotes the partial derivative along $X_k$ in any chart, $B(\partial_{X_1}, \ldots, \partial_{X_{d+d_2}})$ is any differential operator of (classical) order $1$ and $b$ is a (matrix-valued) function.\smallskip

We now plug the expression \eqref{eq horizontal Laplacian locally for double commutator} term by term into the double commutator of \eqref{eq double commutator horizontal Laplacian} and use the linearity of the commutator. Since $b$ commutes with any smooth function $f$, this is zero for the last summand. It is also zero for the second summand of \eqref{eq horizontal Laplacian locally for double commutator}, since for any first order differential operator $B = B(\partial_{X_1}, \ldots, \partial_{X_{d+d_2}})$ we have because of the Leibniz rule
\begin{eqnarray*}
 \left[\left[B,f\right],f\right]\sigma &=& \left[Bf-fB,f\right] \sigma\\
&=& B(f^2 \sigma) - 2f \cdot B(f\sigma) +f^2 \cdot B\sigma\\
&=& 2(Bf)\cdot f\sigma + f^2 B\sigma - 2f \cdot (Bf) \cdot \sigma -2f^2 \cdot B\sigma +f^2 \cdot B\sigma\\
&=& 0
\end{eqnarray*}
locally for any section $\sigma \in \Gamma^\infty(M,E)$.\smallskip

Calculating the first term of \eqref{eq horizontal Laplacian locally for double commutator}, we find that for the double commutator applied to any $\partial_{X_j}^2$ we have
\begin{eqnarray*}
&& \left[\left[\partial_{X_j}^2,f\right],f\right] \sigma \\
&=& \partial_{X_j}^2 (f^2\sigma) - 2f \cdot \partial_{X_j}^2(f\sigma) + f^2 \cdot \partial_{X_j}^2 \sigma\\
&=& \partial_{X_j} \left(2(\partial_{X_j}f) \cdot f\sigma + f^2 \cdot \partial_{X_j}\sigma\right) -2f \cdot \partial_{X_j} \left((\partial_{X_j}f) \cdot \sigma + f \cdot \partial_{X_j} \sigma\right) + f^2 \cdot \partial_{X_j}^2 \sigma \\
&=& 2(\partial_{X_j}^2f) \cdot f\sigma + 2(\partial_{X_j} f)^2 \cdot \sigma + 4 (\partial_{X_j} f) \cdot f \cdot \partial_{X_j} \sigma + f^2 \cdot \partial_{X_j}^2 \sigma\\
& & -2f \cdot (\partial_{X_j}^2 f) \cdot \sigma -4f \cdot (\partial_{X_j} f) \cdot (\partial_{X_j} \sigma) -2f^2 \cdot \partial_{X_j}^2 \sigma + f^2 \cdot \partial_{X_j}^2 \sigma\\
&=& 2(\partial_{X_j} f)^2 \cdot \sigma.
\end{eqnarray*}
Finally, we plug everything together into the double commutator expression \eqref{eq horizontal Laplacian locally for double commutator}, which shows us that
\begin{eqnarray*}
 \frac{1}{2} \left[\left[\Delta,f\right],f\right] \sigma &=& -\frac{1}{2}\sum_{j=1}^d \left[\left[\partial_{X_j}^2,f\right],f\right] \sigma \\
&=& -\sum_{j=1}^d (\partial_{X_j} f)^2 \cdot \sigma\\
&=& - \left\| \grad^H f \right\|^2 \sigma.
\end{eqnarray*}
Hence the statement of the lemma is proved. \eB
\end{lemma} \medskip

From the identity \eqref{eq double commutator horizontal Laplacian} we get an expression depending on a horizontal Laplacian instead of a horizontal Dirac operator which can detect the Carnot-Carath\'{e}odory metric, following the theory from Section 3.3. By Corollary \ref{cor d_CC from D^H via smooth functions}, we have for the Carnot-Carath\'{e}odory metric
\begin{equation} \label{eq detection CC metric from horizontal gradient chapter 7}
 d_{CC}(x,y) = \sup \left\{ |f(x)-f(y)|: f \in \mathrm{Lip}_{CC}(M), \mathrm{ess} \sup_{\xi \in M} \left\| \grad^H f(\xi) \right\| \leq 1\right\}.
\end{equation}
Using this expression, we can prove the following theorem.\medskip

\begin{thm} \label{thm metric detection from horizontal Laplacian}
 Let $\Delta^\mathrm{hor}$ be a horizontal Laplacian, acting on a vector bundle $E$ over a closed Carnot manifold $M$. Then we have for any $x,y \in M$
\begin{equation} \label{eq metric detection horizontal laplacian}
 d_{CC}(x,y) = \sup \left\{ |f(x) - f(y)|: f \in C^\infty(M), \left\| \frac{1}{2} \left[\left[\Delta^{\mathrm{hor}},f\right],f\right] \right\| \leq 1 \right\},
\end{equation}
where $d_{CC}$ is the Carnot-Carath\'{e}odory metric on $M$. \smallskip

\B By Lemma \ref{lemma double commutator from horizontal Laplacian} we have
$$\left\| \frac{1}{2} \left[\left[\Delta^{\mathrm{hor}},f\right],f\right] \right\| = \sup_{x \in M} \left\| \grad^H f(x) \right\|^2.$$
In addition, we have shown in Section 3.3 that
$$\left\| [D^H,f] \right\| = \sup_{x \in M} \left\| \mathrm{grad}^H f (x) \right\|$$
for any smooth function $f \in C^\infty(M)$. Thus \eqref{eq detection CC metric from horizontal gradient chapter 7} provides us
$$ d_{CC}(x,y) = \mathrm{ess} \sup \left\{ |f(x) - f(y)|: f \in C^\infty(M), \ \sup_{\xi \in M} \left\| \grad^H f(\xi) \right\| \leq 1\right\}.$$
Since $\left\| \grad^H f(x) \right\|^2 \leq 1$ if and only if $\left\| \grad^H f(x) \right\| \leq 1$, the statement of the theorem follows from Lemma \ref{lemma double commutator from horizontal Laplacian}.\eB
\end{thm}\medskip

Theorem \ref{thm metric detection from horizontal Laplacian} states that it is possible to detect the Carnot-Carath\'{e}odory metric by arbitrary positive horizontal Laplacians. In the case $D^H$ is a horizontal Dirac operator constructed in Chapter 3, both $D^H$ and $(D^H)^2$ detect the Carnot-Carath\'{e}odory metric via the corresponding formulas.\smallskip

Now it is a natural question whether one can use the metric detection by a positive horizontal Laplacians to find estimates for the Connes metric given by its square root operator. For the rest of this section, we take the algebra $C^\infty(M)$ as the dense sub-algebra of $C(M)$ over which the Connes metric is defined, i.e. we consider the metric
\begin{equation} \label{eq Connes metric from smooth functions}
 d_{D}(x,y) = \sup \left\{ |f(x)-f(y)|: f \in C^\infty(M), \left\| [D,f] \right\| \leq 1\right\}
\end{equation}
for an operator $D$ defining a spectral triple. Then we can note the following criterion.\medskip

\begin{prop} \label{prop metric estimate square root of horizontal Laplacian from below}
Let $\tilde{D}$ be a self-adjoint operator acting on a closed Carnot manifold $M$ and defining a spectral triple such that $\tilde{D}^2 = \Delta^\mathrm{hor}$ is a horizontal Laplacian. Let $d_{\tilde{D}}$ denote the Connes metric defined by $\tilde{D}$ via \eqref{eq Connes metric from smooth functions}. \smallskip

Then if the condition
\begin{equation} \label{eq condition for metric estimate from below}
 \left[\tilde{D},f\right]^2 = \frac{1}{2}\left[[\Delta^\hor,f],f\right]
\end{equation}
is fulfilled, we have the estimate
\begin{equation*} \label{eq metric estimate square root of horizontal Laplacian from above}
 d_{\tilde{D}}(x,y) \leq d_{CC}(x,y)
\end{equation*}
for all points $x,y \in M$.\smallskip

\B This follows immediately from the formulas for the metrics. By \eqref{eq condition for metric estimate from below} we see that
\begin{equation} \label{eq estimate operator norm order 2}
 \left\| \frac{1}{2} \left[[\Delta^\hor,f],f\right] \right\| = \left\| \left[\tilde{D},f\right]^2 \right\| \leq \left\| \left[\tilde{D},f\right]\right\|^2. 
\end{equation}
But this shows that $\|[\tilde{D},f]\| \leq 1$ implies $\|[[\Delta^\hor,f],f]\| \leq 1$, which gives together with Theorem \ref{thm metric detection from horizontal Laplacian} the desired estimates for the corresponding metrics:
\begin{eqnarray*}
 d_{\tilde{D}}(x,y) &=& \sup \left\{ \left|f(x)-f(y)\right|: \left\| \left[\tilde{D},f\right]\right\| \leq 1 \right\}\\
&\leq& \sup \left\{ \left|f(x)-f(y)\right|: \left\| \frac{1}{2}\left[[\Delta^\hor,f],f\right] \right\| \leq 1 \right\}\\
&=& d_{CC}(x,y).
\end{eqnarray*}
  \eB
\end{prop}\medskip

Now in case $D^H$ is a horizontal Dirac operator, condition \eqref{eq condition for metric estimate from below} is fulfilled as we will see in a minute. In this situation we even have equality of the metrics, since estimate \eqref{eq estimate operator norm order 2} in the proof of the proposition is an equality. The reason for this is that $[D^H,f]$ is exactly the horizontal Clifford action $c^H$ of the horizontal gradient of $f$, and that at each point $x\in M$ $c^H: H_xM \rightarrow \en_\C(S^HM)$ is an isometry.\smallskip

We now develop a condition equivalent to \eqref{eq condition for metric estimate from below} from elementary commutator calculations.\medskip

\begin{prop} \label{prop condition [A,B]^2 = 1/2[[A^2,B],B]}
 Let $A,B$ be linear operators on a Hilbert space. Then we have
\begin{equation*} \label{eq [A,B]^2 = 1/2[[A^2,B],B]}
 \left[A,B\right]^2 = \frac{1}{2} \left[\left[A^2,B\right]B\right]
\end{equation*}
if and only if
\begin{equation} \label{eq A[[A,B],B]+[[A,B],B]A=0}
 A\left[[A,B],B\right] + \left[[A,B],B\right]A = 0
\end{equation}\smallskip
\B This is just a simple calculation involving commutator rules:
\begin{eqnarray*}
 \left[\left[A^2,B\right]B\right] &=& \left[A[A,B]+[A,B]A,B\right]\\
&=& \left[A[A,B],B\right] + \left[[A,B]A,B\right]\\
&=& A \left[[A,B],B\right] + [A,B][A,B] + [A,B][A,B] + \left[[A,B],B\right]A.
\end{eqnarray*}
From this equation, we see the equivalence of the two statements. \eB
\end{prop}\medskip

From the aspect of spectral triples, we choose in Proposition \ref{prop condition [A,B]^2 = 1/2[[A^2,B],B]} $A=D$ (the Dirac operator of the triple) and $B=f$ (the representation of the corresponding algebra on the Hilbert space, which is multiplication by functions in the commutative case). Hence it is obvious that \eqref{eq A[[A,B],B]+[[A,B],B]A=0} is fulfilled if $D$ is a first-order differential operator, since in this case we have $[D,f] = 0$. 

More generally, \eqref{eq A[[A,B],B]+[[A,B],B]A=0} is fulfilled if the operator $D$ fulfills the order-one condition for a spectral triple by Alain Connes (see e.g. \cite{Con2}), or \cite{GVF}, Section 10.5), which states in the commutative situation that $[[D,f],g] = 0$ for all $f,g \in C^\infty(M)$. Therefore we can write down the following corollary. \medskip

\begin{cor}
Let $M$ be a closed Carnot manifold. If $(C(M), L^2(\Sigma M), \tilde{D})$ is a spectral triple such that $\tilde{D}^2=\Delta^\mathrm{hor}$ is a horizontal Laplacian, which fulfills the order one condition, then we have for any points $x,y \in M$
$$d_{CC}(x,y) \geq \sup\left\{|f(x)-f(y)|: f \in C^\infty(M), \left\| [\tilde{D},f]\right\| \leq 1 \right\}.$$\smallskip

\B This follows immediately by Proposition \ref{prop condition [A,B]^2 = 1/2[[A^2,B],B]} and Proposition \ref{prop metric estimate square root of horizontal Laplacian from below} together with the above discussion. \eB
\end{cor}\medskip

But despite this corollary, we note that the condition \eqref{eq condition for metric estimate from below} or, equivalently by Proposition \ref{prop condition [A,B]^2 = 1/2[[A^2,B],B]},
$$\tilde{D}\left[[\tilde{D},f],f\right] + \left[[\tilde{D},f],f\right]\tilde{D} = 0$$
seem to be rather strong conditions for the case where $\tilde{D}$ is not a differential operator.\bigskip

\section{Approaches for Approximation of the Metric}

Although the horizontal Dirac operator, which detects the Carnot-Carath\'{e}odory metric, does not furnish a spectral triple, we have seen in Section 7.1 that there are quite a lot of spectral triples arising from Heisenberg calculus which give at least the right dimension. In particular the spectral triples constructed in Corollary \ref{cor disturbed spectral triple on Heisenberg manifold} only differ by a small parameter from the horizontal Dirac operator. But sadly we do not know how their Connes metric behaves with respect to the Carnot-Carath\'{e}odory metric.\smallskip

Now in this final section we want to present a few ideas how one can approximate the Carnot-Carath\'{e}odory metric by a family of spectral triples. Therefore we use the observations made in Section 1.2 concerning arbitrary compact quantum metric spaces to develop criteria for this approximation. First of all, we give a reformulation of Proposition \ref{prop uniformly continuity implies convergence of metrices} and Proposition \ref{prop equivalence condition for metricies} for our setting.\medskip

\begin{prop} \label{prop conditions for metric convergence to d_CC}
Let $M$ be a closed Carnot manifold, where $d_{CC}$ is the Carnot-Carath\'{e}odory metric on $M$. Let further $D^H$ denote a horizontal Dirac operator on $M$, acting on a horizontal Clifford bundle $S^HM$.\smallskip

\begin{enumerate}[(i)]
 \item Let $\tilde{D}$ be an operator on $L^2(S^HM)$ such that $(C(M), L^2(S^HM), \tilde{D})$ is a spectral triple. If there exists a constant $0 < C < 1$ such that 
\begin{equation} \label{eq equivalence approximation spectral triples}
 \left\| [\tilde{D} - D^H,f] \right\| \leq C \left\| [D^H,f] \right\|
\end{equation}
for any $f$ belonging to a suitable sub-algebra of $C(M)$, then the Connes metric $d_{\tilde{D}}$ arising from $\tilde{D}$ is equivalent to $d_{CC}$.

\item Let for $0 < \theta \leq 1$ the family $\tilde{D}_\theta$ be a family of operators on $S^HM$ with the property that
\begin{equation} \label{eq convergence approximation spectral triples}
  \forall \varepsilon > 0 \ \exists \delta >0: \ 0 < \theta < \delta \ \Rightarrow \ \left\| [\tilde{D}_\theta - D^H,f] \right\| < \varepsilon \ \forall f \in \Sigma_0,
\end{equation}
 where $\Sigma_0 := \{ f \in \mathcal{A}': \mathrm{Lip}_{CC}(f) = 1\}$ for a suitable sub-algebra $\mathcal{A}'$ of $C(M)$. Then for every $\varepsilon > 0$ with $\varepsilon < 1$ there is a $\delta > 0$ such that
\begin{equation} \label{eq equivalence and convergence of metricies from D}
 \left(1-\varepsilon\right) d_{\tilde{D}_\theta}(x,y) \leq d_{CC}(x,y) \leq \left(1+\varepsilon\right) d_{\tilde{D}_\theta}(x,y) \ \ \ \forall x,y \in M 
\end{equation}
 for every $0 < \theta < \delta$. 
\end{enumerate} \smallskip

\B We set $L_0(f) := \| [D^H,f] \|$, which is a Lip-norm on $C(M)$ coinciding with the Carnot-Carath\'{e}odory-Lipschitz constant $\mathrm{Lip}_{CC}(f)$ of $f$ by Corollary \ref{cor D^H provides compact quantum metric space}. Further $L_1 := \| [\tilde{D},f] \|$ in (i) and $L_\theta := \| [\tilde{D}_\theta,f]\|$ in (ii) are also Lip-norms on $C(M)$, since they arise from spectral triples. But in this context (i) is just a reformulation of Proposition \ref{prop equivalence condition for metricies} and (ii) is just a reformulation of Proposition \ref{prop uniformly continuity implies convergence of metrices}. \eB
\end{prop}

This proposition implies that we would be able to prove good metric approximations by spectral triples for the Carnot-Carath\'{e}odory metric if the following (far more general) assumption is true.\medskip

\begin{ass} \label{ass estimate for general first-order operator}
 Let $P$ be a first order Heisenberg pseudodifferential operator, acting on a Heisenberg manifold $M$, and let $f \in C^\infty(M)$ or (more generally) $f \in \mathrm{Lip}_{CC}(M)$. Then there exists a $C > 0$ such that
\begin{equation} \label{eq metric estimate for general first order PseuHDO}
 \left\| [P,f] \right\| \leq C \cdot \mathrm{Lip}_{CC}(f),
\end{equation}
where $\mathrm{Lip}_{CC} (f)$ denotes the Carnot-Carath\'{e}odory-Lipschitz constant of $f$ from Definition \ref{def Lip_CC}. \eBsp
\end{ass}\medskip

This assumption is motivated by the fact that something analogous is true for the classical Lipschitz functions $\mathrm{Lip}(M)$ and for classical pseudodifferential operators of order $1$ belonging to the class $S^1_{1,0}$. For this classical situation, the statement can be found in \cite{Tay2}, Section 3.6. But the proof presented in the book by Taylor is quite involved, such that it is far away from being trivial to transfer this proof to the Carnot-Carath\'{e}odory situation.\smallskip

Unfortunately, we have not found a way to prove the estimates \eqref{eq equivalence approximation spectral triples} or \eqref{eq convergence approximation spectral triples}, and the proof of the above assumption (if possible) seems to be an even harder and much involved problem. But we assume that the Connes metrics arising from the spectral triples constructed in Corollary \ref{cor disturbed spectral triple on Heisenberg manifold} are a good approximation of the Carnot-Carath\'{e}odory metric (note that we have in the formulation of the corollary $D^H_0 = |D^H|$ with $D^H$ the horizontal Dirac operator).\smallskip

As an open problem, we present an approximation of the degenerate spectral triple 
$$(C(M),L^2(\Sigma^H_\delta M),D^H)$$
for the case where $M$ is the compact nilmanifold of a Heisenberg group $\Hei^{2m+1}$, and $D^H$ is the horizontal pull-back Dirac operator arising from a spin structure $\delta$ of the $2m$-dimensional torus constructed in Chapter 4. Let $\{X_1, \ldots, X_d\}$ be an orthonormal frame of the horizontal distribution $HM$ of $M$ (with $d=2m$), such that the horizontal pull-back Dirac operator is given by
$$D^H = \sum_{j=1}^d c^H(X_j) \nabla^{\Sigma^H}_{X_j},$$
see Section 4.1 for the explicit construction. We have seen in Section 4.3 that the only space where this horizontal Dirac operator degenerates is its infinite dimensional kernel, but otherwise it has all the properties we ask for (see Corollary \ref{cor D^H gives degenerate spectral triple} and Corollary \ref{cor Hausdorff dimension from D^H}).\smallskip

The idea is now to fix the gap of degeneracy on the kernel of $D^H$ by using the sum-of-squares operator
$$\Delta^H = \sum_{j=1}^d \left(\nabla^{\Sigma^H}_{X_j}\right)^\ast \nabla^{\Sigma^H}_{X_j}.$$
It is known that $\Delta^H$ is hypoelliptic and bounded from below (see also Section 7.1 for this). Now let $P$ denote the orthogonal projection onto the kernel of $D^H$. Since the spectrum of $D^H$ is discrete, this projection operator is a pseudodifferential operator of order $0$, and it is also an order-zero operator in the Heisenberg calculus. Then for a small $\theta > 0$ we define an operator $D^H_\theta$ via
\begin{equation} \label{eq D^H_theta approximating family}
 D^H_\theta := D^H + \theta \cdot P \left(\Delta^H\right)^\frac{1}{2} P. 
\end{equation}
By the above argumentation, this is a hypoelliptic Heisenberg pseudodifferential operator, and hence it does define a spectral triple. In addition it detects the metric dimension of $(M,d_{CC})$ via its eigenvalue asymptotics, since we have the right eigenvalues asymptotics on $(\mathrm{ker} D^H)^\perp$ by Corollary \ref{cor Hausdorff dimension from D^H} and on $\mathrm{ker} D^H)$ because of the eigenvalue asymptotics of $\Delta^H$ (see Theorem \ref{thm asymptotic behaviour of eigenvalues horizontal Laplacian}).\smallskip

Now if Assumption \ref{ass estimate for general first-order operator} is true, we have a constant $C>0$ such that
\begin{equation} \label{eq metric estimate on kernel}
 \left\| [P \left(\Delta^H\right)^\frac{1}{2} P,f]\right\| \leq C \cdot \mathrm{Lip}_{CC}(f).
\end{equation}
Note that in this situation it would be enough to have the estimate \eqref{eq metric estimate on kernel} on the kernel of $D^H$, which may be much easier to prove than the general case of \eqref{eq metric estimate for general first order PseuHDO}. This would lead to the estimate
$$ \left\| [D^H - D^H_\theta,f] \right\| \leq \left\| \theta [P \left(\Delta^H\right)^\frac{1}{2} P,f]\right\| \leq \theta C \cdot \mathrm{Lip}_{CC}(f),$$
because of the definition \eqref{eq D^H_theta approximating family} of $D^H_\theta$, and we could use Proposition \ref{prop conditions for metric convergence to d_CC} to show that $d_{D^H_\theta}(x,y) \rightarrow d_{CC}(x,y)$ for $\theta \rightarrow 0$ for all points $x,y \in M$.

\end{document}